\documentclass[10pt, eqno]{article}
\usepackage{amsmath,amssymb,latexsym}

 \newfont{\Bb}{msbm10 scaled\magstephalf}

\begin{document}

\title{Virtual moduli cycles and Gromov-Witten
 invariants \\ of noncompact symplectic manifolds}

\author{\\ Guangcun Lu
\thanks{The author was supported in part by NNSF 19971045 and 10371007 of China.}
\\
Department of Mathematics,  Beijing Normal University\\
Beijing 100875,   P.R.China\\
(gclu@bnu.edu.cn)}
\date{Received: 8 October 2003\\
Accepted: 16 November 2004} \maketitle \vspace{-0.1in}

\begin{abstract}
This paper constructs and studies the Gromov-Witten invariants
 and their properties for noncompact geometrically bounded symplectic manifolds.
Two localization formulas for GW-invariants are
  also proposed and proved. As applications we get solutions of the generalized string
equation and dilation equation and their variants. The more
solutions of WDVV equation and quantum products on cohomology
groups are also obtained for the symplectic manifolds with
finitely dimensional cohomology groups. To realize these purposes
we further develop the language introduced by Liu-Tian to describe
the virtual moduli cycle (defined by Liu-Tian, Fukaya-Ono,
Li-Tian, Ruan and Siebert).
\end{abstract}
\maketitle


\tableofcontents

\section{Introduction and main results}\label{sec:1}

This paper is a continuation of my paper [Lu1]. In [Lu1] we
defined the Gromov-Witten invariants on semi-positive noncompact
geometrically  bounded symplectic manifolds by combining
Ruan-Tian's method([RT1]) with that of McDuff-Salamon ([McS]), and
also used them to generalize work on the topological rigidity of
Hamiltonian loops on semi-positive closed symplectic manifolds
([LMP]) to the Hamiltonian loops with compact support on this
class of noncompact symplectic manifolds. Not long ago Jun Li and
Gang Tian [LiT1] developed a beautiful algebraic method to
establish general theory of GW-invariants on a smooth projective
algebraic variety ( see [Be] for a different treatment) which
satisfy the quantum cohomology axioms proposed by Kontsevich and
Mannin in [KM] based upon predictions of the general properties of
topological quantum field theory [W2]. Soon after several groups
independently developed different new techniques to construct
Gromov-Witten invariants on any closed symplectic manifolds (cf.
[FuO, LT2, R, Sie]).  Actually the method constructing Floer
homology by Liu-Tian [LiuT1] provided another approach to realize
this goal. These new methods were also used in recent studies of
this field, symplectic topology and Mirror symmetry (cf. [CR,
FuOOO, IP1-2, LiR, LiuT2, LiuT3, LLY1-2, Lu3-4, Mc2] etc.) In this
paper we further develop the virtual moduli cycle techniques
introduced in [LiuT1, LiuT2-3] to generalize work of [Lu1] to
arbitrary noncompact geometrically bounded symplectic manifolds.
It should be noted that unlike the case of closed symplectic
manifolds our invariants do not have  the best invariance
properties as expected. This is natural for noncompact case. But
we give their invariance properties as well as we possibly can.
Actually in many cases one does not need very strong invariants
because the almost complex structures are only auxiliary tools to
realize the purposes in the studies of symplectic topology and
geometry. On the other hand the initial insight came from
physicist Witten's topological $\sigma$ models of two-dimensional
gravity (cf. [W1, W2]). According to Witten [W1] the two
dimensional nonlinear sigma model is a theory of maps from a
Riemannian surface $\Sigma$ to a Riemannian manifold $M$. Gromov's
theory of pseudo-holomorphic curves [Gr2] provided a powerful
mathematical tool for it. The latter requires not only that $M$
should have an almost complex structure, but also that $M$ should
have an additional structure--a symplectic structure. It is
well-known that the existence and classification of the symplectic
structures on a closed almost complex manifold is a difficult
mathematics question. However early in his thesis in 1969 [Gr1]
Gromov completely resolved the corresponding question on an open
almost complex manifold.

Recall that a Riemannian manifold $(M,\mu)$ is called {\it
geometrically bounded} if its sectional curvature is bounded above
and has injectivity radius $i(M, \mu)>0$. (It is weaker than the
usual one for which the absolute value of the sectional curvature
is required to be bounded above.) Such a Riemannian metric is
called geometrically bounded and is always complete. The existence
of such Riemannian metrics were proved in [G]. Let ${\mathcal
R}(M)$ denote the set of all Riemannian metrics on $M$, and
${\mathcal G}{\mathcal R}(M)\subset{\mathcal R}(M)$ be the subset
of those Riemannian metrics on $M$ whose injectivity radius are
more than zero and whose sectional curvatures have the upper
bound. Clearly,${\mathcal G}{\mathcal R}(M)$ might not be
connected even if it is nonempty. Let $(M,\omega)$ be a symplectic
manifold and ${\mathcal J}(M,\omega)$ be the space of all
$\omega$-compatible almost complex structures on $M$. The latter
is contractible and each $J\in {\mathcal J}(M,\omega)$ gives a
  Riemannian metric
  $$g_J(u, v)=\omega(u, Jv),\quad\forall u, v\in TM.$$
A symplectic manifold $(M, \omega)$ without boundary is said to be
 {\it geometrically bounded} if there exist  an almost complex structure
 $J$ and a geometrically bounded Riemannian metric $\mu$ on $M$ such that for
some positive constants $\alpha_0$ and $\beta_0$,
$$
\omega(X, JX)\ge\alpha_0\|X\|^2_\mu\quad {\rm and}\quad |\omega(X,
Y)|\le\beta_0\|X\|_\mu\|Y\|_\mu,\quad\forall X, Y\in TM\eqno(1.1)
$$ (cf. [ALP, Gr2, Sik]). Later we also say such a
$J$ to be $(\omega,\mu)$-{\it geometrically bounded}. In $\S2.1$
we shall prove that (1.1) is equivalent to the fact that $g_J$ is
{\it quasi isometric} to $\mu$, i.e., there exist positive
constants $C, C^\prime$ such that
$$C\|u\|_{\mu}\le\|u\|_{g_J}\le C^\prime\|u\|_{\mu}\quad\forall u\in TM.$$
 Clearly, the quasi isometric metrics have the same completeness and incompleteness.
  Denote by
$${\mathcal J}(M,\omega,\mu)\eqno(1.2)$$
 the set of all $(\omega,\mu)$-geometrically bounded almost complex structures in
 ${\cal J}(M,\omega)$.  In Proposition 2.3 we shall study the connectedness of it.
A  family of Riemannian metrics $(\mu_t)_{t\in[0,1]}$ on $M$ is
said to be {\it uniformly geometrically bounded} if their
injectivity radiuses have a uniform positive lower bound and their
sectional curvatures have a uniform upper bound. A $1$-parameter
family of $(M, \omega_t, J_t, \mu_t)_{t\in [0,1]}$ of
geometrically bounded symplectic manifolds starting at $(M,
\omega_0, J_0, \mu_0)=(M, \omega, J, \mu)$ is called a {\it weak
deformation} of a geometrically bounded symplectic manifold
$(M,\omega, J,\mu)$ if (i) both $(\omega_t)_{t\in [0,1]}$ and
$(J_t)_{t\in [0,1]}$ are smoothly dependent $t$ with respect to
the $C^\infty$-weak topology, (ii) $(\mu_t)_{t\in [0,1]}$ is a
path in ${\mathcal G}{\mathcal R}(M)$  and is also uniformly
geometrically bounded, (iii) there exist uniform positive
constants $\alpha_0$ and $\beta_0$, and a distance $d$ on $M$ such
that for all $t\in [0,1]$ and the distance functions $d_{\mu_t}$
induced by $\mu_t$,
$$\left\{\begin{array}{ll}
d_{\mu_t}\ge d,\\
 \omega_t(X, J_tX)\ge\alpha_0\|X\|^2_{\mu_t}\quad
{\rm and}\\
 |\omega_t(X,
Y)|\le\beta_0\|X\|_{\mu_t}\|Y\|_{\mu_t},\quad\forall X, Y\in TM.
\end{array}\right.\eqno(1.3)$$

 Note that the $C^0$-strong topology ( also  called the {\it Whitney} $C^0$-topology) on
${\mathcal R}(M)$  does not have a countable base at any point for
noncompact $M$  in general. However, all metrics $\mu_t$ of a
$C^0$-strong continuous path $(\mu_t)_{t\in [0,1]}$ in ${\mathcal
R}(M)$ are uniformly quasi isometric, i.e., there exists positive
constants $C_1$ and $C_2$ such that
$$C_1\|u\|_{\mu_0}\le\|u\|_{\mu_t}\le C_2\|u\|_{\mu_0}\;\forall u\in TM\quad{\rm and}\quad
t\in [0,1].$$
 We  define a $C^0$ {\it
super-strong} topology on ${\mathcal G}{\mathcal R}(M)$ as
follows: Let $(S, \ge)$ be a net in ${\mathcal G}{\mathcal R}(M)$,
we say that {\it it converges to $\mu\in{\mathcal G}{\mathcal
R}(M)$ with respect to the $C^0$ super-strong topology} if it
converges to $\mu$ in the $C^0$-strong topology and all metrics in
the net are uniformly geometrically bounded. Clearly, a continuous
path $(\mu_t)_{t\in[0,1]}$ in ${\mathcal G}{\mathcal R}(M)$ with
respect to the $C^0$ super-strong topology is always uniformly
geometrically bounded. However it is still difficult to check that
a path in ${\mathcal G}{\mathcal R}(M)$ is continuous with respect
to the $C^0$ super-strong topology. Observe that the sectional
curvature of a Riemannian metric is uniformly continuous with
respect to the $C^2$-strong topology on ${\mathcal R}(M)$. By the
local coordinate expression of the geodesic equations and the
proof of the fundamental theorem of the linear ordinary
differential equation systems it is not hard to prove that if a
Riemannian metric $g$ on $M$ has  the sectional curvature $K_g$
less than a constant $K$ and the injectivity radius  $i(M,g)>0$
then  there exists a neighborhood ${\mathcal U}(g)$ of $g$ in
${\mathcal R}(M)$ with respect to the $C^\infty$-strong topology
(actually the $C^3$-strong topology seem to be enough) such that
every metric $h\in{\mathcal U}(M)$ has the injectivity radius
$i(M, h)\ge\frac{1}{2}i(M,g)>0$ and the sectional curvature
$K_h\le K+1$.
 Therefore  {\it a continuous
path $(\mu_t)_{t\in [0,1]}$ in ${\mathcal G}{\mathcal R}(M)$ with
respect to the $C^\infty$-strong topology is also a continuous
path with respect to the $C^0$ super-strong topology.}
 A $1$-parameter family of $(M, \omega_t,
J_t, \mu_t)_{t\in [0,1]}$ of geometrically bounded symplectic
manifolds starting at $(M, \omega_0, J_0, \mu_0)=(M, \omega, J,
\mu)$ is called a {\it strong deformation} of a geometrically
bounded symplectic manifold $(M,\omega, J,\mu)$ if (i) both
$(\omega_t)_{t\in [0,1]}$ and $(J_t)_{t\in [0,1]}$ are smoothly
dependent $t$ with respect to the $C^\infty$-strong topology, (ii)
$(\mu_t)_{t\in [0,1]}$ is a continuous path in ${\mathcal
G}{\mathcal R}(M)$ with respect to the $C^0$ super-strong
topology. It is not hard to prove that a  strong deformation of a
geometrically bounded symplectic manifold $(M,\omega, J,\mu)$ must
be a weak deformation of $(M,\omega, J,\mu)$.  Denote by
$$
{\rm Symp}_0^S(M,\omega)\eqno(1.4)
$$
the connected component containing $id_M$ of ${\rm
Symp}_0(M,\omega)$ with respect to the $C^\infty$-strong topology.

 For $\mbox{\Bb K}=\mbox{\Bb C},\mbox{\Bb R}$ and $\mbox{\Bb Q}$ let
$H^\ast(M,\mbox{\Bb K})$ (resp. $H^\ast_c(M,\mbox{\Bb K})$) be the
$\mbox{\Bb K}$-coefficient cohomology
 (resp. $\mbox{\Bb K}$-coefficient cohomology
 with compact support).
Given a geometrically bounded symplectic manifold $(M,\omega,
J,\mu)$ of dimension $2n$,  $A\in H_2(M,\mbox{\Bb Z})$ and
integers $g\ge 0, m>0$ with $2g+m\ge 3$, let $\kappa\in
H_\ast(\overline{\mathcal M}_{g,m},\mbox{\Bb Q})$ and
$\{\alpha_i\}_{1\le i\le m}\subset$ $H^\ast(M,\mbox{\Bb Q})\cup
H_c^\ast(M,\mbox{\Bb Q})$ satisfy
$$\left.\begin{array}{ll}
\sum^m_{i=1}\deg\alpha_i+ {\rm codim}(\kappa)=2c_1(M)(A)+
2(3-n)(g-1)+ 2m.\end{array}\right.\eqno(1.5)$$
 If $\{\alpha_i\}_{1\le i\le m}$ has at least one contained in $H_c^\ast(M,\mbox{\Bb
 Q})$ we can define a  number
$${\mathcal G}{\mathcal W}^{(\omega,\mu,
J)}_{A,g,m}(\kappa;\alpha_1,\cdots,\alpha_m)\in\mbox{\Bb
Q}\eqno(1.6)$$
 in (4.7).   When (1.5) is not satisfied we simply define ${\mathcal G}{\mathcal W}^{(\omega,\mu,
J)}_{A,g,m}(\kappa;\alpha_1,\cdots,\alpha_m)$ to be zero. The
following theorem summarizes up the part results in
$\S4$.\vspace{2mm}

\noindent{\bf Theorem 1.1.}\quad{\it The rational number
${\mathcal G}{\mathcal W}^{(\omega,\mu,
J)}_{A,g,m}(\kappa;\alpha_1,\cdots,\alpha_m)$ in (1.6) has
 the following properties:
 \begin{itemize}
 \item[(i)] It is multilinear and supersymmetric on $\alpha_1,\cdots,
 \alpha_m$.
\item[(ii)] It is independent of
choices of $J\in{\mathcal J}(M,\omega,\mu)$.
\item[(iii)] It is invariant under the weak
 deformation of a geometrically bounded symplectic manifold $(M,\omega, J,\mu)$;
 specially, it only depends on the connected component of $\mu$ in
 ${\mathcal G}{\mathcal R}(M)$ with respect to the $C^\infty$ strong
 topology.
\item[(iv)] For any  $\psi\in{\rm Symp}_0^S(M,\omega)$ it holds that
$${\mathcal G}{\mathcal W}^{(\omega,\psi^\ast\mu,
\psi^\ast J)}_{A,g,m}(\kappa;\alpha_1,\cdots,\alpha_m)={\mathcal
GW}^{(\omega,\mu, J)}_{A,g,m}(\kappa;\alpha_1,\cdots,\alpha_m).
$$
\end{itemize}}

It was proved in [Gr1] that  an open manifold $M$ has a symplectic
structure if and only if $M$ has an almost complex structure and
that every symplectic form $\omega$ on $M$ can be smoothly
homotopy to an exact symplectic form through nondegenerate
$2$-forms on $M$. Therefore our Theorem 1.1 implies: \emph{For a
 noncompact geometrically bounded symplectic manifold $(M,\omega, J,\mu)$
  if there exists a Gromov-Witten invariant ${\mathcal GW}^{(\omega,\mu,
J)}_{A,g,m}(\kappa;\alpha_1,\cdots,\alpha_m)\ne 0$ then it cannot
be deformed to an exact geometrically bounded symplectic manifold
$(M,d\alpha, J_1,\mu_1)$ via the weak deformation.}

In $\S5.1$ two localization formulas (Theorem 5.1 and Theorem 5.3)
for GW-invariants are explicitly proposed and proved. Here the
localization means that a (local) virtual moduli cycle constructed
from a smaller moduli space is sometime enough for computation of
some concrete GW-invariants. They are expected to be able simplify
the computation for  GW-invariants. In particular we use them in
the proof of the composition laws in $\S5.2$.

Let  ${\mathcal F}_m: \overline{\mathcal M}_{g,m}\to
\overline{\mathcal M}_{g,m-1}$ be a map forgetting last marked
point. It is a Lefschetz fibration and the integration along the
fibre induces a map $({\mathcal F}_m)_\sharp$ from
$\Omega^\ast(\overline{\mathcal M}_{g,m})$ to
$\Omega^{\ast-2}(\overline{\mathcal M}_{g,m-1})$ (see (5.25)). It
also induces a ``shriek'' map $({\mathcal F}_m)_!$ from
$H_\ast(\overline{\mathcal M}_{g,m-1};\mbox{\Bb Q})$ to
$H_{\ast+2}(\overline{\mathcal M}_{g,m};\mbox{\Bb Q})$ (see
(5.26)). As expected Theorem 5.8 and Theorem 5.9 give desired
reduction formulas.\vspace{2mm}

\noindent{\bf Theorem 1.2.}\quad{\it If $(g,m)\ne (0,3),(1,1)$,
then for any $\kappa\in H_\ast(\overline{\mathcal
M}_{g,m-1};\mbox{\Bb Q})$, $\alpha_1\in H_c^\ast(M;\mbox{\Bb Q})$,
$\alpha_2, \cdots,\alpha_m\in H^\ast(M;\mbox{\Bb Q})$ with
$\deg\alpha_m=2$ it holds that
\begin{eqnarray*}
&&{\mathcal G}{\mathcal W}^{(\omega, \mu, J)}_{A,g,m}(({\mathcal
F}_m)_!(\kappa);\alpha_1,\cdots,\alpha_m)=\alpha_m(A)\cdot
{\mathcal G}{\mathcal W}^{(\omega, \mu,
J)}_{A,g,m-1}(\kappa;\alpha_1,\cdots,\alpha_{m-1}),\\
&&{\mathcal G}{\mathcal W}^{(\omega, \mu,
J)}_{A,g,m}(\kappa;\alpha_1,\cdots,\alpha_{m-1},{\bf 1})={\mathcal
G}{\mathcal W}^{(\omega, \mu,J)}_{A,g,m}(({\mathcal
F}_m)_\ast(\kappa);\alpha_1,\cdots,\alpha_{m-1}).
\end{eqnarray*}
 Here ${\bf 1}\in H^0(M,\mbox{\Bb Q})$ is the identity.}\vspace{2mm}

Having these two formulas, as in [RT2] we follow [W2] to introduce
the gravitational correlators
$\langle\tau_{d_1,\alpha_1}\tau_{d_2,\alpha_2}\cdots\tau_{d_m,\alpha_m}\rangle_{g,
A}$ in (6.5) and define the free energy function $F^M_g(t^a_r;q)$
in (6.7) and Witten's generating function $F^M(t^a_r;q)$ in (6.8).
It is claimed in our Theorem 6.1 that $F^M(t^a_r;q)$ and
$F^M_g(t^a_r;q)$ respectively satisfy  the expected generalized
string equation and the dilation equation under some reasonable
assumptions (6.1). In addition we also construct a family of new
solutions for a few variants of the generalized string equation
and the dilation equation in Theorem 6.2, which are also new even
for the case of closed symplectic manifolds.

In the case that $\dim H^\ast(M)<\infty$, our Gromov-Witten
invariants also satisfy the composition laws of some form. Let
integers $g_i\ge 0$ and $m_i>0$ satisfy: $2g_i+m_i\ge 3$, $i=1,2$.
 Set $g=g_1+g_2$ and $m=m_1+m_2$ and fix a
decomposition $Q=Q_1\cup Q_2$ of $\{1,\cdots,m\}$ with
$|Q_i|=m_i$. Then one gets a canonical embedding $\varphi_Q:
\overline{\mathcal M}_{g_1, m_1+1}\times\overline{\mathcal
M}_{g_2, m_2+1}\to\overline{\mathcal M}_{g, m}$. Let
$\psi:\overline{\mathcal M}_{g-1, m+2}\to \overline{\mathcal
M}_{g,m}$ be the natural embedding obtained by gluing together the
last two marked points. Take a basis  $\{\beta_i\}$  of
$H^\ast(M)$ and a dual basis $\{\omega_i\}$ of them in
$H^\ast_c(M)$,
 i.e.,
$\langle\omega_j,\beta_i\rangle=\int_M\beta_i\wedge\omega_j=\delta_{ij}$.
Let $\eta^{ij}=\int_M\omega_i\wedge\omega_j$ and
  $c_{ij}=(-1)^{\deg\omega_i\cdot\deg\omega_j}\eta^{ij}$.
 Our Theorem 5.5 and Theorem 5.7 give the composition laws of the following forms.\vspace{2mm}

\noindent{\bf Theorem 1.3.}\quad{\it Assume that  $\dim
H^\ast(M)<\infty$. Let $\kappa\in H_\ast(\overline{\mathcal
M}_{g-1, m+2},\mbox{\Bb Q})$,  and $\alpha_i\in H^\ast(M,\mbox{\Bb
Q})$, $i=1,\cdots, m$. If some
 $\alpha_t\in H_c^\ast(M,\mbox{\Bb Q})$ then
 \begin{eqnarray*}
{\mathcal G}{\mathcal W}^{(\omega, \mu,
J)}_{A,g,m}(\psi_\ast(\kappa);\alpha_1,\cdots,\alpha_m) =
\sum_{i,j} c_{ij}\cdot {\mathcal G}{\mathcal W}^{(\omega, \mu,
J)}_{A, g-1, m+2}(\kappa;\alpha_1,\cdots,\alpha_m,\beta_i,
\beta_j).
\end{eqnarray*}
Moreover,  let $\kappa_i\in H_\ast(\overline{\mathcal M}_{g_i,
m_i},\mbox{\Bb Q})$, $i=1,2$, and  $\alpha_s, \alpha_t\in
H_c^\ast(M,\mbox{\Bb Q})$ for some $s\in Q_1$ and $t\in Q_2$.
 Then
\begin{eqnarray*}
&& \!\!\!{\mathcal G}{\mathcal W}^{(\omega, \mu,
J)}_{A,g,m}(\varphi_{Q\ast}(\kappa_1\times\kappa_2);\alpha_1,\cdots\!,\!\alpha_m)
\!=\epsilon(Q)(-1)^{\deg\kappa_2\sum_{i\in Q_1}\!{\deg\alpha_i}}
\!\!\sum_{A=A_1+A_2}\\
&&\sum_{k,l} \eta^{kl}\cdot {\mathcal G}{\mathcal W}^{(\omega,
\mu, J)}_{A_1,g_1,m_1+1}(\kappa_1;\{\alpha_i\}_{i\in
Q_1},\beta_k)\cdot{\mathcal G}{\mathcal W}^{(\omega, \mu,
J)}_{A_2,g_2,m_2+1}(\kappa_2; \beta_l, \{\alpha_i\}_{i\in Q_2}).
\end{eqnarray*}
Here $g=g_1+ g_2$ and $\epsilon(Q)$ is the sign of permutation
$Q=Q_1\cup Q_2$ of $\{1,\cdots,m\}$ with $|Q_1|=m_1$ and
$|Q_2|=m_2$.}\vspace{2mm}

If $(M,\omega)$ is a closed symplectic manifold it is easily
proved that they are reduced to the ordinary ones. However they
are unsatisfactory because we cannot obtain the desired WDVV
equation and quantum products. Indeed, for a basis
$\{\beta_i\}_{1\le i\le L}$ of $H^\ast(M, \mbox{\Bb Q})$ as in
Theorem 1.3 and $w=\sum t_i\beta_i\in H^\ast(M,\mbox{\Bb C})$,  we
can only define $\underline\alpha$-Gromov-Witten potential
\begin{eqnarray*}
\hspace{36mm}\Phi_{(q,\underline\alpha)}(w)=\!\!\sum_{A\in
H_2(M)}\sum_{m\ge
\max(1, 3-k)}\frac{1}{m!}\; &&\\
\hspace{36mm}{\mathcal G}{\mathcal W}^{(\omega,\mu,
J)}_{A,0,k+m}([\overline{\mathcal M}_{0,k+m}]; \underline\alpha,
w, \cdots, w) q^A.&&\hspace{35mm}(1.7)
\end{eqnarray*}
 and prove that they satisfy  WDVV-equation of our form in Theorem 6.3. Here
$\underline\alpha=\{\alpha_i\}_{1\le i\le k}$ be a collection of
nonzero homogeneous elements in
  $H^\ast_c(M,\mbox{\Bb C})\cup H^\ast(M,\mbox{\Bb C})$ and at
  least one of them belongs to $H^\ast_c(M,\mbox{\Bb C})$.
Similarly, we can only define the quantum product on
$QH^\ast(M,\mbox{\Bb Q})$ of the following form
$$\alpha\star_{\underline\alpha}\beta=\sum_{A\in H_2(M)}\sum_{i,j}
{\mathcal G}{\mathcal W}^{(\omega,\mu,
J)}_{A,0,3+k}([\overline{\mathcal M}_{0,3+k}];\underline\alpha,
\alpha, \beta,\beta_i)\eta_{ij}\beta_j q^A. \eqno(1.8)$$
 and in Theorem 6.5 prove that $QH^\ast(M,\mbox{\Bb Q})$ is a
 supercommutative ring without
identity under the quantum product.

Using the method in this paper many results obtained on closed
symplectic manifolds with GW-invariants theory may be generalized
to this class of noncompact symplectic manifolds, i.e., the
results in [Lu1, Lu4 and LiuT2] (adding a condition $\dim
H^\ast(M)<+\infty$ if necessary) and the relative Gromov-Witten
invariants and symplectic sums formula for them in [IP1, IP2 and
LiR].

This paper is organized as follows. Section 2 will discuss some
properties of the moduli spaces of stable maps in noncompact
geometrical symplectic manifolds and then construct local
uniformizers. The construction of the virtual moduli cycle is
completed in Section 3. In Section 4 we shall prove Theorem 1.1
and some simple properties for GW-invariants. Section 5 contains
two localization formulas for GW-invariants and proofs of Theorem
1.2 and Theorem 1.3. In Section 6 we construct a family of
solutions of the generalized string equation, dilation
  equation and WDVV equation, and also define a family of quantum products.
\vspace{2mm}

{\bf Acknowledgements}.\hspace{2mm}I am very grateful to Dr.
Leonardo Macarini for pointing out an error in the first version
of this paper. I would also like to thank Professor Gang Tian
 for some discussions on the smoothness of virtual moduli cycles and
  concrete revision suggestions. The author thank the referee for pointing
  out an unclear point and a feasible revision. This revised version
was finished during author's visit to ICTP Trieste and IHES Paris.
I would like to express deep gratitude to Professors D. T. L\^e
and J. P. Bourguignon for their invitations and for the financial
support and hospitality of both institutions.

\section{The stable maps in noncompact geometrically  bounded symplectic manifolds}\label{2}

In this section we first study almost complex structures on a
geometrically bounded symplectic manifold and then discuss the
pseudo-holomorphic curves and stable maps in this kind of
manifolds. In $\S2.4$ the strong $L^{k,p}$-topology and local
uniformizers will be defined and constructed respectively.

\subsection{Geometrically  bounded symplectic manifolds}\label{2.1}

We first prove the following proposition.\vspace{2mm}

\noindent{\bf Proposition 2.1.}\hspace{2mm}{\it For a Riemannian
metric $\mu$ on $M$, $J\in{\mathcal J}(M,\omega)$ satisfies (1.1)
if and only if the metric $g_J$ is quasi isometric to $\mu$.}

\noindent{\it Proof.}\quad If (1.1) is satisfied, then for any
$X\in TM$,
$$\alpha_0\|X\|_{\mu}\le\|JX\|_{g_J}\le\beta_0\|X\|_{\mu}.$$
Since $\|JX\|_{g_J}=\|X\|_{g_J}$ for every $X\in TM$, $g_J$ is
quasi isometric to $\mu$.

Conversely, assume that $C_1\|X\|_\mu\le\|X\|_{g_J}\le
C_2\|X\|_\mu$ for some positive constants $C_1, C_2$ and all $X\in
TM$. Then for any $X, Y\in TM$ we have
$$|\omega(X,Y)|=|\omega(JX,JY)|=|g_J(JX,Y)|\le\|JX\|_{g_J}\|Y\|_{g_J}\le C_2^2\|X\|_\mu\|Y\|_\mu.$$
That is, (1.1) is satisfied for $\alpha_0=C_1^2$ and
$\beta_0=C_2^2$.\hfill$\Box$\vspace{2mm}

Therefore a symplectic manifold $(M,\omega)$ is geometrically
bounded if and only if ${\mathcal G}{\mathcal
R}(M)\cap\{g_J\,|\,J\in{\mathcal J}(M,\omega)\}\ne\emptyset$.

Next we need to study the  connectivity of ${\mathcal
J}(M,\omega,\mu)$. The first part of the following lemma is due to
S\'evennec (cf. Prop. 1.1.6 in Chapter II of [ALP] or Lemma 1.2.4
in [IS]). \vspace{2mm}

\noindent{\bf Lemma 2.2.}\hspace{2mm}{\it Let $(V, \omega)$ be a
$2n$-dimensional symplectic vector space and ${\mathcal
J}(V,\omega)$ be the space of compatible linear complex structures
on $(V,\omega)$. Fix a $J_0\in{\mathcal J}(V,\omega)$ and
corresponding inner product $g_0=g_{J_0}=\omega\circ({\rm
Id}\times J_0)$. For $W\in{\rm End}(V)$ let $W^t$ and
$\|W\|_{g_0}$ denote the conjugation and norm of $W$ with respect
to $g_0$. Setting
$${\mathcal W}_{J_0}:=\{W\in{\rm End}(V)\,|\,
WJ_0=-J_0W,\,W^t=W,\,\|W\|_{g_0}<1\}$$ then
$${\bf L}:{\mathcal J}(V,\omega)\to {\mathcal W}_{J_0},\;J\mapsto -(J-J_0)(J+ J_0)^{-1}$$
is a diffeomorphism with the inverse map
 $${\bf K}={\bf L}^{-1}:{\mathcal W}_{J_0}\to {\mathcal J}(V,\omega),\;
 W\mapsto J_0({\rm Id}-W)^{-1}({\rm Id}+ W).$$
 Furthermore, if $J={\bf K}(W)\in{\mathcal J}(V, \omega)$ is such
that for some $\alpha_0\in (0,1)$,
$$g_J(v,v)\ge\alpha_0 g_0(v,v)\;\forall v\in V$$
then the smooth path $(J_t={\bf K}(tW))_{t\in [0,1]}$ in
${\mathcal J}(V,\omega)$ also satisfies
 $$g_{J_t}(v,v)\ge \alpha_0 g_0(v,v)\;\forall v\in V\,{\rm and}\;t\in[0,1].\eqno(2.1)$$}

\noindent{\it Proof.}\quad We only need to prove the second part.
Note that $({\rm Id}-W)^{-1}({\rm Id}+ W)$ is symmetric with
respect to $g_0$ and that
$$g_J(v,v)=g_0(v, ({\rm
Id}-W)^{-1}({\rm Id}+ W)v)\ge\alpha_0 g_0(v,v)\;\forall v\in V$$
 if and only if
 $$\min\{\frac{1+\lambda}{1-\lambda}\,|\,
 \lambda\in\sigma(W)\}\ge\alpha_0,$$
 or equivalently, $1>\lambda\ge \frac{\alpha_0-1}{\alpha_0+1}$
 for any $\lambda\in\sigma(W)$. It is easily checked that
  $1>t\lambda\ge t\frac{\alpha_0-1}{\alpha_0+1}\ge\frac{\alpha_0-1}{\alpha_0+1}$
 for any $t\in[0,1]$ and $\lambda\in\sigma(W)$.  (2.1) follows.
\hfill$\Box$\vspace{2mm}

\noindent{\bf Proposition 2.3.}\hspace{2mm}{\it For any symplectic
manifold $(M,\omega)$  the space ${\mathcal J}(M,\omega)$ is
contractible with respect to the $C^\infty$ weak topology.
Moreover, for any geometrically bounded Riemannian metric $\mu$ on
$M$ with ${\mathcal J}(M,\omega,\mu)\ne\emptyset$, if $J_0,
J_1\in{\mathcal J}(M,\omega, \mu)$ satisfy (1.1) then there exists
a smooth path $(J_t)_{t\in [0,1]}$ in ${\mathcal J}(M,\omega)$
connecting $J_0$ to $J_1$ such that
$$\omega(u, J_tu)=g_{J_t}(u, u)\ge\frac{\alpha_0^3}{\beta^2_0}\|u\|^2_\mu
 \;\forall u\in TM\,{\rm and}\,t\in[0,1].$$
 In particular ${\mathcal J}(M,\omega,\mu)$ is connected with respect to the $C^\infty$ weak topology.}\vspace{2mm}

\noindent{\it Proof.}\quad ${\mathcal J}(M,\omega)$ is always
nonempty(cf. the proof of Proposition 2.50(i) in [McSa2]). Fix a
$J_0\in{\mathcal J}(M,\omega)$. Set $g_{J_0}=\omega\circ({\rm
Id}\times J_0)$ and
$${\mathcal W}_{J_0}:=\{W\in{\rm End}(TM)\,|\,
W(x)J_0(x)=-J_0(x)W(x),\,W^t(x)=W(x),\,\|W(x)\|_{g_0}<1\,\forall
x\in M\},$$
 where $W^t(x)$ is the conjugation  of $W(x)$ with respect
to $g_{J_0}$. It follows from Lemma 2.2 that
$${\mathcal J}(M,\omega)\to {\mathcal W}_{J_0},\;J\mapsto -(J-J_0)(J+ J_0)^{-1}$$
is a diffeomorphism with respect to the $C^\infty$  weak topology
on them. Since ${\mathcal W}_{J_0}$ is contractible, so is
${\mathcal J}(M,\omega)$.

Next we  prove the second claim. Assume that ${\mathcal
J}(M,\omega,\mu)$ is nonempty. Fix a $J_0\in{\mathcal
J}(M,\omega,\mu)$. Then for any $J\in{\mathcal J}(M,\omega,\mu)$
there exists a unique $W\in{\mathcal W}_{J_0}$ such that
$J=J_0({\rm Id}-W)^{-1}({\rm Id}+ W)$. Let $J_0$ and $J$ satisfy
(1.1). We easily deduce that
\begin{eqnarray*}
&&\|J_0u\|_\mu\le\frac{\beta_0}{\alpha_0}\|u\|_\mu,\quad\|u\|_\mu\le\frac{\beta_0}{\alpha_0}\|J_0u\|_\mu\\
&&\|Ju\|_\mu\le\frac{\beta_0}{\alpha_0}\|u\|_\mu,\quad\|u\|_\mu\le\frac{\beta_0}{\alpha_0}\|Ju\|_\mu
\end{eqnarray*}
for all $u\in TM$. Using (1.1) again we get that
$$\min(g_{J_0}(u,u), g_J(u,u))\ge\alpha_0\|u\|^2_\mu\quad{\rm
and}\quad\max(g_{J_0}(u,u), g_J(u,
u))\le\frac{\beta_0^2}{\alpha_0}\|u\|^2_\mu\;\forall u\in TM.$$
 It follows from these that
 $$g_J(u,u)\ge \frac{\alpha_0^2}{\beta_0^2}g_{J_0}(u,u)\;\forall u\in TM.$$
Note that we can assume that $\alpha_0<\beta_0$ and therefore
$0<\frac{\alpha_0^2}{\beta_0^2}<1$. By Lemma 2.2, a smooth path
$J_t=J_0({\rm Id}-tW)^{-1}({\rm Id}+ tW)$, $t\in [0,1]$ in
${\mathcal J}(M,\omega)$ connecting $J_0$ and $J_1=J$, satisfies
$$\omega(u, J_tu)=g_{J_t}(u,u)\ge \frac{\alpha_0^2}{\beta^2_0} g_0(u,u)\ge\frac{\alpha_0^3}{\beta^2_0}\|u\|^2_\mu
 \;\forall u\in TM\;{\rm and}\;t\in[0,1].$$
That is, $J_t\in{\mathcal J}(M,\omega,\mu)$ for all $t\in [0,1]$.
Proposition 2.3 is proved. \hfill$\Box$\vspace{2mm}

Remark that the same reasoning may give the similar results for
${\mathcal J}_\tau(M,\omega)$ and ${\mathcal J}_\tau(M,\omega,
J)$.

\subsection{Pseudo-holomorphic curves in geometrically
bounded symplectic manifolds}\label{2.3}

Let $(M,\omega, J, \mu)$ be a geometrically bounded symplectic
manifold  satisfying (1.1).
 Since $(M, J,\mu)$ is a tame almost complex manifold we have the following
version of Monotonicity principle(cf. [Sik,
Prop.4.3.1(ii)]).\vspace{2mm}

\noindent{\bf Lemma 2.4.}\hspace{2mm}{\it Let $C_0>0$ be an upper
bound of all sectional curvatures $K_\mu$ of $\mu$ and
$r_0=\min\{i(M,\mu), \pi/\sqrt{C_0}\}$. Assume that $h:S\to M$ is
a $J$-holomorphic map from a compact connected Riemannian surface
$S$ with boundary to $M$. Let $p\in M$ and $0<r\le r_0$ be such
that
$$p\in h(S)\subset B_\mu(p, r)\quad{\rm and}\quad h(\partial
S)\subset B_\mu(p, r),$$ then
 $${\rm Area}_\mu(h(S))\ge\frac{\pi\alpha_0}{4\beta_0}r^2.\eqno(2.2)$$}

In fact, by the comments below Definition 4.1.1 of [Sik] we may
choose $C_1=1/\pi$ and $C_2=\beta_0/\alpha_0$ in the present case.
Moreover, $C_6$ in Proposition 4.4.1 of [Sik] is equal to
$\frac{4\beta_0}{\pi\alpha_0^2 r_0}$.\vspace{2mm}

\noindent{\bf Proposition 2.5}([Sik, Prop.4.4.1]).\hspace{2mm}{\it
Let $(M,\omega,\mu, J)$ be as above, and $h:S\to M$ be a
$J$-holomorphic map from a connected compact Riemannian surface
$S$ to $M$. Assume that a compact subset $K\subset M$ is such that
$h(S)\cap K\ne\emptyset$ and $h(\partial S)\subset K$. Then $f(S)$
is contained in
$$U_\mu(K, \frac{4\beta_0}{\pi\alpha_0^2 r_0}{\rm Area}_\mu(h(S)))=\{p\in M\,|\,
d_\mu(p, K)\le\frac{4\beta_0}{\pi\alpha_0^2 r_0}{\rm
Area}_\mu(h(S))\}.$$ In particular, if $S$ is a connected closed
Riemannian surface then
$${\rm diam}_\mu(h(S))\le\frac{4\beta_0}{\pi\alpha_0^2 r_0}
\int_{S}h^\ast\omega.\eqno(2.3)$$}

Lemma 2.4 also leads to\vspace{2mm}

\noindent{\bf Proposition 2.6.}\hspace{2mm}{\it Under the
assumptions of Lemma 2.4, for any nonconstant $J$-holomorphic map
 $f:\Sigma\to M$ from a closed connected Riemannian surface $\Sigma$ to $M$
 it holds that
$$\int_\Sigma f^\ast\omega\ge\frac{\pi\alpha_0}{8\beta_0}r^2_0.\eqno(2.4)$$}

\noindent{\it Proof.}\quad Assume that there is a nonconstant
$J$-holomorphic map  $f:\Sigma\to M$ from a closed connected
Riemannian surface $\Sigma$ to $M$ such that
$$\int_\Sigma f^\ast\omega<\frac{\pi\alpha_0}{8\beta_0}r^2_0.$$
Take $z_0\in\Sigma$ and let  $S$ be the connected component of
$f^{-1}(B_\mu(p, 3r_0/4))$ containing $p=f(z_0)$.  After
perturbing $3r_0/4$ a bit we may assume that $S$ is a Riemannian
surface with smooth boundary. If $\partial S\ne\emptyset$ Lemma
2.4 gives
$$\int_S f^\ast\omega\ge\frac{\pi\alpha_0}{4\beta_0}(3r_0/4)^2>\frac{\pi\alpha_0}{8\beta_0}r^2_0.$$
This contradicts the above assumption. Hence $S=\Sigma$ and
$f(\Sigma)\subset B_\mu(p, 3r_0/4)$. But $3r_0/4<i(M,\mu)$ and
thus $f$ is homotopic to the constant map. It follows that $f$ is
constant because $f$ is $J$-holomorphic. \hfill$\Box$

\subsection{Stable maps in geometrically
bounded symplectic manifolds}\label{2.4}

Following Mumford, by a {\it semistable curve with $m$ marked
points}, one means a pair $(\Sigma;\bar{\bf z})$ of a connected
Hausdorff topological space $\Sigma$ and $m$ different points
$\bar{\bf z}=\{z_1,\cdots,z_m\}$ on it such that there exist a
finite family of smooth Riemann surfaces
$\{\Sigma_s:s\in\Lambda\}$ and continuous maps
$\pi_{\Sigma_s}:\widetilde\Sigma_s\to\Sigma$ satisfying the
following properties: $(i)$ each  $\pi_{\Sigma_s}$ is a local
homeomorphism; $(ii)$ for each $p\in\Sigma$ it holds that
$1\le\sum_s \sharp\pi_{\Sigma_s}^{-1}(p)\le 2$, and all points
which satisfy $\sum_s \sharp\pi_{\Sigma_s}^{-1}(p)=2$ are
isolated; $(iii)$ for each $z_i$, $\sum_s
\sharp\pi_{\Sigma_s}^{-1}(z_i)=1$.

The points in $\Sigma_{\rm sing}:=\{p\in\Sigma: \sum_s
\sharp\pi_{\Sigma_s}^{-1}(p)=2\}$ are called singular points of
$\Sigma$. Each singular point $p$ such that
$\sharp\pi_{\Sigma_s}^{-1}(p)=2$ is called the self-intersecting
point of $\Sigma$. $\Sigma_s:=\pi_{\Sigma_s}(\widetilde\Sigma_s)$
is called the $s^{th}$ components of $\Sigma$, and
$\widetilde\Sigma_s$ is called the smooth resolution of
$\Sigma_s$. Each $z_i$ is called the marked point. The points in
$\pi_{\Sigma_s}^{-1}(\Sigma_{\rm sing})$ (resp.
$\pi_{\Sigma_s}^{-1}(\bar{\bf z})$) are called the singular points
(resp. the marked points) on $\widetilde\Sigma_s$, respectively.
Let $m_s$ be the number of all singular and marked points on
$\widetilde\Sigma_s$ and $g_s$ be the genus of
$\widetilde\Sigma_s$. The {\it genus} $g$ of $(\Sigma;\bar{\bf
z})$ is defined by
$$\left.\begin{array}{ll}
1 + \sum_s g_s + \sharp{\rm Inter}(\Sigma)-\sharp{\rm
Comp}(\Sigma),\end{array}\right.$$
 where $\sharp{\rm Inter}(\Sigma)$ and $\sharp{\rm Comp}(\Sigma)$ stand for the
number of the intersecting points on $\Sigma$ and the number of
the components of $\Sigma$ respectively.
 When $m_s+2g_s\ge 3$ we call the component
$(\widetilde\Sigma_s;\bar{\bf z}_s)$ stable. If all components of
$(\Sigma;\bar{\bf z})$ are stable, $(\Sigma;\bar{\bf z})$ is
called {\it a stable curve of genus $g$ and with $m$ marked
points.}

Two such  genus $g$ semi-stable curves $(\Sigma; z_1,\cdots, z_m)$
and $(\Sigma^\prime; z_1^\prime,\cdots, z_m^\prime)$ are said to
be isomorphic if there is a homeomorphism
$\phi:\Sigma\to\Sigma^\prime$ such that (i) $\phi(z_i)=z_i^\prime,
i=1,\cdots, m$, and (ii) the restriction of $\phi$ to each
component $\Sigma_s$ of it can be lifted to a biholomorphic
isomorphism
$\phi_{st}:\widetilde\Sigma_s\to\widetilde\Sigma_t^\prime$. Denote
by $[\Sigma,\bar{\bf z}]$  the isomorphism class of
$(\Sigma,\bar{\bf z})$ and by $Aut(\Sigma,\bar{\bf z})$ the group
of all automorphisms of $(\Sigma,\bar{\bf z})$. Then
$(\Sigma,\bar{\bf z})$ is stable if and only if
$Aut(\Sigma,\bar{\bf z})$ is a finite group. Let
$\overline{\mathcal M}_{g,m}$ be the set of all isomorphism
classes of stable curves with $m$ marked points and of genus $g$.
It is called the Deligne-Mumford compactifcation of the moduli
space ${\mathcal M}_{g,m}$ of all isomorphism classes of smooth
stable curves with $m$ marked points and of genus $g$.
$\overline{\mathcal M}_{g,m}$ is not only a projective variety but
also complex orbifold of complex dimension $3g-3+m$.

 For the above genus $g$ semi-stable curve $(\Sigma;\bar{\bf
z})$ a continuous map $f:\Sigma\to M$ is called $C^l$-smooth($l\ge
1$) if each $f\circ\pi_{\Sigma_s}$ is so. The homology class of
$f$ is defined by
$f_*([\Sigma])=\sum_s(f\circ\pi_{\Sigma_s})_*[\Sigma]$. Similarly,
for an almost complex structure $J$ on $M$,
 a continuous map $f:\Sigma\to M$ is called $J$-holomorphic
  if each $f\circ\pi_{\Sigma_s}$ is so.
Moreover, for fixed $k\in\mbox{\Bb N}$ and $p\ge 1$ satisfying
$k-\frac{2}{p}\ge 0$ we say a continuous map $f:\Sigma\to M$ to be
$L^{k,p}$-map if each $f\circ\pi_{\Sigma_s}$ is so. \vspace{2mm}

\noindent{\bf Definition 2.7}([KM]).\hspace{2mm}{\it Given a genus
$g$ semi-stable curve $(\Sigma;\bar{\bf z})$ with $m$ marked
points and $J$-holomorphic map $f$ from $\Sigma$ to $M$ a triple
$(f;\Sigma,\bar{\bf z})$ is called a {\bf $m$-pointed  stable
$J$-map} of genus $g$ in $M$ if for each component $\Sigma_s$ of
$\Sigma$ the composition
$f\circ\pi_{\Sigma_s}:\widetilde\Sigma_s\to M$ cannot be constant
map in the case $m_s+2g_s<3$.}\vspace{2mm}

Two such stable maps $(f;\Sigma,\bar{\bf z})$ and
$(f^\prime;\Sigma^\prime,\bar{\bf z}^\prime)$ are called {\it
equivalence} if there is an isomorphism $\phi:(\Sigma,\bar{\bf
z})\to (\Sigma^\prime,\bar{\bf z}^\prime)$ such that
$f^\prime\circ\phi=f$. Denote $[f;\Sigma,\bar{\bf z}]$ by the
equivalence class of $(f;\Sigma,\bar{\bf z})$. The automorphism
group of $(f;\Sigma,\bar{\bf z})$ is defined by
$$Aut(f;\Sigma,\bar{\bf z})=\{\phi\in Aut(\Sigma,\bar z)\,|\, f\circ\phi=f\}.$$
Then $Aut(f;\Sigma,\bar{\bf z})$ is a finite group. For a given
class $A\in H_2(M,\mbox{\Bb Z})$ let us denote
$$\overline{\mathcal M}_{g,m}(M, J, A)$$
by the set of equivalence classes of all  $m$-pointed
 stable $J$-maps of genus $g$ and of class $A$ in  $M$.
There is a natural stratification of $\overline{\mathcal
M}_{g,m}(M, J, A)$, and the number of all stratifications is also
finite provided that $(M,\omega)$ is compact. If $(M,\omega)$ is
not compact one cannot guarantee that the number of the natural
stratifications is finite. In order to understand it better note
that  the following result can follow from (2.3)
directly.\vspace{2mm}

\noindent{\bf Proposition 2.8.}\hspace{2mm}{\it For any
$[f;\Sigma,\bar{\bf z}]\in\overline{\mathcal M}_{g,m}(M, J, A)$ it
holds that
$${\rm diam}_\mu(f(\Sigma))\le\frac{4\beta_0}{\pi\alpha_0^2 r_0}\omega(A).$$}

Fix $k\in\mbox{\Bb N}$ and $p\ge 1$ satisfying $k-\frac{2}{p}\ge
2$. Following [LiuT1, LiuT2, LiuT3] we introduce\vspace{2mm}

\noindent{\bf Definition 2.9.}\hspace{2mm}{\it Given a genus $g$
semi-stable curve $(\Sigma;\bar{\bf z})$ with $m$ marked points
and a $L^{k,p}$-map $f$ from $\Sigma$ to $M$, a triple
$(f;\Sigma,\bar{\bf z})$ is called a {\bf $m$-pointed  stable
$L^{k,p}$-map} of genus $g$ in $M$ if for each component
$\Sigma_s$ of $\Sigma$ satisfying
$(f\circ\pi_{\Sigma_s})_\ast([\widetilde\Sigma_s])=0\in
H_2(M,\mbox{\Bb Z})$ it holds that $m_s+2g_s\ge 3$.}\vspace{2mm}

Notice that the stable $L^{k,p}$-maps are $C^2$-smooth under the
above assumptions. As above we can define the equivalence class
$[f;\Sigma,\bar{\bf z}]$ and the automorphism group
$Aut(f;\Sigma,\bar{\bf z})$  of such a  stable $L^{k,p}$-map
$(f;\Sigma,\bar{\bf z})$. Later, without occurrences of confusions
we abbreviate $(f;\Sigma,\bar{\bf z})$ to ${\bf f}$, and
$Aut(f;\Sigma,\bar{\bf z})$ to $Aut({\bf f})$. Consider the
disjoint union $\widetilde\Sigma:=\cup_s\widetilde\Sigma_s$, which
is called the normalization of $\Sigma$.  Denote by $\tilde f:
=\cup_s f\circ\pi_{\Sigma_s}:\widetilde\Sigma\to M$. We define the
energy of such a stable $L^{k,p}$-map ${\bf f}=(f;\Sigma,\bar{\bf
z})$ by
$$E({\bf f})=\frac{1}{2}\int_{\Sigma}|df|^2=\frac{1}{2}\int_{\widetilde\Sigma}|d\tilde f|^2=
\frac{1}{2}\sum_s\int_{\widetilde\Sigma_s}|d(f\circ\pi_{\Sigma_s})|^2.$$
Here the integral uses a metric $\tau$ in the conformal class
determined by $j_{\widetilde\Sigma}$ and $J$-Hermitian metric
$g_J=\mu$, it is conformally invariant and thus depends only on
$j_{\widetilde\Sigma}$. It is easily proved that $E({\bf f})$ is
independent of choice of representatives in $[{\bf f}]$, so we can
define $E([{\bf f}])=E({\bf f})$. If ${\bf f}$ is $J$-holomorphic
then $E({\bf
f})=\int_{\Sigma}f^\ast\omega=\omega(f_\ast([\Sigma]))$.

Let ${\mathcal B}^M_{A, g,m}$ be the set of equivalence classes of
all $m$-pointed  stable $L^{k,p}$-maps of genus $g$ and of class
$A$ in $M$. Then it is clear that
$$\overline{\mathcal M}_{g,m}(M, J, A)\subset{\mathcal B}^M_{A, g,m}.$$
For each $[{\bf f}]\in{\mathcal B}^M_{A, g,m}$ we set
$${\mathcal E}^M_{A,g,m}=\bigcup_{[{\bf f}]\in{\mathcal B}^M_{A,
g,m}}{\mathcal E}_{[{\bf f}]}\;\;{\rm and}\;\; {\mathcal E}_{[{\bf
f}]}=\Bigl(\bigcup_{{\bf f}\in[{\bf
f}]}L^p_{k-1}(\wedge^{0,1}(f^\ast TM))\Bigr)/\sim.$$
 Here the equivalence relation $\sim$ is defined by the pull-back
 of the sections induced from the identification of the domains.
Following [T2] we associate to each stable map ${\bf
f}=(f;\Sigma,\bar{\bf z})$  a dual graph $\Gamma_{\bf f}$ as
follows. Each irreducible component $\Sigma_s$ of $\Sigma$
corresponds to a vertex $v_s$ in $\Gamma_{\bf f}$ with weight
$(g_s, A_s)$, where $g_s$ is the geometric genus of $\Sigma_s$ and
$A_s=(f\circ\pi_{\Sigma_s})([\widetilde\Sigma_s])$. For each
marked point $z_i$ of $\Sigma_s$ we attach a leg $le_i$ to $v_s$.
For each intersection point of distinct components $\Sigma_s$ and
$\Sigma_t$ we attach an edge $e_{st}$ joining $v_s$ and $v_t$. For
each self-intersection point of $\Sigma_s$ we attach a loop $lo_s$
to $v_s$. $\Gamma_{\bf f}$ is independent of the representatives
in $[{\bf f}]$. We define $\Gamma_{[{\bf f}]}=\Gamma_{\bf f}$. The
genus $g(\Gamma_{[{\bf f}]})$ of $\Gamma_{[{\bf f}]}$ is defined
by the sum of $\sum_s g_s$ and the number of holes in
$\Gamma_{[{\bf f}]}$. The homology class
$A=f_\ast([\Sigma])=\sum_s A_s$ is defined as that of
$\Gamma_{[{\bf f}]}$.

A genus $g$ graph $\Gamma$ with $m$ legs and of homology class $A$
is called {\it effective} if there is a $m$-pointed
$J$-holomorphic stable map ${\bf f}$ of genus $g$ and of class $A$
in $M$ such that its dual graph $\Gamma_{\bf f}=\Gamma$. For $e>0$
let ${\mathcal D}^e_{g,m}$ be the set of those effective genus $g$
graphs $\Gamma$ with $m$ legs and of homology class $A$ satisfying
$\omega(A)\le e$.

Following [LiuT1-LiuT3] $\overline{\mathcal M}_{g,m}(M, J, A)$ is
equipped with the (weak) topology for which it is also Hausdorff.
Notice that Proposition 2.8 implies that for any compact subset
$K\subset M$ the images of all maps in
$$\overline{\mathcal M}_{g,m}(M, J, A; K):=\{[{\bf f}]\in\overline{\mathcal
M}_{g,m}(M, J, A)\,|\, f(\Sigma)\cap K\ne\emptyset\}$$ is
contained in a compact subset of $M$, $\{x\in M\,|\, d_\mu(x,
K)\le 4\omega(A)/(\pi r_0)\}$. We can follow [LiuT1] to
prove\vspace{2mm}

\noindent{\bf Proposition 2.10.}\hspace{2mm}{\it For any compact
subset $K\subset M$, $\overline{\mathcal M}_{g,m}(M, J, A; K)$ is
compact with respect to the weak topology.}\vspace{2mm}

 As in [RT1, \S4 FO, Prop.8.7] it follows from
this result that the following result corresponding to Lemma 3.3
in [Lu1] holds.\vspace{2mm}

\noindent{\bf Proposition 2.11.}\hspace{2mm}{\it For any compact
subset $K\subset M$,
$${\mathcal D}^e_{g,m}(K):=\{\Gamma\in {\mathcal
D}^e_{g,m}\,|\,\exists\; [{\bf f}]\in \overline{\mathcal
M}_{g,m}(M, J, A; K)\;s. t.\; \Gamma_{[{\bf f}]}=\Gamma\}$$ is
finite.}

\subsection{Strong $L^{k,p}$-topology and local uniformizer}\label{2.5}

In this subsection we follow the ideas of [LiuT1-LiuT3] to
construct the strong $L^{k,p}$-topology and local uniformizers on
${\mathcal B}^M_{A,g,m}$. We shall give the main steps even if the
arguments are the same as there. The reason is that we need to
know the size of these local uniformizers when we imitate the
methods in [LiuT1, LiuT2, LiuT3] to construct the desired virtual
moduli cycles in our case.

\subsubsection{Local deformation of a stable curve}
 We review the local deformation of a
stable curve $\sigma\in\overline{\mathcal M}_{g,m}$ in [FO].
Choose a representative $(\Sigma,\bar{\bf z})$ of $\sigma$ and let
$\Sigma=\cup_s\Sigma_s= \cup_s\pi_{\Sigma_s}(\widetilde\Sigma_s)$
be the decomposition of $\Sigma$ to the irreducible components as
in \S2.3. Firstly we describe the deformation of $\sigma$ in the
same stratum. Let $V_{deform}(\sigma)$ be a neighborhood of $0$ in
the product $\prod_s\mbox{\Bb C}^{3g_s-3+m_s}$. The universal
family induces a fiberwise complex structure on the fiber bundle
$V_{deform}(\sigma)\times\Sigma\to V_{deform}(\sigma)$. One may
take the following representative. By unique continuation we can
consider the direct product $V_{deform}(\sigma)\times\Sigma$ and
change the complex structure in a compact set
$K_{deform}(\sigma)\subset\Sigma\setminus(Sing(\Sigma)\cup\bar{\bf
z})$ so that it gives a universal family. One can also take a
family of K\"ahler metric which is constant in $\Sigma\setminus
K_{deform}(\sigma)$. It can also be assumed that
$V_{deform}(\sigma)\times\Sigma$ together with fiberwise complex
structure and K\"ahler metric is equivariant by the diagonal
action of $Aut(\Sigma,\bar{\bf z})$. These show that for every
$u\in V_{deform}$ we have a complex structure $j_u$ on a K\"ahler
metric $\tau_u$ satisfying:
\begin{description}
\item[(i)] $j_u=j_0=j_\Sigma$ outside
$K_{deform}(\sigma)\subset\Sigma\setminus(Sing(\Sigma)\cup\bar{\bf
z})$;
\item[(ii)] $\tau_u=\tau_0$ and is flat outside
$K_{deform}(\sigma)\subset\Sigma\setminus(Sing(\Sigma)\cup\bar{\bf
z})$; \item[(iii)] $j_u$ and $\tau_u$ depend on $u$ smoothly;
\item[(iv)] For $\gamma\in Aut(\Sigma,\bar{\bf z})$ and $u\in
V_{deform}(\sigma)$ it holds that $j_{\gamma\cdot u}=\gamma^\ast
j_u$ and $\tau_{\gamma\cdot u}=\gamma^\ast \tau_u$.
\end{description}
 Let us denote $\Sigma^u$ by $\Sigma$ equipped with
the complex structure $j_u$ and K\"ahler metric $\tau_u$. Then
${\mathcal U}(\sigma):=\{\sigma_u=[\Sigma^u,\bar{\bf z}]\,|\, u\in
V_{deform}(\sigma)\}$ gives a neighborhood of $\sigma$ in the same
stratum.

Next we construct a parametric representation of a neighborhood of
$\sigma$ in $\overline{\mathcal M}_{g,m}$ by gluing. Since
$(\Sigma^u,\bar{\bf z})$ is equal to $(\Sigma,\bar{\bf z})$
outside $K_{deform}(\sigma)$ for every $u\in V_{deform}(\sigma)$
and the gluing only occurs near $Sing(\Sigma)$ we only need to
consider $\Sigma$. Let $V_{resolve}(\sigma)$ be a small
neighborhood of the origin in $\prod_{z\in
Sing(\Sigma)}T_{z_s}\widetilde\Sigma_s\otimes
T_{z_t}\widetilde\Sigma$. Here
$\pi_{\Sigma_s}(z_s)=\pi_{\Sigma_t}(z_t)$ (and $s=t$ is allowed).
Near $Sing(\Sigma)$ the K\"ahler metric $\tau=\tau_0$ induces the
Hermitian metrics on $T_{z_s}\widetilde\Sigma_s$ and
$T_{z_t}\widetilde\Sigma_t$ respectively. They in turn induce one
on the tensor product $T_{z_s}\widetilde\Sigma_s\otimes
T_{z_t}\widetilde\Sigma_t$. For a vector $v=(v_z)\in
V_{resolve}(\sigma)$, if $v_z=0$ nothing is made. If $v_z\ne 0$ we
have a biholomorphism map $\Phi_{v_z}:
T_{z_s}\widetilde\Sigma_s\setminus\{0\}\to
T_{z_t}\widetilde\Sigma_t\setminus\{0\}$ such that
$w\otimes\Phi_{v_z}(w)=v_z$. Setting $|v_z|=R^{-2}$ then for
$|v_z|$ sufficiently small (and thus $R$ sufficiently large) the
$\tau$-metric $2$-disc $D_{z_s}(3R^{-1/2})$ with center $z_s$ and
radius $3R^{-1/2}$ in $\widetilde\Sigma_s$ and
$D_{z_t}(3R^{-1/2})$ with center $z_t$ and radius $3R^{-1/2}$ in
$\widetilde\Sigma_t$ are contained outside
$K_{deform}(\sigma)\cup\bar{\bf z}$. {\bf Later we assume}
$V_{resolve}$ so small that these conclusions hold for all $v\in
V_{resolve}$. Using the exponential maps
$\exp_{z_s}:T_{z_s}\widetilde\Sigma_s\to\widetilde\Sigma_s$ and
$\exp_{z_t}:T_{z_t}\widetilde\Sigma_t\to\widetilde\Sigma_t$ with
respect to the metric $\tau$ we have a biholomorphism
$$\exp_{z_s}^{-1}\circ\Phi_{v_z}\circ\exp_{z_t}^{-1}:
D_{z_s}(R^{-1/2})\setminus D_{z_s}(R^{-3/2})\to
D_{z_t}(R^{-1/2})\setminus D_{z_t}(R^{-3/2}).$$
 Clearly, this biholomorphism maps the internal (resp. outer)
 boundary of the annulus $D_{z_s}(R^{-1/2})\setminus
 D_{z_s}(R^{-3/2})$ to the outer (resp. internal) boundary of the annulus
$D_{z_t}(R^{-1/2})\setminus D_{z_t}(R^{-3/2})$.
 Using it we can glue $\Sigma_s$ and $\Sigma_t$.
 Hereafter we identify $D_{z_s}(R^{-1/2})\setminus
 D_{z_s}(R^{-3/2})$ (resp. $D_{z_t}(R^{-1/2})\setminus
 D_{z_t}(R^{-3/2})$) with $\pi_{\Sigma_s}(D_{z_s}(R^{-1/2})\setminus
 D_{z_s}(R^{-3/2}))$(resp. $\pi_{\Sigma_t}(D_{z_t}(R^{-1/2})\setminus
 D_{z_t}(R^{-3/2}))$).
After performing this construction for each nonzero component
$v_z$ we obtain a $2$-dimensional manifold with possible some
points in $Sing(\Sigma)$ as singular points. Then we  define a
K\"ahler metric (and thus a complex structure) on this new
``manifold". Notice that the construction in [FO] shows that this
new metric is only changed on the glue part
$$\cup_{v_z\ne 0}D_{z_s}(R^{-1/2})\setminus D_{z_s}(R^{-3/2})\equiv
\cup_{v_z\ne 0}D_{z_t}(R^{-1/2})\setminus D_{z_t}(R^{-3/2}).$$
Therefore, if we choose the metric $\tau_u$ on the other part of
this new ``manifold" then a new stable curve
$(\Sigma_{(u,v)},\bar{\bf z})$ associated with $(u,v)\in
V_{deform}(\sigma)\times V_{reslove}(\sigma)$ is obtained.
Moreover the corresponding complex structure $j_{(u,v)}$ and
K\"ahler metric $\tau_{(u,v)}$ only change on
$K_{deform}(\sigma)\cup K_{neck}(\sigma, v)$ and depend on $(u,v)$
smoothly, where
$$K_{neck}(\sigma, v):=\bigcup_{z\in
Sing(\Sigma), v_z\ne 0}D_{z_s}(R^{-1}_z)\cup
D_{z_t}(R^{-1/2}_z),$$
and $R_z=|v_z|^{-2}$. We get a parametric
representation of a neighborhood of $\sigma$ in
$\overline{\mathcal M}_{g,m}$ by
$$V_{deform}(\sigma)\times V_{resolve}(\sigma)\to\overline{\mathcal
M}_{g,m},\; (u,v)\mapsto [\Sigma_{(u,v)},\bar{\bf z}].\eqno(2.5)$$
Notice that $(\Sigma_{(u,0)},\bar{\bf z})=(\Sigma^u,\bar{\bf z})$
for every $u\in V_{deform}(\sigma)$, and that
\begin{eqnarray*}
\hspace{18mm}&&
K_{deform}\cup\bar{\bf z}\subset\\
&&  \Sigma_{(u,v)}\setminus (\bigcup_{z\in Sing(\Sigma), v_z\ne
0}D_{z_s}(3R^{-1}_z)\cup D_{z_t}(3R^{-1/2}_z)\bigcup
Sing(\Sigma_{(u,v)}))\hspace{13mm}(2.6)\end{eqnarray*}
 for all $(u,v)\in
V_{deform}\times V_{resolve}$. It is also clear that
$Sing(\Sigma_{(u,v)})\subset Sing(\Sigma)$ for all such $(u,v)$.

\subsubsection{Local deformation of a stable $J$-map}

Let $[{\bf f}]\in\overline{\mathcal M}_{g,m}(M, J, A)$ and ${\bf
f}=(f;\Sigma,\bar{\bf z})$ be a representative of it. Then
$(\Sigma,\bar{\bf z})$ might not be stable. For each unstable
components $\widetilde\Sigma_s$ we add one or two marked point(s)
to $\Sigma_s$ according to whether $\widetilde\Sigma_s$ contains
two or one distinguished point(s). Since ${\bf f}$ is stable
$f_\ast([\Sigma_s])\ne 0$ for each unstable component
$\widetilde\Sigma_s$. So these added points can be required to be
the smooth points of $\Sigma$, where the differential $df$ is
injective.
 Let $\bar{\bf
y}=\{y_1,\cdots, y_l\}$ be the set of all points added to
$\Sigma$. Then $\sigma:=(\Sigma,\bar{\bf z}\cup\bar{\bf y})$ is a
genus $g$ stable curve with $m+l$ marked points. By (2.5) we have
a local deformation of $\sigma$:
$$V_{deform}(\sigma)\times V_{resolve}(\sigma)\to\overline{\mathcal
M}_{g,m+l},\; (u,v)\mapsto [\Sigma_{(u,v)},\bar{\bf z}\cup\bar{\bf
y} ].\eqno(2.7)$$ Let us define maps $f_{(u,v)}:\Sigma_{(u,v)}\to
M$ as follows.\vspace{2mm}

\noindent{\it Step 1}. For each $u\in V_{deform}(\sigma)$ the
above construction shows $\Sigma_{(u,0)}=\Sigma$ as
$2$-dimensional ``manifolds". We define $f_{(u,0)}=f$. However,
$f_{(u,0)}:\Sigma_{(u,0)}\to M$ might not be $(j_{(u,0)},
J)$-holomorphic.\vspace{2mm}

\noindent{\it Step 2}. For a vector $v=(v_z)\in
V_{resolve}(\sigma)$, if $v_z\ne 0$ for some $z\in Sing(\Sigma)$,
we  assume $|v|$ so small that
$\pi_{\Sigma_s}(D_{z_s}(3R_z^{-1/2}))\cup\pi_{\Sigma_t}(3D_{z_t}(R_z^{-1/2}))$
is contained in $\Sigma\setminus K_{deform}(\sigma)$ where
$\tau_{(u,v)}=\tau$ is flat. Here $R_z=|v_z|^{-1/2}$. Let
\begin{eqnarray*}
&&f\circ\pi_{\Sigma_s}(x)=\exp_{f(z)}(\xi_s(x)),\;{\rm if}\; x\in
D_{z_s}(3R_z^{-1/2}),\\
&&f\circ\pi_{\Sigma_t}(x)=\exp_{f(z)}(\xi_t(x)),\;{\rm if}\; x\in
D_{z_t}(3R_z^{-1/2}),
\end{eqnarray*}
where $\xi_s(x),\xi_t(x)\in T_{f(z)}M$.
 Take a
smooth cut function $\chi:\mbox{\Bb R}\to [0,1]$ such that
\begin{eqnarray*}
\chi(r)=\left\{\begin{array}{ll}0 \;&{\rm as}\; r\le 1\\
1\;&{\rm as}\;r\ge 4\end{array}\right.,\;0\le\chi^\prime(r)\le 1.
\end{eqnarray*}
Denote by $\chi_{v_z}(r)=\chi(R_z r)$ for all $r\ge
0$.
\begin{description}
\item[$\bullet$] If
$x\in\pi_{\Sigma_s}(D_{z_s}(2R_z^{-1/2})\setminus
D_{z_s}(R_z^{-1/2}))=\pi_{\Sigma_s}(D_{z_s}(2R_z^{-1/2}))\setminus
\pi_{\Sigma_s}(D_{z_s}(R_z^{-1/2}))$ we define
$$f_{(u,v)}(x)=\exp_{f(z)}[\chi_{v_z}(|\exp_{z_s}^{-1}(\pi^{-1}_{\Sigma_s}(x))|^2\xi_s(x)].$$
\item[$\bullet$] If
$x\in\pi_{\Sigma_t}(D_{z_t}(2R_z^{-1/2})\setminus
D_{z_t}(R_z^{-1/2}))=\pi_{\Sigma_t}(D_{z_t}(2R_z^{-1/2}))\setminus
\pi_{\Sigma_t}(D_{z_t}(R_z^{-1/2}))$ we define
$$f_{(u,v)}(x)=\exp_{f(z)}[\chi_{v_z}(|\exp_{z_t}^{-1}(\pi^{-1}_{\Sigma_t}(x))|^2\xi_t(x)].$$
\item[$\bullet$] If $x$ belongs to the glue part
$$\pi_{\Sigma_s}(D_{z_s}(R^{-1/2})\setminus D_{z_s}(R^{-3/2}))\equiv
\pi_{\Sigma_t}(D_{z_t}(R^{-1/2})\setminus D_{z_t}(R^{-3/2})),$$ we
define $f_{(u,v)}(x)=f(z)$.
\end{description}
Near other singular points we make similar definitions. If
$x\in\Sigma_{(u,v)}$ is not in
$$\cup_{v_z\ne
0}(\pi_{\Sigma_s}(D_{z_s}(2R_z^{-1/2}))\cup\pi_{\Sigma_t}(D_{z_t}(2R_z^{-1/2})))$$
we define $f_{(u,v)}(x)=f(x)$. Clearly,
$f_{(u,v)}:\Sigma_{(u,v)}\to M$ is smooth in the sense of \S2.4.
It should be noted that the above definition of $f_{(u,v)}$
implies
\begin{eqnarray*}
d_H(f_{(u,v)}(\Sigma_{(u,v)}),
f(\Sigma)):=&&\sup_{x\in\Sigma_{(u,v)}}\inf_{y\in\Sigma}d_M(f_{(u,v)}(x),
 f(y))\\
 \le&& 2\max_{v_z\ne 0}R_z^{-1/2}\\
 =&& 2\max_{v_z\ne 0}|v_z|^{1/4}\le 2|v|^{1/4},
 \end{eqnarray*}
  where $ d_M$ is a fixed Riemannian distance on $M$. It follows
  that $${\rm diam}_\mu(f_{(u,v)}(\Sigma_{(u,v)}))\le {\rm
  diam}_\mu (f(\Sigma))+ 4|v|^{1/4}\quad{\rm and}\eqno(2.8)$$\vspace{-4mm}
  $$f_{(u,v)}(\Sigma_{(u,v)}))\subset\{q\in M\,|\, d_\mu(q,
  f(\Sigma))\le 4|v|^{1/4}\}\eqno(2.9)$$
for every $(u,v)\in V_{deform}(\sigma)\times
V_{resolve}(\sigma)$.\vspace{2mm}

\noindent{\bf Remark 2.12.}\quad For $[{\bf f}]=
[(f;\Sigma,\bar{\bf z})]\in\overline{\mathcal M}_{g,m}(M, J, A)$
and $\sigma=(\Sigma,\bar{\bf z}\cup\bar{\bf y})$ as above,  since
$V_{deform}({\sigma})$ (resp. $V_{resolve}({\sigma})$) is a
neighborhood of the origin in the vector space $\prod_s\mbox{\Bb
C}^{3g_s-3+m_s+l_s}$ (resp. $\prod_{z\in
Sing(\Sigma)}T_{z_s}\widetilde\Sigma_s\otimes
T_{z_t}\widetilde\Sigma_t$) we can fix a $\delta_{\bf f}>0$ such
that the open ball $V_{\delta_{\bf f}}$ with center $0$ and radius
$\delta_{\bf f}$ in the product space $\prod_s\mbox{\Bb
C}^{3g_s-3+m_s+l_s}\times\prod_{z\in
Sing(\Sigma)}T_{z_s}\widetilde\Sigma_s\otimes
T_{z_t}\widetilde\Sigma_t$ is contained in $
V_{deform}({\sigma})\times V_{resolve}({\sigma})$. Moreover, (2.6)
also implies
$$f_{(u,v)}(z_j)=f(z_j)\quad\forall z_j\in\bar{\bf
z}.\eqno(2.10)$$
 Furthermore we require $\delta_{\bf f}>0$ so small that
 each ${\bf f}_{(u,v)}=(f_{(u,v)}, \Sigma_{(u,v)},\bar{\bf
 z})$ is a stable $L^{k,p}$-map and has the homology class $A$ and the energy
 $$E({\bf f}_{(u,v)})\le\omega(A)+1\eqno(2.11)$$
since $E({\bf f})=\omega(A)$.

\subsubsection{Strong $L^{k,p}$-topology}
 For $\epsilon>0$ we denote
$$\widetilde{\bf U}_\epsilon(\delta_{\bf f})$$ by
the set of all tuples $(g_{(u,v)},\Sigma_{(u,v)},\bar{\bf
z}\cup\bar{\bf y})$ satisfying the following conditions:
  $$(u,v)\in V_{\delta_{\bf f}};\eqno(2.12)$$
  $$g_{(u,v)}:\Sigma_{(u,v)}\to M\;{\rm is}\;{\rm a}\; L^{k,p}-{\rm map};\eqno(2.13)$$
\begin{eqnarray*}
 && \|g_{(u,v)}-f_{(u,v)}\|_{k,p}<\epsilon,\;{\rm where}\;{\rm the}\;{\rm
 norm}\; \|\cdot\|_{k,p}\;{\rm is}\;{\rm measured}\;{\rm with}\;{\rm respect}\;{\rm to}\;{\rm the}\\
&&{\rm metrics}\;  \tau_{(u,v)}\;{\rm and}\;
g_J=\mu.\hspace{91mm}(2.14)
 \end{eqnarray*}
 Note that our choices of $(k,p)$ ensure that the norm
 $\|\cdot\|_{k,p}$ is stronger than $C^2$-topology. Therefore
 if $\epsilon>0$ is small enough then each ${\bf g}_{(u,v)}=
(g_{(u,v)},\Sigma_{(u,v)},\bar{\bf z})$ is a stable $L^{k,p}$-map
and has also the homology class $A$ and the energy
$$E({\bf g}_{(u,v)})\le\omega(A)+2\eqno(2.15)$$
because of (2.11). By making $\delta_{\bf f}>0$ smaller we can
assume that these hold for all $0<\epsilon\le\delta_{\bf f}$.
Moreover it follows from (2.8) and (2.14) above that
 $$\max_{(u,v)\in V_{\bf f}}{\rm diam}_\mu(g_{(u,v)}(\Sigma_{(u,v)}))\le {\rm
  diam}_\mu (f(\Sigma))+ 4\delta_{\bf f}^{1/4}+ 2\epsilon$$
    for $(u,v)\in V_{\bf f}$.
We can also assume that $0<\epsilon\le\delta_{\bf f}<1$.  Then
(2.9) and (2.14)  imply
$$\bigcup_{(u,v)\in V_{\bf f}}g_{(u,v)}(\Sigma_{(u,v)})\subset
\{q\in M\,|\, d_\mu(q, f(\Sigma))\le 6\}.\eqno(2.16)$$
 For each $0<\delta\le\delta_{\bf f}$ we denote by
 $$\widetilde{\bf U}_\delta({\bf f}):=\{(g_{(u,v)},\Sigma_{(u,v)},\bar{\bf
z}\cup\bar{\bf y})\in\widetilde{\bf U}_\delta(\delta_{\bf f})\,|\,
(u,v)\in V_{\delta_{\bf f}},\,|(u,v)|<\delta\}.$$ Remark that the
above constructions still hold for any $[{\bf f}]\in {\mathcal
B}^M_{g, m, A}$. Let
$${\bf U}_\delta({\bf f}):=\{[g_{(u,v)},\Sigma_{(u,v)},\bar{\bf
z}]\,|\,(g_{(u,v)},\Sigma_{(u,v)},\bar{\bf z}\cup\bar{\bf
y})\in\widetilde{\bf U}_\delta({\bf f})\}\subset{\mathcal
B}^M_{g,m,A}$$
 and ${\mathcal U}=\{{\bf U}_\delta({\bf f})\,|\, {\bf
 f}\in [{\bf f}]\in{\mathcal B}^M_{g,m,A},\;0<\delta\le\delta_{\bf f} \}$. As in [LiuT1] we can
 prove \vspace{2mm}

 \noindent{\bf Proposition 2.13.}\quad{\it ${\mathcal U}$
 generates a topology, called the (strong) $L^{k,p}$-topology, on
${\mathcal B}^M_{g,m,A}$. The topology is equivalent to the weak
topology on $\overline{\mathcal M}_{g,m}(M,J,A)$.}\vspace{2mm}

Since $\overline{\mathcal M}_{g,m}(M,J,A; K)$ is compact with
respect to the weak topology for any compact subset $K\subset M$
we get\vspace{2mm}

\noindent{\bf Corollary 2.14.}\quad{\it There is an open
neighborhood ${\mathcal W}$ of $\overline{\mathcal
M}_{g,m}(M,J,A)$ in ${\mathcal B}^M_{g,m,A}$ such that ${\mathcal
W}$ is Hausdorff with respect to the $L^{k,p}$-topology.}

\subsubsection{Local uniformizers}

Let ${\bf f}$ be as in \S2.4.2.  For each $y_j\in\bar{\bf y}$ we
choose a codimension two small open disc $H_j\subset M$ such that
(i) $H_j$ intersects $f(\Sigma)$ uniquely and transversely at
$f(y_j)$, (ii) $f^{-1}(H_j)=f^{-1}(f(y_j))$, and (iii)  $H_j$ is
oriented  so that it has positive intersection with $f(\Sigma)$.
Let ${\rm H}=\prod_j H_j$ and
$$\widetilde{\bf U}_\delta({\bf f}, {\rm H}):=\{
(g_{(u,v)},\Sigma_{(u,v)},\bar{\bf z}\cup\bar{\bf y})\in
\widetilde{\bf U}_\delta({\bf f})\,|\, g_{(u,v)}(y_j)\in
H_j\;\forall y_j\in\bar{\bf y}\}$$ for $0<\delta\le\delta_{\bf
f}$. Slightly modifying the proof of Lemma 3.5 in [LiuT2] we have
\vspace{2mm}

 \noindent{\bf Proposition 2.15.}\quad{\it For $\delta>0$ sufficiently small
 there exists a continuous right action of $Aut({\bf f})$ on
$\widetilde{\bf U}_\delta({\bf f}, {\rm H})$ that is smooth on
each open stratum of $\widetilde{\bf U}_\delta({\bf f}, {\rm H})$.
Moreover this action also commutes with the projection
$$\widetilde{\bf U}_\delta({\bf f}, {\rm H})\to{\bf U}_\delta({\bf f}, {\rm H})\subset{\mathcal
B}^M_{g,m,A},\;(g_{(u,v)},\Sigma_{(u,v)},\bar{\bf z}\cup\bar{\bf
y})\mapsto  [g_{(u,v)},\Sigma_{(u,v)},\bar{\bf z}],$$ and the
induced quotient map $\widetilde{\bf U}_\delta({\bf f}, {\rm
H})/Aut({\bf f})\to {\bf U}_\delta({\bf f}, {\rm H})$ is a
homomorphism. Here ${\bf U}_\delta({\bf f}, {\rm
H}):=\{[g_{(u,v)},\Sigma_{(u,v)},\bar{\bf z}]\,|\,
(g_{(u,v)},\Sigma_{(u,v)},\bar{\bf z}\cup\bar{\bf
y})\in\widetilde{\bf U}_\delta({\bf f}, {\rm H})\}$.}
\vspace{2mm}\\
By making $\delta_{\bf f}$ smaller we later always assume that
Proposition 2.15 holds for all $0<\delta\le\delta_{\bf f}$.

\section{Virtual moduli cycles}\label{3}

In this section we shall follow the ideas in [LiuT1-3] to
construct the  virtual moduli cycle. For this goal we give the
necessary reviews of the arguments in [LiuT1-3] and different
points.

\subsection{The local construction}\label{3.1}

There is a natural  stratified Banach bundle $\widetilde{\bf
E}_\delta({\bf f}, {\rm H})\to\widetilde{\bf U}_\delta({\bf f},
{\rm H})$ whose fiber  $\widetilde{\bf E}_\delta({\bf f}, {\rm
H})_{{\bf g}_{(u,v)}}$ at ${\bf
g}_{(u,v)}=(g_{(u,v)},\Sigma_{(u,v)},\bar{\bf z}\cup\bar{\bf y})$
is given by $L^p_{k-1}(\wedge^{0,1}(g^\ast_{(u,v)}TM))$. Here
$\wedge^{0,1}(g^\ast_{(u,v)}TM)$ is the bundle of $(0,1)$-forms on
$\Sigma_{(u,v)}$ with respect to the complex structure $j_{(u,v)}$
on $\Sigma_{(u,v)}$ and the almost complex structure $J$ on $M$,
and the norm $\|\cdot\|_{k-1,p}$ in
$L^p_{k-1}(\wedge^{0,1}(g^\ast_{(u,v)}TM))$  is   with respect to
the metric $\tau_{(u,v)}$ on $\Sigma_{(u,v)}$ and the metric $\mu$
(or $g_J$) on $M$. The action of $Aut({\bf f})$ on $\widetilde{\bf
U}_\delta({\bf f}, {\rm H})$ can be lifted to a linear action on
$\widetilde{\bf E}_\delta({\bf f}, {\rm H})_{{\bf g}_{(u,v)}}$
such that the natural projection $\tilde p_{\bf f}: \widetilde{\bf
E}_\delta({\bf f}, {\rm H})\to\widetilde{\bf U}_\delta({\bf f},
{\rm H})$ is a $Aut({\bf f})$-equivariant and locally trivial
vector bundle when restricted over each stratum of $\widetilde{\bf
U}_\delta({\bf f}, {\rm H})$. Moreover there is a local orbifold
bundle ${\bf E}_\delta([{\bf f}])\to{\bf U}_\delta({\bf f}, {\rm
H})$ with $\widetilde{\bf E}_\delta({\bf f}, {\rm H})$ as the
local uniformizer. That is, $\widetilde{\bf E}_\delta({\bf f},
{\rm H})/Aut({\bf f})={\bf E}_\delta([{\bf f}])$. The fiber ${\bf
E}_\delta([{\bf f}])_{[{\bf h}]}$ at ${[{\bf h}]}\in {\bf
U}_\delta({\bf f}, {\rm H})$ consists of all elements of
$L^p_{k-1}(\wedge^{0,1}(h^\ast TM))$ modulo equivalence relation
defined via pull-back of sections induced by the equivalences of
the domains of different representatives of $[{\bf h}]$. Consider
the section
$$\bar\partial_J:\widetilde{\bf U}_\delta({\bf f}, {\rm H})\to
\widetilde{\bf E}_\delta({\bf f}, {\rm H}),\;{\bf
g}_{(u,v)}\mapsto\bar\partial_J  g_{(u,v)}.\eqno(3.1)
$$
As in [LiuT2, \S5.1] we choose the previous each open disk $H_j$
to be totally geodesic with respect to the metric $g_J$ so that
the tangent space of $\widetilde{\bf U}_\delta({\bf f}, {\rm H})$
at ${\bf g}_{(u,v)}=(g_{(u,v)},\Sigma_{(u,v)}, \bar{\bf
z}\cup\bar{\bf y})$, denoted by $T_{{\bf g}_{(u,v)}}\widetilde{\bf
U}_\delta({\bf f}, {\rm H})$, is equal to
$$L^{k,p}(g^\ast_{(u,v)}TM, T_{{\bf f}(\bar{\bf y})}{\rm H}):=\{
\xi\,|\, \xi\in L^{k,p}(g^\ast_{(u,v)}TM),\;\xi(y_j)\in
T_{y_j}H_j\,\forall y_j\in\bar{\bf y}\}.
$$
For $(u,v)\in \Lambda_\delta({\bf f})=\{(u,v)\in V_{\delta_{\bf
f}}\,|\, |(u,v)|<\delta\}$ let
$$
\widetilde{\bf U}^{(u,v)}_{\delta}({\bf f}, {\rm H}):=\{{\bf
g}_{(u,v)}\in\widetilde{\bf U}_{\delta}({\bf f}, {\rm
H})\}\quad{\rm and}\quad
\widetilde{\bf E}^{(u,v)}_{\delta}({\bf
f}, {\rm H}):=\widetilde{\bf E}_{\delta}({\bf f}, {\rm
H})|_{\widetilde{\bf U}^{(u,v)}_{\delta}({\bf f}, {\rm H})}.
$$
Following [LiuT1-2] one has  the trivialization of $\widetilde{\bf
E}^{(u,v)}_{\delta}({\bf f}, {\rm H})$,
$$\psi^{(u,v)}:\widetilde{\bf U}^{(u,v)}_{\delta}({\bf f}, {\rm
H})\times L^p_{k-1}(f^\ast_{(u,v)}TM)\to\widetilde{\bf
E}^{(u,v)}_{\delta}({\bf f}, {\rm H})\eqno(3.2)
$$
obtained by the $J$-invariant parallel transformation of $(M, J)$.
 Under the trivialization the
restriction of the section $\bar\partial_J$ in (3.1) to
$\widetilde{\bf U}^{(u,v)}_\delta({\bf f}, {\rm H})$ has the
following representative,
$$F_{(u,v)}=\pi_2\circ(\psi^{(u,v)})^{-1}\circ\bar\partial_J:\widetilde{\bf U}^{(u,v)}_{\delta}({\bf
f}, {\rm H})\to L^p_{k-1}(f^\ast_{(u,v)}TM).\eqno(3.3)
$$
 Since
the vertical differential
$$D\bar\partial_J({\bf f}):
L^{k,p}(f^\ast TM, T_{{\bf f}(\bar{\bf y})}{\rm H})\to
\widetilde{\bf E}_\delta({\bf f}, {\rm H})_{\bf f}
$$
is a Fredholm operator its cokernel $R({\bf f})\subset
\widetilde{\bf E}_\delta({\bf f}, {\rm H})_{\bf f}$ is finite
dimensional. So
$$(D\bar\partial_J({\bf f}))\oplus I_{\bf f}:
(T_{\bf f}\widetilde{\bf U}_\delta({\bf f}, {\rm H})) \oplus
R({\bf f})\to \widetilde{\bf E}_\delta({\bf f}, {\rm H})_{\bf f}
$$
is surjective, where $I_{\bf f}$ denotes the inclusion $R({\bf
f})\hookrightarrow \widetilde{\bf E}_\delta({\bf f}, {\rm H})_{\bf
f}$. Note that $Aut_{\bf f}$ is a finite group. $\sum_{\sigma\in
Aut_{\bf f}}\sigma(R({\bf f}))$ is still finite dimensional and
$Aut_{\bf f}$ acts on it. By slightly perturbing $\sum_{\sigma\in
Aut_{\bf f}}\sigma(R({\bf f}))$ we can assume that $Aut_{\bf f}$
acts on it freely. {\it We still denote $R({\bf f})$ by this
perturbation below}.  Take a smooth cut-off function
$\gamma_\epsilon({\bf f})$ supported outside of the
$\epsilon$-neighborhood of double points of
$\Sigma_{(0,0)}=\Sigma$ for a small $\epsilon>0$, and denote by
$R_\epsilon({\bf f})=\{\gamma_\epsilon({\bf f})\cdot\xi\,|\,
\xi\in R({\bf f})\}$. Let us still use $I_{\bf f}$ to denote the
inclusion $R_\epsilon({\bf f})\hookrightarrow \widetilde{\bf
E}_\delta({\bf f}, {\rm H})_{\bf f}$.
 It was proved in [LiuT1] that there exists a sufficiently small
$\epsilon_{\bf f}>0$ such that
$$L_{\bf f}:=(D\bar\partial_J({\bf f}))\oplus I_{\bf f}:
(T_{\bf f}\widetilde{\bf U}_\delta({\bf f},{\rm H})) \oplus
R_\epsilon({\bf f})\to \widetilde{\bf E}_\delta({\bf f}, {\rm
H})_{\bf f}\eqno(3.4)
$$
 is still surjective  for all
$\epsilon\in (0,\epsilon_{\bf f}]$. We may require $\epsilon_{\bf
f}>0$ so small that the $\epsilon_{\bf f}$-neighborhood of double
points of $\Sigma$, denoted by ${\cal N}(Sing(\Sigma),
\epsilon_{\bf f})$, is not intersecting with $K_{deform}(\sigma)$.
\emph{Now we take $0<\delta(\epsilon_{\bf f})\le\delta_{\bf f}$ so
small that when $|(u,v)|\le \delta(\epsilon_{\bf f})$ our
construction for $(f_{(u,v)},\Sigma_{(u,v)})$ in \S2.4.2  only
need to change $f$ and $\Sigma$ in the $\epsilon_{\bf
f}/2$-neighborhood of double points of $\Sigma$.} Let ${\cal
N}(u,v;\epsilon_{\bf f})$ be a part of $\Sigma_{(u,v)}$ which is
obtained after resolving the double points of $\Sigma$ in ${\cal
N}(Sing(\Sigma), \epsilon_{\bf f})$ with gluing parameters $v$,
and such that
$$
\Sigma_{(u,v)}\setminus {\cal N}(u,v;\epsilon_{\bf
f})=\Sigma\setminus{\cal N}(Sing(\Sigma), \epsilon_{\bf
f})\eqno(3.5)
$$
as manifolds. Since $\nu\equiv 0$ in ${\cal N}(Sing(\Sigma),
\epsilon_{\bf f})$ each $\nu\in R_{\epsilon_{\bf f}}({\bf f})$ can
determine a unique element $\nu({\bf f}_{(u,v)})$ of
$L^p_{k-1}(\wedge^{0,1}(f^\ast_{(u,v)}TM))$ for all deformation
parameters $(u,v)$ satisfying $|(u,v)|\le \delta(\epsilon_{\bf
f})$.
$$
\nu({\bf f}_{(u,v)})(z)= \left\{
\begin{array}{ll}
&\!\!\!\!\!\frac{1}{2}(\nu(z)+ J(f(z))\circ\nu(z)\circ
j_{(u,v)})\quad{\rm if}\;z\in
\Sigma_{(u,v)}\setminus {\cal N}(u,v;\epsilon_{\bf f}),\vspace{1mm}\\
&\!\!\!\!\!0 \hspace{50mm}{\rm if}\;z\in
 {\cal N}(u,v;\epsilon_{\bf f}).
 \end{array}
 \right.
$$
Here $j_{(u,v)}$ is the complex structure on $\Sigma_{(u,v)}$
obtained in \S2.4.1. It is clear that
$$
\nu({\bf f}_{(u,v)})=\nu({\bf f}_{(u,0)})\quad{\rm on}\quad
\Sigma_{(u,v)}\setminus {\cal N}(u,v;\epsilon_{\bf
f})=\Sigma\setminus{\cal N}(Sing(\Sigma), \epsilon_{\bf f})
$$
for any $|(u,v)|\le \delta(\epsilon_{\bf f})$. For each ${\bf
g}_{(u,v)}=(g_{(u,v)},\Sigma_{(u,v)},\bar{\bf z}\cup\bar{\bf
y})\in\widetilde{\bf U}_{\delta(\epsilon_{\bf f})}({\bf f})$,
using the $J$-invariant parallel transformation of $(M, J)$ we get
a unique element $\nu({\bf g}_{(u,v)})\in
L^p_{k-1}(\wedge^{0,1}(g^\ast_{(u,v)}TM))$ from $\nu({\bf
f}_{(u,v)})$, which also satisfies:
$$
\nu({\bf g}_{(u,v)})(z)=0\quad{\rm for}\quad z\in
 {\cal N}(u,v;\epsilon_{\bf f}).\eqno(3.6)
 $$
  This gives rise to a section
$$\tilde\nu: \widetilde{\bf U}_{\delta(\epsilon_{\bf f})}({\bf f},
{\rm H})\to \widetilde{\bf E}_{\delta(\epsilon_{\bf f})}({\bf f},
{\rm H}),\;{\bf g}_{(u,v)}\mapsto \nu({\bf g}_{(u,v)}),\eqno(3.7)
$$
 which is continuous and stratawise smooth.  Clearly, $\tilde\nu=0$ as $\nu=0$.
 If $\nu\ne 0$ and $\delta(\epsilon_{\bf f})>0$ sufficiently small then
$\tilde\nu({\bf g}_{(u,v)})\ne 0$ for all $|(u,v)|\le
\delta(\epsilon_{\bf f})$, i.e., $\tilde\nu$ is not equal to the
zero at any point of $\widetilde{\bf U}_{\delta(\epsilon_{\bf
f})}({\bf f}, {\rm H})$. Take a basis $(\nu_1,\cdots,\nu_q)$ of
$R_{\epsilon_{\bf f}}({\bf f})$. We  require each $\nu_j$ so close
to the origin of $R_{\epsilon_{\bf f}}({\bf f})$ that
$$\sup\{\|\tilde\nu_j({\bf g}_{(u,v)})\|\,:\,{\bf g}_{(u,v)}
\in \widetilde{\bf U}_{\delta(\epsilon_{\bf f})}({\bf f}, {\rm
H})\}\le 1\eqno(3.8)
$$
 for $j=1,\cdots, q$.
For convenience we denote by
 $$
 \begin{array}{lccr}
 \widetilde W_{\bf f}:=\widetilde{\bf U}_{\delta(\epsilon_{\bf f})}({\bf f},
{\rm H})&\quad{\rm and}&\quad &W_{\bf f}:={\bf
U}_{\delta(\epsilon_{\bf f})}({\bf f}, {\rm H})=\pi_{\bf
f}(\widetilde W_{\bf
f}),\hspace{4mm}\\
\widetilde E_{\bf f}:=\widetilde{\bf E}_{\delta(\epsilon_{\bf
f})}({\bf f}, {\rm H})&\quad{\rm and}&\quad &E_{\bf f}:={\bf
E}_{\delta(\epsilon_{\bf f})}({\bf f}, {\rm H}).\hspace{22mm}
\end{array}
$$
 Take a $Aut({\bf
f})$-invariant, continuous and stratawise smooth cut-off function
$\beta_{\bf f}$ on $\widetilde W_{\bf f}$ such that it is equal to
$1$ near ${\bf f}$ and that
$$
\widetilde U^{0}_{\bf f}:=\{x\in\widetilde W_{\bf f}\;|\;
\beta_{\bf f}(x)>0\}\subset\subset\widetilde W_{\bf f}.\eqno(3.9)
$$
Hereafter the notation $A\subset\subset B$ denotes that the
closure $Cl(A)$ of $A$ is contained in $B$. Let  $U^0_{\bf
f}=\pi_{\bf f}(\widetilde U_0)$,
$$
\widetilde V^{0}_{\bf f}:={\rm Int}(\{x\in\widetilde W_{\bf
f}\;|\; \beta_{\bf f}(x)=1\}).\eqno(3.10)
$$
and $V^0_{\bf f}=\pi_{\bf f}(\widetilde V^0_{\bf f})$. Then
$\widetilde V^{0}_{\bf f}=(\pi_{\bf f})^{-1}(V^0_{\bf f})$, and it
is also a $Aut({\bf f})$-invariant open neighborhood of ${\bf f}$
in $\widetilde W_{\bf f}$.  Let
  $$\tilde s_j=\beta_{\bf f}\cdot\tilde\nu_j,\; j=1,\cdots,q.\eqno(3.11)$$
 Each $\tilde s_j$ is a stratawise smooth section of the bundle $\widetilde E_{\bf f}\to\widetilde W_{\bf
 f}$, and for $j=1,\cdots, q$,
 $$
 \{{\bf h}\in\widetilde W_{\bf f}\,|\, \tilde s_j({\bf h})\ne 0\}=
 \widetilde U_{\bf f}^{0}\quad{\rm and}\quad
 \sup\{\|\tilde s_j({\bf h})\|\,:\,{\bf h}\in\widetilde W_{\bf f}\}\le 1.
 \eqno(3.12)
 $$
  Here the second claim  follows from (3.8).  For any $|(u,v)|\le \delta(\epsilon_{\bf f})$
 let $\widetilde W_{\bf f}^{(u,v)}=\widetilde{\bf U}^{(u,v)}_{\delta(\epsilon_{\bf f})}({\bf f},
{\rm H})$ and
$$
s_j^{(u,v)}:\widetilde W_{\bf f}^{(u,v)}\to
L^p_{k-1}(f^\ast_{(u,v)}TM)\eqno(3.13)
$$
be the representative of the restriction of the section $\tilde
s_j$ to  $\widetilde W_{\bf f}^{(u,v)}$ under the trivialization
in (3.2), $j=1,\cdots, q$. It is easily checked  that  for
$j=1,\cdots, q$,
$$
\left.
\begin{array}{ll}
 s_j^{(u,v)}({\bf g}_{(u,v)})=\beta_{\bf
f}({\bf g}_{(u,v)}) \cdot\nu_j({\bf f}_{(u,v)}),\;{\rm and}\;{\rm
thus}\vspace{1mm}\\
ds_j^{(u,v)}({\bf g}_{(u,v)})=0 \;\bigr(T_{{\bf
g}_{(u,v)}}\widetilde W_{\bf f}^{(u,v)}\to
L^p_{k-1}(f^\ast_{(u,v)}TM)\bigl)\;{\rm if}\; {\bf g}_{(u,v)}\in
\widetilde V^{0}_{\bf f}.
\end{array}\right.\eqno(3.14)
$$
Moreover, since $Z(\bar\partial_J)\cap Cl(\widetilde U^0_{\bf f})$
is a compact subset in $\widetilde W_{\bf f}$ there exists a
constant $c(\beta_{\bf f})>0$ only depending on $\beta_{\bf f}$
and an open neighborhood ${\cal N}\bigl(Z(\bar\partial_J)\cap
Cl(\widetilde U^0_{\bf f})\bigr)$ of it in $\widetilde W_{\bf f}$
such that
$$
 \|ds_j^{(u,v)}({\bf g}_{(u,v)})\|\le c(\beta_{\bf
f})\;\forall{\bf g}_{(u,v)}\in{\cal N}\bigl(Z(\bar\partial_J)\cap
Cl(\widetilde U_{\bf f})\bigr)\cap\widetilde W_{\bf
f}^{(u,v)}\eqno(3.15)
$$
for $|(u,v)|< \delta(\epsilon_{\bf f})$ if $\delta(\epsilon_{\bf
f})>0$ is small enough. Here $Z(\bar\partial_J)$ is the zero set
of the section $\bar\partial_J:\widetilde W_{\bf f}\to\widetilde
E_{\bf f}$. Note that the map $L_{\bf f}$ in (3.4) being
surjective is equivalent to that the map
\begin{eqnarray*}
&&T_{\bf f}\widetilde{\bf U}^{(0,0)}_{\delta(\epsilon_{\bf
f})}({\bf f}, {\rm
H})\times\mbox{\Bb R}^q\to L^p_{k-1}(f^\ast_{(0,0)}TM),\\
&&(\xi; t_1,\cdots,t_q)\mapsto dF_{(0, 0)}({\bf
f})(\xi)+\sum^q_{j=1}t_j \nu_j
\end{eqnarray*}
is surjective.   By refining the gluing arguments in [McSa1] and
[Liu] the following result was proved in [LiuT1] . \vspace{2mm}

\noindent{\bf Lemma 3.1.}\quad{\it For sufficiently small
$\delta(\epsilon_{\bf f})>0$, $|(u,v)|\le \delta(\epsilon_{\bf
f})$ and ${\bf g}_{(u,v)}\in \widetilde U^{0}_{\bf
f}\cap\widetilde W_{\bf f}^{(u,v)}$ the map
\begin{eqnarray*}
&&T_{{\bf g}_{(u,v)}}(\widetilde U^{0}_{\bf f}\cap\widetilde
W_{\bf
f}^{(u,v)})\times \mbox{\Bb R}^q\to L^p_{k-1}(f^\ast_{(u,v)}TM),\\
&& (\xi;e_1,\cdots,e_q)\mapsto dF_{(u,v)}({\bf
g}_{(u,v)})(\xi)+\sum^q_{j=1}e_j s_j^{(u,v)}({\bf g}_{(u,v)}),
\end{eqnarray*}
is surjective.  In particular the map
\begin{eqnarray*}
&&T_{{\bf g}_{(u,v)}}\widetilde U^{0}_{\bf f}\times\mbox{\Bb R}^q\to L^p_{k-1}(f^\ast_{(u,v)}TM),\\
&&(\xi;e_1,\cdots,e_q)\mapsto dF_{(u,v)}({\bf
g}_{(u,v)})(\xi)+\sum^q_{j=1}e_j s_j^{(u,v)}({\bf g}_{(u,v)})
\end{eqnarray*}
is surjective since $T_{{\bf g}_{(u,v)}}(\widetilde U^{0}_{\bf
f}\cap\widetilde W_{\bf f}^{(u,v)})\subset T_{{\bf
g}_{(u,v)}}\widetilde U^{0}_{\bf f}$.}\vspace{2mm}

 For $\eta>0$ let $\Pi_1$ be the projection to the first factor of
 $\widetilde U^{0}_{\bf f}\times {\bf B}_{\eta}(\mbox{\Bb R}^q)$.
 Consider the sections
$$
\Upsilon^{({\bf f})}:\widetilde U^{0}_{\bf f}\times B_{\eta}(
\mbox{\Bb R}^q)\to \Pi_1^\ast\widetilde E_{\bf f},\;({\bf h},{\bf
t})\mapsto\bar\partial_J(h)+\sum^q_{j=1}t_j\tilde s_j({\bf h}),
$$
 and $\Upsilon^{({\bf f})}_{\bf t}=\Upsilon^{({\bf f})}(\cdot,{\bf t}):
 \widetilde U^{0}_{\bf f}\to\widetilde E_{\bf
 f}$ for ${\bf t}=(t_1,\cdots,t_q)\in {\bf B}_\eta(\mbox{\Bb R}^q)$.
 Clearly, they are the continuous and stratawise smooth.
For any $|(u,v)|\le \delta(\epsilon_{\bf f})$
 let  $\Upsilon^{({\bf f}; u,v)}$ be the restriction of
$\Upsilon^{({\bf f})}$ to $(\widetilde U^{0}_{\bf f}\cap\widetilde
W_{\bf f}^{(u,v)})\times {\bf B}_{\eta}(\mbox{\Bb R}^q)$. As in
(3.3), under the trivialization in (3.2) $\Upsilon^{({\bf f};
u,v)}$ has the representation
\begin{eqnarray*}
&&\widetilde U^{0}_{\bf f}\cap\widetilde W_{\bf f}^{(u,v)}\times
{\bf B}_{\eta}( \mbox{\Bb R}^q)\to
 L^p_{k-1}(f^\ast_{(u,v)}TM),\\
&&({\bf g}_{(u,v)}, {\bf t})\mapsto F_{(u,v)}({\bf
g}_{(u,v)})+\sum^q_{j=1}e_j s_j^{(u,v)}({\bf g}_{(u,v)}).
\end{eqnarray*}
 Its tangent map at $({\bf
g}_{(u,v)}, {\bf t})\in\widetilde V^0_{\bf f}\cap\widetilde
W^{(u,v)}_{\bf f}\times {\bf B}_{\eta}(\mbox{\Bb R}^q)$ is exactly
given by
\begin{eqnarray*}
&&T_{{\bf g}_{(u,v)}}(\widetilde V^{0}_{\bf f}\cap\widetilde
W_{\bf
f}^{(u,v)})\times \mbox{\Bb R}^q\to L^p_{k-1}(f^\ast_{(u,v)}TM),\\
&&(\xi;e_1,\cdots,e_q)\mapsto dF_{(u,v)}({\bf
g}_{(u,v)})(\xi)+\sum^q_{j=1}\bigr(e_j s_j^{(u,v)}({\bf
g}_{(u,v)})+ t_j ds_j^{(u,v)}({\bf g}_{(u,v)})(\xi)\bigl)\\
&&\hspace{25mm}=dF_{(u,v)}({\bf g}_{(u,v)})(\xi)+\sum^q_{j=1}e_j
s_j^{(u,v)}({\bf g}_{(u,v)}).
\end{eqnarray*}
Here the final equality is because of (3.14). By Lemma 3.1 it is
surjective. For any given open subset $\widetilde V_{\bf
f}^{0+}\subset\widetilde U_{\bf f}^{0}$ with $\widetilde V_{\bf
f}^{0}\subset\widetilde V_{\bf f}^{0+}\subset\subset\widetilde
U_{\bf f}^{0}$, it follows from this and (3.6) that if $\eta>0$ is
small enough the map
\begin{eqnarray*}
&&T_{{\bf g}_{(u,v)}}\widetilde V^{0+}_{\bf f}\times\mbox{\Bb R}^q\to L^p_{k-1}(f^\ast_{(u,v)}TM),\\
&&(\xi;e_1,\cdots,e_q)\mapsto dF_{(u,v)}({\bf
g}_{(u,v)})(\xi)+\sum^q_{j=1}\bigr(e_j s_j^{(u,v)}({\bf
g}_{(u,v)})+ t_j ds_j^{(u,v)}({\bf g}_{(u,v)})(\xi)\bigl)
\end{eqnarray*}
is surjective for ${\bf t}\in {\bf B}_\eta(\mbox{\Bb R}^q)$.
 So we get immediately\vspace{2mm}

\noindent{\bf Proposition 3.2.}\quad{\it For any given open subset
$\widetilde V_{\bf f}^{0+}\subset\widetilde U_{\bf f}^{0}$ with
$\widetilde V_{\bf f}^{0}\subset\widetilde V_{\bf
f}^{0+}\subset\subset\widetilde U_{\bf f}^{0}$, there exists a
sufficiently small $\eta=\eta_{\bf f}>0$ such that the restriction
of $\Upsilon^{({\bf f})}$ to each stratum of $\widetilde
V^{0+}_{\bf f}\times {\bf B}_{\eta}(\mbox{\Bb R}^q)$ is
transversal to the zero section. So there exists a residual subset
${\bf B}_{\eta}^{res}(\mbox{\Bb R}^q)$ in ${\bf
B}_{\eta}(\mbox{\Bb R}^q)$ such that for each ${\bf t}\in {\bf
B}_{\eta}^{res}(\mbox{\Bb R}^q)$ the restriction of the section
$\Upsilon^{({\bf f})}_{\bf t}: \widetilde V^{0+}_{\bf f}\to
\widetilde E_{\bf f}$ to each stratum of $\widetilde V^{0+}_{\bf
f}$ is transversal to the zero section. Consequently, for each
${\bf t}\in {\bf B}_{\eta}^{res}(\mbox{\Bb R}^q)$,
$\widetilde{\mathcal M}({\bf f})^{{\bf t}}:=(\Upsilon^{({\bf
f})}_{\bf t}|_{\widetilde V^{0+}_{\bf f}})^{-1}(0)$ is a
stratified and cornered manifold in $\widetilde V^{0+}_{\bf f}$
that has the (even) dimension given by the Index Theorem on all of
its strata. Moreover, any two different ${\bf t}_0,\, {\bf t}_1\in
{\bf B}_{\eta}^{res}(\mbox{\Bb R}^q)$ give cobordant stratified
and cornered manifolds $\widetilde{\mathcal M}({\bf f})^{{\bf
t}_0}$ and $\widetilde{\mathcal M}({\bf f})^{{\bf t}_1}$. }
\vspace{2mm}

\subsection{The global construction}\label{3.2}

The main ideas still follow [LiuT1-3]. But we need to be more
careful and must give suitable modifications for their arguments
because our $\overline{\mathcal M}_{g,m}(M,J,A)$ is not compact.

 For any $[{\bf f}]=[f;\Sigma,\bar{\bf z}]\in\overline{\mathcal
M}_{g,m}(M, J, A)$, let $\widetilde W_{\bf f}=\widetilde{\bf
U}_{\delta(\epsilon_{\bf f})}({\bf f}, {\rm H})$ as above. By
(2.16),
 $$\bigcup_{{\bf g}_{(u,v)} \in\widetilde W_{\bf f}}g_{(u,v)}(\Sigma_{(u,v)})\subset
\{y\in M\,|\, d_\mu(y, f(\Sigma))\le 6\}.\eqno(3.16)
$$
 Fix a compact subset $K_0\subset M$, and denote by
$$K_j:=\{y\in M\,|\, d_\mu(y, K_0)\le
jC\},\;j=1,2,\cdots,\eqno(3.17)
$$
 where
 $$C=C(\alpha_0,\beta_0, C_0, i(M,\mu))=
 6+\frac{4\beta_0}{\pi\alpha^2_0 r_0}\omega(A).\eqno(3.18)
 $$
  It easily follows from Proposition 2.8 and (3.16) that
$$\bigcup_{[{\bf f}]\in\overline{\mathcal M}_{g,m}(M, J,
A;K_j)}\, \bigcup_{{\bf g}_{(u,v)}\in\widetilde W_{\bf
f}}g_{(u,v)}(\Sigma_{(u,v)})\subset
 K_{j+1}\eqno(3.19)
 $$
 for any $j=0,1,2,\cdots$. Since each $\overline{\mathcal M}_{g,m}(M, J,
A;K_j)$ is compact we can choose
 \begin{eqnarray*}
 &&{\bf f}^{(i)}\in\overline{\mathcal M}_{g,m}(M, J,
 A;K_0),\;i=1,\cdots,n_0,\quad{\rm and}\\
&&{\bf f}^{(i)}\in\overline{\mathcal M}_{g,m}(M, J,
 A;K_{j+1})\setminus\overline{\mathcal M}_{g,m}(M, J,
 A;K_j),\;i=n_j+1,\cdots,n_{j+1},
 \end{eqnarray*}
 $j=0,1,2$,
 such that the corresponding $V^0_{{\bf f}^{(i)}}$ as in (3.10) and $W_{{\bf f}^{(i)}}$  satisfy
$$\cup^{n_j}_{i=1}V^0_{{\bf f}^{(i)}}\supset\overline{\mathcal M}_{g,m}(M, J,
A;K_j)\quad{\rm and}\quad\cup^{n_0}_{i=1}W_{{\bf
f}^{(i)}}=\emptyset \eqno(3.20)
$$
 for $j=0,1,2,3$. The second requirement is needed later and
 can always be satisfied because we may add two disjoint open subsets to
$W_{{\bf f}^{(i)}}$, $i=1,\cdots,n_0$.\vspace{2mm}

\noindent{\bf Abbreviation Notations:}\quad For $i=1,2,\cdots,
n_3$, we abbreviate
$$\begin{array}{lcccccr}
 \widetilde W_i=\widetilde W_{{\bf f}^{(i)}},&\quad &W_i=W_{{\bf f}^{(i)}},&\quad
 &\widetilde E_i=\widetilde E_{{\bf f}^{(i)}},&\quad &E_i:=E_{{\bf f}^{(i)}},\\
\widetilde U^0_i=\widetilde U^0_{{\bf f}^{(i)}},&\quad
&U^0_i=U^0_{{\bf f}^{(i)}},&\quad &\widetilde V^0_i=\widetilde
V^0_{{\bf f}^{(i)}},&\quad &V^0_i:=V^0_{{\bf f}^{(i)}},\\
\pi_i=\pi_{{\bf f}^{(i)}}, &\quad &\Gamma_i=Aut({\bf
f}^{(i)}),&\quad &\eta_i:=\eta_{\bf f^{(i)}},&\quad &
\beta_i=\beta_{{\bf f}^{(i)}},\\
\delta_i=\delta(\epsilon_{\bf f}^{(i)}).
\end{array}$$
As in (3.11), for $i=1,2,\cdots,n_3$ let  $\tilde s_{ij}$,
$j=1,\cdots,q_i$, be the corresponding sections of the bundle
$\widetilde E_i\to\widetilde W_i$ such that
$$
\{{\bf h}\in\widetilde W_i\,|\, \tilde s_{ij}({\bf h})\ne
0\}=\widetilde U_i^0,\quad \sup\{\|\tilde s_{ij}({\bf
h})\|\,:\,{\bf h}\in\widetilde W_i\}\le 1, \eqno(3.21)
$$
and that the corresponding  conclusions with those in Lemma 3.1
and Proposition 3.2 hold, i.e.,\vspace{2mm}

\noindent{\bf Lemma 3.3.}\quad{\it For  any $|(u,v)|\le \delta_i$
and ${\bf g}_{(u,v)}\in \widetilde U^{0}_i\cap\widetilde
W_i^{(u,v)}$ the map
$$
\begin{array}{ll}
T_{{\bf g}_{(u,v)}}(\widetilde U^{0}_i\cap\widetilde W_i^{(u,v)})\times
\mbox{\Bb R}^{q_i}\to L^p_{k-1}(f^{(i)\ast}_{(u,v)}TM),\vspace{1mm}\\
 (\xi;e_1,\cdots,e_{q_i})\mapsto dF^{(i)}_{(u,v)}({\bf
g}_{(u,v)})(\xi)+\sum^{q_i}_{j=1}e_j s_{ij}^{(u,v)}({\bf
g}_{(u,v)})
\end{array}$$
and thus
$$
\begin{array}{ll}
T_{{\bf g}_{(u,v)}}\widetilde U^{0}_i\times\mbox{\Bb R}^{q_i}\to L^p_{k-1}(f^{(i)\ast}_{(u,v)}TM),\vspace{1mm}\\
(\xi;e_1,\cdots,e_{q_i})\mapsto dF^{(i)}_{(u,v)}({\bf
g}_{(u,v)})(\xi)+\sum^{q_i}_{j=1}e_j s_{ij}^{(u,v)}({\bf
g}_{(u,v)})
\end{array}
$$
is surjective, $i=1,\cdots, n_3$, where
$$ F^{(i)}_{(u,v)},\, s_{ij}^{(u,v)}:\widetilde W^{(u,v)}_{{\bf
f}^{(i)}}\to L^p_{k-1}(f^{(i)\ast}_{(u,v)}TM),\;j=1,\cdots, q_i,
$$
are the representatives of the sections $\bar\partial_J, \tilde
s_{ij}: \widetilde W_i\to \widetilde E_i$ under the trivialization
as in (3.2). Moreover, for any given open subset $\widetilde
V_i^{0+}\subset\widetilde U_i^{0}$ with $\widetilde
V_i^{0}\subset\widetilde V_i^{0+}\subset\subset\widetilde
U_i^{0}$, there exists a sufficiently small $\eta_i>0$ such that
 for any
$|(u,v)|\le \delta_i$, ${\bf g}_{(u,v)}\in \widetilde V^{0}_i$ and
$(t_{i1},\cdots, t_{iq_i})\in {\bf B}_{\eta_i}(\mbox{\Bb
R}^{q_i})$  the maps
$$
\begin{array}{ll}
T_{{\bf g}_{(u,v)}}\widetilde V^{0+}_i\times \mbox{\Bb R}^{q_i}\to
L^p_{k-1}(f^{(i)\ast}_{(u,v)}TM), \vspace{1mm}\\
 (\xi;e_1,\cdots,e_{q_i})\mapsto dF^{(i)}_{(u,v)}({\bf
g}_{(u,v)})(\xi)+\sum^{q_i}_{j=1}\bigr(e_i s_{ij}^{(u,v)}({\bf
g}_{(u,v)})+ t_{ij} ds_{ij}^{(u,v)}({\bf g}_{(u,v)})(\xi)\bigl)
\end{array}
$$
is surjective, $i=1,\cdots, n_3$. }\vspace{2mm}

\noindent{Denote by}
$${\mathcal W}:=\cup_{i=1}^{n_3}W_i\quad{\rm and}\quad {\mathcal
E}:={\mathcal E}^M_{A,g,m}|_{\mathcal W}.
$$
 The proof of Lemma 4.2 in [LiuT1] or that of
 Theorem 3.4 in [LiuT2] showed:
 \vspace{2mm}

\noindent{\bf Lemma 3.4.}\quad{\it ${\mathcal W}$ is a stratified
orbifold and ${\mathcal E}\to {\mathcal W}$ is a stratified
orbifold bundle with respect to the above local uniformizers.}
\vspace{2mm} \\
Let ${\mathcal  N}$ be the set of all finite
subsets $I=\{i_1,\cdots, i_l\}$ of $\{1,\cdots,n_3\}$ such that
the intersection $W_I:=\cap_{i\in I}W_i$ is nonempty and that
$i_1<i_2<\cdots<i_k$ if $k>1$. Then each $I\in{\mathcal N}$ has
the length $|I|=\sharp(I)<n_3$ because of the second condition in
(3.20). Define $E_I:={\mathcal E}|_{W_I}$.
  For each $I\in {\mathcal  N}$ denote by the group
$\Gamma_I:=\prod_{i\in I}\Gamma_i$ and the fiber product
$$\widetilde W_I^{\Gamma_I}=\Bigl\{(u_i)_{i\in I}\in\prod_{i\in I}
\widetilde W_i\,\bigm|\, \pi_i(u_i)=\pi_j(u_j)\;\forall i, j\in
I\Bigr\}
$$
equipped with the induced topology from $\prod^k_{l=1} \widetilde
W_{i_l}$.
 Then the obvious projection  $\tilde\pi_I:
\widetilde W_I^{\Gamma_I}\to W_I$ has covering group $\Gamma_I$
whose action on $\widetilde W_I^{\Gamma_I}$ is given by
$$
 \phi_I\cdot (\tilde u_i)_{i\in
I}=(\phi_i\cdot \tilde u_i)_{i\in
I},\quad\forall\phi_I=(\phi_i)_{i\in I}\in \Gamma_I. \eqno(3.22)
$$
 (Hereafter we often write $(\tilde
v_{i_1},\cdots, \tilde v_{i_k})$ as $(\tilde v_i)_{i\in I}$.)
Namely, $\tilde\pi_I$ induces a homomorphism from $\widetilde
W_I^{\Gamma_I}/\Gamma_I$ onto $W_I$. Moreover, for $J\subset I\in
{\mathcal  N}$ there are  projections
$$
 \tilde\pi^I_J: \widetilde
W_I^{\Gamma_I}\to \widetilde W_J^{\Gamma_J}, \;\;(\tilde
u_i)_{i\in I}\mapsto (\tilde u_j)_{j\in J}, \eqno(3.23)
$$
 and $\lambda^I_J:\Gamma_I\to\Gamma_J$ given by $(\phi_i)_{i\in
I}\mapsto (\phi_j)_{j\in J}$. Later we also write $\lambda^I_J$ as
$\lambda^i_J$ if $I=\{i\}$. Note that $\pi^I_J$ is not surjective
in general. Repeating the same construction from $\widetilde E_i$
one obtains the bundles $\tilde p_I:\widetilde
E_I^{\Gamma_I}\to\widetilde W_I^{\Gamma_I} $ and the projections
$\widetilde\Pi_I:\widetilde E_I^{\Gamma_I}\to E _I$ and
$\widetilde\Pi^I_J:\widetilde E_I^{\Gamma_I}\to\widetilde
E_J^{\Gamma_J}$ as $J\subset I$. Their properties may be
summarized as:\vspace{2mm}

\noindent{\bf Lemma 3.5.}\quad {\it {\rm (i)}
$\tilde\pi_J\circ\tilde\pi^I_J=\iota^W_{IJ}\circ\tilde\pi_I$ and
$\widetilde\Pi_J\circ\widetilde\Pi^I_J=\iota^E_{IJ}\circ\widetilde\Pi_I$
 for any $J\subset I\in{\cal N}$ and the inclusions

 $\iota_{IJ}^W:W_I\hookrightarrow
 W_J$, $\iota_{IJ}^E:E_I\hookrightarrow E_J$.\\
{\rm (ii)}
$\tilde\pi^I_J\circ\phi_I=\lambda^I_J(\phi_I)\circ\tilde\pi^I_J$
and
$\widetilde\Pi^I_J\circ\phi_I=\lambda^I_J(\phi_I)\circ\widetilde\Pi^I_J$
for any $J\subset I\in{\cal N}$ and   $\phi_I\in\Gamma_I$.\\
{\rm (ii)} The obvious projection $\tilde p_I:\widetilde
E_I^{\Gamma_I}\to\widetilde W_I^{\Gamma_I}$ is equivariant with
respect to the induced actions of

$\Gamma_I$ on them, i.e. $\tilde p_I\circ\psi_I=\psi_I\circ\tilde
p_I$
for any $\psi_I\in\Gamma_I$.\\
{\rm (iii)} $\tilde\pi_I\circ\tilde
p_I=p\circ\widetilde\Pi_I\;\forall\phi_I\in\Gamma_I$, and
$\tilde\pi^I_J\circ\tilde p_I=\tilde
 p_J\circ\widetilde\Pi^I_J$ for any $J\subset I\in{\cal N}$.\\
 {\rm (iv)} The generic
fiber of both $\tilde\pi^I_J$ and $\widetilde\Pi^I_J$ contains
$|\Gamma_I|/|\Gamma_J|$ points.}\vspace{2mm}

Therefore we get a system of bundles
$$(\widetilde{\mathcal E}^\Gamma, \widetilde W^\Gamma)=
\Bigl\{(\widetilde E_I^{\Gamma_I}, \widetilde W_I^{\Gamma_I}),
\tilde\pi_I, \widetilde\Pi_I, \Gamma_I, \tilde\pi^I_J,
\widetilde\Pi^I_J,\,\lambda^I_J\,\bigm|\, J\subset I\in{\mathcal
N}\Bigr\}.
$$
Note that as proved by Liu-Tian each $(\widetilde E_I^{\Gamma_I},
\widetilde W_I^{\Gamma_I})$ is only a pair of stratified Banach
varieties in the sense that locally they are finite union of
stratified Banach manifolds. For reader's convenience let us
explain in details  the arguments in terms of the concept of
``local component" in [LiuT3]. A continuous map from a stratified
manifold to another is called a {\it partially smooth} (
\emph{abbreviated as} PS) if it restricts to a smooth map on each
stratum. Similarly, in the partial smooth category one has the
notions of the PS diffeomorphism, PS embedding and  PS bundle map
(resp. isomorphism).

For $I=\{i_1,\cdots,i_k\}\in{\mathcal N}$ with $k>1$ let $\tilde
u_I\in\widetilde W_I^{\Gamma_I}$ be given. We choose a small
connected open neighborhood $O$ of $u_I=\pi_I(\tilde u_I)$ in
$W_I$ and consider the inverse image $O(\tilde
u_{i_l})=\pi_{i_l}^{-1}(O)$ of $O$ in $\widetilde W_{i_l}$,
$l=1,\cdots,k$. Then each $O(\tilde u_{i_l})$ is an open
neighborhood of $\tilde u_{i_l}$ in $\widetilde W_{i_l}$, and the
\emph{stabilizer subgroup} $\Gamma(\tilde u_{i_l})$  of
$\Gamma_{i_l}$ at $\tilde u_{i_l}$ acts on $O(\tilde u_{i_l})$. If
$O$ is small enough, for any fixed $s\in \{1,\cdots,k\}$ it
follows from the definition of the orbifold that there exist PS
diffeomorphisms
 $\lambda_{i_si_l}: O(\tilde u_{i_s})\to O(\tilde
u_{i_l})$  mapping $\tilde u_{i_s}$ to $\tilde u_{i_l}$, and group
isomorphisms ${\mathcal A}_{i_si_l}:\Gamma(\tilde
u_{i_1})\to\Gamma(\tilde u_{i_l})$  such that
$$\lambda_{i_si_l}\circ\phi={\mathcal A}_{i_si_l}(\phi)\circ\lambda_{i_si_l}\eqno(3.24)
$$
for $\phi\in\Gamma(\tilde u_{i_s})$ and
$l\in\{1,\cdots,k\}\setminus\{s\}$.
 These PS diffeomorphisms
$\lambda_{i_si_l}$ ($s\ne l$) in (3.24) are unique up to
composition with elements in $\Gamma(\tilde u_{i_s})$ and
$\Gamma(\tilde u_{i_l})$. Namely, if we have another PS
diffeomorphism $\lambda'_{i_si_l}: O(\tilde u_{i_s})\to O(\tilde
u_{i_l})$ that maps $\tilde u_{i_s}$ to $\tilde u_{i_l}$, and the
group isomorphism ${\mathcal A}'_{i_si_l}:\Gamma(\tilde
u_{i_s})\to\Gamma(\tilde u_{i_l})$ such that
$$
\lambda'_{i_si_l}\circ\phi={\mathcal
A}'_{i_si_l}(\phi)\circ\lambda'_{i_si_l}
$$
 for any $\phi\in\Gamma(\tilde u_{i_s})$, then there exist
  $\phi_{i_s}\in\Gamma(\tilde u_{i_s})$ and
 $\phi_{i_l}\in\Gamma(\tilde u_{i_l})$ such that
$$
\lambda_{i_si_l}=\phi_{i_l}\circ\lambda'_{i_si_l}\quad{\rm
and}\quad \lambda_{i_si_l}=\lambda'_{i_si_l}\circ\phi_{i_s}.
$$
 In particular these imply that
$$
\lambda_{i_li_s}\circ\lambda_{i_si_l}=\phi_{i_s}\quad{\rm
and}\quad
\lambda_{i_ti_s}\circ\lambda_{i_li_t}\circ\lambda_{i_si_l}=\phi'_{i_s}.
$$
 for some $\phi_{i_s},\phi'_{i_s}\in\Gamma(\tilde
u_{i_s})$ if $s\ne l, t$ and $l\ne t$. For convenience we make the
convention that $\lambda_{i_si_s}=id_{O(\tilde u_{i_s})}$ and
${\mathcal A}_{i_si_s}=1_{\Gamma(\tilde u_{i_s})}$. Then we have
$|\Gamma(\tilde u_{i_s})|^{k-1}$ PS embeddings from $O(\tilde
u_{i_s})$ into $\widetilde W^{\Gamma_I}_I$ given by
$$
 \phi_I\circ\lambda_I^s:\tilde
u_s\mapsto \bigl(\phi_{i}\circ\lambda_{i_si}(\tilde
u_s)\bigr)_{i\in I}, \eqno(3.25)
$$
 where $\lambda_I^s=(\lambda_{i_si})_{i\in I}$ with
 $\lambda_{i_si_s}=id_{O(\tilde u_{i_s})}$,
 and $\phi_I=(\phi_{i})_{i\in I}$ belongs to
$$
 \Gamma(\tilde
u_I)_s:=\Bigl\{(\phi_{i})_{i\in I}\in\prod^k_{l=1}\Gamma(\tilde
u_{i_l})\,\bigm|\, \phi_{i_s}=1\Bigr\}.\eqno(3.26)
$$
 Hereafter saying $\phi_I\circ\lambda_I^s$ to
be a PS embedding from $O(\tilde u_{i_s})$ into $\widetilde
W^{\Gamma_I}_I$ means that $\phi_I\circ\lambda_I^s$ is a PS
embedding from $ O(\tilde u_{i_s})$ into the PS Banach manifold
$\prod^k_{l=1}\widetilde W_{i_l}$ and that the image
$\phi_I\circ\lambda_I^s(O(\tilde u_{i_s}))$ is in $\widetilde
W^{\Gamma_I}_I$ if $I=\{i_1,\cdots, i_k\}$.  Furthermore, $(\tilde
u_{i_1},\cdots,\tilde u_{i_k})\in\prod^k_{l=1}O(\tilde
u_{i_l})\subset \prod^k_{l=1}\widetilde W_{i_l}$ is contained in
$\widetilde W_I^{\Gamma_I}$ if and only if
$$
(\tilde u_{i_1},\cdots, \tilde
u_{i_k})=\phi_I\circ\lambda_I^s(\tilde u_{i_s})$$
 for some
$\phi_I\in\Gamma(\tilde u_I)_s$. So the neighborhood
$$
 \widetilde O(\tilde
u_I):=\pi_I^{-1}(O)=\Bigl(\prod^k_{l=1}O(\tilde
u_{i_l})\Bigr)\cap\widetilde W_I^{\Gamma_I}
$$
 of  $\tilde u_I$ in $\widetilde W_I^{\Gamma_I}$
can be identified with an union of $|\Gamma(\tilde
u_{i_s})|^{k-1}$ copies of $O(\tilde u_{i_s})$,
$$
 \bigcup_{\phi_I\in\Gamma(\tilde
u_I)_s}\phi_I\circ\lambda_I^s \bigl(O(\tilde
u_{i_s})\bigr).\eqno(3.27)
$$
 Clearly, for any two different elements $\phi_I,\phi'_I$ in
$\Gamma(\tilde u_I)_s$ we have
$$
\begin{array}{lc}
\phi_I\circ\lambda_I^s\bigl(O(\tilde
u_{i_s})\bigr)\cap\phi'_I\circ\lambda_I^s\bigl(O(\tilde
u_{i_s})\bigr)=\vspace{2mm}\\
\bigl\{(\phi_i\circ\lambda_{i_si}(\tilde v))_{i\in I}\,\bigm|\,
\tilde v\in O(\tilde u_{i_s})\;\&\;
 (\phi_i\circ\lambda_{i_si}(\tilde
v))_{i\in I}=(\phi'_i\circ\lambda_{i_si}(\tilde v))_{i\in
I}\bigr\}\supseteq\{\tilde u_I\}.
\end{array}
\eqno(3.28)
$$
In order to show that $\widetilde W_I^{\Gamma_I}$ is a stratified
Banach variety, we need to prove that these $|\Gamma(\tilde
u_{i_s})|^{k-1}$ sets, $\phi_I\circ\lambda_I^s(O(\tilde
u_{i_s})),\,\phi_I\in\Gamma(\tilde u_I)_s$ are intrinsic in the
following sense that
$$
\bigl\{\phi_I\circ\lambda_I^s\bigl(O(\tilde
u_{i_s})\bigr)\,\bigm|\, \phi_I\in\Gamma(\tilde
u_I)_s\bigr\}=\bigl\{\phi_I\circ\lambda_I^t\bigl(O(\tilde
u_{i_t})\bigr)\,\bigm|\, \phi_I\in\Gamma(\tilde u_I)_t\bigr\}
\eqno(3.29)
$$
for any two $s, t$ in $\{1,\cdots,k\}$. This is easily proved,
cf., [LuT].

 Similarly, (by shrinking
$O$ if necessary) using the properties of orbifold bundles one has
also the PS bundle isomorphisms
$$
\Lambda_{i_si_l}:\widetilde E_{i_s}|_{O(\tilde
u_{i_s})}\to\widetilde E_{i_l}|_{O(\tilde u_{i_l})}\eqno(3.30)
$$
 that are lifting of $\lambda_{i_si_l}$, i.e. $\lambda_{i_si_l}\cdot\tilde
p_{i_s}=\tilde p_{i_l}\cdot\Lambda_{i_si_l}$, to satisfy
$$
\Lambda_{i_si_l}\circ\Phi=\boxed{{\mathcal
A}_{i_si_l}(\phi)}\circ\Lambda_{i_si_l}\eqno(3.31)
$$
 for any $\phi\in\Gamma(\tilde u_{i_s})$ and $l\in\{1,\cdots, k\}\setminus\{s\}$.
Hereafter
$$
\Phi: \widetilde E_{i_s}|_{O(\tilde u_{i_s})}\to\widetilde
E_{i_s}|_{O(\tilde u_{i_s})}\; \bigl({\rm resp.}\;  \boxed{{\cal
 A}_{i_si_l}(\phi)}: \widetilde
E_{i_l}|_{O(\tilde u_{i_l})}\to\widetilde E_{i_l}|_{O(\tilde
u_{i_l})}\bigr)\eqno(3.32)
$$
is the PS bundle isomorphism lifting of $\phi$ (resp. ${\cal
A}_{i_si_l}(\phi)$) produced in the definition of orbifold
bundles. ({\it Here we understand $\phi_{i_l}\in\Gamma_{i_l}$ to
acts on $\widetilde E_{i_l}$ via
$\phi_{i_l}\cdot\tilde\xi=\Phi_{i_l}(\tilde\xi)$ with
$\Phi_{i_l}=\lambda_{W_{i_l}}(\phi_{i_l})$ the lifting bundle
automorphism of $\phi_{i_l}$ and $\tilde\xi\in\widetilde
E_{i_l}$}.)
 It is easy to see that $\Phi_{i_l}$ maps $\widetilde
E_{i_l}|_{O(\tilde u_{i_l})}$ to $\widetilde E_{i_l}|_{O(\tilde
u_{i_l})}$. In particular the restriction of it to $\widetilde
E_{i_l}|_{O(\tilde u_{i_l})}$ is also a PS bundle automorphism of
$\widetilde E_{i_l}|_{O(\tilde u_{i_l})}$ that is a lift of
$\phi_{i_l}$, still denoted by $\Phi_{i_l}$. So for the above
$\lambda_I^s=(\lambda_{i_si})_{i\in I}$ with
 $\lambda_{i_si_s}=id_{O(\tilde u_{i_s})}$,
 and $\phi_I=(\phi_{i})_{i\in I}\in\Gamma(\tilde u_I)_s$ we get
 the corresponding
$$
 \Phi_I:=(\Phi_{i})_{i\in
I}\quad{\rm and}\quad\Lambda_I^s:=(\Lambda_{i_si})_{i\in
I},\eqno(3.33)
$$
 where we have made the convention that
$$
 \Phi_{i_s}=1\quad{\rm and}\quad
\Lambda_{i_si_s}=id_{\widetilde E_{i_s}|_{O(\tilde
u_{i_s})}}.\eqno(3.34)
$$
 Corresponding with (3.25) we obtain $|\Gamma(\tilde
u_{i_s})|^{k-1}$ PS bundle embeddings
$$
 \Phi_I\circ\Lambda_I^s:\widetilde
E_{i_s}|_{O(\tilde u_{i_s})}\to\widetilde E_I^{\Gamma_I},\; \tilde
\xi\mapsto \bigl(\Phi_i\circ\Lambda_{i_si}(\tilde \xi)\bigr)_{i\in
I}\eqno(3.35)
$$
for $\phi_I\in\Gamma(\tilde u_I)_s$.
  Clearly, each map
\begin{eqnarray*}
&&  \Phi_I\circ\Lambda_I^s(\widetilde E_{i_s}|_{O(\tilde
u_{i_s})})\to \phi_I\circ\lambda_I^s(O(\tilde u_{i_s})),\\
&&\bigl(\Phi_i\circ\Lambda_{i_si}(\tilde \xi)\bigr)_{i\in
I}\mapsto \bigl(\phi_i\circ\lambda_{i_si}(\tilde p_{i_s}(\tilde
\xi))\bigr)_{i\in I}
 \end{eqnarray*}
 gives a stratified Banach vector bundle over the stratified Banach
manifold $\phi_I\circ\lambda_I^s(O(\tilde u_{i_s}))$, and
 $$
 \bigcup_{\phi_I\in\Gamma(\tilde u_I)_s}
 \Phi_I\circ\Lambda_I^s(\widetilde E_{i_s}|_{O(\tilde u_{i_s})})=\widetilde
E_I^{\Gamma_I}\Bigm|_{\bigcup_{\phi_I\in\Gamma(\tilde
u_I)_s}\phi_I\circ\lambda_I^s(O(\tilde u_{i_s}))}.
$$
   Corresponding with (3.28) and
 (3.29) we have also:
$$
\begin{array}{lc}
\Phi_I\circ\Lambda_I^s(\widetilde E_{i_s}|_{O(\tilde
u_{i_s})})\bigcap\Phi'_I\circ\Lambda_I^s(\widetilde
E_{i_s}|_{O(\tilde u_{i_s})})=\vspace{2mm}\\
\bigl\{\bigl(\Phi_i\circ\Lambda_{i_si}(\tilde \xi)\bigr)_{i\in
I}\,\bigm|\, \tilde \xi\in\widetilde E_{i_s}|_{O(\tilde
u_{i_s})}\;\&\; \bigl(\Phi_i\circ\Lambda_{i_si}(\tilde
\xi)\bigr)_{i\in I}=\bigl(\Phi'_i\circ\Lambda_{i_si}(\tilde
\xi)\bigr)_{i\in I}\bigr\}.
\end{array}\eqno(3.36)
$$
for any two different $\phi_I,\psi_I$ in $\Gamma(\tilde u_I)_s$,
and
$$
\bigl\{\Phi_I\circ\Lambda_I^s(\widetilde E_{i_s}|_{O(\tilde
u_{i_s})})\,\bigm|\, \phi_I\in\Gamma(\tilde
u_I)_s\bigr\}=\bigl\{\Phi_I\circ\Lambda_I^t(\widetilde
E_{i_t}|_{O(\tilde u_{i_t})})\,\bigm|\, \phi_I\in\Gamma(\tilde
u_I)_t\bigr\}\eqno(3.37)
 $$
 for any  $s, t\in\{1,\cdots,k\}$.
{\it Therefore $\widetilde E_I^{\Gamma_I}\to\widetilde
W_I^{\Gamma_I}$ is
 a union of $|\Gamma(\tilde u_{i_s})|^{k-1}$ stratified Banach
 vector bundles near $\tilde u_I$.}\vspace{2mm}

\noindent{\bf Remark 3.6.}\quad By shrinking $O$ we may also
assume that  each $\lambda_{i_si_l}$ in (3.24) is not only a PS
diffeomorphism from $O(\tilde u_{i_s})$ to $O(\tilde u_{i_l})$,
but also one from a neighborhood of the closure of $O(\tilde
u_{i_s})$ to that of $O(\tilde u_{i_l})$. In this case the union
in (3.27) may be required to satisfy:
$$
 \left\{\begin{array}{ll} {\rm
If}\; \tilde v\in\bigcup_{\phi_I\in\Gamma(\tilde
u_I)_s}\phi_I\circ\lambda_I^s(O(\tilde u_{i_s}))\;{\rm
and}\;\tilde v\notin \phi_I\circ\lambda_I^t(O(\tilde u_{i_t}))\;{\rm for}\;{\rm some}\,t,\\
{\rm then} \; \tilde v\notin
Cl\bigr(\phi_I\circ\lambda_I^t(O(\tilde u_{i_t}))\bigl).
\end{array}\right.\eqno(3.38)
$$
 In fact, if $\tilde v\in
Cl\bigr(\phi_I\circ\lambda_I^t(O(\tilde
u_{i_t}))\bigl)\setminus\phi_I\circ\lambda_I^t(O(\tilde u_{i_t}))$
then there exists a unique $\tilde x\in Cl(O(\tilde u_{i_t}))$,
which must sit in $Cl(O(\tilde u_{i_t}))\setminus O(\tilde
u_{i_t})$, such that $\phi_I\circ\lambda_I^t(\tilde x)=\tilde v$.
Moreover, since $\tilde v$ belongs to some
$\phi'_I\circ\lambda_I^s(O(\tilde u_{i_s}))$ one has a unique
$\tilde y\in Cl(O(\tilde u_{i_s}))$, which must sit in $O(\tilde
u_{i_s})$, such that $\phi'_i\circ\lambda_I^s(\tilde y)=\tilde v$.
So $\phi_I\circ\lambda_I^t(\tilde
x)=\phi'_I\circ\lambda_I^s(\tilde y)$. For the sake of clearness
we assume $t=1$ and $s=2$ then
$$
\bigl(\tilde x, \phi_{i_2}\circ\lambda_{i_1i_2}(\tilde x),\cdots,
\phi_{i_k}\circ\lambda_{i_1i_k}(\tilde x)\bigr)=
\bigl(\phi'_{i_1}\circ\lambda_{i_2i_1}(\tilde y), \tilde
y,\phi'_{i_3}\circ\lambda_{i_3i_1}(\tilde y),\cdots,
\phi'_{i_k}\circ\lambda_{i_2i_k}(\tilde x)\bigr).
$$
Therefore $\tilde x=\phi'_{i_1}\circ\lambda_{i_2i_1}(\tilde y)$.
Note that both PS diffeomorphisms $\phi'_{i_1}$ and
$\lambda_{i_2i_1}$ map $O(\tilde u_{i_2})$ to $O(\tilde u_{i_1})$.
We deduce that $\tilde x=\phi'_{i_1}\circ\lambda_{i_2i_1}(\tilde
y)\in O(\tilde u_{i_1})$. This contradicts that $\tilde x\in
Cl(O(\tilde u_{i_1}))\setminus O(\tilde u_{i_1})$. Equation (3.38)
is proved. Similarly, we also require that
$$
\left\{\begin{array}{ll} {\rm If}\; \tilde
\xi\in\bigcup_{\Phi_I\in\Gamma(\tilde
u_I)_s}\Phi_I\circ\Lambda_I^s(O(\tilde u_{i_s}))\;{\rm
and}\;\tilde \xi\notin \Phi_I\circ\Lambda_I^t(O(\tilde
u_{i_t}))\;{\rm for}\;{\rm some}\,t,\\
{\rm then} \; \tilde \xi\notin
Cl\bigr(\Phi_I\circ\Lambda_I^t(O(\tilde u_{i_t}))\bigl).
\end{array}\right.\eqno(3.39)
$$
 Later we always assume that (3.38) and
(3.39) hold without special statements.

In terms of [LiuT3] we introduce:\vspace{2mm}

\noindent{\bf Definition 3.7.}\quad{\it The family of the PS
embeddings  given by (3.25),
 $\bigl\{\phi_I\circ\lambda_I^s\,|\, \phi_I\in\Gamma(\tilde u_I)_s\bigr\}$,
  is called a
{\bf local coordinate chart of} $\widetilde W^{\Gamma_I}_I$ over a
neighborhood $\widetilde O(\tilde u_I)$ of $\tilde u_I$, each
$\phi_I\circ\lambda_I^s$ a {\bf component} of this chart.
Similarly, we call the family of the PS bundle embeddings given by
(3.35), $\bigl\{\Phi_I\circ\Lambda_I^s\,|\, \phi_I\in\Gamma(\tilde
u_I)_s\bigr\}$,  a {\bf local bundle coordinate chart of
$\widetilde E^{\Gamma_I}_I$ over $\widetilde O(\tilde u_I)$}, each
$\Phi_I\circ\Lambda_I^s$ a {\bf component} of it. For $\tilde
x_I\in\phi_I\circ\lambda_I^s(O(\tilde u_{i_s}))$, a connected
relative open subset $W\subset\phi_I\circ\lambda_I^s(O(\tilde
u_{i_s}))$ containing $\tilde u_I$ is called a {\bf local
component of $\widetilde W^{\Gamma_I}_I$ near} $\tilde x_I$. In
particular, $\phi_I\circ\lambda_I^s(O(\tilde u_{i_s}))$ is a local
component of $\widetilde W^{\Gamma_I}_I$ near $\tilde u_I$. The
restriction of $\Phi_I\circ\Lambda_I^s(\widetilde
E_{i_s}|_{O(\tilde u_{i_s})})$ to a  local component near $\tilde
u_I$ is called a {\bf local component of $\widetilde
E^{\Gamma_I}_I$ near $\tilde u_I$}. Two local components of
$\widetilde W^{\Gamma_I}_I$ (or $\widetilde E^{\Gamma_I}_I$) near
a point $\tilde u_I$ of $\widetilde W^{\Gamma_I}_I$ are said to
belong to different kinds if the intersection of both is not a
local component of $\widetilde W^{\Gamma_I}_I$ near $\tilde u_I$.
Let ${\mathcal A}$ be a family of local components near a point
$\tilde u_I\in \widetilde W^{\Gamma_I}_I$.
 If the union of sets in ${\mathcal A}$ forms an open neighborhood
${\mathcal O}(\tilde u)$ of $\tilde u_I$ in $\widetilde
W^{\Gamma_I}_I$ we call ${\mathcal A}$ {\bf a complete family of
local components of $\widetilde W^{\Gamma_I}_I$ over ${\mathcal
O}(\tilde u_I)$}.}\vspace{2mm}

Equations (3.29) and (3.37) show that we construct the same
complete families of local components starting from two different
$s$ and $t$. Moreover, suppose that $Q\subset O$ is another small
open neighborhood of $u_I$ in $W_I$. Then all above constructions
work if we replace $O$ with $Q$, and the corresponding local
components of $\widetilde W^{\Gamma_I}_I$ near $\tilde u_I$ and
those of $\widetilde E^{\Gamma_I}_I$ near $(\widetilde
E^{\Gamma_I}_I)_{\tilde u_I}$ are, respectively, given by
$$
\bigl\{\phi_I\circ\lambda_I^s( Q(\tilde u_{i_s}))\,\bigm|\,
\phi_I\in\Gamma(\tilde u_I)_s\bigr\}\quad{\rm and}\quad
\bigl\{\Phi_I\circ\Lambda_I^s(\widetilde E_{i_s}|_{Q(\tilde
u_{i_s})})\,\bigm|\, \phi_I\in\Gamma(\tilde
u_I)_s\bigr\}.\eqno(3.40)
$$
 Here  $Q(\tilde u_{i_l})=\pi_{i_l}^{-1}(Q)$ are
the inverse images of $Q$ in $\widetilde W_{i_l}$, $l=1,\cdots,k$,
and $\phi_I\circ\lambda_I^s$ and $\Phi_I\circ\Lambda_I^s$ as in
(3.25) and (3.35). As stated in Lemma 4.3 of [LiuT3] we have:
\vspace{2mm}

\noindent{\bf Proposition 3.8.}\quad{\it The notion of local
component is functorial with respect to restrictions and
projections, that is:
\begin{description}
\item[(i)] If $\{\phi_I\circ\lambda_I^s\,|\,\phi_I\in\Gamma(\tilde
u_I)_s\}$ is a local coordinate chart near $\tilde
u_I\in\widetilde W_I^{\Gamma_I}$, then for each point $\tilde
y_I\in\cup_{\phi_I\in\Gamma(\tilde
u_I)_s}\phi_I\circ\lambda_I^s(O(\tilde u_{i_s}))$ there exists a
neighborhood ${\mathcal Q}$ of $\tilde y_I$ in $\widetilde
W_I^{\Gamma_I}$ such that $\{{\mathcal Q}\cap
\phi_I\circ\lambda_I^s(O(\tilde
u_{i_s}))\,|\,\phi_I\in\Gamma(\tilde u_I)_s\}$ is a complete
family of local components of $\widetilde W^{\Gamma_I}_I$ over
${\mathcal Q}$.

\item[(ii)] If
$\{\phi_I\circ\lambda_I^s\,|\,\phi_I\in\Gamma(\tilde u_I)_s\}$ is
a local coordinate chart near $\tilde u_I\in\widetilde
W_I^{\Gamma_I}$, $i_s\in J\hookrightarrow I$, then
 $\{\tilde\pi^I_J\circ\phi_I\circ\lambda_I^s\,|\,\phi_I\in\Gamma(\tilde
u_I)_s\}$ is a local coordinate chart near $\tilde\pi^I_J(\tilde
u_I)\in\widetilde W_J^{\Gamma_J}$ (after  the repeating maps
$\tilde\pi^I_J\circ\phi_I\circ\lambda_I^s$ is only taken a time).
In particular, not only the projection $\tilde\pi^I_J:\widetilde
W_I^{\Gamma_I}\to\widetilde W_J^{\Gamma_J}$ maps  the local
component $\phi_I\circ\lambda_I^s(O(\tilde u_{i_s}))$ near $\tilde
u_I$ to the local component
$\tilde\pi^I_J(\phi_I\circ\lambda_I^s(O(\tilde u_{i_s})))$ near
$\tilde\pi^I_J(\tilde u_I)$ in $\widetilde W_J^{\Gamma_J}$ but
also the restriction
$$\tilde\pi^I_J|_{\phi_I\circ\lambda_I^s(O(\tilde u_{i_s}))}:
\phi_I\circ\lambda_I^s(O(\tilde u_{i_s}))\to
\tilde\pi^I_J(\phi_I\circ\lambda_I^s(O(\tilde u_{i_s})))
$$
is a PS map. Actually one has
$\tilde\pi^I_J\circ\phi_I\circ\lambda_I^s=\lambda^I_J(\phi_I)\circ\lambda_J^s
=\phi_J\circ\lambda_J^s$, where
$\lambda_J^s=(\lambda_{i_sl})_{l\in J}$.

\item[(iii)] Correspondingly, the restriction of the projection
$\widetilde\Pi^I_J:\widetilde E_I^{\Gamma_I}\to\widetilde
E_J^{\Gamma_J}$ to $\widetilde
E_I^{\Gamma_I}|_{\phi_I\circ\lambda_I^s(O(\tilde u_{i_s}))}$ is a
PS bundle map from  $\widetilde
E_I^{\Gamma_I}|_{\phi_I\circ\lambda_I^s(O(\tilde u_{i_s}))}$ to
$\widetilde
E_J^{\Gamma_J}|_{\pi^I_J(\phi_I\circ\lambda_I^s(O(\tilde
u_{i_s})))}$.
\end{description}}

\noindent{\it Proof.} (i) By (3.38) we may choose a connected open
neighborhood ${\mathcal Q}\subset\pi_I^{-1}(O)$ of $\tilde y_I$ in
$\widetilde W_I^{\Gamma_I}$ such that for any
$\phi_I\in\Gamma(\tilde u_I)_s$,
$${\mathcal Q}\cap Cl(\phi_I\circ\lambda_I^s(O(\tilde
u_{i_s})))=\emptyset\quad {\rm as}\quad \tilde
y_I\notin\phi_I\circ\lambda_I^s(O(\tilde u_{i_s})).$$
 If $y_I\in\phi_I\circ\lambda_I^s(O(\tilde
u_{i_s}))$ then the connected relative
 open subset ${\mathcal Q}\cap \phi_I\circ\lambda_I^s(O(\tilde
u_{i_s}))$ in $\phi_i\circ\lambda_I^s(O(\tilde u_{i_s}))$ is a
local component of $\widetilde W^{\Gamma_I}_I$ near $\tilde v_I$.
It is easy to see that the union of sets in $\{{\mathcal Q}\cap
\phi_I\circ\lambda_I^s(O(\tilde
u_{i_s}))\,|\,\phi_I\in\Gamma(\tilde u_I)_s\}$ is equal to
${\mathcal Q}$. So (i) holds.

(ii) We may assume that $I=\{i_1,\cdots,i_k\}$ and
$J=\{i_1,\cdots,i_{k'}\}$ with $k'<k$. By (3.39) we may also take
$s=1$. Then by (3.26) we have
\begin{eqnarray*}
&&\Gamma(\tilde
u_I)_1=\bigl\{(1,\phi_{i_2},\cdots,\phi_{i_k})\,\bigm|\,
\phi_{i_l}\in\Gamma(\tilde u_{i_l}), l=2,\cdots,k\bigr\},\nonumber\\
&&\Gamma(\tilde\pi^I_J(\tilde
u_I))_1=\bigl\{(1,\phi_{i_2},\cdots,\phi_{i_{k'}})\,\bigm|\,
\phi_{i_l}\in\Gamma(\tilde u_{i_l}), l=2,\cdots,k'\bigr\}.
\end{eqnarray*}
That is, $\lambda^I_J(\Gamma(\tilde
u_I)_1)=\Gamma(\tilde\pi^I_J(\tilde u_I))_1$.  So for a given
local component
$$
\phi_I\circ\lambda_I^1(O(\tilde u_{i_1}))=\bigl\{ \bigl(\tilde x,
\phi_{i_2}\circ\lambda_{i_1i_2}(\tilde x),\cdots,
\phi_{i_k}\circ\lambda_{i_1i_k}(\tilde x)\bigr)\,\bigm|\,\tilde
x\in O(\tilde u_{i_1})\bigr\}
$$
of $\widetilde W_I^{\Gamma_I}$ near $\tilde u_I$, we have a
corresponding local component of $\widetilde W_J^{\Gamma_J}$ near
$\tilde\pi^I_J(\tilde u_I)$ as follow:
$$
 \phi_J\circ\lambda_J^1(O(\tilde u_{i_1}))=\bigl\{ \bigl(\tilde
x_1, \phi_{i_2}\circ\lambda_{i_1i_2}(\tilde x),\cdots,
\phi_{i_k'}\circ\lambda_{i_1i_{k'}}(\tilde
x)\bigr)\,\bigm|\,\tilde x\in O(\tilde u_{i_1})\bigr\}.
$$
Here
$\phi_J=(1,\phi_{i_2},\cdots,\phi_{i_{k'}})=\lambda^I_J(\phi_I)\in\Gamma(\tilde\pi^I_J(\tilde
u_I))_1$ and $\lambda_J^1=(\phi_{i_1l})_{l\in J}$.
 Clearly,
 $\tilde\pi^I_J\circ\phi_I\circ\lambda_I^1=\lambda^I_J(\phi_I)\circ\lambda_J^1=\phi_J\circ\lambda_J^1$,
 and so   $\tilde\pi^I_J\bigl(\phi_I\circ\lambda_I^1(O(\tilde
u_{i_1}))\bigr)=\phi_J\circ\lambda_J^1(O(\tilde u_{i_1}))$.
Moreover the restriction
$$
\tilde\pi^I_J|_{\phi_I\circ\lambda_I^1(O(\tilde
u_{i_1}))}:\phi_I\circ\lambda_I^1(O(\tilde u_{i_1}))\to
\phi_J\circ\lambda_J^1(O(\tilde u_{i_1}))
$$
given by
$$\bigl(\tilde x,
\phi_{i_2}\circ\lambda_{i_1i_2}(\tilde x),\cdots,
\phi_{i_k}\circ\lambda_{i_1i_k}(\tilde x)\bigr)\mapsto
\bigl(\tilde x, \phi_{i_2}\circ\lambda_{i_1i_2}(\tilde x),\cdots,
\phi_{i_{k'}}\circ\lambda_{i_1i_{k'}}(\tilde x)\bigr)
$$
is a PS map. (ii) is proved.
\hfill$\Box$\vspace{2mm}

The arguments below Definition 3.7 and Proposition 3.8 show that
the notion of the local component is intrinsic. Now for each
$I=\{i_1,\cdots,i_k\}\in{\mathcal N}$ with $|I|>1$, we can
desingularize $\widetilde W_I^{\Gamma_I}$ to get a true stratified
Banach manifold $\widehat W_I^{\Gamma_I}$.  Consider the
stratified Banach manifold
$$
\boxed{\widetilde W_I^{\Gamma_I}\!}:= \coprod_{\tilde
u_I\in\widetilde W_I^{\Gamma_I}}\coprod_{O\in{\mathcal
U}(u_I)}\coprod_{1\le s\le |I|}\coprod_{\phi_I\in\Gamma(\tilde
u_I)_s}\phi_I\circ\lambda_I^s(O(\tilde u_{i_s}))\eqno(3.41)
$$
where ${\mathcal U}(u_I)$ is the germ of small connected open
neighborhood of $u_I=\pi_I(\tilde u_I)$ in $W_I$  and $O(\tilde
u_{i_s})=\pi^{-1}_{i_s}(O)$, $s=1,\cdots,|I|$. We define the
equivalence relation $\sim$ in $\boxed{\widetilde
W_I^{\Gamma_I}\!}$ as follows. For $\tilde u_I,\tilde
u'_I\in\widetilde W_I^{\Gamma_I}$,
 and $\phi_I\in\Gamma(\tilde u_I)_s$, $\phi'_I\in\Gamma(\tilde u'_I)_t$,
 and $\tilde y_I\in\phi_I\circ\lambda_I^s(O(\tilde u_{i_s}))$
 and $\tilde y_I^\prime\in\phi_I^\prime\circ\lambda_I^{t\prime}(O'(\tilde u_{i_t}))$
 we define
$$(\tilde u_I,\phi_I, \tilde y_I)\sim(\tilde u'_I,\phi'_I,
 \tilde y'_I)\eqno(3.42)$$
 if and only if the following two conditions are satisfied:
\begin{description}
\item[(i)] $\tilde y_I=\tilde y'_I$ as points in $ \widetilde
W_I^{\Gamma_I}$;

\item[(ii)] $\phi_I\circ\lambda_I^s(O(\tilde u_{i_s}))
\bigcap\phi'_I\circ\lambda_I^{\prime t}(O'(\tilde u_{i_t}))$ is a
local component of $\widetilde W_I^{\Gamma_I}$ near $\tilde
y_I=\tilde y'_I$.
\end{description}
Clearly, any two different points in $\{(\tilde
u_I,\phi_I)\}\times\phi_I\circ\lambda_I^s(O(\tilde u_{i_s}))$ are
not equivalent with respect to $\sim$. Denote by $[\tilde
u_I,\phi_I, \tilde y_I]$ the equivalence class of $(\tilde
u_I,\phi_I, \tilde y_I)$, and by
$$
\widehat W_I^{\Gamma_I}=\boxed{\widetilde
W_I^{\Gamma_I}\!}\bigg/\sim
$$
 if $|I|>1$, and by $\widehat W_I^{\Gamma_I}=\widetilde W_i$ if
 $I=\{i\}\in{\mathcal N}$. Then $\widehat
 W_I^{\Gamma_I}$ is a stratified Banach manifold. Note that (3.26)
 implies that each point in $\{(\tilde
u_I,\phi_I)\}\times\phi_I\circ\lambda_I^s(O(\tilde u_{i_s}))$ with
$s>1$ must be equivalent to some point in $\{(\tilde
u_I,\phi'_I)\}\times\phi'_I\circ\lambda_I^1(O(\tilde u_{i_1}))$
for some $\phi'_I\in\Gamma(\tilde u_I)_1$ with respect to $\sim$.
So
$$
\widehat W_I^{\Gamma_I}=\coprod_{\tilde u_I\in\widetilde
W_I^{\Gamma_I}}\coprod_{O\in{\mathcal
U}(u_I)}\coprod_{\phi_I\in\Gamma(\tilde
u_I)_1}\phi_I\circ\lambda_I^1(O(\tilde u_{i_1}))\bigg/\sim
\eqno(3.43)
$$
 if $|I|>1$. Now for each $s=1,\cdots,k$,
$$
 \widehat{\phi_I\circ\lambda_I^s}: O(\tilde u_{i_s})\to \widehat
W_I^{\Gamma_I},\;\tilde x\mapsto [\tilde u_I,\phi_I,
\phi_I\circ\lambda_I^s(\tilde x)],\eqno(3.44)
$$
gives a PS open embedding, called {\it a local coordinate chart of
$\widehat W_I^{\Gamma_I}$ near $[\tilde u_I,\phi_I,\tilde u_I]$}.
The projections $\tilde\pi_I$ and $\tilde\pi^I_J$ induce natural
ones
$$
\hat\pi_I: \widehat W_I^{\Gamma_I}\to W_I,\quad [\tilde
u_I,\phi_I, \tilde y_I]\mapsto \tilde\pi_I(\tilde y_I),$$
$$
\hat\pi^I_J:\widehat W_I^{\Gamma_I}\to \widehat
W_J^{\Gamma_J},\;[\tilde u_I,\phi_I, \tilde y_I]\mapsto [\tilde
u_J,\phi_J, \tilde y_J],\eqno(3.45)
$$
where $\phi_J=\lambda^I_J(\phi_I)=(\phi_l)_{l\in J}$, $\tilde
u_J=\tilde\pi^I_J(\tilde u_I)=(\tilde u_l)_{l\in J}$ and $\tilde
y_J= \tilde\pi^I_J(\tilde y_I)=(\tilde y_l)_{l\in J}$.
 By Proposition 3.8(ii) the map $\hat\pi^I_J$ is
well-defined, and in the charts $\widehat{\phi_I\circ\lambda_I^1}$
and
$$ \widehat{\phi_J\circ\lambda_J^1}: O(\tilde u_{j_1})\to
\widehat W_J^{\Gamma_J},\;\tilde x\mapsto [\tilde u_J,\phi_J,
\phi_J\circ\lambda_J^1(\tilde x)],
$$
 the projection $\hat\pi^I_J$ may be represented by
$$
  O(\tilde u_{i_1})\to O(\tilde u_{j_1}),\; \tilde x\mapsto
\lambda_{i_1j_1}(\tilde x)
$$
because $j_1\in J\subset I$. Here $\lambda_{i_1j_1}=1_{ O(\tilde
u_{i_1})}$ if $j_1=i_1$.
 This shows that $\hat\pi^I_J$ {\it is a PS map}. As
 $\tilde\pi^I_J$ in (3.23),  $\hat\pi^I_J$ is not surjective in general if $I\ne
 J$.

Note that the action of $\Gamma_I$ on $\widetilde W_I^{\Gamma_I}$
in (3.22) induces a natural one on  $\boxed{\widetilde
W_I^{\Gamma_I}\!}\,$. For any $\psi_I=(\psi_l)_{l\in
I}\in\Gamma_I$ and $\tilde y_I=\phi_I\circ\lambda_I^1(\tilde x)\in
\phi_I\circ\lambda_I^1(O(\tilde u_{i_1}))$ we have
$$
\psi_I\cdot\tilde y_I= \psi_I(\tilde
y_I)=\phi^\ast_I\circ\lambda_I^{1\ast}(\psi_{i_1}(\tilde x))\in
\phi^\ast_I\circ\lambda_I^{1\ast}(\psi_{i_1}(O(\tilde u_{i_1}))),
$$
where $\phi^\ast_I=(\phi^\ast_l)_{l\in I}\in\Gamma(\psi_I(\tilde
u_I))_1 $, $\phi^\ast_l= \psi_l\circ\phi_l\circ\psi_l^{-1}$ for
$l\in I$,  and $\lambda^{1\ast}_I=(\lambda^\ast_{i_1l})_{l\in I}$,
$\lambda^\ast_{i_1l}=\psi_l\circ\lambda_{i_1l}\circ\psi_{i_1}^{-1}$
for any $l\in I$. So one may define
$$
\psi_I\cdot (\tilde u_I,\phi_I, \tilde y_I)=\bigl(\psi_I(\tilde
u_I), \phi^\ast_I, \psi_I(\tilde y_I)\bigr).\eqno(3.46)
$$
Clearly, this gives a PS action of $\Gamma_I$ on the space
$\boxed{\widetilde W_I^{\Gamma_I}\!}\,$. Note that the above
$\phi^\ast_I$ belongs to $\Gamma(\psi_I(\tilde u_I))_s$ as
$\phi_I\in\Gamma(\tilde u_I)_s$. Moreover, for
$s\in\{1,\cdots,k\}$ and $l\in I$, if we define
$\lambda^{s\ast}_I=(\lambda^\ast_{i_sl})_{l\in I}$,
$\lambda^\ast_{i_sl}=\psi_l\circ\lambda_{i_sl}\circ\psi_{i_s}^{-1}$,
then
$$
\phi^\ast_I\circ\lambda_I^{s\ast}:\psi_{i_s}(O(\tilde u_{i_s}))\to
\widetilde W_I^{\Gamma_I},\; \tilde z_s\mapsto
\bigl(\phi_l^\ast\circ\lambda^\ast_{i_sl}(\tilde z_s)\bigr)_{l\in
I}
$$
is a PS open embedding, and
$\{\phi_I^\ast\circ\lambda_I^{s\ast}\,|\,\phi_I\in\Gamma(\tilde
u_I)_s\}$ is a local coordinate chart of $\widetilde
W_I^{\Gamma_I}$ over the neighborhood $\psi_I\bigl( O(\tilde
u_I)\bigr)$ of $\psi_I(\tilde u_I)$ in the sense of Definition
3.7. In particular, we get a local coordinate chart of $\widehat
W_I^{\Gamma_I}$ near $[\psi_I(\tilde
u_I),\phi^\ast_I,\psi_I(\tilde u_I)]$,
$$
\widehat{\phi^\ast_I\circ\lambda_I^{\ast 1}}:\psi_{i_s}(O(\tilde
u_{i_s}))\to \widehat W_I^{\Gamma_I},\;\tilde x\mapsto
[\psi_I(\tilde u_I),\phi^\ast_I, \phi^\ast_I\circ\lambda_I^{\ast
s}(\tilde x)].\eqno(3.47)
$$
From these and (3.43) it follows that the action in (3.46)
preserves the equivalence relation $\sim$ in $\boxed{\widetilde
W_I^{\Gamma_I}\!}\,$, i.e.
 $$(\tilde u_I,\phi_I, \tilde y_I)\sim(\tilde
u'_I,\phi'_I,
 \tilde y'_I)\Longleftrightarrow
\psi_I\cdot (\tilde u_I,\phi_I, \tilde y_I)\sim\psi_I\cdot(\tilde
u'_I,\phi'_I,  \tilde y'_I).
$$
 Consequently, (3.46) induces a
natural  action of $\Gamma_I$ on $\widehat W_I^{\Gamma_I}$:
$$
 \psi_I\cdot[\tilde u_I,\phi_I,
\tilde y_I]=[\psi_I\cdot(\tilde u_I,\phi_I, \tilde
y_I)]\;\forall\psi_I\in\Gamma_I.\eqno(3.48)
$$
 The action is also partially smooth because in the charts in (3.44)
and (3.47)
$$[\tilde u_I,\phi_I, \tilde y_I]=[\tilde u_I,\phi_I, \phi_I\circ\lambda_I^s(\tilde
x)]\mapsto\psi_I\cdot[\tilde u_I,\phi_I, \tilde y_I]$$
  is given by
$$
O(\tilde u_{i_s})\to\psi_{i_s}(O(\tilde u_{i_s})),\quad\tilde
x\mapsto \psi_{i_s}(\tilde x).
$$
There is also a continuous surjective map
$$
 q_I: \widehat W_I^{\Gamma_I}\to \widetilde W_I,\;[\tilde u_I,\phi_I, \tilde
 y_I]\mapsto \tilde y_I\eqno(3.49)
$$
such that $\tilde\pi_I\circ q_I=\hat\pi_I$. By (3.46) and (3.48),
$q_I$ commutates with the actions on $\widehat W_I^{\Gamma_I}$ and
$\widetilde W_I$, i.e. $q_I(\psi_I\cdot[\tilde u_I,\phi_I, \tilde
y_I])=\psi_I\cdot (q_I([\tilde u_I,\phi_I, \tilde  y_I]))$ for
each $[\tilde u_I,\phi_I, \tilde  y_I]\in \widehat W_I^{\Gamma_I}$
and $\psi_I\in\Gamma_I$. These imply that $q_I$ induces a
continuous surjective map $\hat q_I$ from $\widehat
 W_I^{\Gamma_I}/\Gamma_I$ to $\widetilde W_I^{\Gamma_I}/\Gamma_I$.
 Hence $\hat\pi_I$ induces a continuous
surjective map from $\widehat W_I^{\Gamma_I}/\Gamma_I$ to $W_I$.
Summarizing the above arguments we get:\vspace{2mm}

\noindent{\bf Proposition 3.9}\quad {\it {\rm (i)} Each $\widehat
 W_I^{\Gamma_I}$ is a stratified Banach manifold and $\hat\pi^I_J$
is a PS map.\\
{\rm (ii)} There exists a natural PS action of $\Gamma_I$ on
$\widehat
 W_I^{\Gamma_I}$ under which $\hat\pi_I$ is invariant and thus
induces a continuous surjective map from $\widehat
 W_I^{\Gamma_I}/\Gamma_I$ to $W_I$.\\
{\rm (iii)} $\hat\pi_J\circ\hat\pi_J^I=\iota^W_{IJ}\circ\hat\pi_I$
for any $J\subset I\in {\cal N}$ and the inclusion
$\iota^W_{IJ}: W_I\hookrightarrow W_J$.\\
 {\rm (iv)} $\hat\pi^I_J\circ\psi_I=\lambda^I_J(\psi_I)\circ\hat\pi^I_J$ for
any $J\subset I\in {\cal N}$ and $\psi_I\in\Gamma_I$.}\vspace{2mm}

 Similarly we can define a desingularization  $\widehat E_I^{\Gamma_I}$ of
  $\widetilde E_I^{\Gamma_I}$ and the bundle projection $\hat p_I:\widehat E_I^{\Gamma_I}
  \to\widehat W_I^{\Gamma_I}$ so that the restriction of it to
  each stratum of $\widehat W_I^{\Gamma_I}$ is a Banach vector
  bundle. We shall give necessary details.  Consider the stratified Banach manifold
$$
\boxed{\widetilde E_I^{\Gamma_I}\!}:= \coprod_{\tilde
u_I\in\widetilde W_I^{\Gamma_I}}\coprod_{O\in{\mathcal
U}(u_I)}\coprod_{1\le s\le |I|}\coprod_{\phi_I\in\Gamma(\tilde
u_I)_s}\{(\tilde
u_I,\phi_I)\}\times\Phi_I\circ\Lambda_I^s(\widetilde
E_{i_s}|_{O(\tilde u_{i_s})}).\eqno(3.50)
$$
Clearly, it is  a stratified Banach bundle over $\boxed{\widetilde
W_I^{\Gamma_I}\!}\,$. Here ${\mathcal U}(u_I)$ and $O(\tilde
u_{i_s})=\pi^{-1}_{i_s}(O)$ are as in (3.41), and $\Phi_I$ and
$\Lambda_I^s$ as in (3.33). Define the equivalence relation
$\overset{e}{\sim}$ in $\boxed{\widetilde E_I^{\Gamma_I}\!}$ as
follows:
 For $\tilde u_I,\tilde u'_I\in\widetilde W_I^{\Gamma_I}$,
 and $\phi_I\in\Gamma(\tilde u_I)_s$, $\phi'_I\in\Gamma(\tilde u'_I)_t$,
 and  $\tilde\xi_I\in\Phi_I\circ\Lambda_I^s(\widetilde E_{i_s}|_{O(\tilde
u_{i_s})})$
 and $\tilde\xi_I^\prime\in\Phi_I^\prime\circ\Lambda_I^{\prime t}(\widetilde E_{i_t}|_{O(\tilde
u'_{i_t})})$
  we define
$$(\tilde u_I,\phi_I, \tilde \xi_I)\overset{e}{\sim}(\tilde
u'_I,\phi'_I,  \tilde\xi'_I)\eqno(3.51)
$$
if and only if the following two conditions hold:
\begin{description}
\item[(A)] $\tilde \xi_I=\tilde \xi'_I$ as points in $\widetilde
W_I^{\Gamma_I}$;

\item[(B)] $\Phi_I\circ\Lambda_I^s(\widetilde E_{i_s}|_{O(\tilde
u_{i_s})}) \bigcap\Phi'_I\circ\Lambda_I^{\prime t}(\widetilde
E_{i_t}|_{O(\tilde u'_{i_t})})$ is a local component of
$\widetilde E_I^{\Gamma_I}\;{\rm near}$ the fibre at $\tilde
p_I(\tilde \xi_I)=\tilde p_I(\tilde \xi'_I)$, where $\tilde
p_I:\widetilde E_I^{\Gamma_I}\to\widetilde W_I^{\Gamma_I}$ is the
obvious projection.
\end{description}
It is easy to see that the relation $\overset{e}{\sim}$ is
compatible with  $\sim$ in (3.42). Denote by $\langle\tilde
u_I,\phi_I, \tilde \xi_I\rangle$ the equivalence class of $(\tilde
u_I,\phi_I, \tilde \xi_I)$, and by
$$
 \widehat E_I^{\Gamma_I}:=\boxed{\widetilde
E_I^{\Gamma_I}\!}\bigg/\overset{e}{\sim}
$$
 if $|I|>1$, and by $\widehat E_I^{\Gamma_I}=\widetilde E_i$ if
 $I=\{i\}\in{\mathcal N}$. As in (3.43) it holds that
$$
\widehat E_I^{\Gamma_I}=\coprod_{\tilde u_I\in\widetilde
W_I^{\Gamma_I}}\coprod_{O\in{\mathcal
U}(u_I)}\coprod_{\phi_I\in\Gamma(\tilde u_I)_1}\{(\tilde
u_I,\phi_I)\}\times\Phi_I\circ\Lambda_I^1(\widetilde
E_{i_1}|_{O(\tilde u_{i_1})})\bigg/\overset{e}{\sim}.
$$
Clearly, we get a stratified Banach bundle,
$$
 \hat p_I: \widehat
 E_I^{\Gamma_I}\to\widehat
 W_I^{\Gamma_I},\;\langle\tilde u_I,\phi_I, \tilde
 \xi_I\rangle\mapsto [\tilde u_I,\phi_I, \tilde p_I(\tilde \xi_I)].\eqno(3.52)
$$
 For each $s=1,\cdots,k$, we have also the PS bundle open embedding,
$$
\widehat{\Phi_I\circ\Lambda_I^s}:\widetilde E_{i_s}|_{\widetilde
O_{i_s}}\to \widehat E_I^{\Gamma_I},\;\tilde \xi\mapsto
\langle\tilde u_I,\phi_I, \Phi_I\circ\Lambda_I^s(\tilde
\xi)\rangle,\eqno(3.53)
$$
 called {\it a local bundle
coordinate chart of $\widehat E_I^{\Gamma_I}$ near $\langle\tilde
u_I,\phi_I, \tilde o_I(\tilde u_I)\rangle$}. Here $\tilde
o_I:\widetilde W_I^{\Gamma_I}\to\widetilde E_I^{\Gamma_I}$ is the
zero section. By Proposition 3.8(iii) the projections
$\widetilde\Pi_I$ and $\widetilde\Pi^I_J$ induce natural ones
$$\widehat\Pi_I: \widehat
E_I^{\Gamma_I}\to E_I,\;\langle\tilde u_I,\phi_I, \tilde
\xi_I\rangle\mapsto \widetilde\Pi_I(\tilde \xi_I),
$$
$$
\widehat\Pi^I_J:\widehat E_I^{\Gamma_I}\to \widehat
E_J^{\Gamma_J},\;\langle\tilde u_I,\phi_I, \tilde
\xi_I\rangle\mapsto \langle\tilde u_J,\phi_J, \tilde
\xi_J\rangle,\eqno(3.54)
$$
  where
$\phi_J=\lambda^I_J(\phi_I)=(\phi_l)_{l\in J}$, $\tilde
u_J=\tilde\pi^I_J(\tilde u_I)=(\tilde u_l)_{l\in J}$ and $\tilde
\xi_J=\widetilde\Pi^I_J(\tilde \xi_I)=(\tilde \xi_l)_{l\in J}$. In
the bundle charts $\widehat{\Phi_I\circ\Lambda_I^1}$ and
$$ \widehat{\Phi_J\circ\Lambda_J^1}:
\widetilde E_{j_1}|_{O(\tilde u_{j_1})}\to \widehat
E_J^{\Gamma_J},\;\tilde\eta\mapsto \langle\tilde u_J,\phi_J,
\Phi_J\circ\Lambda_J^1(\tilde \eta)\rangle,\eqno(3.55)
$$
 the projection $\widehat\Pi^I_J$ may be represented by
$$
 \widetilde E_{i_1}|_{O(\tilde u_{i_1})}\to \widetilde E_{j_1}|_{O(\tilde u_{j_1})},
 \; \tilde \xi\mapsto
\Lambda_{i_1j_1}(\tilde \xi)\eqno(3.56)
$$
because $j_1\in J\subset I$. Here $\Lambda_{i_1j_1}=1_{\widetilde
E_{i_1}|_{O(\tilde u_{i_1})}}$ if $j_1=i_1$. Hence
$\widehat\Pi^I_J$ is a {\it PS bundle map} and also {\it
isomorphically maps} the fibre of $\widehat E_I^{\Gamma_I}$ at
$\hat x_I\in\widehat W_I^{\Gamma_I}$ to that of $\widehat
E_J^{\Gamma_J}$ at $\hat\pi^I_J(\hat x_I)\in\widehat
W_J^{\Gamma_J}$. By (3.51) we may define a natural norm in each
fibre of $\widehat E_I^{\Gamma_I}$ by
$$
 \|\langle\tilde u_I, \phi_I,
\tilde\xi_I\rangle\|:=\sum_{l\in I}\|\tilde\xi_l\|\eqno(3.57)
$$
 for any vector
$\langle\tilde u_I, \phi_I, \tilde\xi_I\rangle$ in the fibre of
$\widehat E_I^{\Gamma_I}$ at $[\tilde u_I, \phi_I,\tilde x_I)]\in
\widehat W_I^{\Gamma_I}$, where $\|\tilde\xi_l\|$ denotes the norm
of $\tilde\xi_l$ in $(\tilde E_l)_{\tilde x_l}$.
 The action of $\Gamma_I$ on $\widetilde E_I^{\Gamma_I}$ induces
a natural one on $\boxed{\widetilde E_I^{\Gamma_I}\!}\,$,
$$
 \psi_I\cdot (\tilde u_I,\phi_I,
\tilde \xi_I)=\bigl(\psi_I(\tilde u_I), \phi^\ast_I,
\psi_I\cdot\tilde \xi_I\bigr)\eqno(3.58)
$$
for any $\psi_I=(\psi_l)_{l\in I}\in\Gamma_I$. Here
$\phi^\ast_I=(\phi^\ast_l)_{l\in
I}=(\psi_l\circ\phi_l\circ\psi_l^{-1})_{l\in I}
\in\Gamma(\psi_I(\tilde u_I))_1 $ as in (3.45), and
$$
\psi_I\cdot\tilde\xi_I=\Psi_I(\tilde\xi_I)=(\Psi_l(\tilde\xi_l))_{l\in
I}=\Phi^\ast_I\circ\Lambda_I^{1\ast}(\Psi_{i_1}(\tilde \eta))\in
\Phi^\ast_I\circ\Lambda_I^{1\ast}\bigr(\Psi_{i_1}(\widetilde
E_{i_1}|_{O(\tilde u_{i_1})})\bigl)
$$
for $ \tilde \xi_I=\Phi_I\circ\Lambda_I^1(\tilde \eta)\in
\Phi_I\circ\Lambda_I^1(\widetilde E_{i_1}|_{O(\tilde u_{i_1})})$
and $\Phi^\ast_I=(\Phi^\ast_l)_{l\in I}=
(\Psi_l\circ\Phi_l\circ\Psi_l^{-1})_{l\in I}$,  and
$\Lambda^{1\ast}_I=(\Lambda^\ast_{i_1l})_{l\in I}$,
$\Lambda^\ast_{i_1l}=\Psi_l\circ\Lambda_{i_1l}\circ\Psi_{i_1}^{-1}$
for any $l\in I$.
 Equation (3.58) gives a
PS action of $\Gamma_I$ on $\boxed{\widetilde
E_I^{\Gamma_I}\!}\,$. For each $s\in\{1,\cdots,k\}$, as above if
$\phi_I\in\Gamma(\tilde u_I)_s$ then $\phi^\ast_I$ belongs to
$\Gamma(\psi_I(\tilde u_I))_s$. Let us define
$\Lambda^{s\ast}_I=(\Lambda^\ast_{i_sl})_{l\in I}$,
$\Lambda^\ast_{i_sl}=\Psi_l\circ\Lambda_{i_sl}\circ\Psi_{i_s}^{-1}$
for any $l\in I$, then
$$
\Phi^\ast_I\circ\Lambda_I^{s\ast}:\Psi_{i_s}(\widetilde
E_{i_s}|_{O(\tilde u_{i_s})})\to \widetilde E_I^{\Gamma_I},\;
\tilde \xi_s\mapsto
\bigl(\Phi_l^\ast\circ\Lambda^\ast_{i_sl}(\tilde
\xi_s)\bigr)_{l\in I}
$$
is a PS bundle open embedding, and
$\{\Phi_I^\ast\circ\Lambda_I^{s\ast}\,|\,\phi_I\in\Gamma(\tilde
u_I)_s\}$ is a local bundle coordinate chart of $\widetilde
E_I^{\Gamma_I}$ over the neighborhood $\psi_I\bigl(\widetilde
O(\tilde u_I)\bigr)$ of $\psi_I(\tilde u_I)$.  In particular, we
get a local bundle coordinate chart of $\widehat E_I^{\Gamma_I}$
near $\langle\psi_I(\tilde u_I),\phi^\ast_I,\tilde
o_I(\psi_I(\tilde u_I))\rangle$,
$$
\widehat{\Phi^\ast_I\circ\Lambda_I^{\ast 1}}:\widetilde
E_{i_1}|_{\psi_{i_1}(O(\tilde u_{i_1}))}\to \widehat
E_I^{\Gamma_I},\;\tilde \xi\mapsto \langle\psi_I(\tilde
u_I),\phi^\ast_I, \Phi^\ast_I\circ\Lambda_I^{\ast 1}(\tilde
\xi)\rangle.\eqno(3.59)
$$
Moreover one easily checks that the action in (3.58) preserves the
equivalence relation $\overset{e}{\sim}$ in $\boxed{\widetilde
E_I^{\Gamma_I}\!}\,$ and thus
 induces a natural action of
$\Gamma_I$ on $\widehat E_I^{\Gamma_I}$:
$$
 \psi_I\cdot\langle\tilde
u_I,\phi_I, \tilde \xi_I\rangle=\langle\psi_I\cdot(\tilde
u_I,\phi_I, \tilde\xi_I)\rangle\eqno(3.60)
$$
for any $\psi_I\in\Gamma_I$. Note that in the charts
$\widehat{\Phi_I\circ\Lambda_I^1}$ in (3.53) and
$\widehat{\Phi^\ast_I\circ\Lambda_I^{\ast 1}}$ in (3.59)
$$\langle\tilde u_I,\phi_I, \tilde\xi_I\rangle=\langle\tilde u_I,\phi_I,
\Phi_I\circ\Lambda_I^1(\tilde\eta)\rangle\mapsto\psi_I\cdot\langle\tilde
u_I,\phi_I, \tilde\xi_I\rangle$$
 is given by
$$
 \widetilde E_{i_1}|_{O(\tilde u_{i_1})}\to\widetilde E_{i_1}|_{\psi_{i_1}(O(\tilde
u_{i_1}))},\quad\tilde \eta\mapsto \Psi_{i_1}(\tilde\eta).
$$
So (3.60) defines a {\it PS action of }$\Gamma_I$ on $\widehat
E_I^{\Gamma_I}$.
 Moreover, the clear continuous surjective map
$$
 Q_I: \widehat E_I^{\Gamma_I}\to \widetilde E_I,\;\langle\tilde u_I,\phi_I, \tilde
 \xi_I\rangle\mapsto \tilde\xi_I
$$
 commutates with the above actions of $\Gamma_I$ on
$\widehat E_I^{\Gamma_I}$ and $\widetilde E_I$ and also satisfies
$\widetilde\Pi_I\circ Q_I=\widehat\Pi_I$. It follows that  the
projection $\widehat\Pi_I$ is invariant under the
 $\Gamma_I$-action and induces a continuous surjective map
$$
  \widehat E_I^{\Gamma_I}/\Gamma_I\to  E_I,\;\langle\tilde u_I,\phi_I, \tilde
 \xi_I\rangle\mapsto \widetilde\Pi_I(\tilde\xi_I).
$$
Summarizing up the above arguments we get:\vspace{2mm}

\noindent{\bf Proposition 3.10.}\quad{\it {\rm (i)} Each $\hat
p_I:\widehat
 E_I^{\Gamma_I}\to\widehat
 W_I^{\Gamma_I}$ is a stratified Banach bundle and $\widehat\Pi^I_J$
is a PS

bundle map.\\
{\rm (ii)} There exists a PS action of $\Gamma_I$ on $\widehat
 E_I^{\Gamma_I}$ (given by
(3.60)) under which $\widehat\Pi_I$ is invariant and

induces a continuous surjective map from $\widehat
 E_I^{\Gamma_I}/\Gamma_I$ to $E_I$.\\
 {\rm (iii)} The projection $\hat p_I:\widehat
E_I^{\Gamma_I}\to\widehat W_I^{\Gamma_I}$ is equivariant with
respect to the actions of $\Gamma_I$ on them, i.e.

 $\hat
p_I\circ\psi_I=\psi_I\circ\hat p_I$
for any $\psi_I\in\Gamma_I$.\\
 {\rm (iv)} $\hat\Pi_J\circ\hat\Pi_J^I=\iota^E_{IJ}\circ\hat\Pi_I$
   for any $J\subset I\in{\cal N}$ and the inclusion $\iota^E_{IJ}:
E_I\hookrightarrow E_J$.\\
  {\rm (v)} $\widehat\Pi^I_J\circ\psi_I=\lambda^I_J(\psi_I)\circ\widehat\Pi^I_J$
for any $J\subset I\in{\cal N}$ and  $\psi_I\in\Gamma_I$.\\
{\rm (vi)} $\hat\pi_I\circ\hat p_I=p\circ\widehat\Pi_I$ for any
 $\phi_I\in\Gamma_I$.\\
 {\rm (vii)} $\hat p_J\circ\widehat\Pi^I_J=\hat\pi^I_J\circ\hat p_I$ for any $J\subset I\in{\cal
 N}$.\\
 {\rm (viii)} For any $L, J, I\in{\cal N}$ with $L\subset J\subset
 I$ it holds that $\hat\pi^J_L\circ\hat\pi^I_J=\hat\pi^I_L$ and
 $\widehat\Pi^J_L\circ\widehat\Pi^I_J=\widehat\Pi^I_L$.}\vspace{2mm}

The final (vii) may easily follow from (3.45) and (3.54). We get a
system of PS Banach bundles
$$
\bigl(\widehat{\mathcal E}^\Gamma, \widehat W^\Gamma\bigr)=
\Bigl\{\bigl(\widehat E_I^{\Gamma_I}, \widehat
W_I^{\Gamma_I}\bigr), \hat\pi_I, \widehat\Pi_I, \Gamma_I,
\hat\pi^I_J,\widehat\Pi^I_J, \lambda^I_J\,\bigm|\, J\subset
I\in{\mathcal N}\Bigr\},\eqno(3.61)
$$
which is called a {\it desingularization} of the bundle system
$\bigl(\widetilde{\mathcal E}^\Gamma, \widetilde W^\Gamma\bigr)$.

As in [LiuT1] we may take the pairs of $\Gamma_i$-invariant open
sets $W^j_i\subset\subset U^j_i$, $j=1, 2, \cdots,n_3-1$, such
that
$$U^1_i\subset\subset W^2_i\subset\subset U^2_i\cdots
\subset\subset W^{n_3-1}_i\subset\subset U^{n_3-1}_i\subset\subset
W_i^{n_3}=W_i.$$
 Note that each $I\in{\mathcal N}$ has length $|I|<n_3$. For each $I\in{\mathcal N}$  we define
$$V_I=(\cap_{i\in I}W_i^k)\setminus Cl(\bigcup_{J:|J|>k}(\cap_{j\in
J}Cl(U^k_j))).\eqno(3.62)
$$
Then the second condition in (3.20) implies
$$
V_I=W_I\quad\forall I\in{\cal N}\;\;{\rm with}\; |I|=n_3-1.
$$
 For $\{V_I\,|\, I\in{\cal
N}\}$, unlike $\{W_I\,|\, I\in{\cal N}\}$ we cannot guarantee that
$V_I\subset V_J$ even if $J\subset I$. However we have still the
corresponding result with Lemma 4.3 in [LiuT1].\vspace{2mm}

 \noindent{\bf Lemma 3.11.}\quad{\it
 $\{V_I\,|\,I\in{\mathcal
N}\}$ is an open covering of $\cup^{n_3}_{i=1}W_i^1$  and
satisfies:
\begin{description}
\item[(i)] $V_I\subset W_I$  for any $I\in {\mathcal N}$.

\item[(ii)] $Cl(V_I)\cap V_J\ne\emptyset$  only if $I\subset J$ or
$J\subset I$. (Actually  $Cl(V_I)\cap V_J\ne\emptyset$ implies
$V_I\cap V_J\ne\emptyset$, and thus $Cl(V_I)\cap V_J=\emptyset$
for any $I, J\in{\cal N}$ with $|I|=|J|=n_3-1$.)

\item[(iii)] For any nonempty open subset ${\cal
W}^\ast\subset\cup^{n_3}_{i=1}W_i^1$ and $I\in{\cal N}$ let
$W_I^\ast=W_I\cap{\cal W}^\ast$ and $V_I^\ast=V_I\cap{\cal
W}^\ast$. Then $\{V_I^\ast\,|\, I\in{\cal N}\}$ is an open
covering of ${\cal W}^\ast$, and also satisfies the above
corresponding properties (i)-(ii).
\end{description}}

\noindent{\it Proof.}\quad {\it Step 1}. Let ${\cal
N}_k=\{I\in{\cal N}\,|\, |I|=k\}$, $k=1,\cdots,n_3-1$. We first
prove that
 $\{V_I\,|\,I\in{\mathcal N}\}$  is an open covering of
$\cup^{n_3}_{i=1}W_i^1$. It is easy to see that for any
$I\in{\mathcal N}$ with $|I|=k$,
$$
V_I=\Bigl(\bigcap_{i\in I}W_i^k\Bigr)\setminus
\Bigl(\bigcup_{J\in{\cal N}_{k+1}}\bigl(\bigcap_{j\in
J}Cl(U_j^k)\bigr)\Bigr).\eqno(3.63)
 $$
 For $x\in W_i$, set $I_1=\{i\}$.
 If $x\in V_{I_1}$ nothing is done. Otherwise, there is
 $J_1\in{\mathcal N}$ with $|J_1|=2$  such that
 $x\in\cap_{j\in J_1}Cl(U^1_j)\subset W^2_{J_1}$. Set $I_2=I_1\cup J_1$ then
 $x\in W^1_{I_1}\cap W^2_{J_1}\subset W^2_{I_1}\cap W^2_{J_1} $ and
 thus
 $$x\in W^2_{I_2}\subset W^{|I_2|}_{I_2},\;  |I_2|\ge |J_1|=|I_1|+1=2.$$
If $x\in V_{I_2}$ then nothing is done. Otherwise, because of
(3.62) there is $J_2\in{\mathcal N}$ with $|J_2|=|I_2|+1$ such
that $x\in \cap_{j\in J_2}Cl(U_j^{|I_2|})\subset W_{J_2}^{|J_2|}$.
Set $I_3=I_2\cup J_2$ then
 $x\in  W^2_{I_2}\cap W_{J_2}^{|J_2|}\subset W^{|I_3|}_{I_2}\cap W_{J_2}^{|I_3|}$ and
 thus
$$
x\in W_{I_3}^{|I_3|},\;  |I_3|\ge |J_2|=|I_2|+1\ge 3.
$$
After repeating finite times this process there must exist some
$I_k\in{\mathcal N}$ such that $x\in V_{I_k}$ because $
V_I=W_I^{n_3-1}$ for any $I\in{\cal N}$ with $|I|=n_3-1$.
 The desired conclusion is proved.\vspace{2mm}

{\it Step 2.}  We prove (i) and (ii). (i) is obvious.  We only
need to prove (ii). Let $ Cl(V_I)\cap V_J\ne\emptyset$ for two
different  $I, J\in{\mathcal N}$ with $|I|=k$ and $|J|=l$. Since
$I=\{1,\cdots, n_3\}\notin{\cal N}$ we can assume that $n_3>k\ge
l$ below. Suppose that there exists a $r\in J\setminus I$. Let
$x\in Cl(V_I)\cap V_J$. Take a sequence $\{x_k\}\subset V_I$ such
that $x_k\to x$. Since $V_J$ is open then $x_k\in V_J$ for $k$
sufficiently large. So  $V_I\cap V_J\ne\emptyset$ and we may
assume that $x\in V_I\cap V_J$. Now one hand $x\in V_J$ implies
that
$$x\in Cl\bigl(\cap_{i\in J}W^l_i\bigr)\subset Cl\bigl(\cap_{i\in J}W^k_i\bigr)\subset Cl(W^k_r).$$
 On the other hand $x\in V_I$
 implies that $x$ does not belong  to
 \begin{eqnarray*}
 \cap_{i\in I}Cl(U^k_i)\cap Cl(U^k_r)\!\!\!\!\!\!\!&&\supseteq\cap_{i\in I}Cl(W^k_i)\cap
 Cl(W^k_r)\\
&&\supseteq Cl(\cap_{i\in I}W^k_i)\cap Cl(W^k_r).
\end{eqnarray*}
So $x\notin Cl(W^k_r)$ because $x\in V_I\subset Cl\bigl(\cap_{i\in
I}W^k_i\bigr)$. This contradiction shows that $J\subset
I$.\vspace{2mm}

{\it Step 3.} Note that $Cl(V_I^\ast)\cap V_J^\ast\subset
Cl(V_I)\cap V_J\cap{\cal W}^\ast$ for any $I, J\in{\cal N}$. The
desired conclusions follow immediately.
   \hfill$\Box$\vspace{2mm}

  Set
$$ \widehat V_I=(\hat\pi_I)^{-1}(V_I)\quad{\rm and}\quad\widehat
E_I=(\widehat\Pi_I)^{-1}(E_I|_{V_I}),\eqno(3.64)
$$
For $J\subset
I\in{\cal N}$, though $\hat\pi^I_J(\widehat
W_I^{\Gamma_I})\subset\widehat W_J^{\Gamma_J}$ we cannot guarantee
that $\hat\pi^I_J(\widehat V_I)\subset\widehat V_J$. So
$\hat\pi^I_J$ only defines a PS map from the open subset
$(\hat\pi^I_J)^{-1}(\widehat V_J)\cap\widehat V_I$
 to $\widehat V_J$. However, we still denote by
$\hat\pi^I_J$ the restriction of $\hat\pi^I_J$ to
$(\hat\pi^I_J)^{-1}(\widehat V_J)\cap\widehat V_I$. In this case
it holds that
$$\hat\pi_J\circ\hat\pi^I_J=\iota^V_{IJ}\circ\hat\pi_I\eqno(3.65)
$$
for any $J\subset I\in {\cal N}$, where $\iota^V_{IJ}:V_I\cap
V_J\hookrightarrow V_J$ is the inclusion. Similarly,
$\widehat\Pi^I_J$ denote the restriction of $\widehat\Pi^I_J$ to
$(\widehat\Pi^I_J)^{-1}(\widehat E_J)\cap\widehat E_I$.
 Since the bundle projection
$\hat p_I:\widehat E_I^{\Gamma_I}\to\widehat W_I^{\Gamma_I}$ in
(3.52) maps $\widehat E_I$ to $\widehat V_I$  we still use $\hat
p_I$ to denote the bundle projection $\widehat E_I\to\widehat
V_I$. The system of stratified Banach bundles
$$
(\widehat {\mathcal E}, \widehat V)=\bigl \{(\widehat E_I,
\widehat V_I), \pi_I,
 \pi^I_J, \widehat\Pi_I, \widehat\Pi^I_J, p_I,\Gamma_I\,\bigm|\, J\subset I\in{\mathcal
 N}\bigr\},\eqno(3.66)
$$
 is called {\it renormalization} of
$(\widehat {\mathcal E}^{\Gamma}, \widehat W^\Gamma)$.

 A {\it global section} of the bundle system $(\widehat {\mathcal
E}^\Gamma, \widehat W^\Gamma)$ (resp. $(\widehat {\mathcal E},
\widehat V)$ ) is a compatible collection $S=\{S_I\,|\,
I\in{\mathcal N}\}$ of sections $S_I$ of $\widehat
E_I^{\Gamma_I}\to\widehat W_I^{\Gamma_I}$ (resp.  $\widehat
E_I\to\widehat V_I$ ) in the sense that $(\widehat\Pi^I_J)^\ast
S_J=S_I$, i.e., $S_J(\hat\pi^I_J(x))=\widehat\Pi^I_J(S_I(x))$
  for any $x\in
\widehat W_I^{\Gamma_I}$ (resp. $(\hat\pi^I_J)^{-1}(\widehat
V_J)\cap\widehat V_I$)) and $J\subset I\in{\cal N}$. As before $S$
is said to be {\it partially smooth} if each $S_I$ is continuous
and stratawise smooth. Such a section $S$ is called {\it
transversal} to the zero section if each $S_I$ is so.\vspace{2mm}

\noindent{\bf Lemma 3.12.}\quad{\it Let $I_l=\{l\}$ for $1\le l\le
n_3$. For $I\in{\cal N}$ with $|I|>1$, and $l\in I$ each PS
section   $\tilde\sigma_l:\widetilde W_l\to\widetilde E_l$ can
define a PS
  section of the bundle $\hat p_I:\widehat
  E_I^{\Gamma_I}\to\widehat W_I^{\Gamma_I}$ by
 $$
 \hat\sigma_{lI}:\widehat W_I^{\Gamma_I}\to\widehat
E_I^{\Gamma_I},\; \hat x\mapsto
\bigl((\widehat\Pi^I_{I_l})^\ast\tilde\sigma_l\bigr)(\hat x)=:
   (\widehat\Pi^I_{I_l}|_{(\widehat E_I)_{\hat x}})^{-1}\bigl(\tilde\sigma_l(\hat \pi^I_{I_l}(\hat
   x))\bigr)\eqno(3.67)
   $$
which is equal to $\tilde\sigma_l$ if $I=I_{l}$. In particular, if
$\tilde\sigma_l:\widetilde W_l\to\widetilde E_l$ has the support
contained in $\widetilde W_l$ (i.e., $\pi_l(\{\tilde
x\in\widetilde W_l\,|\,\tilde\sigma_l(x)\ne 0\})\subset\subset
W_l$) then for each $I\in{\cal N}$ it may determine a PS section
$(\hat\sigma_l)_I$ of the bundle $\widehat E_I\to\widehat V_I$
even if $l\notin I$, and these PS sections give rise to a global
PS section $\hat \sigma_l=\{(\hat\sigma_l)_I: I\in{\cal N}\}$ of
the bundle system $(\widehat {\mathcal E}, \widehat
V)$.}\vspace{2mm}

\noindent{\it Proof.}\quad To see that the section $\hat
\sigma_{lI}$ defined by (3.67)  is partially smooth let us write
it in the charts in (3.44) and (3.53) with $s=1$. Let $\hat
x=\widehat{\phi_I\circ\lambda_I^1}(\tilde x_{i_1})=[\tilde
u_I,\phi_I, \phi_I\circ\lambda^1_I(\tilde x_{i_1})]$ for $\tilde
x_{i_1}\in O(\tilde u_{i_1})$. Then $\hat\pi^I_{I_l}(\hat
x)=[\tilde u_l,\phi_l,\phi_l\circ\lambda_{i_1l}(\tilde
x_{i_1})]=\phi_l\circ\lambda_{i_1l}(\tilde x_{i_1})$, and in
particular $\hat\pi^I_{I_l}(\hat x)=\tilde x_{i_1}$ if $l=i_1$. So
$\tilde\sigma_l(\phi_l\circ\lambda_{i_1l}(\tilde
x_{i_1}))=\langle\tilde u_l,\phi_l,
\tilde\sigma_l(\phi_l\circ\lambda_{i_1l}(\tilde x_{i_1}))\rangle$.
 Under the bundle charts
$\widehat{\Phi_I\circ\Lambda^1_I}$ in (3.53) and
$$
\widehat{\Phi_{I_l}\circ\Lambda^1_{I_l}}: \widetilde
E_{l}|_{O(\tilde u_{l})}\to \widehat
E_{I_l}^{\Gamma_{I_l}},\;\tilde\eta\mapsto \langle\tilde
u_{I_l},\phi_{I_l}, \Phi_{I_l}\circ\Lambda^1_{I_l}(\tilde
\eta)\rangle=\langle\tilde u_l,\phi_l, \tilde
\eta\rangle=\tilde\eta,
$$
we can derive from (3.55) and (3.56) that the projection
$\hat\Pi^I_{I_l}$ may be represented by
$$
 \widetilde E_{i_1}|_{O(\tilde
u_{i_1})}\to \widetilde E_{l}|_{O(\tilde u_l)},
 \; \tilde \xi_{i_1}\mapsto
\Phi_l\circ\Lambda_{i_1l}(\tilde \xi_{i_1}).
 $$
 It follows that
$\hat\sigma_{lI}$ may locally be represented as
$$
 O(\tilde u_{i_1})\to \widetilde E_{i_1}|_{O(\tilde
u_{i_1})},  \; \tilde x_{i_1}\mapsto
(\Phi_l\circ\Lambda_{i_1l})^{-1}\bigl(\tilde\sigma_l(\phi_l\circ\lambda_{i_1l}(\tilde
x_{i_1}))\bigr).\eqno(3.68)
$$
 This shows that $\hat\sigma_{lI}$ is partially smooth.

Next we prove the second claim.  For $I=\{i_1,\cdots,i_k\}\in
{\mathcal N}$, if $l\notin I$, (3.63) implies that
\begin{eqnarray*} V_I\!\!\!\!\!&&=\Bigl(\bigcap_{i\in
I}W_i^k\Bigr)\setminus
\Bigl(\bigcup_{J:|J|=k+1}\bigl(\bigcap_{j\in
J}Cl(U^k_j)\bigr)\Bigr)\\
&&\subset\Bigl(\bigcap_{i\in I}W_i^k\Bigr)\setminus
\Bigl(\bigl(\bigcap_{i\in I}Cl(U^k_i)\bigr)\cap
Cl(U^k_l)\Bigr)\\
&&\subset\Bigl(\bigcap_{i\in I}W_i^k\Bigr)\setminus
\Bigl(\bigl(\bigcap_{i\in I}Cl(U^k_i)\bigr)\cap
Cl(W^1_l)\Bigr)\\
&&=\Bigl(\bigcap_{i\in I}W_i^k\Bigr)\setminus
 Cl(W^1_l)\\
&&\subset\Bigl(\bigcap_{i\in I}Cl(W_i^k)\Bigr)\setminus
 W^1_l
\end{eqnarray*}
where the second inclusion is because $W^1_l\subset U^1_l\subset
U^k_l$, and the second equality comes from the fact that
$\cap_{i\in I}W_i^k\subset\cap_{i\in I}U^k_i$.

It follows that
$$
Cl(V_I)\subset\Bigl(\bigcap_{i\in I}Cl(W_i^k)\Bigr)\setminus
 W^1_l
$$
since the left side is closed.  Moreover, $\pi_l({\rm
supp}(\tilde\sigma_l))\subset W_l^1$.

 We get that
$$
\pi_l({\rm supp}(\tilde\sigma_l))\cap
Cl(V_I)=\emptyset\quad\forall l\notin I,\,1\le l\le n_3.
\eqno(3.69)
 $$
 Note that $\widehat V_{I_l}\subset\widehat W_{I_l}=\widetilde W_l$.
 For $I\in{\cal N}$ we define
$$
\left\{
\begin{array}{ll}
(\hat\sigma_l)_{I_l}=\tilde\sigma_l|_{\widehat V_{I_l}},
(\hat\sigma_l)_I=0\;{\rm if}\; l\notin I,\,{\rm and}\,\\
(\hat\sigma_l)_I:=\hat\sigma_{lI}|_{\widehat V_I}\;{\rm if}\;l\in
I\;{\rm and}\;|I|>1. \end{array}\right.
$$
  By the local
expression in (3.68) they are partially smooth. Clearly, the
collection $\{(\hat\sigma_l)_I: I\in{\cal N}\}$ is compatible in
the sense that $(\widehat\Pi^I_J)^\ast (\hat
\sigma_l)_J=(\hat\sigma_l)_I$ for any $J\subset I\in{\cal N}$.
Hence $\hat \sigma_l=\{(\hat\sigma_l)_I: I\in{\cal N}\}$ is a PS
global section of the bundle system $(\widehat{\cal E},\widehat
V)$.
 \hfill$\Box$\vspace{2mm}

\noindent{\bf Remark 3.13.}\quad For each $i=1,\cdots,n_3$ let
$W^+_i$ be a slight larger open set than $W_i$ with
$W_i\subset\subset W_i^+$. For each $i=1,\cdots,n_3$ let us take
pairs of open subsets $W^{+j}_i\subset\subset U^{+j}_i$,
$j=1,\cdots, n_3-1$ such that
$$
W^j_i\subset\subset W^{+j}_i\subset\subset U^{+j}_i\subset\subset
U^j_i,\;j=1,\cdots, n_3-1.
$$
Then for each $I\in{\cal N}$ with $|I|=k$ we follow (3.62) to
define
$$
 V^+_I:=\Bigl(\bigcap_{i\in I}W^{+k}_i\Bigr)\setminus
\Bigl(\bigcup_{J:|J|>k}\bigl(\bigcap_{j\in
J}Cl(U^{+k}_j)\bigr)\Bigr)
$$
 and thus get a system of stratified Banach bundles  $(\widehat
{\mathcal E}^{\Gamma^+}, \widehat W^{\Gamma^+})$ and its
renormalization $(\widehat {\mathcal E}^+, \widehat V^+)$.  For
any section $\tilde\sigma_l:\widetilde W_l\to\widetilde E_l$ with
$\pi_l({\rm supp}(\tilde\sigma_l))\subset W_l$, which is naturally
viewed a section $\tilde\sigma_l:\widetilde W_l^+\to\widetilde
E_l^+$, the same reason as above may yield a global PS section
$\hat\sigma_l^+=\{(\hat\sigma_l)^+_I: I\in{\cal N}\}$ of the
bundle system $(\widehat {\mathcal E}^+, \widehat V^+)$ which
restricts to $\hat\sigma_l$ on $(\widehat {\mathcal E}, \widehat
V)$. Actually, $(\hat\sigma_l)^+_I$ vanishes outside $\widehat
V_I\subset\widehat V_I^+$. Note that $Cl(V_I)\subset V_I^+$ for
any $I\in{\cal N}$. We may consider the closure of $(\widehat
{\mathcal E}, \widehat V)$ in $(\widehat {\mathcal E}^{+},
\widehat V^{+})$,
$$
(Cl(\widehat {\mathcal E}), Cl(\widehat V))=\bigl \{(Cl(\widehat
E_I), Cl(\widehat V_I)), \hat\pi_I,
 \hat\pi^I_J, \widehat\Pi_I, \widehat\Pi^I_J, \hat p_I,\Gamma_I\,\bigm|\, J\subset I\in{\mathcal N}\bigr\}
$$
and define its global section by requiring that
$S_J(\hat\pi^I_J(x))=\widehat\Pi^I_J(S_I(x))$  for $x\in
(\hat\pi^I_J)^{-1}(Cl(\widehat V_J))\cap Cl(\widehat V_I)$ and
$J\subset I\in{\cal N}$. In this case the above section
$\hat\sigma_l$ produced by a section $\tilde\sigma_l:\widetilde
W_l\to\widetilde E_l$ with $\pi_l({\rm
supp}(\tilde\sigma_l))\subset W_l$ can naturally be extended into
a global section of $(Cl(\widehat {\mathcal E}), Cl(\widehat V))$
(by the zero extension), still denoted by $\hat \sigma_l$ without
confusions. Later when we need to use $(Cl(\widehat {\mathcal E}),
Cl(\widehat V))$ we assume that we have taken a $(\widehat
{\mathcal E}^+, \widehat W^+)$ to contain it. \vspace{2mm}

 It is not difficult to check
that the Cauchy-Riemann operator may define a global PS section
$\hat\partial_J=\{(\hat\partial_J)_I\,|\,I\in{\mathcal N}\}$ of
$(\widehat{\cal E}^\Gamma,\widehat W^\Gamma)$ (resp.
$(\widehat{\cal E},\widehat V)$ and $(Cl(\widehat{\cal E}),
Cl(\widehat V)$)). In particular, from (3.45) and (3.54) it
follows that
$$
(\hat\partial_J)_L\circ\hat\pi^I_L=\widehat\Pi^I_L\circ(\hat\partial_J)_I\eqno(3.70)
 $$
 for any $L\subset I\in{\cal N}$.
In the charts of (3.44) and (3.53) the section
$(\hat\partial_J)_I$ may be expressed as
$$
 O(\tilde u_{i_s})\to\widetilde
E_{i_s}|_{O(\tilde u_{i_s})},\; {\bf f}=(f,\Sigma,\bar z)\mapsto
\bar\partial_Jf.\eqno(3.71)
$$
So $(\hat\partial_J)_I$ is  partially smooth, and its zero set is
given by
$$
{\cal Z}((\hat\partial_J)_I)=(\hat
\pi_I)^{-1}\bigl(Z(\bar\partial_J)\cap W_I\bigr).
$$
Since $Z(\bar\partial_J)\cap Cl(V_I)$ is compact (as a subset of
 $V_I^+$ by Remark 3.13), so is
$$
 {\cal VZ}((\hat\partial_J)_I):=(\hat
\pi_I)^{-1}\bigl(Z(\bar\partial_J)\cap
Cl(V_I)\bigr)\subset\widehat V_I^+.
$$
 Let $\tilde s_{ij}$, $j=1,\cdots,q_i$, $i=1,\cdots,n_3$, be
the sections satisfying (3.21) and Lemma 3.3. By Lemma 3.12 each
$\tilde s_{ij}$ determines a global PS section
 $\hat s_{ij}=\{(\hat s_{ij})_I\,|\,
I\in{\cal N}\}$ of the  bundle system $(\widehat{\cal E},\widehat
V)$. It follows from Remark 3.13 that $\hat s_{ij}$ can also be
viewed a global PS section of $(\widehat {\mathcal E}^+, \widehat
V^+)$. This and finiteness of ${\cal N}$ imply that there exists a
constant $C>0$ and an open neighborhood $U({\cal
VZ}((\hat\partial_J)_I))$  of ${\cal VZ}((\hat\partial_J)_I)$ in
 $\widehat W_I^+$ such that
$$\|(\hat s_{ij})_I(x)\|\le C\;\forall x\in U\bigr({\cal
VZ}((\hat\partial_J)_I)\bigl)\cap \widehat V_I^+. \eqno(3.72)
$$
Here the norm $\|\cdot\|$ is defined by (3.57). Note that
$U\bigr({\cal VZ}((\hat\partial_J)_I)\bigl)\cap\widehat V_I$ is an
open subset in $\widehat V_I$ and that the union $\cup_{I\in{\cal
N}}\hat\pi_I(U\bigr({\cal
VZ}((\hat\partial_J)_I))\bigl)\cap\widehat V_I)$
 is still an open
neighborhood of $\overline{\cal M}_{g,m}(M, J, A; K_3)$ in ${\cal
W}$. Take an open neighborhood ${\cal W}^\ast$ of $\overline{\cal
M}_{g,m}(M, J, A; K_3)$ in ${\cal W}$ so that it is contained in
the intersection of $\cup_{I\in{\cal N}}\hat\pi_I(U\bigr({\cal
VZ}((\hat\partial_J)_I))\bigl)\cap\widehat V_I)$ and
$\cup^{n_3}_{i=1}V_i^0$ in (3.20). Set
$$
V^\ast_I:= V_I\cap {\cal W}^\ast,\quad \widehat
V^\ast_I=(\hat\pi_I)^{-1}(V^\ast_I)\quad{\rm and}\quad\widehat
E^\ast_I=(\widehat\Pi_I)^{-1}(E_I|_{V^\ast_I}).\eqno(3.73)
$$
As above we get a stratified Banach bundle system
$$
(\widehat {\mathcal E}^\ast, \widehat V^\ast)=\bigl \{(\widehat
E_I^\ast, \widehat V_I^\ast), \hat\pi_I,
 \hat\pi^I_J, \widehat\Pi_I, \widehat\Pi^I_J, \hat p_I,\Gamma_I\,\bigm|\,
 J\subset I\in{\mathcal N}\bigr\},\eqno(3.74)
$$
which is called the {\it restriction} of $(\widehat {\mathcal E},
\widehat V)$ to the open subset ${\cal W}^\ast\subset {\cal W}$.
Similarly, we can define its global section. Clearly, each global
section $\hat \sigma =\{(\hat\sigma)_I\,|\,I\in{\cal N}\}$ of
$(\widehat{\cal E},\widehat V)$ may restrict to a global section
of $(\widehat {\mathcal E}^\ast, \widehat V^\ast)$,
$\hat\sigma|_{\widehat V^\ast} =\{(\hat\sigma)_I|_{\widehat
V_I^\ast}\,|\, I\in{\cal N}\}$. As in Remark 3.13 we may also
consider the closure of $(\widehat {\mathcal E}^\ast, \widehat
V^\ast)$ in  $(\widehat {\mathcal E}^{\Gamma^+}, \widehat
W^{\Gamma^+})$,
$$
(Cl(\widehat {\mathcal E}^\ast), Cl(\widehat V^\ast))=\bigl
\{(Cl(\widehat E^\ast_I), Cl(\widehat V^\ast_I)), \hat\pi_I, \hat
\pi^I_J, \widehat\Pi_I, \widehat\Pi^I_J, p_I,\Gamma_I\,\bigm|\,
J\subset I\in{\mathcal N}\bigr\}
$$
and define its global section.

 Denote by $q=q_1+\cdots+ q_{n_3}$. Taking
$0<\eta\le\min\{\eta_i\,|\,1\le i\le n_3\}$,
 we have the obvious
pullback stratified Banach bundle system
$$
\bigl({\bf P}_1^\ast\widehat {\mathcal E}^\ast, \widehat
V^\ast\times {\bf B}_\eta({\mbox{\Bb R}}^q)\bigr) =\bigl \{({\bf
P}_1^\ast\widehat E_I^\ast, \widehat V_I^\ast\times {\bf
B}_\eta({\mbox{\Bb R}}^q)), \pi_I, \pi^I_J, \Pi_I,\Pi^I_J,
p_I,\Gamma_I\,\bigm|\, J\subset I\in{\mathcal N}\bigr\},
\eqno(3.75)
$$
where ${\bf P}_1$ are the projections to the first factor, and
$\hat\pi_I, \hat\pi^I_J, \hat\Pi_I,\hat\Pi^I_J, \hat p_I$ are
naturally pullbacks of those projections in (3.74). Consider its
 global section
$\Psi=\{\Psi_I\,|\, I\in{\cal N}\}$,
$$\Psi_I: \widehat
V_I^\ast\times {\bf B}_\eta({\mbox{\Bb R}}^q)\to {\bf
P}_1^\ast\widehat E_I^\ast,\; (\hat x_I, {\bf t})\mapsto
(\hat\partial_J)_I(\hat x_I)+
 \sum^{n_3}_{i=1}\sum^{q_i}_{j=1}t_{ij}(\hat s_{ij})_I(\hat
 x_I),\eqno(3.76)
$$
where ${\bf t}=\{t_{ij}| 1\le j\le q_i,\,1\le i\le n_3\}\in
{\mbox{\Bb R}}^q$.
  Clearly, $\Psi_I(\hat x_I, 0)=0$ for any zero $\hat x_I$ of
 $(\hat\partial_J)_I$ in $\widehat V_I$.
 \vspace{2mm}

\noindent{\bf Theorem 3.14.}\quad{\it If the above $\eta>0$ is
small enough, then the global section $\Psi=\bigl\{\Psi_I\,|\,
I\in{\cal N}\bigr\}$ is transversal the zero section at each zero
$(\hat x_I,{\bf t})\in \widehat V_I^\ast\times {\bf
B}_\eta(\mbox{\Bb R}^q)$. So for each $I\in{\cal N}$ the set
$$
\widehat\Omega_I(K_0, \eta):=\{(\hat x_I, {\bf t})\in\widehat
V_I^\ast\times {\bf B}_\eta({\mbox{\Bb R}}^q)\,|\, \Psi_I(\hat
x_I,{\bf t})=0\}
$$
 is a  stratified manifold of top dimension $2m+ 2c_1(A)+ 2(3-n)(g-1)+q$. Moreover, for any
$J\subset I\in{\cal N}$ the projection
$$
 \hat\pi^I_J: (\hat
\pi^I_J)^{-1}\bigl(\widehat\Omega_I(K_0, \eta)\bigr)\to {\rm
Im}(\hat\pi^I_J)\subset\widehat\Omega_J(K_0, \eta)
$$
is a $|\Gamma_I|/|\Gamma_J|$-fold PS covering.}\vspace{2mm}

\noindent{\it Proof.}\quad  We may extend $\Psi_I$ on $Cl(\widehat
V_I^\ast)\times {\bf B}_\eta({\mbox{\Bb R}}^q)$ naturally, still
denoted by $\Psi_I$. Let $(\hat u_I, {\bf t}^0)\in Cl(\widehat
V_I^\ast)\times {\bf B}_\eta({\mbox{\Bb R}}^q)$ be a zero of
$\Psi_I$. Then $u_I=\hat\pi_I(\hat u_I)\in\cup^{n_3}_{i=1}V_i^0$,
and by the choice of ${\cal W}^\ast$ above it
 sits in  $V^0_{i_s}\subset W_{i_s}$ for some $i_s\in I$.
By (3.49) let $\tilde u_I=q_I(\hat u_I)$ and denote by $\tilde
u_{i_s}$ the s-th component of $\tilde u_I$. Then in the chart
$\widehat{\phi_I\circ\lambda_I^s}$ in (3.44) we have
 $\hat u_I=\widehat{\phi_I\circ\lambda_I^s}(\tilde u_{i_s})$.
Note that
 $$\Psi_I(\hat x_I, {\bf t})=(\hat\partial_J)_I(\hat x_I)+
 \sum_{i\in I}\sum^{q_i}_{j=1}t_{ij}(\hat s_{ij})_I(\hat
 x_I)
$$
for any $(\hat x_I, {\bf t})\in\widehat V_I\times {\bf
B}_\eta({\mbox{\Bb R}}^q)$. In the charts of (3.44) and (3.53),
(3.71) gives the local expression of the section
$(\hat\partial_J)_I$, and as in (3.68) for each $i\in I$ and
$j=1,\cdots,q_i$ we have also  one of $(\hat s_{ij})_I$ ,
 $$
 O(\tilde u_{i_s})\to
\widetilde E_{i_s}|_{O(\tilde u_{i_s})},  \; \tilde x\mapsto
\tau_{ij}(\tilde x):=(\Phi_i\circ\Lambda_{i_si})^{-1}\bigl(\tilde
s_{ij}(\phi_i\circ\lambda_{i_si}(\tilde x))\bigr).
$$
Hereafter we need to shrink $O(\tilde u_{i_s})$ so that $O(\tilde
u_{i_s})\subset (\widehat{\phi_I\circ\lambda_I^s})^{-1}(\widehat
V_I)\subset\widetilde U^0_{i_s}$. Using these we can get the
following local expression of $\Psi_I$  in the natural pullback
charts,
\begin{eqnarray*}
\hspace{25mm}&&
 O(\tilde u_{i_s})\times {\bf B}_\eta({\mbox{\Bb R}}^q)\to{\bf
P}_1^\ast(\widetilde E_{i_s}|_{O(\tilde u_{i_s})}),\hspace{50mm}(3.77)\\
&&(\tilde x, {\bf t}) \mapsto\bar\partial_J(\tilde x) + \sum_{i\in
I}\sum^{q_i}_{j=1}t_{ij}\tau_{ij}(\tilde
x)\\
&&=\bar\partial_J(\tilde x) +\sum^{q_{i_s}}_{j=1}t_{i_sj} \tilde
s_{i_sj}(\tilde x)+ \sum_{i\in
I\setminus\{i_s\}}\sum^{q_i}_{j=1}t_{ij}\tau_{ij}(\tilde x).
\end{eqnarray*}
Let $\tilde u_{i_s}\in\widetilde V^0_{i_s}\cap\widetilde
W^{(u,v)}_{i_s}$ for some $|(u,v)|\le \delta(\epsilon_{\bf
f}^{(i_s)})$. By Lemma 3.3 we have the trivialization
representative of the restriction of the section in (3.77) to
$O(\tilde u_{i_s})\cap\widetilde W^{(u,v)}_{i_s}\times {\bf
B}_\eta(\mbox{\Bb R}^q)$,
\begin{eqnarray*}
\hspace{18mm}&&
 (O(\tilde
u_{i_s})\cap\widetilde W^{(u,v)}_{i_s})\times
{\bf B}_\eta(\mbox{\Bb R}^q)\to L^p_{k-1}(f^{(i_s)\ast}_{(u,v)}TM),\hspace{38mm}(3.78)\\
&&(\tilde x;{\bf t})\mapsto F^{(i_s)}_{(u,v)}(\tilde
x)+\sum^{q_{i_s}}_{j=1}t_{i_sj} s_{i_sj}^{(u,v)}(\tilde
x)+\sum_{i\in
I\setminus\{i_s\}}\sum^{q_i}_{j=1}t_{ij}\tau^{(u,v)}_{ij}(\tilde
x),
\end{eqnarray*}
where $\tau^{(u,v)}_{ij}$ are the corresponding trivialization
representatives of the sections $\tau_{ij}$ under (3.2).
 Now the tangent map of the map in
(3.78) at $(\tilde u_{i_s},{\bf t}^0)$ is given by
\begin{eqnarray*}
\hspace{13mm}&&
 T_{\tilde u_{i_s}}(O(\tilde
u_{i_s})\cap\widetilde W^{(u,v)}_{i_s})\times{\mbox{\Bb R}}^q\to
L^p_{k-1}(f^{(i_s)\ast}_{(u,v)}TM),\hspace{42mm}(3.79)\\
  && (\tilde \xi,{\bf v})\mapsto dF^{(i_s)}_{(u,v)}(\tilde
u_{i_s})(\tilde\xi)+ \sum^{q_{i_s}}_{j=1}v_{i_sj}
s^{(u,v)}_{i_sj}(\tilde u_{i_s})+ \sum^{q_{i_s}}_{j=1}t^0_{i_sj}
ds^{(u,v)}_{i_sj}(\tilde u_{i_s})(\tilde\xi)\\
 &&\hspace{15mm}+
\sum_{i\in
I\setminus\{i_s\}}\sum^{q_i}_{j=1}v_{ij}\tau^{(u,v)}_{ij}(\tilde
u_{i_s})+\sum_{i\in I\setminus\{i_s\}}
\sum^{q_i}_{j=1}t^0_{ij}d\tau^{(u,v)}_{ij}(\tilde
u_{i_s})(\tilde\xi).
\end{eqnarray*}
Since $Cl(\widehat V_I^\ast)\subset\widehat V_I^+\subset\widehat
W_I^+$  it follows from (3.76) that
$$
 {\cal Z}\bigl((\hat\partial_J)_I|_{Cl(\widehat V_I^\ast)}\bigr)=(\hat
\pi_I)^{-1}\bigl(Z(\bar\partial_J)\cap Cl(V_I^\ast)\bigr)
$$
is compact in $\widehat W_I^+$  for each $I\in{\cal N}$. Note that
$ds^{(u,v)}_{i_sj}(\tilde u_{i_s})=0$, $j=1,\cdots, q_{i_s}$. As
in Lemma 3.1 we can use this fact, (3.6), (3.14) and (3.15) to
prove that the map in (3.79) is surjective for $\eta>0$
 and ${\cal
W}^\ast$ sufficiently small. Then the standard arguments can
complete the proof.
 \hfill$\Box$\vspace{2mm}

 Since ${\cal N}$ is a finite set and each $\widehat V_I^\ast$ (and thus
$\widehat\Omega_I$) has only finitely many strata, Theorem 3.14
and the Sard-Smale theorem immediately gives:\vspace{2mm}

\noindent{\bf Corollary 3.15.}\quad{\it There exists a residual
subset ${\bf B}^{res}_\eta(\mbox{\Bb R}^q)\subset {\bf
B}_\eta(\mbox{\Bb R}^q)$ such that for each ${\bf t}\in {\bf
B}^{res}_\eta(\mbox{\Bb R}^q)$ the global section
 $\Psi^{({\bf t})}=\{\Psi^{({\bf t})}_I\,|\, I\in{\cal
N}\}$ of the PS bundle system $\bigl(\widehat {\mathcal E}^\ast,
\widehat V^\ast)\bigr)$ is transversal the zero section,  where $
\Psi^{({\bf t})}_I: \widehat V_I^\ast\to \widehat E_I^\ast,\; \hat
x_I \mapsto \Psi_I(\hat x_I,{\bf t})$. Therefore the set
$\widehat{\cal M}^{\bf t}_I(K_0):= (\Psi^{({\bf t})}_I)^{-1}(0)$
is a  stratified Banach manifold of dimension $2m+2c_1(A)+
2(3-n)(g-1)$.   It also holds that
\begin{description}
\item[(i)] The stratified Banach manifolds $\widehat{\mathcal
M}_I^{{\bf t}}(K_0)$ has no strata of codimension odd, and each
stratum of $\widehat{\cal M}^{\bf t}_I(K_0)$ of codimension $r$ is
exactly the intersection of $\widehat{\cal M}^{\bf t}_I(K_0)$ and
the stratum of $\widehat V_I^\ast$ of codimension $r$ for
$r=0,\cdots, 2m+ 2c_1(A)+ 2(3-n)(g-1)$.

\item[(ii)] The family $\widehat{\cal M}^{\bf
t}(K_0)=\bigl\{\widehat{\cal M}^{\bf t}_I(K_0)\,|\, I\in{\cal
N}\bigr\}$ is compatible in the sense that for any $J\subset
I\in{\cal N}$,
$$
 \hat\pi^I_J: (\hat \pi^I_J)^{-1}\bigl(\widehat{\cal M}^{\bf
t}_I(K_0)\bigr)\to {\rm Im}(\hat\pi^I_J)\subset\widehat{\cal
M}^{\bf t}_J(K_0)\eqno(3.80)
 $$
 is a $|\Gamma_I|/|\Gamma_J|$-fold partially smooth covering.

\item[(iii)]  For each $I\in{\cal N}$ and any two ${\bf t},\, {\bf
t}'\in {\bf B}^{res}_\eta(\mbox{\Bb R}^q)$, the cornered
stratified Banach manifolds $\widehat{\mathcal M}_I^{{\bf
t}}(K_0)$ and $\widehat{\mathcal M}_I^{{\bf t}'}(K_0)$ are
cobordant, and thus maps $\hat\pi_I: \widehat{\mathcal M}_I^{{\bf
t}}(K_0)\to {\mathcal W}^\ast$ and $\hat\pi_I: \widehat{\mathcal
M}_I^{{\bf t}'}(K_0)\to {\mathcal W}^\ast$ are also cobordant.
\end{description}}\vspace{2mm}

\noindent{We} can also show that each $\widehat{\cal M}^{\bf
t}_I(K_0)$ carries a natural orientation and all $\hat\pi^I_J$
perverse the orientations, cf., [LuT] for details. By the
definition of $V_I^\ast$ in (3.73) it follows from (3.72) that
$$
\|(\hat s_{ij})_I(x)\|\le C\;\forall x\in \widehat
V_I^\ast\eqno(3.81)
$$
and $i=1,\cdots,n_3$, $j=1,\cdots,q_i$. To understand the family
$\widehat{\cal M}^{\bf t}(K_0)=\bigl\{\widehat{\cal M}^{\bf
t}_I(K_0)\,|\, I\in{\cal N}\bigr\}$ we first consider the case
that $M$ is a closed manifold. \vspace{2mm}

\noindent{\bf  Proposition 3.16.}\quad{\it If $M$ is a closed
symplectic manifold, we may take $K_0=M$ and $i=1,\cdots,n_0$ in
the construction above.  Using some ideas in the proof of Theorem
4.1 in [LiuT1] we can derive from (3.81) and the compactness of
$$
\overline{\mathcal M}_{g,m}(M, J, A)= \bigcup_{I\in{\cal
N}}\hat\pi_I({\cal Z}\bigl((\hat\partial_J)_I|_{\widehat
V_I^\ast}\bigr))\eqno(3.82)
$$
that for any given small neighborhood ${\mathcal U}$ of
$\overline{\mathcal M}_{g,m}(M,J,A)$ there exists a small
$\delta>0$ such that the closure of $\cup_{I\in{\mathcal N}}
\hat\pi_I(\widehat{\mathcal M}_I^{{\bf t}}(K_0))$ may be contained
in ${\cal U}$ for any ${\bf t}\in {\bf B}_\delta(\mbox{\Bb R}^q)$.
 In particular there exists a positive number $\varepsilon\le\eta$ such that $\cup_{I\in{\mathcal N}}
\hat\pi_I(\widehat{\mathcal M}_I^{{\bf t}}(K_0))$ is compact for
each ${\bf t}\in {\bf B}_\varepsilon(\mbox{\Bb R}^q)$.
Consequently, the family $\widehat{\cal M}^{\bf
t}(K_0)=\bigl\{\widehat{\cal M}^{\bf t}_I(K_0)\,|\, I\in{\cal
N}\bigr\}$ is ``like'' an open cover of an oriented closed
manifold.  Especially, each $\widehat{\mathcal M}_I^{{\bf
t}}(K_0)$ is a finite set provided that $2m+ 2c_1(A)+
2(3-n)(g-1)=0$. The formal summation
$${\mathcal C}^{\bf t}:=\sum_{I\in{\mathcal
 N}}\frac{1}{|\Gamma_I|}\{\hat\pi_I:{\mathcal M}^{\bf t}_I(K_0)\to{\mathcal
 W}\}\quad\forall{\bf t}\in {\bf B}_\varepsilon^{res}(\mbox{\Bb R}^q),
 $$
gives rise to a family of  cobordant singular cycles in
 ${\mathcal W}$,  called the \emph{virtual moduli cycles} in ${\mathcal
 W}$ (cf. [LiuT1, LiuT2, LiuT3]). As explained  on page 65 of [LiuT1] the summation
precisely means that on the overlap of two pieces ${\mathcal
C}^{\bf t}(K_0)$ we only count them once. Hereafter we omit the
wide hat $\widehat{\quad}$ over $\widehat{\mathcal M}_I^{{\bf t}}$
and the dependence marks on $\omega,J,\mu$ and $A$ without
occurring of confusions.  }\vspace{2mm}

\noindent{\it Proof.}\quad Since $\cup_{I\in{\cal N}}V_I^\ast$ is
an open neighborhood of $\overline{\mathcal M}_{g,m}(M, J, A)$ in
${\cal W}$ we can take a positive $\varepsilon\le\eta$ such that
the closure of $\cup_{I\in{\mathcal N}}
\hat\pi_I(\widehat{\mathcal M}_I^{{\bf t}}(K_0))$ is contained in
$\cup_{I\in{\cal N}}V_I^\ast$ for any ${\bf t}\in {\bf
B}_\varepsilon(\mbox{\Bb R}^q)$. We shall prove that
$\cup_{I\in{\mathcal N}} \hat\pi_I(\widehat{\mathcal M}_I^{{\bf
t}}(K_0))$ is compact for any ${\bf t}\in {\bf
B}_\varepsilon(\mbox{\Bb R}^q)$. By the proof of Theorem 4.1 in
[LiuT1] each $\hat\pi_I(\widehat{\mathcal M}_I^{{\bf t}}(K_0))$
has the compact closure in ${\cal W}$. So the compact subset
$Cl\bigl(\cup_{I\in{\mathcal N}} \hat\pi_I(\widehat{\mathcal
M}_I^{{\bf t}}(K_0))\bigr)=\cup_{I\in{\mathcal N}}
Cl\bigl(\hat\pi_I(\widehat{\mathcal M}_I^{{\bf t}}(K_0))\bigr)$ is
contained in the open subset $\cup_{I\in{\cal N}}V_I^\ast$. Note
that any given point $x\in Cl\bigl(\hat\pi_I(\widehat{\mathcal
M}_I^{{\bf t}}(K_0))\bigr)\setminus \hat\pi_I(\widehat{\mathcal
M}_I^{{\bf t}}(K_0))\subset Cl(V_I^\ast)$ may be contained in some
$V_L^\ast$. So $Cl(V_I)\cap V_L\ne\emptyset$. It follows from (ii)
in Lemma 3.11 that $I\subset L$ or $L\subset I$. Let $I\subset L$
and $x=\hat\pi_L(\hat x_L)$ for some $\hat x_L\in\widehat
V_L^\ast$. By Remark 3.13 we may assume $x=\hat\pi_I(\hat x_I)$
for some $\hat x_I\in Cl(\widehat V_I^\ast)$. Since $\Psi^{({\bf
t})}=\{\Psi^{({\bf t})}_I\,|\, I\in{\cal N}\}$ may be extended
into a global section of the PS bundle system $\bigl(Cl(\widehat
{\mathcal E}^\ast), Cl(\widehat V^\ast)\bigr)$ naturally, we get
that $\Psi^{({\bf t})}_I(\hat x_I)=0$. Since $\hat\pi_I^L(\hat
x_L)=\hat x_I$, $\Psi^{({\bf t})}_L(\hat x_L)=0$ and thus
$x\in\hat\pi_L(\widehat{\mathcal M}_L^{{\bf t}}(K_0))$. This shows
that $Cl\bigl(\hat\pi_I(\widehat{\mathcal M}_I^{{\bf
t}}(K_0))\bigr)\subset\cup_{J\in{\mathcal N}}
\hat\pi_J(\widehat{\mathcal M}_J^{{\bf t}}(K_0))$. For the case
$L\subset I$ we have $\hat\pi^I_L(\hat x_I)=\hat x_L$. It follows
from $\Psi^{({\bf t})}_I(\hat x_I)=0$ that $\Psi^{({\bf
t})}_L(\hat x_L)=0$. Hence it also holds that $x=\hat\pi_L(\hat
x_L)\in \hat\pi_L(\widehat{\mathcal M}_L^{{\bf t}}(K_0))$.
 In summary we obtain that $\cup_{I\in{\mathcal
N}}Cl\bigr(\hat\pi_I(\widehat{\mathcal M}_I^{{\bf
t}}(K_0))\bigr)\subset\cup_{I\in{\mathcal N}}
\hat\pi_I(\widehat{\mathcal M}_I^{{\bf t}}(K_0))$ and hence
$\cup_{I\in{\mathcal N}} \hat\pi_I(\widehat{\mathcal M}_I^{{\bf
t}}(K_0))=\cup_{I\in{\mathcal
N}}Cl\bigr(\hat\pi_I(\widehat{\mathcal M}_I^{{\bf t}}(K_0))\bigr)$
is compact.

If $2m+ 2c_1(A)+ 2(3-n)(g-1)=0$ then  each $\widehat{\mathcal
M}_I^{{\bf t}}(K_0)$ is a manifold of dimension zero. Assume that
it contains infinitely many points $\hat x^{(k)}$, $k=1,
2,\cdots$. We may assume that $\{\hat\pi_I(\hat x^{(k)})\}$
converges to some $[{\bf f}]\in\hat\pi_L\bigl(\widehat{\mathcal
M}_L^{{\bf t}}(K_0)\bigr)$ because $\cup_{I\in{\mathcal N}}
\hat\pi_I(\widehat{\mathcal M}_I^{{\bf t}}(K_0))$ is compact. Note
that inverse image of each point by $\hat\pi_I$ contains at most
$|\Gamma_I|$ points and that $\{\hat\pi_I(\hat x^{(k)})\}$ are
contained in the closed subset $Cl(\widehat
V_I^\ast)\subset\widehat W_I^{\Gamma_I}$. After passing a
subsequence we may assume that $\{\hat x^{(k)}\}$ converges to
some $\hat x_I\in Cl\bigl(\widehat{\mathcal M}_I^{{\bf
t}}(K_0)\bigr)$. Then $\Psi^{({\bf t})}_I(\hat x_I)=0$. As showed
at the end of proof of Theorem 3.14 the section $\Psi^{({\bf
t})}_I$ is still transversal to the zero section near $\hat x_I$
in $Cl(\widehat V_I^\ast)$. This destroys the manifold structure
of ${\cal Z}\bigl((\hat\partial_J)_I|_{Cl(\widehat
V_I^\ast)}\bigr)$. The desired conclusion is proved.
\hfill$\Box$\vspace{2mm}

 In the case (3.82) does not hold. However, (3.19) and (3.20)
 imply that
\begin{eqnarray*}
&&\overline{\mathcal M}_{g,m}(M, J, A;K_3)=\Bigl\{[{\bf f}]\in
\bigcup_{I\in{\cal N}}\hat\pi_I({\cal
Z}\bigl((\hat\partial_J)_I|_{\widehat V_I^\ast}\bigr))\,\Bigm|\,
{\it Im}(f)\cap K_3\ne\emptyset\Bigr\}\quad{\rm and} \\
&& {\it Im}(f)\subset K_4\;\forall [{\bf f}]\in \bigcup_{I\in{\cal
N}}\hat\pi_I({\cal Z}\bigl((\hat\partial_J)_I|_{\widehat
V_I^\ast}\bigr)).
\end{eqnarray*}
Using (3.81) and the compactness of $\overline{\mathcal
M}_{g,m}(M, J, A;K_3)$ we can, as above, prove:\vspace{2mm}

\noindent{\bf Proposition 3.17.}\quad{\it For any given small
neighborhood ${\mathcal U}$ of $\overline{\mathcal M}_{g,m}(M,J,A,
K_3)$ there exists a small $\delta>0$ such that for any ${\bf
t}\in {\bf B}_\delta(\mbox{\Bb R}^q)$,
$$
\overline{\mathcal M}_{g,m}^{\bf t}(M, J, A;K_3):=\Bigl\{[{\bf
f}]\in\cup_{I\in{\mathcal N}} \hat\pi_I(\widehat{\mathcal
M}_I^{{\bf t}}(K_0))\,\Bigm|\, {\it Im}(f)\cap
K_3\ne\emptyset\Bigr\}\subset{\cal U}.
$$
 In particular there exists a positive constant $\varepsilon\le\eta$  such
that for any ${\bf t}\in{\bf B}_\varepsilon(\mbox{\Bb R}^q)$ the
following hold:\\
(i) The closure of $\overline{\mathcal M}_{g,m}^{\bf t}(M, J,
A;K_3)$ is contained in $\cup_{I\in{\cal N}}V_I^\ast$.\\
(ii) For any $[{\bf f}]\in Cl\bigl((\hat\pi_I(\widehat{\mathcal
M}_I^{{\bf t}}(K_0))\bigr)$ either $[{\bf f}]\in\overline{\mathcal
M}_{g,m}^{\bf t}(M, J, A;K_3)$ or ${\it Im}(f)\cap K_3=\emptyset$.
 Especially $\overline{\mathcal M}_{g,m}^{\bf t}(M, J, A;K_3)$ is compact. }\vspace{2mm}

In the present case the family $\widehat{\cal M}^{\bf
t}(K_0)=\bigl\{\widehat{\cal M}^{\bf t}_I(K_0)\,|\, I\in{\cal
N}\bigr\}$ is ``like'' an open cover of an open submanifold of an
oriented  geometrically bounded manifold. But the corresponding
formal summation
$${\mathcal C}^{\bf t}(K_0):=\sum_{I\in{\mathcal
 N}}\frac{1}{|\Gamma_I|}\{\hat\pi_I:{\mathcal M}^{\bf t}_I(K_0)\to{\mathcal
 W}^\ast\}\quad\forall{\bf t}\in {\bf B}_\varepsilon^{res}(\mbox{\Bb R}^q),\eqno(3.83)
 $$
gives rise to a family of  cobordant singular chains in
 ${\mathcal W}^\ast$  in general.
 We call it a {\it virtual moduli chain} in ${\mathcal W}$.
 However, it can still be used to define our GW-invariants in
next section.

\section{Gromov-Witten invariants}

In the first subsection we shall define our Gromov-Witten
invariants and prove their simple properties. In $\S4.2-\S4.6$ we
shall prove in details  the independence of our invariants for
various choices.

\subsection{Definition and simple properties}

From now on we always assume that $(M,\omega)$ is noncompact and
geometrically bounded. For $\mbox{\Bb K}=\mbox{\Bb C}, \mbox{\Bb
R}$ and $\mbox{\Bb Q}$ let $H^{\rm II}_\ast(M, \mbox{\Bb K})$ be
the {\it singular homology of the second kind} (or the homology
based on infinite locally finite singular chains), and $H^\ast(M,
\mbox{\Bb K})$ (resp. $H^\ast_c(M, \mbox{\Bb K})$) denote the
singular cohomology (resp. the singular cohomology with compact
supports) with coefficients in $\mbox{\Bb K}$. (see [Sc, Sk]). By
the deRham
 theorem we may also consider $H^\ast(M,\mbox{\Bb K})$ (resp. $H^\ast_c(M,\mbox{\Bb K})$)
as the deRham cohomology (resp. the deRham cohomology with compact
support) with coefficients in $\mbox{\Bb K}$. In particular
$H^\ast(M,\mbox{\Bb Q})$ (resp. $H^\ast_c(M,\mbox{\Bb Q})$)
consists of all deRham cohomology classes in $H^\ast(M,\mbox{\Bb
R})$ (resp. $H^\ast_c(M,\mbox{\Bb R})$) which take rational values
over all integral cycles.
 One has the Poincar\'e
duality isomorphisms
$$
PD: H_p(M, \mbox{\Bb K})\to H_c^{n-p}(M, \mbox{\Bb K})\quad{\rm
and}\quad PD^{\rm II}: H_p^{\rm II}(M, \mbox{\Bb K})\to H^{n-p}(M,
\mbox{\Bb K}).
$$
 The orientation determines  a class $[M]\in H^{\rm II}_{2n}(M,\mbox{\Bb Z})$, which is
 Poinca\'re dual to ${\bf 1}\in H^0(M,\mbox{\Bb Z})$, also
 called the \emph{fundamental cohomology class}.

 Consider the evaluation maps
$$
{\rm EV}_{g,m}=(\prod_{i=1}^m{\rm ev}_i):{\mathcal B}^M_{A,g,m}\to
M^m\quad{\rm and}\quad \Pi_{g,m}:{\mathcal
B}^M_{A,g,m}\to\overline{\mathcal
 M}_{g,m},
 $$
where ${\rm ev}_i([f,\Sigma,\bar{\bf z}])=f(z_i)$, $i=1,\cdots,m$
and $\Pi_{g,m}([f,\Sigma,\bar{\bf z}])=[\Sigma^\prime,\bar{\bf
 z}^\prime]$ is obtained by collapsing  components of $(\Sigma,\bar{\bf z})$
with genus $0$ and at most two special points.

Let $\{\alpha_i\}_{1\le i\le m}\subset H_c^\ast(M,\mbox{\Bb
Q})\cup H^\ast(M,\mbox{\Bb Q})$, and at least one of them, saying
$\alpha_1$, belong to $H_c^\ast(M,\mbox{\Bb Q})$. Later, if
$\alpha_i\in H_c^\ast(M,\mbox{\Bb Q})$ (resp. $H^\ast(M,\mbox{\Bb
Q})$) we understand the Poinca\'re duality of it to sit in
$H_\ast(M, \mbox{\Bb Q})$ (resp. $H_\ast^{\rm II}(M, \mbox{\Bb
Q})$.) Let $\bar{\alpha_i}: C(\alpha_i)\to M$ be the
representative cycle of the Poinca\'re duality of $\alpha_i$. For
the sake of simplicity we can assume that each $C(\alpha_i)$ is a
closed smooth manifold. Take a compact set $K_0$ in $M$ such that
$$
\bar{\alpha_1}(C(\alpha_1))\subset K_0.\eqno(4.1)
$$
Note that  each map ${\rm EV}_{g,m}\circ\hat\pi_I:{\mathcal
M}^{\bf t}_I(K_0)\to  M^m$ is partially smooth. For the virtual
moduli chain ${\mathcal C}^{\bf t}(K_0)$ in (3.83), we can jiggle
$\prod_{i=1}^m\bar\alpha_i:\prod_{i=1}^mC(\alpha_i)\to M^m$ to
make it intersect each ${\rm EV}_{g,m}\circ\hat\pi_I$
transversally. In this case we also say it to be transversal to
$$
{\rm EV}_{g,m}\circ{\mathcal C}^{\bf t}(K_0):=\sum_{I\in{\mathcal
 N}}\frac{1}{|\Gamma_I|}\{{\rm EV}_{g,m}\circ\hat\pi_I:{\mathcal M}^{\bf t}_I(K_0)\to
  M^m\}.\eqno(4.2)
$$

Define $C(K_0;\{\bar\alpha_i\}^m_{i=1};\omega,\mu,J)_I^{\bf t}$ to
be the fiber product
$$
\{(x,u)\in \prod_{i=1}^m C(\alpha_i)\times {\mathcal M}^{\bf
t}_I(K_0)\;|\; \prod_{i=1}^m\bar\alpha_i(x)={\rm
EV}_{g,m}\circ\hat\pi_I(u) \}\eqno(4.3)
$$
It is either empty or a stratified Banach manifold of dimension
$r=2c_1(M)(A)+ 2(3-n)(g-1)+ 2m-\sum^m_{i=1}\deg\alpha_i$ in the
case $r\ge 0$. \emph{The key is that (3.16)-(3.18) and the choices
of ${\bf f}^{(i)}$ imply:}
$$
C(K_0;\{\bar\alpha_i\}^m_{i=1};\omega,\mu,J)_I^{\bf
t}\ne\emptyset\Longrightarrow i\le n_1\;{\rm for}\;{\rm any}\;i\in
I.\eqno(4.4)
$$
 Denote by
$$
\Pi^I_{g,m}:C(K_0;\{\bar\alpha_i\}^m_{i=1};\omega,\mu,J,A)_I^{\bf t}\to \overline{\mathcal
 M}_{g,m},\;(x,u)\mapsto \Pi_{g,m}\circ\hat\pi_I(u).
 $$
For any $L\subset I\in{\cal N}$, the map (3.80) naturally induces
a $|\Gamma_I|/|\Gamma_L|$-fold partially smooth covering,
$$
 P^I_L: (P^I_L)^{-1}\bigl( C(K_0;\{\bar\alpha_i\}^m_{i=1};\omega,\mu,J)_I^{\bf
t}\bigr)\to {\rm Im}(P^I_J)\subset
C(K_0;\{\bar\alpha_i\}^m_{i=1};\omega,\mu,J)_L^{\bf t}
 $$
This is, the family
$C(K_0;\{\bar\alpha_i\}^m_{i=1};\omega,\mu,J)^{\bf
t}=\bigl\{C(K_0;\{\bar\alpha_i\}^m_{i=1};\omega,\mu,J)_I^{\bf
t}\,|\, I\in{\cal N}\bigr\}$ is compatible.\vspace{2mm}

\noindent{\bf Proposition 4.1.}\quad{\it The set $\cup_{I\in{\cal
N}}\Pi^I_{g,m}(C(K_0;\{\bar\alpha_i\}^m_{i=1};\omega,\mu,J)_I^{\bf
t})$ is compact in $\overline{\mathcal M}_{g,m}$. As in
Proposition 3.16, for any ${\bf t}\in{\bf
B}^{res}_\varepsilon(\mbox{\Bb R}^q)$ the formal summation
\begin{eqnarray*}
\hspace{22mm}&&\quad
 {\mathcal C}_g^{\bf
t}(K_0;\{\bar\alpha_i\}^m_{i=1};\omega,\mu,J,A)\hspace{68mm}(4.5)\\
&&=\sum_{I\in{\mathcal
 N}}\frac{1}{|\Gamma_I|}\{\Pi^I_{g,m}:C(K_0;
 \{\bar\alpha_i\}^m_{i=1};\omega,\mu,J,A)_I^{\bf t}\to
 \overline{\mathcal  M}_{g,m}\}
 \end{eqnarray*}
is a rational \emph{cycle} in $\overline{\mathcal M}_{g,m}$ of
dimension $r=2c_1(M)(A)+ 2(3-n)(g-1)+ 2m-\sum^m_{i=1}\deg\alpha_i$
if $r\ge 0$.  Moreover, the cobordant virtual moduli chains in
(3.83) yield the cobordant virtual moduli cycles in (4.5). As in
proof of Proposition 3.16 we can use (4.4) to prove that
$\cup_{I\in{\cal N}}C(K_0;
 \{\bar\alpha_i\}^m_{i=1};\omega,\mu,J,A)_I^{\bf t}$ is
 a finite set for $r=0$.}\vspace{2mm}

\noindent{\it Proof.}\quad We only need to prove that
$\cup_{I\in{\cal
N}}\Pi^I_{g,m}(C(K_0;\{\bar\alpha_i\}^m_{i=1};\omega,\mu,J)_I^{\bf
t})$ is compact. Since both $\prod_{i=1}^m C(\alpha_i)$ and
$Cl\bigl(\hat\pi_I({\mathcal M}^{\bf t}_I(K_0))\bigr)$ are
compact,
$Cl\bigl(\Pi^I_{g,m}(C(K_0;\{\bar\alpha_i\}^m_{i=1};\omega,\mu,J)_I^{\bf
t})\bigr)$ is also compact. Let $\Pi^I_{g,m}(x_k, u_k)$ be a
sequence converging to $w\in
Cl\bigl(\Pi^I_{g,m}(C(K_0;\{\bar\alpha_i\}^m_{i=1};\omega,\mu,J)_I^{\bf
t})\bigr)$. Here $(x_k, u_k)\in
C(K_0;\{\bar\alpha_i\}^m_{i=1};\omega,\mu,J)_I^{\bf t}$ for
$k=1,2\cdots$. Passing to subsequences we may assume that $x_k\to
x$ and $\hat\pi_I(u_k)\to [{\bf f}]\in Cl\bigl(\hat\pi_I({\mathcal
M}^{\bf t}_I(K_0))\bigr)$. By (4.4), $i\le n_1$ for any $i\in I$.
It follows from (3.19) that ${\it Im}(f)\subset K_2\subset K_3$.
By Proposition 3.17(ii), $[{\bf f}]\in\overline{\mathcal
M}_{g,m}^{\bf t}(M, J, A;K_3)\subset\cup_{L\in{\mathcal N}}
\hat\pi_L(\widehat{\mathcal M}_L^{{\bf t}}(K_0))$. So there exists
a $u$ in some $\widehat{\mathcal M}_J^{{\bf t}}(K_0)$ such that
$[{\bf f}]=\hat\pi_J(u)$. By (4.3),
$$
\prod_{i=1}^m\bar\alpha_i(x_k)={\rm
EV}_{g,m}(\hat\pi_I(u_k)),\;k=1,2,\cdots.
$$
Setting $k\to\infty$ we get
$$
\prod_{i=1}^m\bar\alpha_i(x)={\rm EV}_{g,m}([{\bf f}])={\rm
EV}_{g,m}(\hat\pi_L(u)).
$$
So $(x,u)\in C(K_0;\{\bar\alpha_i\}^m_{i=1};\omega,\mu,J)_L^{\bf
t}$ and $\Pi^L_{g,m}(x,u)=\Pi_{g,m}(\hat\pi_L(u))=\Pi_{g,m}([{\bf
f}])$. Moreover,
$$
w=\lim_{k\to\infty}\Pi^I_{g,m}(x_k,u_k)=
\lim_{k\to\infty}\Pi_{g,m}(\hat\pi_I(u_k))=\Pi_{g,m}([{\bf
f}]).
$$
Hence
$w=\Pi^L_{g,m}(x,u)\in\Pi^L_{g,m}(C(K_0;\{\bar\alpha_i\}^m_{i=1};\omega,\mu,J)_L^{\bf
t})$. This implies that
$$Cl\bigl(\Pi^I_{g,m}(C(K_0;\{\bar\alpha_i\}^m_{i=1};\omega,\mu,J)_I^{\bf
t})\bigr)\subset\cup_{I\in{\cal
N}}\Pi^I_{g,m}(C(K_0;\{\bar\alpha_i\}^m_{i=1};\omega,\mu,J)_I^{\bf
t})
$$
and thus that
$$\cup_{I\in{\cal
N}}\Pi^I_{g,m}(C(K_0;\{\bar\alpha_i\}^m_{i=1};\omega,\mu,J)_I^{\bf
t})=\cup_{I\in{\cal
N}}Cl\bigl(\Pi^I_{g,m}(C(K_0;\{\bar\alpha_i\}^m_{i=1};\omega,\mu,J)_I^{\bf
t})\bigr)$$
 is compact.\hfill$\Box$\vspace{2mm}

 Clearly,  two
different ${\bf t}, {\bf t}^\prime\in{\bf
B}_{\varepsilon}^{res}(\mbox{\Bb R}^q)$ give cobordant rational
cycles ${\mathcal C}_g^{\bf
t}(K_0;\{\bar\alpha_i\}^m_{i=1};\omega,\mu,J,A)$ and ${\mathcal
C}_g^{{\bf
t}^\prime}(K_0;\{\bar\alpha_i\}^m_{i=1};\omega,\mu,J,A)$.
Moreover, the choices of cycle representatives
$\{\bar\alpha_i\}^m_{i=1}$ satisfying (4.1) also yield cobordant
rational cycles ${\mathcal C}_g^{\bf
t}(K_0;\{\bar\alpha_i\}^m_{i=1};\omega,\mu,J,A)$. So they
determine a unique homology class in $H_r(\overline{\mathcal
M}_{g,m}, \mbox{\Bb Q})$, denoted by $[{\mathcal C}_g^{\bf
t}(K_0;\{\bar\alpha_i\}^m_{i=1};\omega,\mu,J,A)]$ or
$$
[{\mathcal C}_g^{\bf t}(K_0;\{\alpha_i\}^m_{i=1};\omega,\mu,J,A)]
 \eqno(4.6)
$$
because it is independent of choices of ${\bf t}\in{\bf
B}_{\varepsilon}^{res}(\mbox{\Bb R}^q)$ and of ones of cycle
representatives $\bar\alpha$ of $PD(\alpha)$ satisfying (4.1).
 One can easily prove that
shrinking ${\cal W}^\ast$ and thus $\varepsilon>0$ gives the same
class as in (4.6). More generally, it only depends on cobordant
class of virtual moduli chains in (3.83). Furthermore we can prove
that this class does not depend on all related choices. Namely,
\begin{itemize}
\item it is independent of choices of section $s_i$;

\item it is independent of choices of
 $W_i$;
\item it is independent of choices of
 $K_0$;
\item it is independent of choices of  almost complex structures
$J\in{\mathcal J}(M,\omega,\mu)$;

\item it is independent of the weak deformation of the
geometrically bounded symplectic manifold $(M,\omega, J,\mu)$.
\end{itemize}
We first admit them and put off their proofs to $\S4.2-\S4.6$. Let
$2g+m\ge 3$. If $r\ge 0$ we call the Poincar\'e dual of the
homology class in (4.6),
$$
I^{(M,\mu, J)}_{A,g,m}(\alpha_1\otimes\cdots\otimes\alpha_m):=PD(
[{\mathcal C}_g^{\bf
t}(K_0;\{\alpha_i\}^m_{i=1};\omega,\mu,J,A)])
$$
 the {\it Gromov-Witten class} in terms of [KM].
It is clear that we can understand and write \linebreak
$[{\mathcal C}_g^{\bf
t}(K_0;\{\alpha_i\}^m_{i=1};\omega,\mu,J,A)]$ as $[{\mathcal
C}_g^{\bf t}(K_0; \otimes^m_{i=1}\alpha_i;\omega,\mu,J,A)]$ or
\begin{eqnarray*}
&&[{\mathcal C}_g^{\bf
t}(K_0;\{\alpha_i\}^{m_1}_{i=1},\{\alpha_i\}^{m_2}_{i=m_1+1},\cdots,
\{\alpha_i\}^{m_k}_{i=m_{k-1}+1};\omega,\mu,J,A)]\quad{\rm or}\\
&& [{\mathcal C}_g^{\bf
t}(K_0;\otimes^{m_1}_{i=1}\alpha_i,\otimes^{m_2}_{i=m_1+1}\alpha_i,
 \cdots, \otimes^{m_k}_{i=m_{k-1}+1}\alpha_i;\omega,\mu,J,A)]
\end{eqnarray*}
 for $1\le m_1<m_2\cdots <m_k=m$.

Let $2g+m\ge 3$. If (1.5) is satisfied we define the Gromov-Witten
invariants
 $$\left.\begin{array}{ll}
{\mathcal G}{\mathcal W}^{(\omega,\mu,
J)}_{A,g,m}(\kappa;\alpha_1,\cdots,\alpha_m) :=\langle
PD(\kappa),[{\mathcal C}_g^{\bf
t}(K_0;\{\alpha_i\}^m_{i=1};\omega,\mu,J,A)]\rangle,
\end{array}\right.\eqno(4.7)
$$
and define ${\mathcal G}{\mathcal W}^{(\omega,\mu,
J)}_{A,g,m}(\kappa;\alpha_1,\cdots,\alpha_m)=0$ if (1.5) is not
satisfied. Here $\langle\cdot,\cdot\rangle$ is the usual pair
between homology and cohomology.

Note that the space $\overline{\mathcal M}_{0,m}$ reduces to a
single point for $m<3$ and that it is a smooth oriented compact
manifold of dimension $2m-6$ for $m\ge 3$ ([HS2]). Then for $m\ge
3$ it is easy to see that the Gromov-Witten invariants ${\mathcal
G}{\mathcal W}^{(\omega,\mu,
J)}_{A,0,m}(\kappa;\alpha_1,\cdots,\alpha_m)$ is equal to the
intersection number of $(\Pi_{0,m}\times{\rm
EV}_{0,m})\circ{\mathcal C}^{\bf t}(K_0)$ with the class
$\kappa\times\prod^m_{i=1}PD(\alpha_i)$ in
$H_\ast(\overline{\mathcal M}_{0,m}\times M^m,\mbox{\Bb Q})$,
i.e.,
\begin{eqnarray*}
\hspace{22mm}&&\quad{\mathcal G}{\mathcal W}^{(\omega,\mu,
J)}_{A,0,m}(\kappa;\alpha_1,\cdots,\alpha_m)\hspace{65mm}(4.8)\\
&&=(\Pi_{0,m}\times{\rm EV}_{0,m})\circ{\mathcal C}^{\bf
t}(K_0)\cdot_{\overline{\mathcal M}_{0,m}\times
M^m}(\kappa\times\prod^m_{i=1}PD(\alpha_i)).
\end{eqnarray*}
These suggest an equivalent definition of the invariants in (4.7)
by the integration theory at the end of the subsection.

 Let the rational number in (1.6)
be chosen as one given by (4.7) then these give the claims
(i)-(iii) in Theorem 1.1 immediately. As to (iv), its proof is as
follows:

 By (1.4) we can
choose a continuous path $(\psi_t)_{t\in [0,1]}$ in ${\rm
Symp}_0(M,\omega)$ with respect to the $C^\infty$-strong topology
to connect $id_M=\psi_0$ to $\psi_1=\psi$. Then $(\psi_t^\ast
J)_{t\in [0,1]}$ and $(\psi_t^\ast\mu)_{t\in [0,1]}$ are
respectively continuous paths in ${\mathcal J}(M,\omega)$ and
${\mathcal R}(M)$ with respect to the $C^\infty$-strong topology.
It easily follows that each $(M,\omega, \psi^\ast_t
J,\psi^\ast_t\mu)$ is geometrically bounded and that $(M,\omega,
\psi^\ast_t J,\psi^\ast_t\mu)_{t\in [0,1]}$ is a strong
deformation of the geometrically bounded symplectic manifold
$(M,\omega, J_0,\mu_0)=(M,\omega, J,\mu)$. Therefore Theorem
1.1(iii) leads to Theorem 1.1(iv).
 \vspace{2mm}

 Since $M$ is a connected, oriented
and noncompact manifold of dimension $2n$ one always has
$$
H^0_c(M, \mbox{\Bb R})=0=H^{2n}(M, \mbox{\Bb R})\quad{\rm
and}\quad H^0(M, \mbox{\Bb R})\cong \mbox{\Bb R}\cong H^{2n}_c(M,
\mbox{\Bb R}).\eqno(4.9)
$$
 Later without special statements we always denote by ${\bf 1}\in H^0(M,\mbox{\Bb Q})$
the fundamental class.  From the definition we easily obtain the
following simple properties.\vspace{2mm}

\noindent{\bf Theorem 4.2.}\quad{\it Let $\alpha_1\in
H_c^\ast(M,\mbox{\Bb Q})$ and $\alpha_2,\cdots,\alpha_m\in
H^\ast(M,\mbox{\Bb Q})$. Then
$$
{\mathcal G}{\mathcal W}^{(\omega, \mu,
J)}_{A,g,m}([\overline{\mathcal
M}_{g,m}];\alpha_1,\cdots,\alpha_{m-1}, {\bf 1})=0\; {\rm if}\;
m\ge 4\; {\rm or}\;A\ne 0\,{\rm and}\, m\ge 1.
$$
$$
{\mathcal G}{\mathcal W}^{(\omega, \mu,
J)}_{A,0,m}([\overline{\mathcal
M}_{0,m}];\alpha_1,\cdots,\alpha_m)=\left\{\begin{array}{lll}
\int_M\alpha^\ast_1\wedge\alpha_2^\ast\wedge\alpha_3^\ast\;\!\!\!&&{\rm
for}\;m=3, A=0\\
0\;\!\!\!&&{\rm for}\;m>3, A=0
\end{array}\right.
$$
Here $\alpha_i^\ast$ are closed form representatives of $\alpha_i$
and $\alpha_1$ has a compact set.}\vspace{2mm}

\noindent{\it Proof.}\quad The first result is a direct
consequence of the reduction formulas. We put it here so as to
compare it with the second conclusion. Theorem 5.9 gives rise to
\begin{eqnarray*}
&&\quad {\mathcal G}{\mathcal W}^{(\omega, \mu,
J)}_{A,g,m}([\overline{\mathcal
M}_{0,m}];\alpha_1,\cdots,\alpha_{m-1}, {\bf 1})\\
&&={\mathcal G}{\mathcal W}^{(\omega, \mu,
J)}_{A,g,m-1}(({\mathcal F}_m)_\ast([\overline{\mathcal
M}_{g,m}]);\alpha_1,\cdots,\alpha_{m-1})\\
&&={\mathcal G}{\mathcal W}^{(\omega, \mu,
J)}_{A,g,m-1}(0;\alpha_1,\cdots,\alpha_{m-1})=0
\end{eqnarray*}
because $({\mathcal F}_m)_\ast([\overline{\mathcal M}_{g,m}])\in
H_{6g-6+2m}(\overline{\mathcal M}_{g,m-1};\mbox{\Bb Q})$ must be
zero.

Next we prove the second result. Assume that $m\ge 3$ and $A=0$.
One easily checks that $\overline{\mathcal M}_{0,m}(M,J,
0)=\overline{\mathcal M}_{0,m}\times M$ and that the differential
$D\bar\partial_J$ at any point $[{\bf f}]\in\overline{\mathcal
M}_{0,m}(M,J, 0)$ is surjective. So for any compact submanifold
$K_0\subset M$ of codimension $0$ the associated virtual moduli
chain with $\overline{\mathcal M}_{0,m}(M,J, 0;K_3)$ (constructed
by Liu-Tian's method as above) can be chosen as
$$
{\mathcal C}(\overline{\mathcal
M}_{0,m}\times K_3)=\{\hat\pi:\overline{\mathcal M}_{0,m}\times
{\rm Int}(K_3)\to {\mathcal W}\},
$$
where $\hat\pi(([\Sigma,\bar{\bf z}], x))=[\Sigma,\bar{\bf z},
x]$. Note that ${\rm EV}_{0,m}\circ\hat\pi(([\Sigma,\bar{\bf z}],
x))=(x,\cdots,x)\in M^m$ and that
$\Pi_{0,m}\circ\hat\pi(([\Sigma,\bar{\bf z}],x))=[\Sigma,\bar{\bf
z}]$. By (4.8) the GW-invariant ${\mathcal G}{\mathcal
W}^{(\omega, \mu, J)}_{0,0,m}([\overline{\mathcal
M}_{0,m}];\alpha_1,\cdots,\alpha_m)$ is exactly the intersection
number of $(\Pi_{0,m}\times{\rm EV}_{0,m})\circ\hat\pi$ and the
product map
$$
id\times\prod^m_{i=1}\bar\alpha_i:\overline{\mathcal M}_{0,m}\times
C(\alpha_i)\to\overline{\mathcal M}_{0,m}\times M^m,
$$
which is equal to the intersection number of
the map $\prod^m_{i=1}\bar\alpha_i: C(\alpha_i)\to M^m$ with
$\triangle(M^m)=\{(x,\cdots,x)\in M^m\,|\, x\in M\}\subset M^m$.
Therefore
$$
{\mathcal G}{\mathcal W}^{(\omega, \mu,
J)}_{0,0,m}([\overline{\mathcal
M}_{0,m}];\alpha_1,\cdots,\alpha_m)=\int_M\wedge^m_{i=1}\alpha_i^\ast.
$$
Here $\alpha_i^\ast$ are closed form representatives of $\alpha_i$
and $\alpha_1$ has a compact support contained in $K_0$. By the
dimension condition (1.5) we have
$\sum^m_{i=1}\deg\alpha_i^\ast=2m-6+2n$. It follows from (4.8)
that
$$
{\mathcal
G}{\mathcal W}^{(\omega, \mu, J)}_{0,0,m}([\overline{\mathcal
M}_{0,m}];\alpha_1,\cdots,\alpha_m)=\left\{\begin{array}{ll}0\hspace{29mm}
{\rm if}\;m>3\\
\int_M\alpha_1^\ast\wedge\alpha_2^\ast\wedge\alpha_3^\ast\quad{\rm
if}\;m=3.
\end{array}\right.
$$
\hfill$\Box$\vspace{2mm}

For any $0\le p\le 2n$ there  exists a natural homomorphism
   $$
   {\bf I}^p: H^p_c(M, \mbox{\Bb C})\to H^p(M, \mbox{\Bb  C}),\eqno(4.10)
   $$
  which is given by $[\eta]_c\mapsto [\eta]$ when $H^\ast_c(M, \mbox{\Bb C})$ (resp.
  $H^\ast(M, \mbox{\Bb  C})$) is identified with the de-Rham cohomology with compact
  supports (resp. the de-Rham cohomology.) Here $\eta$ is a closed $p$-form on $M$ with
  compact support.  Clearly, ${\bf I}^p$ maps $H^p_c(M, \mbox{\Bb R})$
  (resp. $H^p_c(M, \mbox{\Bb Q})$)    to $H^p(M, \mbox{\Bb R})$ (resp.
   $H^p(M, \mbox{\Bb Q})$).
By (4.10), both ${\bf I}^0$ and ${\bf I}^{2n}$ are zero
homomorphisms.  It is easy to prove that $\alpha\in{\rm Ker}({\bf
I}^p)$ if and only if $\alpha$ has a closed form representative
(and thus all closed form representatives) that is exact on $M$.
For any $0<p<2n$, we denote by
$$
\widehat H^p_c(M,\mbox{\Bb C}):=H^p_c(M,\mbox{\Bb C})/{\rm
Ker}({\bf I}^p)
$$
 and by  $\widehat\alpha=\alpha+ {\rm Ker}({\rm
I}^p)$ the equivalence  class of $\alpha$ in the quotient space
$\widehat H^p_c(M,\mbox{\Bb C})$. By (4.6) it is easily seen that
$$
{\mathcal G}{\mathcal W}^{(\omega, \mu,
J)}_{A,g,m}(\kappa;\alpha_1,{\bf
I}^\ast(\alpha_2),\cdots,\alpha_m)={\mathcal G}{\mathcal
W}^{(\omega,
\mu,J)}_{A,g,m}(\kappa;\alpha_1,\alpha_2,\cdots,\alpha_m)\eqno(4.11)
$$
if $\alpha_1,\alpha_2\in H^\ast_c(M, \mbox{\Bb Q})$. However, even
if $\alpha_2\in H^q_c(M,\mbox{\Bb C})$, $0<q<2n$,  but no one of
$\{\alpha_i\}_{2<i\le m}$
 belongs to $H^\ast_c(M,\mbox{\Bb C})$, we cannot   define
${\mathcal G}{\mathcal W}^{(\omega, \mu, J)}_{A,g,m}(\kappa; {\bf
I}^p( \alpha_1),{\bf I}^q(\alpha_2),\cdots,\alpha_m)$ yet. Let
$\alpha_1$ and $\alpha_1^\prime$ be two representatives of
$\widehat\alpha_1\in\widehat H^p_c(M,\mbox{\Bb C})$. Then
$\alpha_1-\alpha_1^\prime\in {\rm Ker}({\bf I}^p)$. Since
${\mathcal G}{\mathcal W}^{(\omega, \mu,
J)}_{A,g,m}(\kappa;\alpha_1,\cdots,\alpha_m)$ is multilinear on
$\alpha_1,\cdots,\alpha_m,\kappa$ it follows from (4.11) that
\vspace{2mm}

\noindent{\bf Definition and Proposition 4.3.}\quad{\it Assume
that one of the classes $\alpha_2,\cdots,\alpha_m$ also belongs to
$H_c^\ast(M,\mbox{\Bb Q})$. Then
$$
{\mathcal G}{\mathcal W}^{(\omega, \mu,
J)}_{A,g,m}(\kappa;\alpha_1,\alpha_2,\cdots,\alpha_m)={\mathcal
G}{\mathcal W}^{(\omega, \mu,
J)}_{A,g,m}(\kappa;\alpha_1^\prime,\alpha_2,\cdots,\alpha_m).
$$
In particular ${\mathcal G}{\mathcal W}^{(\omega, \mu,
J)}_{A,g,m}(\kappa;\{\alpha_i\}_{1\le i\le m})=0$ if
$\alpha_1\in{\rm Ker}({\bf I}^p)$.
 Thus if $m\ge 2$ one may define for any representatives $\alpha_i$ of $\widehat\alpha_i$,
$i=1,\cdots,k$,
$${\mathcal G}{\mathcal W}^{(\omega, \mu,
J)}_{A,g,m}(\kappa;\{\widehat\alpha_i\}_{1\le i\le k},
,\{\alpha_i\}_{k+1\le i\le m}):={\mathcal G}{\mathcal W}^{(\omega,
\mu, J)}_{A,g,m}(\kappa;\{\alpha_i\}_{1\le i\le m}).
$$
Consequently, our Gromov-Witten invariants may
descend to $\widehat H_c^\ast(M,\mbox{\Bb Q})^{\otimes m}$ as long
as $m\ge 2$. }\vspace{2mm}

\noindent{\bf An equivalent definition.}\quad It is well-known
that the intersection theory is closely related to  integration.
The latter is more convenient in computation (see \S5.3). In the
original version of this paper we use integration theory to define
our invariants. By de Rham's theorem for manifolds (cf. Th3.1 and
the last remark (b) in Appendix of [Ma]),
 for every class  $\alpha\in H_c^\ast(M,\mbox{\Bb Q})\cup H^\ast(M,\mbox{\Bb Q})$,
  we may always choose a closed representative form  $\alpha^\ast$
 on $M$  such that its integral over
 every smooth integral cycle is an rational number.
Similarly, since the orbifold is an rational (co)homology
manifold, using Satake's de Rham's theorem for orbifolds (cf.
[Sat]) we may also take for a given $\kappa\in
H_\ast(\overline{\mathcal
 M}_{g,m},\mbox{\Bb Q})$  a closed representative form $\kappa^\ast$ on
$\overline{\mathcal  M}_{g,m}$ P\'oincare dual to $\kappa$
 such that the integral of $\kappa^\ast$ over
 every smooth integral cycle is a rational number.
Let $\{\alpha_i\}_{1\le i\le m}\subset H_c^\ast(M,\mbox{\Bb
Q})\cup H^\ast(M,\mbox{\Bb Q})$, and at least one of them, saying
$\alpha_1$, belong to $H_c^\ast(M,\mbox{\Bb Q})$. Then we can
choose their closed representative forms $\alpha_i^\ast$,
$i=1,\cdots, m$ and a compact set $K_0$ in $M$ such that
$$
{\rm supp}(\wedge^m_{i=1}\alpha^\ast_i)\subset K_0.\eqno(4.12)
$$
By the construction of $C(K_0;
\{\bar\alpha_i\}^m_{i=1};\omega,\mu,J,A)_I^{\bf t}$ in (4.3)
 it is not hard to check that
 \begin{eqnarray*}
\hspace{28mm}&&\quad\int_{C_g(K_0;\{\alpha_i\}^m_{i=1};\omega,\mu,J,
  A)_I^{\bf t}}(\Pi^I_{g,m})^\ast\kappa^\ast\hspace{50mm}(4.13)\\
 &&=\int_{{\mathcal M}^{\bf t}_I(K_0)}\bigl((\Pi_{g,m}\times {\rm
EV}_{g,m})\circ\hat\pi_I\bigl)^\ast(\kappa^\ast\times\prod^m_{i=1}\alpha_i^\ast)
\end{eqnarray*}
for each $I\in{\mathcal N}$. Here because of (1.5) both
integrations are taken on the top strata of two stratified
manifolds in an obvious sense and exist since
$C_g(K_0;\{\alpha_i\}^m_{i=1};\omega,\mu,J,
  A)_I^{\bf t}$ and ${\mathcal M}^{\bf t}_I(K_0)$ have compact
  closures.
 By (4.7)
 \begin{eqnarray*}
&&\quad {\mathcal G}{\mathcal W}^{(\omega,\mu,
J)}_{A,g,m}(\kappa;\alpha_1,\cdots,\alpha_m)\\
&&=\langle
PD(\kappa),[{\mathcal C}_g^{\bf
t}(K_0;\{\alpha_i\}^m_{i=1};\omega,\mu,J,A)]\rangle\\
&&=\int_{{\mathcal C}_g^{\bf
t}(K_0;\{\alpha_i\}^m_{i=1};\omega,\mu,J,A)}\kappa^\ast\\
&&``="\sum_{I\in{\mathcal
 N}}\frac{1}{|\Gamma_I|}\int_{C(K_0;
 \{\bar\alpha_i\}^m_{i=1};\omega,\mu,J,A)_I^{\bf
 t}}(\Pi^I_{g,m})^\ast\kappa^\ast.
\end{eqnarray*}
This and (4.13) gives
\begin{eqnarray*}
\hspace{40mm}&&\quad
 {\mathcal G}{\mathcal
W}^{(\omega,\mu, J)}_{A,g,m}(\kappa;\alpha_1,\cdots,\alpha_m)\hspace{48mm}(4.14)\\
&&``="\sum_{I\in{\mathcal
 N}}\frac{1}{|\Gamma_I|}
\int_{{\mathcal M}^{\bf t}_I(K_0)}\bigl((\Pi_{g,m}\times {\rm
EV}_{g,m})\circ\hat\pi_I\bigl)^\ast(\kappa^\ast\times\prod^m_{i=1}\alpha_i^\ast).
\end{eqnarray*}
As explained for the formal summation in (4.5) we here write
$``\sum_{I\in{\mathcal N}}"$ to means that we only count these
overlap parts once in the summation of right side of (4.14).
Precisely, for $I, J\in{\mathcal N}$
 with $J\hookrightarrow I$ it is claimed in Corollary 3.15 that
 $$\hat\pi^I_J: (\hat\pi^I_J)^{-1}({\mathcal M}_J^{{\bf
 t}}(K_0))\subset{\mathcal M}_I^{{\bf t}}(K_0)
 \to {\rm Im}(\hat\pi^I_J)\subset{\mathcal M}_J^{{\bf t}}(K_0)$$
 is a $|\Gamma_I|/|\Gamma_J|$-fold stratawise smooth covering. So
 \begin{eqnarray*}
&\qquad\quad\int_{(\hat\pi^I_J)^{-1}({\mathcal M}_J^{{\bf
 t}}(K_0))}((\Pi_{g,m}\times {\rm
EV}_{g,m})\circ\hat\pi_I\bigl)^\ast(\kappa^\ast\times\prod^m_{i=1}\alpha_i^\ast)\\
&=\frac{|\Gamma_I|}{|\Gamma_J|}\int_{{\rm
Im}(\hat\pi^I_J)}((\Pi_{g,m}\times {\rm
EV}_{g,m})\circ\hat\pi_J\bigl)^\ast(\kappa^\ast\times\prod^m_{i=1}\alpha_i^\ast).
\end{eqnarray*}
 This shows that the part
$$
\frac{1}{|\Gamma_I|}\int_{(\hat\pi^I_J)^{-1}({\mathcal M}_J^{{\bf
 t}}(K_0))}((\Pi_{g,m}\times {\rm
EV}_{g,m})\circ\hat\pi_I\bigl)^\ast\Bigl(\kappa^\ast\times\prod^m_{i=1}\alpha_i^\ast\Bigr)$$
 in $\frac{1}{|\Gamma_I|}\int_{{\mathcal M}^{\bf
t}_I(K_0)}((\Pi_{g,m}\times {\rm
EV}_{g,m})\circ\hat\pi_I\bigl)^\ast(\kappa^\ast\times\prod^m_{i=1}\alpha_i^\ast)$
is equal to one
$$\frac{1}{|\Gamma_J|}\int_{{\rm Im}(\hat\pi^I_J)}((\Pi_{g,m}\times {\rm
EV}_{g,m})\circ\hat\pi_I\bigl)^\ast\Bigl(\kappa^\ast\times\prod^m_{i=1}\alpha_i^\ast\Bigr)
$$
in $\frac{1}{|\Gamma_J|}\int_{{\mathcal M}^{\bf
t}_J(K_0)}((\Pi_{g,m}\times {\rm
EV}_{g,m})\circ\hat\pi_I\bigl)^\ast(\kappa^\ast\times\prod^m_{i=1}\alpha_i^\ast)$.
For more details the readers may refer to our expository paper
[LuT]. For simplicity  we later shall omit the double quotation
marks in $``\sum_{I\in{\mathcal N}}"$ without occurring of
confusions. Hence we may use (4.14) to define the same invariants
as (4.7).

\subsection{The homology class in (4.6 ) being independent of choices of
sections $s_i$.}

 We only need to prove the conclusion in the case
increasing a section. Assume that a global smooth section $\hat
s=\{\hat s_I\}_{I\in{\mathcal N}}$  of the bundle system
$(\widehat{\mathcal E}^\ast,\widehat V^\ast)=\{
(\widehat{E}_I^\ast,\widehat V_I^\ast)\,|\,I\in{\mathcal N}\}$ as
$\hat s_{ij}=\{(\hat s_{ij})_I\}_{I\in{\cal N}}$ in (3.76) is such
that the global section $\Psi'=\{\Psi'_I\,|\, I\in{\mathcal N}\}$
of the bundle system
$$
\bigl({\bf P}_1^\ast\widehat{\mathcal E}^\ast, \widehat
V^\ast\times {\bf B}_{\eta'}(\mbox{\Bb R}^{q+1})\bigr)=
\bigl\{\bigl({\bf P}_1^\ast\widehat E_I^\ast, \widehat
V_I^\ast\times {\bf B}_{\eta'}(\mbox{\Bb R}^{q+1})\bigr)\,\bigm|\,
I\in {\mathcal N}\bigr\}
$$
is transversal to the zero section. Here ${\bf P}_1$ is as in
(3.75) and

\begin{eqnarray*}
&&\Psi'_I: \widehat V_I^\ast\times {\bf B}_\tau({\mbox{\Bb
R}}^{q+1})\to {\bf P}_1^\ast\widehat E_I^\ast,\\
&&
 (\hat x_I; {\bf t},t)\mapsto (\bar\partial_J)_I(\hat x_I) +
\sum^{n_3}_{i=1}\sum^{q_i}_{j=1}t_{ij}(\hat s_{ij})_I(\hat
 x_I)+ t\hat s_I(\hat x_I).
\end{eqnarray*}
 As in Corollary 3.15 and (3.75) we get a family of
cobordant singular chains in
 ${\mathcal W}^\ast$,
$$
{\mathcal C}^{({\bf t},t)}(K_0)^\prime:=\sum_{I\in{\mathcal
 N}}\frac{1}{|\Gamma_I|}\{\hat\pi_I:
{\mathcal M}^{({\bf t},t)}_I(K_0)^\prime\to{\mathcal
 W}^\ast\}\quad\forall ({\bf t},t)\in{\bf B}_{\varepsilon'}^{res}(\mbox{\Bb R}^{q+1}),
 $$
where $0<\varepsilon'<\eta'$ and ${\bf
B}_{\varepsilon'}^{res}(\mbox{\Bb R}^{q+1})$ is a residual subset
of ${\bf B}_{\varepsilon_1}(\mbox{\Bb R}^{q+1})$. For
$0<\tau\le\min\{\varepsilon,\varepsilon'\}$
 let us consider the bundle system
$$
\bigl({\bf P}_1^\ast\widehat {\mathcal E}^\ast, \widehat
V^\ast\times {\bf B}_{\tau}(\mbox{\Bb R}^{q+1})\times[0,1]\bigr)=
\bigl\{\bigr({\bf P}_1^\ast\widehat E_I^\ast, \widehat
V_I^\ast\times {\bf B}_{\tau}(\mbox{\Bb R}^{q+1})\times[0,1])\,|\,
I\in{\mathcal N}\bigr\}
$$
 and its smooth section $\Phi=\{\Phi_I\,|\, I\in{\mathcal
N}\}$,
$$\Phi_I:\, (\hat x_I; {\bf t}, t; r)\mapsto
(\bar\partial_J)_I(\hat x_I) +
\sum^{n_3}_{i=1}\sum^{q_i}_{j=1}t_{ij}(\hat s_{ij})_I(\hat
 x_I)+ rt\hat s_I(\hat x_I)
$$
 for any $(\hat x_I; {\bf t}, t; r)\in\widehat V_I^\ast\times{\bf
B}_{\varepsilon'}(\mbox{\Bb R}^{q+1})\times[0,1]$, where
 ${\bf P}_1$ is still the projection to
the first factor of $\widehat V_I^\ast\times {\bf
B}_{\tau}(\mbox{\Bb R}^{q+1})\times[0,1]$. If $\tau>0$ is small
enough, $\Phi$ is transversal to the zero section. It follows that
for each $({\bf t},t)$ in some residual subset ${\bf
B}_{\tau}^{res}(\mbox{\Bb R}^{q+1})$  of ${\bf B}_{\tau}(\mbox{\Bb
R}^{q+1})$ and each $I\in{\cal N}$ the section
$$
\Phi_I^{({\bf t},t)}:\,
(\hat x_I; r)\mapsto \Phi_I(\hat x, {\bf t},t;
r)
$$
of the bundle  $({\bf P}_1^\ast\widehat E^\ast_I, \widehat
V^\ast_I\times[0,1])$ is transversal to the zero section. Any
$({\bf t},t)\in {\bf B}_{\tau}^{res}(\mbox{\Bb R}^{q+1})\cap {\bf
B}_{\varepsilon'}^{res}(\mbox{\Bb R}^{q+1})$ yields a virtual
module chain in ${\mathcal W}\times [0,1]$
$$
\left.\begin{array}{ll}
{\mathcal C}^{({\bf t},t)}:=\sum_{I\in{\mathcal
 N}}\frac{1}{|\Gamma_I|}\{\hat\pi_I^{\prime\prime}:
(\Phi_I^{({\bf t},t)})^{-1}(0)\to{\mathcal
 W}\}\end{array}\right.
  $$
with $\partial{\mathcal C}^{({\bf t},t)}=({\mathcal C}^{\bf
t}(K_0)\times\{0\})\cup(-{\mathcal C}^{({\bf
t},t)}(K_0)^\prime\times\{1\})$. ($\partial{\mathcal C}^{({\bf
t},t)}$ is the part of ${\mathcal C}^{({\bf t},t)}$ of codimension
one.) Here ${\bf t}$ can be required to belong to ${\bf
B}_{\varepsilon}^{res}(\mbox{\Bb R}^{q})$ because of the Claim
A.1.11 in [LeO]. \vspace{2mm}

\noindent{\bf Lemma 4.4}([LeO]).\hspace{2mm}{\it Let $X$ and $Y$
be metric spaces which satisfy the second axiom of countability.
Suppose that $S$ is a countable intersection of open dense subsets
in the product space $X\times Y$. Consider the space $X_S$
consisting of those $x\in X$ such that $S\cap\{x\}\times Y$ is a
countable intersection of open dense subsets in $\{x\}\times Y$.
Then $X_S$ is a countable intersection of open dense subsets in
$X$.}\vspace{2mm}

So ${\mathcal C}^{\bf t}(K_0)$ and ${\mathcal C}^{({\bf t},
t)}(K_0)^\prime$ are virtual module chain cobordant. It implies
that  the corresponding ${\rm EV}_{g,m}\circ{\mathcal C}^{\bf
t}(K_0)$ and ${\rm EV}_{g,m}\circ{\mathcal C}^{{\bf
t},t)}(K_0)^\prime$, and thus the rational cycles ${\mathcal
C}_g^{\bf t}(K_0;\{\bar\alpha_i\}^m_{i=1};\omega,\mu,J,A)$ and
${\mathcal C}_g^{({\bf
t},t)}(K_0;\{\bar\alpha_i\}^m_{i=1};\omega,\mu,J,A)^\prime$ in
$\overline{\mathcal M}_{g,m}$ ( as given by (4.5)) are cobordant.
The desired conclusion follows.

\subsection{The homology class in (4.6) being independent of choices of
 $W_i$.}

Firstly, we prove the case that {\it a $W_{n_3+1}$ and the
corresponding sections $\tilde
s_{(n_3+1)j}=\beta_{(n_3+1)j}\cdot\tilde\nu_{(n_3+1)j}$,
$j=1,\cdots, q_{n_3+1}$, are added.}  In the case let us denote by
$${\mathcal W}_0=\cup_{i=1}^{n_3}W_i,\;{\mathcal E}_0={\mathcal
E}^M_{A,g,m}|_{{\mathcal W}_0},\;{\mathcal
W}_1=\cup^{n_3+1}_{i=1}W_i,\;{\mathcal E}_1={\mathcal
E}^M_{A,g,m}|_{{\mathcal W}_1}.$$
 Then ${\mathcal W}_0$ is an open stratified suborbifold of
 ${\mathcal W}_1$ and ${\mathcal E}_0={\mathcal E}_1|_{{\mathcal
 W}_0}$.  Consider the
stratified orbifold ${\mathcal W}_1\times [0,1]$ with boundary
$\partial ({\mathcal W}_1\times [0, 1])=({\mathcal
W}_1\times\{0\})\cup({\mathcal W}_1\times\{1\})$, and the
pull-back stratified orbifold bundle ${\rm pr}^\ast{\mathcal
E}_1\to{\mathcal W}_1\times [0,1]$ via the projection ${\rm
pr}:{\mathcal W}_1\times[0,1]\to {\mathcal W}_1$. Let us construct
a family of uniformizers of an open covering of it. For
$i=1,\cdots, n_3$ and $j=1,\cdots, n_3+1$ set
$$\begin{array}{lcccccr}
 \widetilde W_i^\prime=\widetilde W_i\times[0, \frac{3}{4}),
 \;\!\!\! &W_i^\prime=W_i\times[0, \frac{3}{4}),\;
 &\pi_i^\prime=\pi_i\times{\bf 1},\\
\widetilde W_{n_3+j}^\prime=\widetilde W_j\times(\frac{1}{4},
1],\; &W_{n_3+j}^\prime=W_j\times(\frac{1}{4}, 1],\;
 &\pi_{n_3+j}^\prime=\pi_j\times{\bf 1},\\
\Gamma_i^\prime=\{g\times{\bf 1}\,|\,g\in\Gamma_i\},
&\Gamma_{n_3+j}^\prime=\{g\times{\bf 1}\,|\,
g\in\Gamma_j\},\\
\widetilde E_i^\prime=p_i^\ast\widetilde E_i,\;\!\!\!
&E_i^\prime=({\rm pr}^\ast{\mathcal E})|_{W_i^\prime},\;
 &p_i: \widetilde W_i\times[0, \frac{3}{4})\to\widetilde W_i,\\
\widetilde E_{n_3+j}^\prime=p_{n_3+j}^\ast\widetilde E_j,\;\!\!\!
&E_{n_3+j}^\prime=({\rm pr}^\ast{\mathcal
E})|_{W_{n_3+j}^\prime},\;
 &\;p_{n_3+j}: \widetilde W_j\times(\frac{1}{4}, 1]\to\widetilde
 W_j.
\end{array}$$
Here all $p_i, p_{n_3+j}$ are the projections to the first
factors. Setting
$$\underline{\mathcal W}=\Bigl(\bigcup^{n_3}_{i=1}W_i\times [0,
\frac{3}{4})\Bigr)\cup\Bigl(\bigcup^{n_3+1}_{i=1}W_i\times
(\frac{1}{4}, 1]\Bigr)\quad{\rm and}\quad \underline{\mathcal
E}=({\rm pr}^\ast{\mathcal E}_1)|_{\underline{\mathcal W}}.
$$
Then $\underline{\mathcal W}$ is an open stratified suborbifold of
${\mathcal W}_1\times [0,1]$, and has also the boundary $\partial
\underline{\mathcal W}=({\mathcal W}_0\times\{0\})\cup({\mathcal
W}_1\times\{1\})$, and $\underline{\mathcal
E}\to\underline{\mathcal W}$ is a stratified orbifold bundle.

Take $\Gamma_k^\prime$-invariant, continuous and stratawise smooth
cut-off functions $\beta_k^\prime$ on $\widetilde W_k^\prime$,
$k=1,\cdots, 2n_3+1$ such that:
\begin{itemize}
\item $\beta_k^\prime(x,0)=\beta_k(x)$, $\forall (x,
0)\in\widetilde W_k^\prime$ and $k=1,\cdots, n_3$,

\item $\beta_k^\prime(x,1)=\beta_{k-n_3}(x)$, $\forall (x,
1)\in\widetilde W_k^\prime$ and $k=n_3+1,\cdots, 2n_3+1$,

\item $\widetilde U_k^{\prime 0}:=\{e\in\widetilde W_k^\prime\;|\;
\beta_k^\prime(e)>0\}\subset\subset\widetilde W_k^\prime$,
$k=1,\cdots, 2n_3+1$,

\item $\cup_{k=1}^{2n_3+1} V_k^{\prime 0}\supset\overline{\mathcal
M}_{g,m}(M,J,A, K_3)\times[0,1]$, where $V_k^{\prime
0}=\pi_I'(\widetilde V_k^{\prime 0})$ and $\widetilde V_k^{\prime
0}:=\{e\in\widetilde W_k^\prime\;|\; \beta_k^\prime(e)=1\}$,
$k=1,\cdots, 2n_3+1$.
\end{itemize}
Then we have:
\begin{eqnarray*}
 && U_k^{\prime 0}\cap
 (W_k\times\{0\})=U_k^0\times\{0\},\quad k=1,\cdots,n_3,\\
&& U_{n_3+ k}^{\prime 0}\cap
(W_k\times\{1\})=U_k^0\times\{1\},\;k=1,\cdots,n_3+1,\\
&& V_k^{\prime 0}\cap
  (W_k\times\{0\})=V_k^0\times\{0\},\quad k=1,\cdots,n_3,\\
&& V_{n_3+ k}^{\prime 0}\cap
(W_k\times\{1\})=V_k^0\times\{1\},\;k=1,\cdots,n_3+1.
\end{eqnarray*}

Let $\tilde s_{ij}=\beta_i\cdot\tilde\nu_{ij}$ be the sections of
the bundle $\widetilde E_i\to\widetilde W_i$ satisfying (3.21),
$i=1,\cdots, n_3+1$  and $j=1,\cdots,q_i$. For $k=1,\cdots, n_3$
we choose the sections of the bundle $\widetilde
E_k^\prime\to\widetilde W_k^\prime$,
$$
\tilde
s^\prime_{kj}(x,t)=\beta_k^\prime(x,t)\cdot\tilde\nu_{kj}(x)\quad\forall\;(x,t)\in\widetilde
W_k^\prime,\; j=1,\cdots, q_k,
$$
 and for $k=1,\cdots,n_3+1$,  ones of the bundle
$\widetilde E_{n_3+k}^\prime\to\widetilde W_{n_3+k}^\prime$,
$$
\tilde
s^\prime_{(n_3+k)j}(x,t)=\beta_{n_3+k}^\prime(x,t)\cdot\tilde\nu_{kj}(x)\quad\forall
(x,t)\in\widetilde W_{n_3+k}^\prime,\; j=1,\cdots, q_k.
$$
 These $q+q'$
 sections satisfy the corresponding
 properties with  (3.12). In this subsection
  $q=\sum^{n_3}_{i=1}q_i$ and
 $q'=\sum^{n_3+1}_{i=1}q_i$.
 Let ${\mathcal N}^\prime$ be the set of all
 subsets $I\subset\{1,\cdots,2n_3+1\}$ with $W_I^\prime=\cap_{i\in
 I}W_i^\prime\ne\emptyset$.  Note that $\{W_i'\,|\, 1\le i\le 2n_3+1\}$ still satisfies
 the corresponding condition with the second one in (3.20).
For the original pairs of $\Gamma_i$-invariant open sets
$W^j_i\subset\subset U^j_i$ with
$$U^1_i\subset\subset W_i^2\subset\subset\cdots
\subset\subset W^{n_3}_i\subset\subset U^{n_3}_i\subset\subset
W_i^{n_3+1}= W_i, $$
 $j=1, 2, \cdots,n_3$ and $i=1,\cdots,n_3+1$ we take
pairs of $\Gamma_k^\prime$-invariant open sets $W^{\prime
l}_k\subset\subset U^{\prime l}_k$, $k=1,\cdots,2n_3+1$ and
$l=1,\cdots,2n_3$ such that the following hold:
$$ U^{\prime 1}_k\subset\subset W^{\prime 2}_k\subset\subset\cdots
\subset\subset W^{\prime 2n_3}_k\subset\subset U^{\prime
2n_3}_k\subset\subset W_k^{\prime 2n_3+1}=W_k;$$
$$\left.\begin{array}{ll}
W_k^{\prime l}\cap(W_k\times\{0\})=W_k^l\quad{\rm and}\quad
U_k^{\prime l}\cap(W_k\times\{0\})=U_k^l,\\
\quad k=1,\cdots,n_3\quad{\rm and}\quad
l=1,\cdots,n_3-1;\end{array}\right.$$
$$\left.\begin{array}{ll}
W_{n_3+k}^{\prime n_3+l}\cap(W_k\times\{1\})=W_k^l\quad{\rm
and}\quad U_{n_2+k}^{\prime n_3+l}\cap(W_k\times\{1\})=U_k^l,\\
\quad k=1,\cdots,n_3+1\quad{\rm and}\quad
l=1,\cdots,n_3.\end{array}\right.
$$
Then as in (3.62), for each $I\in{\mathcal N}^\prime$ with $|I|=k$
we define
$$
 V_I^\prime=(\cap_{i\in I}W_i^{\prime k})\setminus
Cl(\bigcup_{J:|J|>k}(\cap_{j\in J}Cl(U^{\prime k}_j))).
$$
 Let ${\mathcal
N}^{\prime\prime}=\{I\subset\{1,\cdots,n_3+1\}\,|\, W_I=\cap_{i\in
I}W_i\ne\emptyset\}$. It is easy to see that
 $$
{\mathcal N}^{\prime\prime}=\{ I-n_3\,|\, I\in{\mathcal
 N}^\prime,\;\max(I)>n_3\}.
$$
 For each
$I\in{\mathcal N}^{''}$ with $|I|=k$ we also define
$$
V_I^{\prime\prime}=(\cap_{i\in
I}W_i^k)\setminus \bigcup_{J:|J|>k}(\cap_{j\in J}Cl(U^k_j)).
$$
Let $\max I$ (resp. $\min I$) be the largest (resp. smallest)
number, $I-n_3=\{i-n_3\,|\,i\in I\}$ if $\min I>n_3$. The keys are
$$
V_I'\cap\bigl(\cup_{k=1}^{n_3}W_k\times\{0\}\bigr)=
\left\{\begin{array}{ll} V_I\times\{0\},\quad\max I\le n_3,\\
\emptyset,\hspace{15mm}\max I>n_3,\end{array}\right.\eqno(4.15)
$$

$$
V_I'\cap\bigl(\cup_{k=1}^{n_3+1}W_k\times\{1\}\bigr)=
\left\{\begin{array}{ll} V''_{I-n_3}\times\{1\},\quad\min I>n_3,\\
\emptyset,\hspace{20mm}\min I\le n_3.\end{array}\right.
\eqno(4.16)
$$
As in (3.73) we also have the corresponding $V_I'^\ast$, which
satisfies the similar properties to (4.15) and (4.16).
 Repeating the previous arguments we get a
(compatible) system of bundles
$$
\bigl({\bf P}_1^\ast\widehat {\mathcal E}^{\prime\ast}, \widehat
V^{\prime\ast}\times {\bf B}_{\epsilon}(\mbox{\Bb R}^q\times
\mbox{\Bb R}^{q'})\bigr)= \bigl\{\bigl({\bf P}_1^\ast\widehat
E_I^{\prime\ast}, \widehat V_I^{\prime\ast}\times {\bf
B}_{\eta}(\mbox{\Bb R}^q\times \mbox{\Bb R}^{q'})\bigr)\,|\,
I\in{\mathcal N}^\prime\bigl\}
$$
 and a global PS
  section ${\cal S}=\{{\cal S}_I\,|\, I\in{\mathcal
N}^\prime\}$ of it given by
$${\cal S}_I:\, (\hat x_I; {\bf t},{\bf t}')\mapsto
(\bar\partial_J)_I(\hat x_I) +
\sum^{n_3}_{i=1}\sum^{q_i}_{j=1}t_{ij}(\hat s'_{ij})_I(\hat
x_I)+\sum^{n_3+1}_{i=1}\sum^{q_i}_{j=1}t'_{ij}(\hat
s'_{(i+n_3)j})_I(\hat x_I)
$$
 for any $(\hat x_I; {\bf t},{\bf t}')\in\widehat V_I^{\prime\ast}\times{\bf
B}_{\epsilon}(\mbox{\Bb R}^q\times\mbox{\Bb R}^{q'})$, where ${\bf
t}=\{t_{ij}\}\in\mbox{\Bb R}^q$ and ${\bf
t}'=\{t'_{ij}\}\in\mbox{\Bb R}^{q'}$.

For $\epsilon>0$ small enough, as in Theorem 3.14 we can use it to
construct a family of cobordant virtual moduli chains with
boundary and of dimension $2m+ 2c_1(A)+ 2(3-n)(g-1)+1$ in
$\underline{\mathcal W}$,
 $${\mathcal C}^{({\bf t},{\bf t}')}(K_0):=\sum_{I\in{\mathcal
 N}^\prime}\frac{1}{|\Gamma_I^\prime|}\{\hat\pi_I^\prime:
 ({\cal S}_I^{({\bf t},{\bf t}')})^{-1}(0)\to\underline{\mathcal
 W}\}
$$
for $({\bf t},{\bf t}')\in{\bf B}_{\epsilon}^{res}(\mbox{\Bb
R}^q\times\mbox{\Bb R}^{q'})$. Here ${\cal S}_I^{({\bf t},{\bf
t}')}(\hat x_I)={\cal S}_I(\hat x_I; {\bf t},{\bf t}')$. Note that
the stratified smooth manifold $({\cal S}_I^{({\bf t},{\bf
t}')})^{-1}(0)$ has no boundary if and only if $I\in{\mathcal
 N}^\prime$ contains at least two numbers, one of them is no more than
 $n_3$ and other is no less than $n_3$. So its two boundary components,
 denoted by ${\mathcal C}^{({\bf t},{\bf t}')}(K_0)_0$ and
  ${\mathcal C}^{({\bf t},{\bf t}')}(K_0)_1$,
 are respectively virtual moduli chains  in ${\mathcal
 W}_0\times\{0\}$ and ${\mathcal W}_1\times\{1\}$  given by
$$
{\mathcal C}^{({\bf t},{\bf t}')}(K_0)_0=\!\!\sum_{I\in{\mathcal
 N}^\prime,\,\max(I)\le n_3}\!\!\frac{1}{|\Gamma_I^\prime|}
 \{\hat\pi_I^\prime:({\cal S}_I^{({\bf t},{\bf t}')})^{-1}(0)\cap(\hat\pi_I^\prime)^{-1}
 ({\mathcal W}_0\times\{0\})\to{\mathcal
 W}_0\times\{0\}\},$$
$${\mathcal C}^{({\bf t},{\bf t}')}(K_0)_1=\!\!\sum_{I\in{\mathcal
 N}^\prime,\,\min(I)>n_3}\!\!\frac{1}{|\Gamma_I^\prime|}
 \{\hat\pi_I^\prime:({\cal S}_I^{({\bf t},{\bf t}')})^{-1}(0)\cap(\hat\pi_I^\prime)^{-1}
 ({\mathcal W}_1\times\{1\})\to{\mathcal
 W}_1\times\{1\}\}.$$
Consider the evaluation  $\overline{{\rm
EV}}_{g,m}:\underline{\mathcal W}\to M^m,\;(x,t)\mapsto{\rm
EV}_{g,m}(x)$, and
  denote by $\overline{{\rm EV}}^0_{g,m}$ and
$\overline{{\rm EV}}^1_{g,m}$ the restrictions of it to ${\mathcal
W}_0\times\{0\}$ and ${\mathcal W}_1\times\{1\}$ respectively.
Using $\overline{{\rm EV}}^0_{g,m}\circ{\mathcal C}^{({\bf t},{\bf
t}')}(K_0)_0$ and $\overline{{\rm EV}}^1_{g,m}\circ {\mathcal
C}^{({\bf t},{\bf t}')}(K_0)_1$, and as in Proposition 4.1  we can
derive two rational cycles
 $$
 {\mathcal C}_g^{({\bf
t},{\bf
t}')}(K_0;\{\bar\alpha_i\}^m_{i=1};\omega,\mu,J,A)_0\quad{\rm
and}\quad {\mathcal C}_g^{({\bf t},{\bf
t}')}(K_0;\{\bar\alpha_i\}^m_{i=1};\omega,\mu,J,A)_1
$$
in $\overline{\mathcal M}_{g,m}$ of dimension $2c_1(M)(A)+
2(3-n)(g-1)+ 2m-\sum^m_{i=1}\deg\alpha_i$, respectively. Note that
$\overline{{\rm EV}}_{g,m}\circ{\mathcal C}^{({\bf t},{\bf
t}')}(K_0)$ has exactly two boundary components, $\overline{{\rm
EV}}^0_{g,m}\circ{\mathcal C}^{({\bf t},{\bf t}')}(K_0)_0$ and
$\overline{{\rm EV}}^1_{g,m}\circ{\mathcal C}^{({\bf t},{\bf
t}')}(K_0)_1$. As usual we can use it to construct a rational
chain ${\mathcal C}_g^{({\bf t},{\bf
t}')}(K_0;\{\bar\alpha_i\}^m_{i=1};\omega,\mu,J,A)$ in
$\overline{\mathcal M}_{g,m}$ of dimension $2c_1(M)(A)+
2(3-n)(g-1)+ 2m+1-\sum^m_{i=1}\deg\alpha_i$ with exact two
boundary components ${\mathcal C}_g^{({\bf t},{\bf
t}')}(K_0;\{\bar\alpha_i\}^m_{i=1};\omega,\mu,J,A)_0$ and
${\mathcal C}_g^{({\bf t},{\bf
t}')}(K_0;\{\bar\alpha_i\}^m_{i=1};\omega,\mu,J,A)_1$. So we get
 $$
 [{\mathcal C}_g^{({\bf t},{\bf
t}')}(K_0;\{\bar\alpha_i\}^m_{i=1};\omega,\mu,J,A)_0]+[{\mathcal
C}_g^{({\bf t},{\bf
t}')}(K_0;\{\bar\alpha_i\}^m_{i=1};\omega,\mu,J,A)_1]=0.\eqno(4.17)
$$
For $({\bf t},{\bf t}')\in{\bf B}_{\epsilon}^{res}(\mbox{\Bb
R}^q\times\mbox{\Bb R}^{q'})$,
 by Lemma 4.4 we can require
${\bf t}\in{\bf B}_{\varepsilon}^{res}(\mbox{\Bb R}^q)$ and ${\bf
t}'\in{\bf B}_{\varepsilon'}^{res}(\mbox{\Bb R}^{q'})$.
 Note that
$\{I\in{\mathcal N}^\prime\,|\,\max(I)\le n_3\}={\mathcal N}$ and
that (4.15) implies that ${\mathcal C}^{({\bf t},{\bf
t}')}(K_0)_0$ is actually independent of ${\bf t}'$. Therefore
${\mathcal C}_g^{({\bf t},{\bf
t}')}(K_0;\{\bar\alpha_i\}^m_{i=1};\omega,\mu,J,A)_0$ is exactly
identified with the virtual module cycle  ${{\mathcal C}_g^{\bf
t}(K_0;\{\bar\alpha_i\}^m_{i=1};\omega,\mu,J,A)}$ in (4.5). So
$$
[{\mathcal C}_g^{({\bf t},{\bf
t}')}(K_0;\{\bar\alpha_i\}^m_{i=1};\omega,\mu,J,A)_0]=[{\mathcal
C}_g^{{\bf
t}}(K_0;\{\bar\alpha_i\}^m_{i=1};\omega,\mu,J,A)].\eqno(4.18)
$$
Using the bundles $\widetilde E_i\to\widetilde W_i$ and their
sections $\tilde s_{ij}=\beta_i\cdot\tilde\nu_{ij},\; i=1,\cdots,
n_3+1,\;j=1,\cdots,q_i$,
 we can construct as before a bundle system
$$
({\bf P}_1^\ast\widehat {\mathcal E}''^\ast, \widehat
V''^\ast\times {\bf B}_{\varepsilon''}(\mbox{\Bb R}^{q'}))=\bigl
\{({\bf P}_1^\ast\widehat E_I''^\ast, \widehat V_I''^\ast\times
{\bf B}_{\varepsilon''}(\mbox{\Bb R}^{q'}))\,|\, I\in{\mathcal
N}^{\prime\prime}\bigr\}
$$
and its global smooth section $\Psi=\{\Psi_I\,|\, I\in{\mathcal
N}^{\prime\prime}\}$ given by
$$\Psi_I:\, (\hat x_I, {\bf t})\mapsto
(\bar\partial_J)_I(\hat x_I) +
\sum^{n_3+1}_{i=1}\sum^{q_i}_{j=1}t_{ij}\cdot(\hat
s''_{ij})_I(\hat x_I)
$$
 for any $(\hat x_I, {\bf t})\in\widehat V_I''^\ast\times{\bf
B}_{\varepsilon''}(\mbox{\Bb R}^{q'})$. Here
$$
V_I''^\ast=V_I''\cap {\cal W}^{\ast},\quad \widehat
V_I''^\ast=(\hat\pi_I'')^{-1}(V_I''^\ast),\quad \widehat
E_I^{\prime\prime\ast}=\widehat
E_I^{\prime\prime\Gamma_I}|_{\widehat V_I''^\ast}.
$$
Denote by
$${\mathcal
C}_g^{{\bf
t}}(K_0;\{\bar\alpha_i\}^m_{i=1};\omega,\mu,J,A)'',\;{\bf t}\in
{\bf B}_{\varepsilon''}^{res}(\mbox{\Bb R}^{q'})
$$
 the corresponding
 cobordant rational cycles in $\overline{\cal M}_{g,m}$
 constructed from them.  As above we can take ${\bf t}={\bf t}'$ and identify
$-{\mathcal C}_g^{({\bf t},{\bf
t}')}(K_0;\{\bar\alpha_i\}^m_{i=1};\omega,\mu,J,A)_1$ with
${\mathcal C}_g^{{\bf
t}'}(K_0;\{\bar\alpha_i\}^m_{i=1};\omega,\mu,J,A)''$. Hence
$$
-[{\mathcal C}_g^{({\bf t},{\bf
t}')}(K_0;\{\bar\alpha_i\}^m_{i=1};\omega,\mu,J,A)_1]=[{\mathcal
C}_g^{{\bf
t}'}(K_0;\{\bar\alpha_i\}^m_{i=1};\omega,\mu,J,A)''].\eqno(4.19)
$$
It follows from (4.17), (4.18) and (4.19) that
$$[{\mathcal C}_g^{\bf
t}(K_0;\{\bar\alpha_i\}^m_{i=1};\omega,\mu,J,A)]=[{\mathcal
C}_g^{{\bf t}'}(K_0;\{\bar\alpha_i\}^m_{i=1};\omega,\mu,J,A)''].$$
 The desired conclusion is proved.

The general case easily follows from  \S4.2 and the above
arguments.\vspace{2mm}

\noindent{\bf Remark 4.5.}\quad If $(M,\omega)$ is a closed
symplectic manifold we may take $K_0=M$ and $W_i$,
$i=1,\cdots,n_0$ to construct a virtual moduli cycle ${\mathcal
C}^{\bf t}$ as in Proposition 3.16.  It gives a homology class
$[{\mathcal C}^{\bf t}]\in H_r({\mathcal B}^M_{A,g,m};\mbox{\Bb
Q})$. Here $r=2m+2c_1(A)+2(3-n)(g-1)$. The arguments in last two
subsections show that it is independent of the choices of
neighborhoods and sections. Suitably modifying the arguments in
this subsection one can easily prove that it is also independent
of the choice of almost complex structures $J\in{\mathcal
J}(M,\omega)$ and in fact only depend on the symplectic
deformation class of the symplectic form $\omega$.
 Using $\Pi_{g,m}$ and ${\rm
EV}_{g,m}$ we get a  homology class
$$
[\overline{\mathcal M}_{g,m}(M, \langle\omega\rangle, A)]=
(\Pi_{g,m}\times{\rm EV}_{g,m})_\ast([{\mathcal C}^{\bf t}])\in
H_r(\overline{\mathcal M}_{g,m}\times M^m;\mbox{\Bb Q}),
$$
 called the {\it virtual fundamental class},  where
$\langle\omega\rangle$ denotes the deformation class of the
symplectic form $\omega$. It can be used to define the
Gromov-Witten invariants for closed symplectic manifolds
$$
{\mathcal G}{\mathcal W}^{M,\langle\omega\rangle}_{A,g,m}:
H_\ast(\overline{\mathcal M}_{g,m};\mbox{\Bb Q})\times
H^\ast(M;\mbox{\Bb Q})^{\otimes m}\to \mbox{\Bb Q}
$$
 by ${\mathcal G}{\mathcal
W}^{M,\langle\omega\rangle}_{A,g,m}(\kappa;\alpha_1,\cdots,\alpha_m)
=\langle[\overline{\mathcal M}_{g,m}(M, \langle\omega\rangle, A)],
PD(\kappa)\oplus\wedge^m_{i=1}\alpha_i\rangle$. This case had been
completed by different authors(cf. [FuO, LiT2, R, Sie]) though we
do not know whether or not our result is same as one of them.
\vspace{2mm}

\subsection{The homology class in (4.6   ) being independent of choices of
 $K_0$.} We shall prove it in two steps.\vspace{2mm}

\noindent{\it Step 1.}\quad Take a positive integer $r>3$. We
furthermore choose
 $$
 {\bf f}^{(i)}\in\overline{\mathcal M}_{g,m}(M, J,
 A;K_r)\setminus\overline{\mathcal M}_{g,m}(M, J,
 A;K_3),\;i=n_3+1,\cdots,n_r
 $$
such that the corresponding $V^0_{{\bf f}^{(i)}}$ and $W_{{\bf
f}^{(i)}}$, and the original ones satisfy
$$
\cup^{n_r}_{i=1}V^0_{{\bf f}^{(i)}}\supset\overline{\mathcal
M}_{g,m}(M, J, A;K_r).
$$
  Note that (3.16)-(3.19) imply
 $$
  W_i\cap W_j=\emptyset,\;\forall i\le
 n_1,\,j>n_3.\eqno(4.20)
 $$
  We still use the abbreviation notations below (3.20). For $i=n_3+1,\cdots, n_r$,
  we furthermore take the sections $\tilde s_{ij}$ of the bundle $\widetilde
  E_i\to\widetilde W_i$, $j=1,\cdots, q_i$ satisfying
  (3.21). We replace ${\cal W}$ and ${\cal E}$ with
  $${\mathcal W}^\prime=\cup^{n_r}_{i=1}W_i\quad{\rm and}\quad
  {\mathcal E}^\prime=\cup^{n_r}_{i=1}{\mathcal E}_i
  $$
 and denote by ${\mathcal N}_r=\{I\subset\{1,\cdots, n_r\}\,|\,
W_I=\cap_{i\in I} W_i\ne\emptyset\}$. On the basis of the original
pairs of $\Gamma_i$-invariant open sets $W^j_i\subset\subset
U^j_i$, $j=1, 2, \cdots,n_3-1$, $i=1,\cdots, n_3$  we add some
pairs of $\Gamma_i$-invariant open sets $W^j_i\subset\subset
U^j_i$ such that for $i=1,\cdots,n_r$,
$$
U^1_i\subset\subset\cdots \subset\subset W^{n_3}_i\subset\subset
U^{n_3}_i\cdots\subset\subset W_i^{n_r-1}\subset\subset
U_i^{n_r-1} \subset\subset W_i^{n_r}=W_i.
$$
Then as before, for each $I\in{\mathcal N}_r$ with $|I|=k$  we
define
$$
V_I^\prime=(\cap_{i\in I}W^k_i) \setminus
(\bigcup_{J:|J|>k}(\cap_{j\in J}Cl(U^k_j)))\eqno(4.21)
$$
and get the corresponding $V_I^{\prime\ast}$. Let
$q'=\sum^{n_r}_{i=1}q_i$. Repeating the previous arguments we get
a (compatible) system of bundles
$$
({\bf P}_1^\ast\widehat {\mathcal E}^{\prime\ast},
\widehat V^{\prime\ast}\times {\bf B}_{\eta_r}(\mbox{\Bb
R}^{q'}))= \{({\bf P}_1^\ast\widehat E_I^{\prime\ast}, \widehat
V_I^{\prime\ast}\times {\bf B}_{\eta_r}(\mbox{\Bb R}^{q'}))\,|\,
I\in{\mathcal N}_r\}
$$
 and a  global
smooth section $\Psi^\prime=\{\Psi_I^\prime\,|\, I\in{\mathcal
N}_r\}$ of it given by
$$\Psi_I^\prime:\, (\hat x_I, {\bf t})\mapsto
(\bar\partial_J)_I(\hat x_I) +
\sum^{n_r}_{i=1}\sum^{q_i}_{j=1}t_{ij}\cdot(\hat
s_{ij}^\prime)_I(\hat x_I)\eqno(4.22)
$$
 for any $(\hat x_I, {\bf t})\in\widehat V_I^{\prime\ast}\times{\bf
B}_{\eta_r}(\mbox{\Bb R}^{q'})$. As in Corollary 3.14 we have a
positive number $\varepsilon_r\in (0,\eta_r)$ and a family of
cobordant virtual moduli chains in  ${\mathcal
W}^{\prime\ast}\supset{\mathcal W}^\ast$,
$$
{\mathcal C}^{\bf t}(K_0)^\prime:=\sum_{I\in{\mathcal
 N}_r}\frac{1}{|\Gamma_I|}\{\hat\pi_I^\prime:
{\mathcal M}^{\bf t}_I(K_0)^\prime\to{\mathcal
 W}^\ast\}\quad\forall{\bf t}\in{\bf B}_{\varepsilon_r}^{res}(\mbox{\Bb R}^{q'}).
 \eqno(4.23)
$$

Now we need to compare these virtual moduli cycles with ones in
(3.83). To this goal  recall that a subset of a topological space
is often called {\it residual} if it is the countable intersection
of open dense sets. Therefore Lemma 4.4 asserts that if $S$ is
residual then $X_S=\{x\in X\,|\, S\cap(\{x\}\times Y)\;{\rm
is}\;{\rm residual}\;{\rm in}\;\{x\}\times Y\}$ is residual in
$X$.
 For $r>3$ and
 ${\bf B}_{\varepsilon_r}(\mbox{\Bb R}^{q'})$ we define
$$\rho_r=\sup\{\rho>0\,|\,{\bf B}_{\rho}(\mbox{\Bb R}^{q})\times
{\bf B}_{\rho}(\mbox{\Bb R}^{q'-q})\subset{\bf
B}_{\varepsilon_r}(\mbox{\Bb R}^{q'})\}.$$
 Then ${\bf B}_{\rho_r}(\mbox{\Bb R}^{q})\times {\bf
B}_{\rho_r}(\mbox{\Bb R}^{q'-q})\subset{\bf
B}_{\varepsilon_r}(\mbox{\Bb R}^{q'})$. Each ${\bf t}\in{\bf
B}_{\varepsilon_r}(\mbox{\Bb R}^{q'})$ can be written as $({\bf
t}^{(1)}, {\bf t}^{(2)})$, where ${\bf t}^{(1)}\in{\bf
B}_{\rho_r}(\mbox{\Bb R}^{q})$ and ${\bf t}^{(2)}\in{\bf
B}_{\rho_r}(\mbox{\Bb R}^{q'-q})$. Denote by ${\bf
B}^{reg}_{\rho_r}(\mbox{\Bb R}^{q})$ be the set of ${\bf
t}^{(1)}\in{\bf B}_{\rho_r}(\mbox{\Bb R}^{q})$ such that ${\bf
B}^{res}_{\varepsilon_r}(\mbox{\Bb R}^{q'})\cap(\{{\bf
t}^{(1)}\}\times {\bf B}_{\rho_r}(\mbox{\Bb R}^{q'-q}))$ is
residual
 in ${\bf t}^{(1)}\times {\bf B}_{\rho_r}(\mbox{\Bb R}^{q'-q})$.
 By Lemma 4.4 it is residual in ${\bf B}_{\rho_r}(\mbox{\Bb R}^{q})$.
Thus the intersection
$$
{\bf B}^{res}_{\varepsilon_r}(\mbox{\Bb R}^{q})\cap{\bf
B}^{reg}_{\rho_3}(\mbox{\Bb R}^{q}) \subset {\bf
B}^{res}_{\varepsilon_r}(\mbox{\Bb R}^{q})
$$
 is also residual in ${\bf B}_{\rho_r}(\mbox{\Bb R}^{q})$.
 In particular it is nonempty. This shows that we can always take
 ${\bf t}=({\bf t}^{(1)}, {\bf t}^{(2)})\in{\bf B}^{res}_{\varepsilon_r}(\mbox{\Bb
R}^{q'})$ such that ${\bf t}^{(1)}\in{\bf
B}^{res}_{\varepsilon_r}(\mbox{\Bb R}^{q})$. Here we can also
assume that $\varepsilon_r$ is less than $\varepsilon$ in (3.83).

Let us compare the virtual moduli cycle ${\mathcal C}^{\bf
t}(K_0)^\prime$ with one ${\mathcal C}^{{\bf t}^{(1)}}(K_0)$
obtained by (3.83). By (4.20), if $I\in{\mathcal N}$ contains some
number $i\in\{n_1+1,\cdots, n_3\}$, then $W_I\subset W_i$. By
(3.16)-(3.19) and choice of ${\bf f}^{(i)}$ it is not hard to see
that each map in $W_i$ has the image set contained outside $K_0$.
Therefore the image of the map ${\rm EV}_{g,m}\circ\hat\pi_I:
{\mathcal M}^{\bf t}_I(K_0)\to M^m$ is contained outside the
subset $K_0^m\subset \overline{\mathcal M}_{g,m}\times M^m$. It
follows that the fibre product
$$
C(K_0;\{\bar\alpha_i\}^m_{i=1};\omega,\mu,J,A)_I^{{\bf t}^{(1)}}=\emptyset
$$
and thus that
\begin{eqnarray*}
\hspace{20mm}&&\quad{\mathcal C}_g^{{\bf
t}^{(1)}}(K_0;\{\bar\alpha_i\}^m_{i=1};\omega,\mu,J,A)\hspace{63mm}(4.24)\\
&&=\sum_{I\in{\mathcal
 N}}\frac{1}{|\Gamma_I|}\{\Pi^I_{g,m}:C(K_0;\{\bar\alpha_i\}^m_{i=1};
 \omega,\mu,J,A)_I^{{\bf t}^{(1)}}\to
 \overline{\mathcal  M}_{g,m}\}\\
&&=\sum_{I\in{\mathcal
 N}_1}\frac{1}{|\Gamma_I|}\{\Pi^I_{g,m}:C(K_0;
 \{\bar\alpha_i\}^m_{i=1};\omega,\mu,J,A)_I^{{\bf t}^{(1)}}\to
 \overline{\mathcal  M}_{g,m}\}.
\end{eqnarray*}
Here ${\mathcal N}_1=\{I\in{\mathcal N}\,|\, \max(I)\le n_1\}$.
The same reason leads to
\begin{eqnarray*}
\hspace{20mm}&&\quad {\mathcal C}_g^{{\bf
t}}(K_0;\{\bar\alpha_i\}^m_{i=1};\omega,\mu,J,A)^\prime\hspace{63mm}(4.25)\\
&&=\sum_{I\in{\mathcal
 N}_r}\frac{1}{|\Gamma_I|}\{\Pi^I_{g,m}:C(K_0;\{\bar\alpha_i\}^m_{i=1};
 \omega,\mu,J,A)_I'^{{\bf t}}\to
 \overline{\mathcal  M}_{g,m}\}\\
&&=\sum_{I\in{\mathcal
 N}_1}\frac{1}{|\Gamma_I|}\{\Pi^I_{g,m}:C(K_0;\{\bar\alpha_i\}^m_{i=1};
 \omega,\mu,J,A)_I'^{{\bf t}}\to
 \overline{\mathcal  M}_{g,m}\}.
 \end{eqnarray*}
 Here $C(K_0;\{\bar\alpha_i\}^m_{i=1};\omega,\mu,J,A)_I'^{{\bf t}}$ is
 the fibre product obtained  in (4.3) by replacing ${\mathcal
M}^{\bf t}_I(K_0)$ and $\hat\pi_I$ by ${\mathcal M}^{\bf
t}_I(K_0)^\prime$ and $\hat\pi_I^\prime$ respectively. For
$I\in{\mathcal N}_1$ with $|I|=k$, (3.62), (3.63) and (4.21) imply
$$V_I=(\cap_{i\in I}W^k_i)\setminus
 (\bigcup_{J\in{\mathcal N}:|J|=k+1}(\cap_{j\in J}Cl(U^k_j))),$$
$$V_I^\prime=(\cap_{i\in I}W^k_i)
\setminus (\bigcup_{J\in{\mathcal N}_r:|J|=k+1}(\cap_{j\in
J}Cl(U^k_j))).
$$
Note that $\cap_{j\in J}Cl(U^k_j)$ and $\cap_{i\in I}W^k_i$ are
disjoint if $J$ contains a number greater than $n_3$. So
\begin{eqnarray*}
 (\cap_{i\in I}W^k_i)\setminus Cl(\bigcup_{J\in{\mathcal
N}_r:|J|=k+1}(\cap_{j\in J}Cl(U^k_j)))\\
 = (\cap_{i\in
I}W^k_i)\setminus Cl(\bigcup_{J\in{\mathcal N}:|J|=k+1}(\cap_{j\in
J}Cl(U^k_j))).
\end{eqnarray*}
This implies that $V_I=V_I^\prime,\;\forall I\in{\mathcal N}_1$.
It follows that
\begin{eqnarray*}
 && \widehat V_I^\ast=\widehat
V_I^{\prime\ast},\;\widehat E_I^\ast=\widehat E_I^{\prime\ast},
\;\widehat \pi_I=\widehat \pi_I^\prime, (\hat s_{ij})_I=(\hat
s_{ij}')_I,\\
&&\forall I\in{\mathcal N}_1,\;  i=1,\cdots,n_3,
j=1,\cdots,q_i.\end{eqnarray*}
 But $(\hat s_{ij}')_I=0$ for any $i>n_3$. It follows from
 (4.22) and the definition of $\Psi_I$ in (3.76) that
$${\mathcal M}^{\bf t}_I(K_0)^\prime={\mathcal M}^{{\bf
t}^{(1)}}_I(K_0)\times\{{\bf t}^{(2)}\},\quad\forall I\in{\mathcal
N}_1.
$$
 From this, (4.24) and (4.25) it follows that
 ${\mathcal C}_g^{{\bf t}^{(1)}}(K_0;\{\bar\alpha_i\}^m_{i=1};\omega,\mu,J,A)$
 can be identified
with ${\mathcal C}_g^{{\bf
t}}(K_0;\{\bar\alpha_i\}^m_{i=1};\omega,\mu,J,A)^\prime$. We get
the desired result
$$[{\mathcal C}_g^{{\bf t}^{(1)}}(K_0;\{\bar\alpha_i\}^m_{i=1};\omega,\mu,J,A)]=
[{\mathcal C}_g^{{\bf
t}}(K_0;\{\bar\alpha_i\}^m_{i=1};\omega,\mu,J,A)^\prime].
$$

\noindent{\it Step 2.}\quad Let $K_0^\prime\subset M$ be another
compact subset containing ${\rm suppt}(\alpha_1)\supseteq{\rm
suppt}(\wedge^m_{i=1}\alpha_i)$. Since $\cup^\infty_{r=1}K_r=M$ we
can take a positive integer $r>6$ such that
$$K_3^\prime:=\{x\in M\,|\, d(x, K_0^\prime)\le 3C\}\subset
K_{r-3}.
$$
 Then
 $\overline{\mathcal M}_{g,m}(M, J,
 A;K_3^\prime)\subset\overline{\mathcal M}_{g,m}(M, J,
 A;K_{r-3})$.
 Assume as above  that ${\bf f}_i^\prime\in\overline{\mathcal M}_{g,m}(M, J,
 A;K_3^\prime)$, $i=1,\cdots,l_3$, and their open neighborhoods
 $W_i^\prime$, $i=1,\cdots$, $l_3$, and the sections
 $\tilde \sigma_{ij},\,i=1,\cdots,k_3$ and $j=1,\cdots, p_i$
can used to construct a family of cobordant virtual moduli chains
$${\bar{\mathcal C}}^{{\bf t}^{(1)}},\;
 \forall {\bf t}^{(1)}\in{\bf B}^{res}_{\delta}(\mbox{\Bb R}^p).
$$
 Here $p=p_1+\cdots + p_{k_3}$.
Furthermore we choose points ${\bf
f}_i^\prime\in\overline{\mathcal M}_{g,m}(M, J,
 A;K_r)\setminus\overline{\mathcal M}_{g,m}(M, J,
 A;K_3^\prime)$, $i=l_3+1,\cdots,l_r$, and their open neighborhoods
 $W_i^\prime$, $i=l_3+1,\cdots,l_r$, and the sections
 $\tilde \sigma_{ij},\,i=l_3+1,\cdots,l_r$ and $j=1,\cdots,p_i$
such that they and the open neighborhoods
 $\widetilde W_i^\prime$, $i=1,\cdots,l_3$ and the sections
 $\tilde \sigma_{ij},\,i=1,\cdots,k_3, j=1,\cdots,p_i$ together can be used to form
 another family of cobordant virtual moduli chains
$$
{\bar{\mathcal C}}^{{\bf t}},\;
 \forall {\bf t}\in{\bf B}^{res}_{\delta_r}(\mbox{\Bb R}^{p'}),
 $$
 where $p'=p_1+\cdots + p_{k_r}$.
By Step 1 we know that using the virtual moduli  chains
${\bar{\mathcal C}}^{{\bf t}^{(1)}}$ and ${\bar{\mathcal C}}^{{\bf
t}}$ yields the same homology class in (4.6). On the other hand,
the above construction of ${\bar{\mathcal C}}^{{\bf t}}$ shows
that the virtual moduli chains ${\mathcal C}^{{\bf
t}}(K_0)^\prime$ of (4.23) and ${\bar{\mathcal C}}^{{\bf
t}^\prime}$  also gives the same homology class in (4.6) for any
${\bf t}\in {\bf B}^{res}_{\varepsilon_r}(\mbox{\Bb R}^{q'})$ and
${\bf t}^\prime\in{\bf B}^{res}_{\delta_r}(\mbox{\Bb R}^{p'})$.
Therefore the statement of this subsection is proved.

 \noindent{\bf Remark 4.6.}\quad In the arguments above
we have chosen the constant $C$ given by (3.18).   That is, $C=6+
\frac{4\beta_0}{\pi\alpha_0^2 r_0}\omega(A)$. Actually the methods
of this subsection may be used to prove that the left side of
(4.11) is independent of choices of a fixed constant $C\ge 6+
\frac{4\beta_0}{\pi\alpha_0^2 r_0}\omega(A)$. Indeed, for a given
constant $\widehat C>6+ \frac{4\beta_0}{\pi\alpha_0^2
r_0}\omega(A)$ and  the compact subsets
$$
\widehat K_j:=\{q\in M\,|\, d_\mu(q, K_0)\le
j\widehat C\},\;j=1,2,\cdots,
$$
 we can choose a positive integer number $r>3$ so large that
$\widehat K_3\subset K_r$ and thus
$$
\overline{\mathcal M}_{g,m}(M, J,
 A;\widehat K_3)\subset\overline{\mathcal M}_{g,m}(M, J,
 A; K_r).
$$
 Here $K_3$ is given by (3.17).
Then we can use the methods of this subsection to construct two
virtual moduli chains ${\mathcal C}_0$ and ${\mathcal C}_1$ from
$\overline{\mathcal M}_{g,m}(M, J,
 A;\widehat K_3)$ and $\overline{\mathcal M}_{g,m}(M, J,
 A; K_r)$ such that using ${\mathcal C}_0$ and ${\mathcal C}_1$ yields the same
 homology class in (4.6).

\subsection{The homology class in (4.6) being independent of choices of
 almost complex structures
$J\in{\mathcal J}(M,\omega,\mu)$.}

We want to prove that
$$
[{\mathcal C}_g^{\bf
t}(K_0;\{\bar\alpha_i\}^m_{i=1};\omega,\mu,J_0,A)]= [{\mathcal
C}_g^{\bf
t}(K_0;\{\bar\alpha_i\}^m_{i=1};\omega,\mu,J_1,A)]\eqno(4.26)
$$
 for any two $J_0, J_1\in{\mathcal J}(M,\omega,\mu)$.
 Without loss of generality we assume that (1.1) holds for them. By Proposition 2.3
we have a smooth path $(J_t)_{t\in [0,1]}$ in ${\mathcal
J}(M,\omega,\mu)$ connecting $J_0$ to $J_1$ such that for
$\gamma_0=\alpha_0^3/\beta_0^2$($<\alpha_0$),
$$
\omega(X, J_tX)\ge\gamma_0\|X\|^2_\mu,\; \forall X\in TM\,{\rm and}\,t\in[0,1].
$$
We still assume that  ${\rm suppt}(\wedge^m_{i=1}\alpha_i)$ is
contained in $K_0$. Take a positive constant  $C$ satisfying
 $$
C>6+ \frac{4\beta_0}{\pi\gamma_0^2
r_0}\omega(A).
$$
 Then for $K_j:=\{q\in M\,|\, d_\mu(q, K_0)\le
jC\},\;j=1,2,\cdots$, we have
$$
{\rm Im}({\bf f})\subset K_4,\quad\forall
{\bf f}\in\overline{\mathcal M}_{g,m}(M, J_t,
 A; K_3)\;{\rm and}\;t\in [0,1].
$$
Using $\cup_{t\in [0,1]}\overline{\mathcal M}_{g,m}(M, J_t,
 A; K_3)$ we can construct a virtual moduli chain cobordsim between
 the virtual moduli chains associated with $\overline{\mathcal M}_{g,m}(M, J_0,
 A; K_3)$ and one with $\overline{\mathcal M}_{g,m}(M, J_1,
 A; K_3)$. It may yield the cobordant rational cycles
 ${\mathcal C}_g^{\bf
 t}(K_0;\{\bar\alpha_i\}^m_{i=1};\omega,\mu,J_0,A)$ and
${\mathcal C}_g^{{\bf
t}'}(K_0;\{\bar\alpha_i\}^m_{i=1};\omega,\mu,J_1,A)$ in
$\overline{\mathcal M}_{g,m}$. Equation (4.26) is proved.

\subsection{The dependence of the homology class in (4.6) on Riemannian metrics
and symplectic forms}

In this subsection we shall prove that the homology class in (4.6)
is invariant under the weak deformation of $(M,\omega, J,\mu)$. As
before  we  assume that ${\rm suppt}(\wedge^m_{i=1}\alpha_i)$ is
contained in $K_0$. Let
 $(M, \omega_t, J_t, \mu_t)_{t\in [0,1]}$ be a weak deformation of
the  geometrically bounded symplectic manifold $(M, \omega, J,
\mu)$. Take a positive constant  $C$ satisfying the following
property:
 $$C>6+ \frac{4\beta_0}{\pi\gamma_0^2
r_0}\omega(A).
$$
 Then for $K_j:=\{q\in M\,|\, d(q, K_0)\le
jC\},\;j=1,2,\cdots$, by the third condition in (1.3) we have
$$
\cup_{t\in [0,1]}\{{\rm Im}({\bf f})\,|\,{\bf f}\in\overline{\mathcal M}_{g,m}(M, J_t,
 A; K_3)\}\subset K_4.
$$
Since $K_4$ is compact, for any $t\in (0, 1]$ we can use
$\overline{\mathcal M}_{g,m}(M, J_s, A; K_3),\, s\in [0,t]$ to
construct a virtual moduli chain cobordsim between  the virtual
moduli chain associated with $\overline{\mathcal M}_{g,m}(M, J_0,
A; K_3)$ and one with $\overline{\mathcal M}_{g,m}(M, J_t,
 A; K_3)$. As in \S4.5 we can easily show that
 $$[{\mathcal C}_g^{\bf
 t}(K_0;\{\alpha_i\}^m_{i=1};\omega_t,\mu_t,
J_t,A)]= [{\mathcal C}_g^{\bf
t}(K_0;\{\alpha_i\}^m_{i=1};\omega,\mu,J,A)]\quad\forall t\in
[0,1].
$$

\noindent{\bf Remark 4.7.}\quad In the above arguments all almost
complex structures are assumed to belong to ${\mathcal
J}(M,\omega,\mu)$. Actually, all arguments from $\S2.1$ to  this
subsection are still valid for all almost complex structures in
${\mathcal J}_\tau(M,\omega,\mu)$.

\subsection{The case of a compact symplectic manifold with
contact type boundaries.}

Similar to the case of closed symplectic manifolds in Remark 4.5
we may follow a way  given in [Lu2] to use the above techniques
 to construct a virtual fundamental class $[\overline{\mathcal M}_{g,m}(M,
\langle\omega\rangle_c, A)]\in H_r(\overline{\mathcal
M}_{g,m}\times M^m;\mbox{\Bb Q})$ which may be used to define
GW-invariants. Hereafter $\langle\omega\rangle_c$ always denote
the deformation class of $\omega$ with contact type boundary. We
below will give it as a consequence of Theorem 1.1 again. Recall
that a compact $2n$-dimensional manifold $(M,\omega)$ with
boundary is said to have contact type boundary $\partial M$ if
there exists a positive contact form $\lambda$ on $\partial M$
such that $d\lambda=\omega|_{\partial M}$. The positivity of
$\lambda$ means that $\lambda\wedge (d\lambda)^{n-1}$ gives the
natural orientation on $\partial M$ induced by $\omega^n$.  Note
that for such a manifold $(M,\omega)$ one may always associate a
noncompact symplectic manifold, denoted $(\widetilde
M,\widetilde\omega)$ such that
$$\widetilde M=M\cup_{\partial M\times\{1\}}(\partial M\times
[1,+\infty)\quad{\rm and}$$\vspace{-4mm}
$$\widetilde\omega=\left\{\hspace{-2mm}\begin{array}{ll}
&\omega\quad\quad\quad{\rm on}\;M\\
&d(z\lambda)\quad{\rm on}\;\partial M\times [1,+\infty),
\end{array}\right.\eqno(4.27)
$$
where $z$ is the projection on the second factor. There is a class
of important almost complex structures on $\partial M\times
[1,+\infty)$. Let $\xi_\lambda={\rm Ker}(\lambda)$ and $X_\lambda$
be the Reeb vector field determined by $i_{X_\lambda}\lambda=1$
and $i_{X_\lambda}(d\lambda)=0$. Using the natural splitting
$T(\partial M)=\mbox{\Bb R}X_\lambda\oplus\xi_\lambda$ and a given
complex structure $\hat J\in{\mathcal J}(\xi_\lambda,
d\lambda|_{\xi_\lambda})$ on the bundle $\xi_\lambda\to\partial M$
we get a $\widetilde\omega$-compatible almost complex structure on
$\partial M\times [1,+\infty)$ as follows:
$$\bar J(z,x)(h,k)=(-\lambda(x)(k), \hat
J(x)(k-\lambda(x)(k)X_\lambda(x))+hX_\lambda(x))\eqno(4.28)
$$
 for $h\in\mbox{\Bb R}\cong T_z\mbox{\Bb R}$ and $k\in
 T_x(\partial M)$. The corresponding compatible  metric
 $g_{\bar J}$ is given by
 \begin{eqnarray*}
 g_{\bar J}(z,x)((h_1,k_1),(h_2,k_2))&=& g_{\hat
J}(x)(k_1-\lambda(x)(k_1)X_\lambda(x),
k_2-\lambda(x)(k_2)X_\lambda(x))\\
&+& h_1\cdot h_2+\lambda(x)(k_1)\cdot\lambda(x)(k_2)
\end{eqnarray*}
 for $h_1, h_2\in\mbox{\Bb R}$ and $k_1, k_2\in T_x(\partial M)$.
Therefore we get
$$\left\{\begin{array}{ll}
 |\widetilde
\omega((h_1,k_1),(h_2,k_2))|\le
\|(h_1,k_1)\|_{g_{\bar J}}\cdot\|(h_2,k_2)\|_{g_{\bar J}},\\
\|(h_1,k_1)\|^2_{g_{\bar
J}}=\|k_1-\lambda(x)(k_1)X_\lambda(x)\|^2_{g_{\hat J}}+ h_1^2+
(\lambda(x)(k_1))^2.
\end{array}\right.\eqno(4.29)
$$
Fix a Riemannian metric $\mu_1$ on $\partial M$. The product of it
with the standard metric $\tau$ on $[1, +\infty)$ may be extended
to a Riemannian metric on $\widetilde M$, denoted by $\mu$. By
Lemma 2.2 of [Lu1] it is easily proved that $\mu$ is a
geometrically bounded Riemannian metric on $\widetilde M$. Since
$M$ is compact there exist constants $0<\gamma_1(\lambda,\hat
J)\le 1\le\gamma_2(\lambda,\hat J)$ such that for any $v\in
T(\partial M)$,
$$
\gamma_1(\lambda,\hat J)\|v\|_{\mu_1}\le\sqrt{\|v-\lambda(v)X_\lambda\|_{g_{\hat J}}+
|\lambda(v)|^2} \le \gamma_2(\lambda,\hat
J)\|v\|_{\mu_1}
$$
 This and (4.29) give
$$\left\{\begin{array}{ll}
\gamma_1(\lambda,\hat J)\|X_1\|_\mu\le\|X_1\|_{g_{\bar
J}}\le\gamma_2(\lambda,\hat J)\|X_1\|_\mu,\\
|\widetilde\omega(X_1, X_2)|\le \gamma_2(\lambda,\hat
J)^2\|X_1\|_\mu\cdot\|X_2\|_\mu
\end{array}\right.\eqno(4.30)
$$
for $X_1=(h_1,k_1), X_2=(h_2, k_2)\in T\widetilde M$.

  By Proposition 1.2.1 in [IS] each
 $J\in{\mathcal J}(M,\omega)$ can be extended to a
 $\widetilde\omega$-compatible almost complex structure
 $\tilde J$ on $\widetilde M$ such that $\tilde J=\bar J$ on
 $\partial M\times [2,+\infty)$ and thus $g_{\tilde
 J}=g_{\bar J}$ on  $\partial M\times [2,+\infty)$.
 Clearly, by decreasing $\gamma_1(\lambda,\hat J)$ and increasing
$\gamma_2(\lambda,\hat J)$ we can assume
 $$\gamma_1(\lambda,\hat
J)\|X_1\|_\mu\le\|X_1\|_{g_{\tilde J}}\le\gamma_2(\lambda,\hat
J)\|X_1\|_\mu,\;\forall X_1\in T\widetilde M.\eqno(4.31)
$$
 This and (4.30) show that $\tilde J\in{\mathcal J}(\widetilde
 M,\widetilde\omega,\mu)$.

Since $M\subset\widetilde M$ is a proper deformation retract there
exists a natural isomorphism  $H_c^\ast(M;\mbox{\Bb Q})\cong
H_c^\ast(\widetilde M;\mbox{\Bb Q})$. So the cohomology classes
$\alpha_1,\cdots,\alpha_m\in H_c^\ast(M;\mbox{\Bb
 Q})$ may naturally be viewed as classes in $H_c^\ast(\widetilde  M;\mbox{\Bb Q})$.
 Let $\kappa\in H_\ast(\overline{\mathcal
M}_{g,m},\mbox{\Bb Q})$; by Theorem 1.1 we have a rational number
$${\mathcal G}{\mathcal
W}^{(\widetilde\omega, \mu, \tilde
J)}_{A,g,m}(\kappa;\alpha_1,\cdots,\alpha_m).\eqno(4.32)
$$
 Seemingly, it depends on $\omega$, $\lambda$, $J$ and $\mu_1$, $\hat J$.
Actually this number is independent of their suitable choices.
Firstly, for two Riemannian metrics $\mu_1^{(0)}$ and
$\mu_1^{(1)}$ on $\partial M$ let $\mu^{(0)}$ and $\mu^{(1)}$ be
the corresponding metrics on $\widetilde M$ obtained by
$\mu_1^{(0)}$ and $\mu_1^{(1)}$ as above. Since the set of all
Riemannian metrics on a manifold is convex we get a family of
Riemannian metrics on $\widetilde M$,
$\mu^{(t)}=(1-t)\mu^{(0)}+t\mu^{(1)}$ with $t\in [0,1]$. Note that
$\mu^{(t)}=((1-t)\mu_1^{(0)}+t\mu_1^{(1)})\times\tau$ outside
$M\subset\widetilde M$. It is not hard to prove that
$(\mu^{(t)})_{t\in [0,1]}$ is a continuous path in ${\mathcal
G}{\mathcal R}(\widetilde M)$ with respect to the
$C^\infty$-strong topology. Therefore Theorem 1.1(iii) implies
that (4.32) does not depend on $\mu_1$ (and thus $\mu$).

Next for a smooth path $(J_t)_{t\in [0,1]}$ in ${\mathcal
J}(M,\omega)$ and a smooth path $(\hat J_t)_{t\in [0,1]}$ in
${\mathcal J}(\xi_\lambda, d\lambda|_{\xi_\lambda})$ we may extend
$(J_t)_{t\in [0,1]}$ into a smooth path $(\tilde J_t)_{t\in
[0,1]}$ in ${\mathcal J}(\widetilde M,\widetilde\omega)$ such that
$\tilde J_t=\bar J_t$ on $\partial M\times [2, +\infty)$ for any
$t\in [0,1]$. Here $\bar J_t$ is obtained by replacing $\hat J$
with $\hat J_t$ in (4.28). One easily proves that (4.31) uniformly
hold for all $\hat J_t$, $\tilde J_t$ and suitable positive
constants $\gamma_1(\lambda,\{\hat J_t\})$ and
$\gamma_2(\lambda,\{\hat J_t\})$. As in the proof of Theorem
1.1(ii) we can show that (4.32) is independent of choices of $\hat
J$ and $J$.

Thirdly, we show that (4.32) is independent of choices of
$\lambda$. Let us  denote by
$${\rm Cont}_+(M,\partial M,\omega)
$$
the set of the positive contact form $\lambda$ on $\partial M$
such that $d\lambda=\omega|_{\partial M}$. It was proved in Lemma
2.1 of [Lu2] that ${\rm Cont}_+(M,\partial M,\omega)$ is a convex
set. Thus any two points $\lambda_0,\lambda_1$ in ${\rm
Cont}_+(M,\partial M,\omega)$ may be connected with a smooth path
$(\lambda_t)_{t\in [0,1]}$ in ${\rm Cont}_+(M,\partial M,\omega)$.
Using (4.27) and (4.28) we  obtain a  path
$(\widetilde\omega_t)_{t\in[0,1]}$ of symplectic forms on
$\widetilde M$ and a path $(\tilde J_t)_{t\in [0,1]}$ in
${\mathcal J}(\widetilde M,\widetilde\omega)$ that are continuous
with respect to the $C^\infty$-strong topology.  Moreover, each
$(\widetilde M,\widetilde\omega, \tilde J,\mu)$ is geometrically
bounded. By Theorem 1.1(iii) we have
$$
{\mathcal G}{\mathcal
W}^{(\widetilde\omega_1, \mu, \tilde
J_1)}_{A,g,m}(\kappa;\alpha_1,\cdots,\alpha_m)={\mathcal
G}{\mathcal W}^{(\widetilde\omega_0, \mu, \tilde
J_0)}_{A,g,m}(\kappa;\alpha_1,\cdots,\alpha_m).
$$

Finally, we study the dependence of (4.32) on $\omega$. Note that
there exists an arbitrary small deformation of a symplectic form
$\omega$ for which $\partial M$ is not of contact type.  Thus we
must restrict to a class of deformations of the form $\omega$. A
deformation of the symplectic form $\omega$ on $M$,
$(\omega_t)_{t\in[0,1]}$ starting at $\omega_0=\omega$, is called
a {\it deformation with contact type boundary} if each
$(M,\omega_t)$ is a symplectic manifold with contact type
boundary. Since each $\omega_t$ is exact near $\partial M$ we may
construct a smooth family of $1$-forms $(\lambda_t)_{t\in [0,1]}$
such that each $\lambda_t$ belongs to ${\rm Cont}_+(M,\partial
M,\omega_t)$. They give a deformation
$(\widetilde\omega_t)_{t\in[0,1]}$ of symplectic forms on
$\widetilde M$ by (4.27). Using these  and a given metric $\mu_1$
on $\partial M$ we may construct a continuous path $(\hat
J_t)_{t\in [0,1]}$ and thus a continuous path $(\tilde J_t)_{t\in
[0,1]}$ of almost complex structures on $\widetilde M$ (with
respect to the $C^\infty$ strong topology) with $\tilde
J_t\in{\mathcal J}(\widetilde M,\widetilde\omega_t,\mu)$. By
Theorem 1.1(iii) we get that
$${\mathcal G}{\mathcal
W}^{(\widetilde\omega_t, \mu, \tilde
J_t)}_{A,g,m}(\kappa;\alpha_1,\cdots,\alpha_m)={\mathcal
G}{\mathcal W}^{(\widetilde\omega_0, \mu, \tilde
J_0)}_{A,g,m}(\kappa;\alpha_1,\cdots,\alpha_m).
$$
 for any $t\in [0,1]$.
Summing up the previous arguments yields\vspace{2mm}

\noindent{\bf Corollary 4.7.}\quad{\it Given a compact
$2n$-dimensional symplectic manifold $(M,\omega)$ with contact
type boundary,  and integers $g\ge 0$, $m>0$ with $2g+m\ge 3$, and
a class $A\in H_2(M,\mbox{\Bb Z})$ the map
$${\mathcal G}{\mathcal W}^{M,\langle\omega\rangle_c}_{A,g,m}:
H_\ast(\overline{\mathcal M}_{g,m};\mbox{\Bb Q})\times
H_c^\ast(M;\mbox{\Bb Q})^{\otimes m}\to \mbox{\Bb
Q}
$$
 given by
 ${\mathcal G}{\mathcal W}^{M,\langle\omega\rangle_c}_{A,g,m}(\kappa;
 \alpha_1,\cdots,\alpha_m):=
{\mathcal G}{\mathcal W}^{(\widetilde\omega, \mu, \tilde
J)}_{A,g,m}(\kappa;\alpha_1,\cdots,\alpha_m)$
 is well-defined and independent of all related choices. }\vspace{2mm}

\section{The properties of Gromov-Witten invariants}

This section will discuss the deep properties of the invariants
constructed in $\S4$.

\subsection{The localization formulas}

The following two localization formulas may be used to simplify
the computation for the Gromov-Witten invariants. They show that
it is sometime enough to use the smaller moduli space to construct
a (local) virtual moduli cycle for computation of some concrete
GW-invariants. Some special forms of them  perhaps appeared in
past literatures. We here explicitly propose and prove them in
generality.

Let $\Lambda\subset\overline{\mathcal M}_{g,m}(M, A, J)$ be a
nonempty compact subset. For given classes $\alpha_1\in
H_c^\ast(M,\mbox{\Bb Q})$ and  $\alpha_2,\cdots,\alpha_m\in
H^\ast(M,\mbox{\Bb Q})$, as before we fix a compact subset
$K_0\subset M$ containing ${\rm supp}(\wedge^m_{i=1}\alpha_i)$
such that $\Lambda\subset\overline{\mathcal M}_{g,m}(M, A,
 J;K_3)$. Now replacing $\overline{\mathcal M}_{g,m}(M, A,
J;K_3)$ with $\Lambda$ we
 take finitely many  points
  $[{\bf f}_i]\in\Lambda$,
$i=1,\cdots, k$, the corresponding uniformizers $\pi_i: \widetilde
W_i=\widetilde W_{{\bf f}_i}\to W_i$  and $\Gamma_i={\rm Aut}({\bf
f}_i)$-invariant partially smooth cut-off function $\beta_i$ on
$\widetilde W_i$ to satisfy
$$\cap^k_{i=1}W_i=\emptyset\quad{\rm
and}\quad\cup^k_{i=1}V_i^0\supset\Lambda\eqno(5.1)
$$
 for $\widetilde V_i^0:={\rm Int}(\{x\in\widetilde W_i\,|\,\beta_i(x)=1\})$ and
 $V_i^0=\pi_i(\widetilde V_i^0)$, $i=1,\cdots,k$.
As in $\S3$ let $\tilde s_{i1},\cdots,\tilde s_{iq_i}$ be the
chosen desired  sections of the bundles $\widetilde
E_i\to\widetilde W_i$,  $i=1,\cdots,k$.  Using these data
 and repeating the arguments in $\S3$ we can construct a system of bundles
 $(\widehat{\mathcal E}^\ast,\widehat V^\ast)=\{(\widehat E_I^\ast, \widehat
V_I^\ast),\;\hat\pi_I\,|\, I\in{\mathcal N}_k\}$ and a family of
cobordant virtual moduli chains in ${\mathcal
W}^\ast(\Lambda)=\cup_{I\in{\cal
N}_k}V_I^\ast\subset\cup^k_{i=1}V_i^0$ of dimension $2m+ 2c_1(A)+
2(3-n)(g-1)$
$${\mathcal C}^{\bf t}(K_0, \Lambda):=\sum_{I\in{\mathcal
 N}_k}\frac{1}{|\Gamma_I|}\{\hat\pi_I:{\mathcal M}^{\bf t}_I(K_0, \Lambda)\to{\mathcal
 W}^\ast(\Lambda)\},
$$
for any ${\bf t}\in{\bf B}_{\varepsilon}^{res}(\mbox{\Bb R}^{q})$
associated with $\Lambda$. Here $q=q_1+\cdots +q_k$ and ${\mathcal
N}_k$ the collection of the subset $I$ of $\{1,\cdots,k\}$ with
$W_I=\cap_{i\in I}W_i\ne\emptyset$.
 Note that ${\mathcal W}(\Lambda)$ is only a neighborhood of
 $\Lambda$ in ${\mathcal B}^M_{A,g,m}$, which may be required to be
arbitrarily close to $\Lambda$.

Now replacing ${\rm EV}_{g,m}\circ{\mathcal C}^{\bf t}(K_0)$ in
(4.2) by ${\rm EV}_{g,m}\circ{\mathcal C}^{\bf t}(K_0, \Lambda)$
we, as in $\S4.1$, have the corresponding fibre product
$C(\Lambda, K_0;\{\bar\alpha_i\}^m_{i=1};\omega,\mu,J,A)_I^{\bf
t}$ to (4.3) and a rational cycle
\begin{eqnarray*}
\hspace{10mm} &&\quad {\mathcal C}_g^{\bf
t}(\Lambda, K_0;\{\bar\alpha_i\}^m_{i=1};\omega,\mu,J,A)\hspace{75mm}(5.2)\\
&&=\sum_{I\in{\mathcal
 N}_k}\frac{1}{|\Gamma_I|}\bigl\{\Pi^I_{g,m}:C(\Lambda, K_0;
 \{\bar\alpha_i\}^m_{i=1};\omega,\mu,J,A)_I^{\bf t}\to
 \overline{\mathcal  M}_{g,m}\bigr\}
 \end{eqnarray*}
in $\overline{\mathcal M}_{g,m}$ of dimension $r=2c_1(M)(A)+
2(3-n)(g-1)+ 2m-\sum^m_{i=1}\deg\alpha_i$, which  corresponds to
(4.5). If $r<0$ it is an empty set. In particular we get a
homology class
$$[{\mathcal C}_g^{\bf t}(\Lambda,K_0;\{\alpha_i\}^m_{i=1};\omega,\mu,J,A)]
\in H_r(\overline{\mathcal  M}_{g,m}, \mbox{\Bb Q}), \eqno(5.3)
$$
 corresponding with (4.6). As in $\S4$  we can
prove that it is independent of all related choices. We call this
class or the Poincar\'e dual of it  the {\it Gromov-Witten class
relative to} $\Lambda$. Our first localization formula
is\vspace{2mm}

\noindent{\bf Theorem 5.1.}\quad{\it Let $2g+m\ge 3$ and ${\rm
supp}({\rm EV}_{g,m}^\ast(\prod^m_{i=1}\alpha_i))\subset\Lambda$.
Then
$$[{\mathcal C}_g^{\bf
t}(\Lambda,K_0;\{\alpha_i\}^m_{i=1};\omega,\mu,J,A)]=[{\mathcal
C}_g^{\bf t}(K_0;\{\alpha_i\}^m_{i=1};\omega,\mu,J,A)] \eqno(5.4)
$$
and thus for any homology class $\kappa\in
H_{\ast}(\overline{\mathcal M}_{g,m};\mbox{\Bb Q})$ satisfying
(1.5) it holds that
$${\mathcal G}{\mathcal W}^{(\omega,\mu,
J)}_{A,g,m}(\kappa;\alpha_1,\cdots,\alpha_m)=\langle PD(\kappa),
[{\mathcal C}_g^{\bf
t}(\Lambda,K_0;\{\alpha_i\}^m_{i=1};\omega,\mu,J,A)]\rangle.
$$}

\noindent{\bf Proof.}\quad On the basis of the previous choices we
furthermore take points
$$[{\bf f}_i]\in\overline{\mathcal M}_{g,m}(M, A,
J;K_3)\setminus\Lambda,\; i=k+1,\cdots, l,$$
  the  corresponding uniformizers $\widetilde W_i$ and the $\Gamma_i$-invariant functions
  $\beta_i$, $i=k+1,\cdots,l$. Let $\tilde s_{i1}\cdots,\tilde
  s_{iq_i}$ be the chosen desired sections of all bundles $\widetilde
  E_i\to\widetilde W_i$, $i=k+1,\cdots, l$. Since ${\mathcal
  W}^\ast(\Lambda)$ is a neighborhood of
$\Lambda$ in ${\mathcal B}^M_{A,g,m}$ we can assume
  $$\Lambda\cap
  W_i=\emptyset,\;i=k+1,\cdots,l.\eqno(5.5)
$$
Denote by ${\mathcal N}_l$ the collection of the subset $I$ of
$\{1,\cdots,l\}$ with $W_I=\cap_{i\in I}W_i\ne\emptyset$. Starting
from $\widetilde W_i$, $\tilde s_{i1},\cdots,\tilde s_{iq_i}$,
$i=1,\cdots,l$
 we can repeat the arguments in $\S3$ to
 construct a system of bundles $(\widehat{\mathcal
E}^{\ast\prime},\widehat V^{\ast\prime})=\{(\widehat
E_I^{\ast\prime}, \widehat
V_I^{\ast\prime}),\,\hat\pi_I^\prime\,|\, I\in{\mathcal N}_l\}$
and a family of cobordant virtual moduli chains in ${\cal
W}^{\ast\prime}=\cup_{I\in{\cal N}_l}V_I^{\ast\prime}$,
$${\mathcal C}^{\bf t}(K_0):=\sum_{I\in{\mathcal
 N}_l}\frac{1}{|\Gamma_I|}\{\hat\pi_I^\prime:{\mathcal M}^{\bf t}_I(K_0)\to{\mathcal
 W}^{\ast\prime}\}$$
for any ${\bf t}\in{\bf B}_{\varepsilon'}^{res}(\mbox{\Bb
R}^{q'})$  associated with $\Lambda$. Here $q'=q_1+\cdots + q_l$.
Let
$$
{\mathcal C}_g^{\bf
t}(K_0;\{\bar\alpha_i\}^m_{i=1};\omega,\mu,J,A)
=\sum_{I\in{\mathcal
 N}_l}\frac{1}{|\Gamma_I|}\{\Pi^I_{g,m}:C(K_0;
 \{\bar\alpha_i\}^m_{i=1};\omega,\mu,J,A)_I^{\bf t}\to
 \overline{\mathcal  M}_{g,m}\}
 $$
be the rational cycle in $\overline{\mathcal M}_{g,m}$ obtained
from  ${\mathcal C}^{\bf t}(K_0)$ as in (4.5). Note that
${\mathcal N}_k\subset{\mathcal N}_l$ and that
$\Gamma_I^\prime=\Gamma_I$ for any $I\in{\mathcal N}_k$. By (5.5)
we have
\begin{eqnarray*}
\hspace{22mm}&&\quad {\mathcal C}_g^{\bf
t}(K_0;\{\bar\alpha_i\}^m_{i=1};\omega,\mu,J,A)\hspace{65mm}(5.6)\\
&&=\sum_{I\in{\mathcal
 N}_k}\frac{1}{|\Gamma_I|}\{\Pi^I_{g,m}:C(K_0;
 \{\bar\alpha_i\}^m_{i=1};\omega,\mu,J,A)_I^{\bf t}\to
 \overline{\mathcal  M}_{g,m}\}.
 \end{eqnarray*}
 Note that $\widehat
V_I^{\ast\prime}\subset\widehat V_I^\ast$ for $I\in{\mathcal N}_k$
and that ${\rm suppt}({\rm
EV}_{g,m}^\ast(\wedge^m_{i=1}\alpha_i))$ is contained in
$\Lambda\subset\cup_{I\in{\mathcal N}_k}V_I^{\ast\prime}$.
 By the arguments in
$\S4.1$ it is not hard to see that when restricted on
$\cup_{I\in{\mathcal N}_k}V_I^{\ast\prime}$ the virtual moduli
chains ${\mathcal C}^{\bf t}(K_0)$ and ${\mathcal C}^{{\bf
t}^{(1)}}(K_0, \Lambda)$ are cobordant for sufficiently small
${\bf t}=({\bf t}^{(1)}, {\bf t}^{(2)})\in{\bf
B}_{\varepsilon'}^{res}(\mbox{\Bb R}^{q'})$ with ${\bf
t}^{(1)}\in{\bf B}_{\varepsilon}^{res}(\mbox{\Bb R}^{q})$. It
follows that the rational cycles as in (5.2) and (5.6), ${\mathcal
C}_g^{{\bf t}^{(1)}}(\Lambda,
K_0;\{\bar\alpha_i\}^m_{i=1};\omega,\mu,J,A)$ and ${\mathcal
C}_g^{\bf t}(K_0;\{\bar\alpha_i\}^m_{i=1};\omega,\mu,J,A)$
 are cobordant. This leads to (5.4).\hfill$\Box$\vspace{2mm}

 Note that for each nonempty compact subset
$G\subset\overline{\mathcal M}_{g,m}$  the set
$$\overline{\mathcal M}_{g,m}(M, A, J;K_3, G):=
\{[{\bf f}]\in\overline{\mathcal M}_{g,m}(M, A, J;K_3)\,|\,
\Pi_{g,m}([{\bf f}])\subset G\}$$
 is a compact subset of $\overline{\mathcal M}_{g,m}(M, A,
 J;K_3)$. Take $\Lambda=\overline{\mathcal M}_{g,m}(M, A,
J;K_3, G)$ we have \vspace{2mm}

\noindent{\bf Corollary 5.2.}\quad{\it Let $2g+m\ge 3$ and (1.5)
hold. If a homology class $\kappa\in H_{\ast}(\overline{\mathcal
M}_{g,m};\mbox{\Bb Q})$ has a cycle representative with the image
set contained in $G$ then
 \begin{eqnarray*}
&&{\mathcal G}{\mathcal W}^{(\omega,\mu,
J)}_{A,g,m}(\kappa;\alpha_1,\cdots,\alpha_m) =\langle PD(\kappa),
[{\mathcal C}_g^{\bf
t}(K_0;\{\alpha_i\}^m_{i=1};\omega,\mu,J,A)]\rangle\\
&&\hspace{40mm}=\langle PD(\kappa), [{\mathcal C}_g^{\bf
t}(\Lambda,K_0;\{\alpha_i\}^m_{i=1};\omega,\mu,J,A)]\rangle.
\end{eqnarray*}}

 Now we are a position to give  the second localization formula.
Denote by ${\mathcal B}^{M, S}_{A,g,m}:=(\Pi_{g,m})^{-1}(S)$ and
${\mathcal E}^{M,S}_{A,g,m}:={\mathcal E}^M_{A,g,m}|_{{\mathcal
B}^{M, S}_{A,g,m}}$, and by $\Pi^s_{g,m}$ and ${\rm EV}_{g,m}^s$
the restrictions of $\Pi_{g,m}$ and ${\rm EV}_{g,m}$ to ${\mathcal
B}^{M,S}_{A,g,m}$ respectively. Since $S$ is a compact suborbifold
of $\overline{\mathcal M}_{g,m}$,  $\overline{\mathcal M}_{g,m}(M,
A, J;K_3, S)$ is a compact subset of ${\mathcal B}^{M,
S}_{A,g,m}$. For each point $[{\bf f}]\in \overline{\mathcal
M}_{g,m}(M, A, J;K_3, S)$ we can, as before, construct a
neighborhood $W^s_{\bf f}$ of $[{\bf f}]$ in ${\mathcal B}^{M,
S}_{A,g,m}$ and a uniformizer $\pi^s_{\bf f}: \widetilde W^s_{\bf
f}\to W^s_{\bf f}$ with uniformization group $\Gamma_{\bf f}^s$.
Repeating the arguments in $\S3$ and $\S4$ we can construct a
virtual moduli chain of dimension $2k+ 2c_1(A)+ 2n-2ng$
$$
{\mathcal C}^{\bf t}_S(K_0):=\sum_{I\in{\mathcal
 N}^s}\frac{1}{|\Gamma_I^s|}\{\hat\pi_I^s:{\mathcal M}^{\bf t}_{SI}(K_0)\to{\mathcal
 W}^\ast_S\}
$$
 in some neighborhood ${\mathcal W}^\ast_S$ of $\overline{\mathcal M}_{g,m}(M, A, J;K_3, S)$
 in ${\mathcal B}^{M,S}_{A,g,m}$  associated with
$\overline{\mathcal M}_{g,m}(M, A, J;K_3, S)$.

As above we replace ${\rm EV}_{g,m}\circ{\mathcal C}^{\bf t}(K_0)$
in (4.2) by ${\rm EV}_{g,m}\circ{\mathcal C}^{\bf t}_S(K_0)$ and
get the corresponding fibre product
$C_S(K_0;\{\bar\alpha_i\}^m_{i=1};\omega,\mu,J,A)_I^{\bf t}$ with
(4.3) and a rational cycle
$$
{\mathcal C}_g^{\bf
t}(K_0;\{\bar\alpha_i\}^m_{i=1};\omega,\mu,J,A)_S
=\sum_{I\in{\mathcal
 N}}\frac{1}{|\Gamma_I|}\{\Pi^I_{g,m}:C_S(K_0;
 \{\bar\alpha_i\}^m_{i=1};\omega,\mu,J,A)_I^{\bf t}\to
 \overline{\mathcal  M}_{g,m}\}.
$$
in $S$ of dimension $r_s=2c_1(M)(A)+ 2k+ 2n(1-g)
-\sum^m_{i=1}\deg\alpha_i$. The latter  corresponds to (4.5). If
$r_s<0$ it is an empty set. In particular we get a homology class
$$
[{\mathcal C}_g^{\bf
t}(K_0;\{\alpha_i\}^m_{i=1};\omega,\mu,J,A)_S] \in H_{r_s}(S,
\mbox{\Bb Q}),
$$
corresponding with (4.6). As in $\S4$  we can prove that it is
independent of all related choices. We call this class or the
Poincar\'e dual of it the {\it Gromov-Witten class with value in}
$H_\ast(S, \mbox{\Bb Q})$.

Note that ${\rm codim}(\kappa)$ in (1.5) is equal to $6g-6+2m-j$.
But the codimension  ${\rm codim}_S(\kappa)$ of $\kappa$ in $S$ is
equal to $2k-j$. Therefore (1.5) implies
$$\sum_{i=1}^m\deg\alpha_i+{\rm
codim}_S(\kappa)=2c_1(M)(A)+ 2n- 2ng+ 2k. \eqno(5.7)$$
   We define
 $$
{\mathcal G}{\mathcal W}^{(\omega,\mu, J)}_{A,g,m, S
}(\kappa;\alpha_1,\cdots,\alpha_m) :=\langle PD_S(\kappa),
[{\mathcal C}_g^{\bf
t}(K_0;\{\alpha_i\}^m_{i=1};\omega,\mu,J,A)_S]\rangle, \eqno(5.8)
$$
and ${\mathcal G}{\mathcal W}^{(\omega,\mu,
J)}_{A,g,m,S}(\kappa;\alpha_1,\cdots,\alpha_m)=0$ if (1.5) ( or
(5.7)) is not satisfied. As in $\S4$ we may prove that ${\mathcal
G}{\mathcal W}^{(\omega,\mu,
J)}_{A,g,m,S}(\kappa;\alpha_1,\cdots,\alpha_m)$ has the similar
properties to Theorem 1.1.  Our second localization formula
is\vspace{2mm}

\noindent{\bf Theorem 5.3.}\quad{\it Under the above assumptions
it holds that
$${\mathcal G}{\mathcal W}^{(\omega,\mu,
J)}_{A,g,m, S }(\kappa;\alpha_1,\cdots,\alpha_m)={\mathcal
G}{\mathcal W}^{(\omega,\mu,
J)}_{A,g,m}(\kappa;\alpha_1,\cdots,\alpha_m).\eqno(5.9)
$$}

\noindent{\it Proof.}\quad From the arguments above Proposition
5.1 we can take finitely many points
  $[{\bf f}_i]\in\Lambda:=\overline{\mathcal M}_{g,m}(M,  A, J;K_3, S)$,
$i=1,\cdots, k$, to construct a system of bundles
 $(\widehat{\mathcal E}^\ast,\widehat V^\ast)=\{(\widehat E_I^\ast, \widehat
V_I^\ast),\,\hat\pi_I\,|\, I\in{\mathcal N}_k\}$, and then choose
sections $\tilde s_{ij}:\widetilde W_i\to\widetilde E_i$,
$i=1,\cdots, k$ and $j=1,\cdots, q_i$ to construct a family of
cobordant virtual moduli chains in ${\mathcal
W}^\ast(\Lambda)\subset{\mathcal B}^M_{A,g,m}$ of dimension $2m+
2c_1(A)+ 2(3-n)(g-1)$
$${\mathcal C}^{\bf t}(K_0, \Lambda)=\sum_{I\in{\mathcal
 N}_k}\frac{1}{|\Gamma_I|}\{\hat\pi_I:{\mathcal M}^{\bf t}_I(K_0, \Lambda)\to{\mathcal
 W}^\ast(\Lambda)\}
$$
for any ${\bf t}\in{\bf B}_{\varepsilon}^{res}(\mbox{\Bb R}^{q})$
associated with $\overline{\mathcal M}_{g,m}(M, A, J;K_3, S)$.
Here $q=q_1+\cdots + q_k$. They determine a rational homology
class as in (5.3). By Corollary 5.2 we have
 $$\left.\begin{array}{ll}
{\mathcal G}{\mathcal W}^{(\omega,\mu,
J)}_{A,g,m}(\kappa;\alpha_1,\cdots,\alpha_m) =\langle PD(\kappa),
[{\mathcal C}_g^{\bf
t}(\Lambda,K_0;\{\alpha_i\}^m_{i=1};\omega,\mu,J,A)]\rangle
\end{array}\right.\eqno(5.10)
$$
Let $R_S$ be a retraction from a tubular neighborhood of $S$ in
$\overline{\mathcal M}_{g,m}$ onto $S$, and $PD_S(\kappa)$ (resp.
$PD(S)$) be the Poincar\'e dual of $\kappa$ in $H^\ast(S,\mbox{\Bb
Q})$ (resp. $S$ in $H^\ast(\overline{\mathcal M}_{g,m},\mbox{\Bb
Q})$). Since one can take a closed form representative $S^\ast$ of
$PD(S)$ with support in  the tubular neighborhood, for any closed
form $\kappa_S^\ast$ representing $PD_S(\kappa)$ we get such a
closed representative form
$\kappa^\ast:=R_S^\ast(\kappa_S^\ast)\cup S^\ast$ of the
Poincar\'e dual $PD(\kappa)$ in $H^{6g-6+2m-j}(\overline{\mathcal
M}_{g,m},\mbox{\Bb Q})$ of $\kappa\in H_j(\overline{\mathcal
M}_{g,m},\mbox{\Bb Q})$.
 So (5.10)
becomes
\begin{eqnarray*}
 \hspace{10mm}\quad{\mathcal G}{\mathcal W}^{(\omega,\mu,
J)}_{A,g,m}(\kappa;\alpha_1,\cdots,\alpha_m)\!\!
&=&\!\!\int_{{\mathcal C}_g^{\bf
t}(\Lambda,K_0;\{\alpha_i\}^m_{i=1};\omega,\mu,J,A)}R_S^\ast(\kappa_S^\ast)\cup S^\ast\hspace{20mm}(5.11)\\
\!\!&=&\!\!\int_{S\cap{\mathcal C}_g^{\bf
t}(\Lambda,K_0;\{\alpha_i\}^m_{i=1};\omega,\mu,J,A)}\kappa_S^\ast.
\end{eqnarray*}

 Let $W_i^s=W_i\cap{\mathcal B}^{M,S}_{A,g,m}$,
$\widetilde W_i^s=\pi_i^{-1}(W_i^s)$,
$\beta_i^s=\beta_i|_{\widetilde W_i^s}$ and $\tilde
s^s_{ij}=\tilde s_{ij}|_{\widetilde W_i^s}$,
$\pi_i^s=\pi_i|_{\widetilde W_i^s}$, $i=1,\cdots,k$, and
$j=1,\cdots, q_i$. As before we assume that
$$W_i^1\subset\subset U_i^1\subset\subset W_i^2\subset\subset
U_i^2\cdots\subset\subset W_i^{k-1}\subset\subset
U_i^{k-1}\subset\subset W_i^k=W_i,$$
 $i=1,\cdots,k$ are the pairs of open subsets constructing ${\mathcal C}^{\bf t}(K_0, \Lambda)$ above.
 Set
 $$W_i^{js}=W_i^j\cap{\mathcal B}^{M,S}_{A,g,m}\quad{\rm and}\quad
 U_i^{js}=U_i^j\cap{\mathcal B}^{M,S}_{A,g,m},$$
 then $\widetilde W_i^{js}=(\pi_i^s)^{-1}(W_i^{js})=\pi_i^{-1}(W_i^{js})$ and $\widetilde
U_i^{js}=(\pi_i^s)^{-1}(U_i^{js})=\pi_i^{-1}(U_i^{js})$. As in
(3.62), for each $I\in{\mathcal N}_k$ with $|I|=r$ we define
$$
 V^s_I:=\Bigl(\bigcap_{i\in I}W^{rs}_i\Bigr)\setminus
\Bigl(\bigcup_{J:|J|>r}\bigl(\bigcap_{j\in
J}Cl(U^{rs}_j)\bigr)\Bigr)
$$
and $V_I^{s\ast}=V_I^\ast\cap V_I^s$. Then it is not hard to check
that $V_I^s$ and $V_I^{s\ast}$ are equal to $V_I\cap{\mathcal
B}^{M,S}_{A,g,m}$ and $V_I^\ast\cap{\mathcal B}^{M,S}_{A,g,m}$
respectively. It follows that the corresponding system of bundles
produced from them $(\widehat{\mathcal
E}^{s\ast},\widehat{\mathcal V}^{s\ast})=\{(\widehat
E^{s\ast}_I,\widehat V_I^{s\ast}),\,\hat\pi_I^s\,|\,I\in{\mathcal
N}_k\}$ is equal to the restriction of the system of bundles
$(\widehat{\mathcal E}^\ast,\widehat V^\ast)=\{(\widehat
E_I^\ast,\widehat V_I^\ast),\,\hat\pi_I\,|\,I\in{\mathcal N}_k\}$
to ${\mathcal B}^{M,S}_{A,g,m}$. By increasing the number of the
sections $\tilde s_{ij}$ if necessary we may assume that the
restriction of the section $\Psi_I^{(\bf t)}$ in Corollary 3.15 to
$(\widehat{\mathcal E}^{s\ast},\widehat{\mathcal V}^{s\ast})$ is
transversal to the zero section. It follows from that the rational
cycle ${\mathcal C}_g^{\bf
t}(K_0;\{\bar\alpha_i\}^m_{i=1};\omega,\mu,J,A)_S$ in $S$
associated with it is equal to the intersection of $S$ with the
rational cycle ${\mathcal C}_g^{\bf
t}(\Lambda,K_0;\{\bar\alpha_i\}^m_{i=1};\omega,\mu,J,A)$ in
$\overline{\mathcal M}_{g,m}$. So we have
\begin{eqnarray*}
\langle PD_S(\kappa), [{\mathcal C}_g^{\bf
t}(K_0;\{\bar\alpha_i\}^m_{i=1};\omega,\mu,J,A)_S]\rangle
\!\!\!\!&=&\!\!\!\int_{{\mathcal C}_g^{\bf
t}(K_0;\{\bar\alpha_i\}^m_{i=1};\omega,\mu,J,A)_S}\kappa_S^\ast\\
&=&\!\!\!\int_{S\cap{\mathcal C}_g^{\bf
t}(K_0;\{\bar\alpha_i\}^m_{i=1};\omega,\mu,J,A)}\kappa_S^\ast.
\end{eqnarray*}
 This and (5.11) together give rise to
\begin{eqnarray*}
&&\quad{\mathcal G}{\mathcal W}^{(\omega,\mu,
J)}_{A,g,m}(\kappa;\alpha_1,\cdots,\alpha_m)\\
&& =\langle PD_S(\kappa), [{\mathcal C}_g^{\bf
t}(K_0;\{\alpha_i\}^m_{i=1};\omega,\mu,J,A)_S]\rangle\\
&&={\mathcal G}{\mathcal W}^{(\omega,\mu,
J)}_{A,g,m,S}(\kappa;\alpha_1,\cdots,\alpha_m).
\end{eqnarray*}
\hfill$\Box$\vspace{2mm}

A special case of the formula (5.9) is the corresponding form of
the formula (4.49) in [R]. For the invariant ${\mathcal
G}{\mathcal W}^{(\omega,\mu, J)}_{A,g,m, S }$ we may also obtain a
similar localization formula to Theorem 5.1.

\subsection{The composition laws}

The composition laws for the Gromov-Witten invariants are the most
useful properties. To state them we assume that $g=g_1+g_2$ and
$m=m_1+m_2$ with $2g_i+m_i\ge 3$. Following [KM] we fix a
decomposition $Q=Q_1\cup Q_2$ of $\{1,\cdots,m\}$ with $|Q_i|=m_i$
and then get a canonical embedding $\varphi_Q: \overline{\mathcal
M}_{g_1, m_1+1}\times\overline{\mathcal M}_{g_2,
m_2+1}\to\overline{\mathcal M}_{g, m}$, which assigns marked
curves $(\Sigma_i, z_1^{(i)},\cdots, z_{m_i+1}^{(i)})$( $i=1,2$)
to their union $\Sigma_1\cup\Sigma_2$ with $z^{(1)}_{m_1+1}$
identified to $z_1^{(2)}$ and remaining points renumbered by
$\{1,\cdots,m\}$ in such a way that their relative order is kept
intact and points on $\Sigma_i$ are numbered by $S_i$. There is
also another natural embedding  $\psi:\overline{\mathcal M}_{g-1,
m+2}\to \overline{\mathcal M}_{g,m}$ obtained by gluing together
the last two marked points. So $\varphi_Q(\overline{\mathcal
M}_{g_1, m_1+1}\times\overline{\mathcal M}_{g_2, m_2+1})$ is a
compact complex suborbifold of $\overline{\mathcal M}_{g, m}$ of
complex dimension $3g-4+m$, and  $\psi(\overline{\mathcal M}_{g-1,
m+2})$ is a compact complex suborbifold of $\overline{\mathcal
M}_{g,m}$  complex dimension $3g-4+m$.
 \vspace{2mm}

\noindent{\bf Lemma 5.4.}\quad{\it Assume that $\dim
H^\ast(M)<\infty$ and that  $\{\beta_i\}$ is  a basis of the
vector space $H^\ast(M,\mbox{\Bb Q})$. By Remark 5.7 of [BoTu] one
has the Poincar\'e duality: $H^p(M)\cong(H_c^{n-p}(M))^\ast$ for
any integer $p$. Thus  $\dim H_c^\ast(M,\mbox{\Bb Q})=\dim
H^\ast(M)<\infty$ and we may choose a dual basis $\{\omega_i\}$ of
$\{\beta_i\}$ in $H^\ast_c(M)$, i.e.,
$\langle\omega_i,\beta_j\rangle=\langle\omega_i\wedge\beta_j,
[M]\rangle=\int_M\omega_i\wedge\beta_j=\delta_{ij}$.
 Then the Poincar\'e dual $PD^{\rm II}([\Delta_M])$ of
the class $[\Delta_M]\in H^{\rm II}_{2n}(M\times M,\mbox{\Bb Q})$
of the diagonal $\Delta_M$ in $M\times M$ are given by
$$PD^{\rm II}([\Delta_M])=\sum_{\deg\beta_i+\deg\beta_j=2n}c_{ij}
\rho_1^\ast\beta_i\wedge\rho_2^\ast\beta_j.$$ where $\rho_k$ are
the projections of $M\times M$ to the k-th factor, $k=1,2$, and
$$\left.\begin{array}{ll}
c_{ij}=(-1)^{\deg\omega_i\cdot\deg\omega_j}\eta^{ij}\quad{\rm
and}\quad \eta^{ij}=\int_M\omega_i\wedge\omega_j.
\end{array}\right.$$
Moreover, if $\beta_i^\ast$  are the closed form representatives
of $\beta_i$, then $PD^{\rm II}([\Delta_M])$  has the closed
representing form
$$\Delta_M^\ast=\sum_{\deg\beta_i+\deg\beta_j=2n}
c_{ij}\rho_1^\ast\beta_i^\ast\wedge\rho_2^\ast\beta_j^\ast.$$}\vspace{2mm}

  By (4.9), $H^0(M,\mbox{\Bb Q})=\mbox{\Bb Q}$ and
$H^{2n}(M,\mbox{\Bb Q})=0$. So each
$2n-\deg\beta_i=\deg\omega_i>0$. We always assume $\beta_1={\bf
1}$. Note that in general  the matrix $(\eta^{ij})$ is degenerate
for noncompact $M$. If $M$ is a closed manifold then $(\eta^{ij})$
is invertible and $(\eta^{ij})^{-1}=(\eta_{ij})$, where
$\eta_{ij}=\int_M\beta_i\wedge\beta_j$. In this case
$\omega_i=\sum_j\eta^{ij}\beta_j$.\vspace{2mm}

\noindent{\it Proof of Lemma 5.4.}\quad Since $H^\ast(M)$ is
finitely dimensional it follows from the K\"unneth formula (cf.
Proposition 9.12 in [BoTu]) that  $\{\beta_i\otimes\beta_j\}$ form
a basis $H^\ast(M\times M,\mbox{\Bb Q})$. (This is a unique place
where $\dim H^\ast(M)<\infty$ is used.) So $PD^{\rm
II}([\Delta_M])\in H^{2n}(M\times M)$ and $[\Delta_M]\in
H_{2n}^{\rm II}(M\times M)$ can be written as
$$\left.\begin{array}{ll}
\sum_{i,j}c_{ij}\beta_i\otimes\beta_j\quad{\rm and}\quad
\sum_{i,j}c_{ij}(PD^{\rm II})^{-1}(\beta_i)\otimes(PD^{\rm
II})^{-1}(\beta_j)\end{array}\right.$$
 respectively.  Now on one hand
\begin{eqnarray*}
&& \langle [\Delta_M],
\rho_1^\ast[\omega_k]\wedge\rho_2^\ast[\omega_l]\rangle=\int_{\Delta_M}
\rho_1^\ast\omega_k\wedge\rho_2^\ast\omega_l\\
&&\hspace{35mm}=\int_M\imath^\ast\rho_1^\ast\omega_k\wedge\imath^\ast\rho_2^\ast\omega_l=
\int_M\omega_k\wedge\omega_l,
\end{eqnarray*}
 where $\imath: M\to\Delta_M\subset M\times M$ is the diagonal map. On
the other hand
\begin{eqnarray*}
 &&\quad\langle [\Delta_M],
\rho_1^\ast[\omega_k]\wedge\rho_2^\ast[\omega_l]\rangle\\
&&=\sum_{i,j}c_{ij}\langle(PD^{\rm
II})^{-1}(\beta_i)\otimes(PD^{\rm
II})^{-1}(\beta_j),\rho_1^\ast[\omega_k]\wedge\rho_2^\ast[\omega_l]\rangle\\
&&=\sum_{i,j}c_{ij}(-1)^{\deg\beta_i\cdot\deg\beta_j}\langle
[\omega_k],(PD^{\rm II})^{-1}(\beta_i)\rangle \langle[\omega_l],
(PD^{\rm II})^{-1}(\beta_j)\rangle.
\end{eqnarray*}
 So
$c_{ij}=(-1)^{\deg\beta_i\cdot\deg\beta_j}
\int_M\omega_i\wedge\omega_j=(-1)^{\deg\omega_i\cdot\deg\omega_j}
\int_M\omega_i\wedge\omega_j$ because $\deg(PD^{\rm II})^{-1}(
\beta_i)=2n-\deg\beta_i=\deg\omega_i$. By the definition of the
Poincar\'e dual of a closed oriented submanifold on the page 51 of
[BoTu], using the similar reason one may readily prove the second
statement because $c_{ij}\ne 0$ implies that
$\deg\omega_i+\deg\omega_j=\dim M=2n$ and thus that
$\deg\beta_i\deg\omega_i+\deg\beta_j\deg\omega_j$ is
even.\hfill$\Box$\vspace{2mm}

\noindent{\bf Theorem 5.5.}\quad{\it Assume that $\dim
H^\ast(M)<\infty$. Let $\kappa\in H_\ast(\overline{\mathcal
M}_{g-1, m+2},\mbox{\Bb Q})$,  and $\alpha_i\in H^\ast(M,\mbox{\Bb
Q})$, $i=1,\cdots, m$. Suppose that some
 $\alpha_t\in H_c^\ast(M,\mbox{\Bb Q})$.
 Then
 \begin{eqnarray*}
{\mathcal G}{\mathcal W}^{(\omega, \mu,
J)}_{A,g,m}(\psi_\ast(\kappa);\alpha_1,\cdots,\alpha_m) =
\sum_{i,j} c_{ij}\cdot {\mathcal G}{\mathcal W}^{(\omega, \mu,
J)}_{A, g-1, m+2}(\kappa;\alpha_1,\cdots,\alpha_m,\beta_i,
\beta_j)
\end{eqnarray*}}

\noindent{\it Proof.}\quad Take a compact subset $K_0\subset M$
containing ${\rm supp}(\wedge^m_{i=1}\alpha_i)$. Since $\psi$ is
an embedding $S:=\psi(\overline{\mathcal M}_{g-1, m+2})$ is a
compact complex suborbifold of  $\overline{\mathcal M}_{g,m}$
complex dimension $3g-4+m$. Consider the evaluation
$$\Xi:{\mathcal B}^M_{A, g-1, m+2}\to M\times M
$$
 given by $[\Sigma,\bar{\bf z}, f]\mapsto
(f(z_{m+1}), f(z_{m+2})$. There is a map $\Psi:
\Xi^{-1}(\Delta_M)\to{\mathcal B}^M_{A, g, m}$, which is a lifting
of $\psi$, such that the following commutative diagram holds:
\begin{center}\setlength{\unitlength}{1mm}
\begin{picture}(80,30)
\thinlines \put(15,25){$\Xi^{-1}(\Delta_M)$} \put(35,25){\vector
(1,0){34}} \put(45,27){$\Psi$} \put(25,22){\vector(0,-1){16}}
\put(8,13){$\Pi_{g-1, m+2}^r$} \put(75,25){${\mathcal B}^M_{A, g,
m}$} \put(80,22){\vector(0,-1){16}} \put(71,13){$\Pi_{g,m}$}
\put(15,0){$\overline{\mathcal M}_{g-1,m+2}$}
\put(35,1){\vector(1,0){34}} \put(45,3){$\psi$}
\put(75,0){$\overline{\mathcal M}_{g,m}$}
\end{picture}\end{center}
Here $\Pi_{g-1, m+2}^r$ is the restriction of $\Pi_{g-1, m+2}$ to
$\Xi^{-1}(\Delta_M)$. For  $[{\bf f}]\in{\mathcal B}^M_{A, g-1,
m+2}$ let $\pi_{\bf f}:\widetilde W_{\bf f}\to W_{\bf f}$ be a
uniformizer  of a neighborhood of $[{\bf f}]$ in ${\mathcal
B}^M_{A, g-1, m+2}$ constructed in $\S2.4$. One has a natural
lifting $\widetilde\Xi_{\bf f}:\widetilde W_{\bf f}\to M\times M$
of $\Xi$ given by $(\Sigma,\bar{\bf z}, f)\mapsto (f(z_{m+1}),
f(z_{m+2})$. It is easily checked that $\widetilde\Xi_{\bf f}$ is
transversal to the diagonal $\triangle_M$ and that
$(\widetilde\Xi_{\bf f})^{-1}(\Delta_M)$ is invariant under the
action of group $\Gamma_{\bf f}={\rm Aut}({\bf f})$. Therefore
$(\widetilde\Xi_{\bf f})^{-1}(\Delta_M)$ is a $\Gamma_{\bf
f}$-invariant stratified Banach submanifold of $\widetilde W_{\bf
f}$. It may follow that $\Xi^{-1}(\Delta_M)$ is a stratified
Banach suborbifold of ${\mathcal B}^M_{A, g-1, m+2}$ and that
$\pi_{\bf f}:(\widetilde\Xi_{\bf f})^{-1}(\Delta_M)\to W_{\bf
f}\cap\Xi^{-1}(\Delta_M)$ gives a uniformizer of a neighborhood of
$[{\bf f}]$ in $\Xi^{-1}(\Delta_M)$.  Note that
$\overline{\mathcal M}_{g,m}(M,J,A;K_3)\ne\emptyset$ implies that
$\overline{\mathcal
M}_{g-1,m+2}(M,J,A;K_3)\cap\Xi^{-1}(\Delta_M)\ne\emptyset$. As in
the proof of Theorem 5.3 we first choose finitely many points
$[{\bf f}_1],\cdots,[{\bf f}_k]$ in $\overline{\mathcal
M}_{g-1,m+2}(M,J,A;K_3)\cap\Xi^{-1}(\Delta_M)$ and then choose
finitely many points $[{\bf f}_{k+1}],\cdots,[{\bf f}_l]$ in
$\overline{\mathcal
M}_{g-1,m+2}(M,J,A;K_3)\setminus\overline{\mathcal
M}_{g-1,m+2}(M,J,A;K_3)\cap\Xi^{-1}(\Delta_M)$ to construct a
virtual moduli chain
$${\mathcal C}^{\bf t}(K_0)=\sum_{I\in{\mathcal
 N}_l}\frac{1}{|\Gamma_I|}\{\hat\pi_I:{\mathcal M}^{\bf t}_I(K_0)\to{\mathcal
 W}^\ast\}
$$
 of dimension $2(m+2)+
2c_1(A)+ 2(3-n)(g-2)$ in ${\cal W}^\ast=\cup_{I\in{\cal
N}_r}V_I^\ast\subset{\cal W}$ associated with $\overline{\mathcal
M}_{g-1,m+2}(M,J,A;K_3)$. As in $\S4.1$ we can use it to get a
rational cycle
$${\mathcal C}_{g-1}^{\bf t}(K_0;\{\bar\alpha_i\}^m_{i=1},
PD^{II}(\triangle_M);\omega,\mu,J,A)\eqno(5.12)
$$
 in $\overline{\mathcal M}_{g-1, m+2}$ and its homology class
 $$[{\mathcal
C}_{g-1}^{\bf t}(K_0;\{\alpha_i\}^m_{i=1},
PD^{II}(\triangle_M);\omega,\mu,J,A)]
$$
in $H_{r_1}(\overline{\mathcal M}_{g-1, m+2},\mbox{\Bb Q})$ for
$r_1=2c_1(A)+2m+2(3-n)(g-1)-2-\sum^m_{i=1}\deg\alpha_i$.

 Since $f_i(z^{(i)}_{m+1})\ne
f_i(z^{(i)}_{m+2})$ for $i=k+1,\cdots,l$ we can require
$$(\cup^l_{i=k+1}W_i)\cap\Xi^{-1}(\Delta_M)=\emptyset.
$$
 Recall that the constructions of the
virtual moduli chains in $\S3$. It is not hard using our above
arguments to show that for ${\mathcal N}_k:=\{I\in{\mathcal
N}_l\,|\,\max(I)\le k\}$,
$${\mathcal C}^{\bf t}_{\Xi^{-1}(\Delta_M)}(K_0):=\sum_{I\in{\mathcal
 N}_k}\frac{1}{|\Gamma_I|}\{\hat\pi_{Ir}:{\mathcal M}^{\bf t}_I(K_0)_r\to{\mathcal
 W}^\ast\cap\Xi^{-1}(\Delta_M)\}
$$
is a virtual moduli chain  of dimension
 $2(m+2)+ 2c_1(A)+ 2(3-n)(g-2)-2n$  in ${\mathcal
 W}^\ast\cap\Xi^{-1}(\Delta_M)$ associated with
 $\overline{\mathcal M}_{g-1,m+2}(M,J,A;K_3)\cap\Xi^{-1}(\Delta_M)$,
 where ${\mathcal M}^{\bf t}_I(K_0)_r=(\Xi\circ\hat\pi_I)^{-1}(\Delta_M)$
 and $\hat\pi_{Ir}$ denotes the restriction of $\hat\pi_I$ to ${\mathcal M}^{\bf
t}_I(K_0)_r$. (This actually implies that $\Xi\circ\hat\pi_I:
{\mathcal M}^{\bf t}_I(K_0)\to M\times M$ is transversal to
$\Delta_M$.) As in $\S4.1$ we use
$(\prod^m_{i=1}ev_i)\circ{\mathcal C}^{\bf
t}_{\Xi^{-1}(\Delta_M)}(K_0)$ to construct a rational cycle
$${\mathcal C}_{g-1}^{\bf t}(K_0;\{\bar\alpha_i\}^m_{i=1}
;\omega,\mu,J,A;\Xi^{-1}(\Delta_M))
$$
 in $\overline{\mathcal M}_{g-1, m+2}$ and its homology class
 $$[{\mathcal
C}_{g-1}^{\bf
t}(K_0;\{\alpha_i\}^m_{i=1};\omega,\mu,J,A;\Xi^{-1}(\Delta_M))]\in
H_{r_1}(\overline{\mathcal M}_{g-1, m+2},\mbox{\Bb Q})
$$
 for $r_1=2c_1(A)+2m+2(3-n)(g-1)-2-\sum^m_{i=1}\deg\alpha_i$.
It is easily checked that
\begin{eqnarray*}
\hspace{30mm}&& \quad [{\mathcal C}_{g-1}^{\bf
t}(K_0;\{\alpha_i\}^m_{i=1},
PD^{II}(\triangle_M);\omega,\mu,J,A)]\hspace{33mm}(5.13)\\
&&=[{\mathcal C}_{g-1}^{\bf t}(K_0;\{\alpha_i\}^m_{i=1}
;\omega,\mu,J,A;\Xi^{-1}(\Delta_M))].
\end{eqnarray*}

Note that $\Xi^{-1}(\Delta_M)$ can be identified with ${\mathcal
B}^{M,S}_{A,g,m}$ via $\Psi$. We may consider ${\mathcal C}^{\bf
t}_{\Xi^{-1}(\Delta_M)}(K_0)$  as a virtual moduli chain
 in ${\mathcal B}^{M,S}_{A,g,m}$ associated with $\overline{\mathcal
M}_{g,m}(M,J,A;K_3, S)$. In this case we have
$$
\psi_\ast[{\mathcal C}_{g-1}^{\bf t}(K_0;\{\alpha_i\}^m_{i=1}
;\omega,\mu,J,A;\Xi^{-1}(\Delta_M))]= [{\mathcal C}_g^{\bf
t}(K_0;\{\alpha_i\}^m_{i=1} ;\omega,\mu,J,A)_S].\eqno(5.14)
$$
 Then by Theorem 5.3, (5.14) and (5.13) we get
 \begin{eqnarray*}
 &&\quad{\mathcal G}{\mathcal W}^{(\omega, \mu,
J)}_{A,g,m}(\psi_{\ast}(\kappa);\alpha_1,\cdots,\alpha_m)\\
&&={\mathcal G}{\mathcal W}^{(\omega, \mu,
J)}_{A,g,m,S}(\psi_{\ast}(\kappa);\alpha_1,\cdots,\alpha_m)\\
&&=\langle PD_S(\psi_{\ast}(\kappa)), [{\mathcal C}_g^{\bf
t}(K_0;\{\alpha_i\}^m_{i=1};\omega,\mu,J,A)_S]\rangle\\
&&=\langle PD(\kappa), [{\mathcal C}_{g-1}^{\bf
t}(K_0;\{\alpha_i\}^m_{i=1};\omega,\mu,J,A;\Xi^{-1}(\Delta_M))]\rangle\\
&&=\langle PD(\kappa), [{\mathcal C}_{g-1}^{\bf
t}(K_0;\{\alpha_i\}^m_{i=1},
PD^{II}(\triangle_M);\omega,\mu,J,A)]\rangle\\
&&=\sum_{i,j}c_{ij}\langle PD(\kappa), [{\mathcal C}_{g-1}^{\bf
t}(K_0;\{\alpha_i\}^m_{i=1},
\beta_i\otimes\beta_j;\omega,\mu,J,A)]\rangle\\
&&=\sum_{i,j}c_{ij}{\mathcal G}{\mathcal W}^{(\omega, \mu,
J,A)}_{A,g-1,m+2}(\kappa;\alpha_1,\cdots,\alpha_m,\beta_i,\beta_j).
\end{eqnarray*}
Here the fifth equality comes from Lemma 5.4. Theorem 5.5 is
proved.
 \hfill$\Box$\vspace{2mm}

For $i=1,2$, let $\kappa_i\in H_\ast(\overline{\mathcal M}_{g_i,
m_i},\mbox{\Bb Q})$, $A_i\in H_2(M,\mbox{\Bb Z})$ and integers
$g_i\ge 0$, $m_i>0$ satisfy $2g_i+m_i\ge 3$. Moreover assume that
$\alpha_i,\gamma_k\in H^\ast(M,\mbox{\Bb Q})$ for $i=1,\cdots,m_1$
and $k=1,\cdots, m_2$. If
 $\alpha_s, \gamma_t\in H_c^\ast(M,\mbox{\Bb Q})$ for some $s$ and $t$ we may take a
  compact  subset $K_0\subset M$ containing ${\rm supp}(\alpha_s)$ and
  ${\rm supp}(\gamma_t)$.
Consider the product stratified Banach orbifold ${\mathcal
B}^M_{A_1,g_1,m_1}\times{\mathcal B}^M_{A_2,g_2,m_2}$. Repeating
the arguments in $\S3$ we can construct a virtual moduli chain
$$
{\mathcal C}^{\bf t}(K_0;\{A_i,g_i,m_i\})=\sum_{I\in{\mathcal
 N}_l}\frac{1}{|\Gamma_I|}\{\hat\pi_I:{\mathcal M}^{\bf t}_I(K_0;
  \{A_i,g_i,m_i\} )\to{\mathcal  W}\}
  $$
  of dimension $2(m_1+m_2)+2c_1(A_1+A_2)+2(3-n)(g_1+g_2-2)$ in ${\mathcal
B}^M_{A_1,g_1,m_1}\times{\mathcal B}^M_{A_2,g_2,m_2}$ associated
with $\overline{\mathcal M}_{g_1,m_1}(M, A_1, J_1;
K_3)\times\overline{\mathcal M}_{g_2,m_2}(M, A_2, J_2; K_3)$.
Using the same ideas as the proof of (4.16) of Proposition 4.3  in
[Lu3] we can derive that the virtual moduli chain ${\mathcal
C}^{\bf t}(K_0;\{A_i,g_i,m_i\})$ induces a virtual moduli chain
${\mathcal C}^{\bf t}(K_0;A_1,g_1,m_1)$ for $\overline{\mathcal
M}_{g_1,m_1}(M, A_1, J_1; K_3)$ and that ${\mathcal C}^{\bf
t}(K_0;A_2,g_2,m_2)$ for $\overline{\mathcal M}_{g_2,m_2}(M, A_2,
J_2; K_3)$, such that
$${\mathcal C}^{\bf
t}(K_0;\{A_i,g_i,m_i\})={\mathcal C}^{\bf
t}(K_0;A_1,g_1,m_1)\times{\mathcal C}^{\bf t}(K_0;A_2,g_2,m_2).
$$
 As in $\S4.1$ we can use them, $\alpha_i$, $i=1,\cdots,m_1$
 and $\gamma_j$, $j=1,\cdots,m_2$ to construct a rational cycle
 in $\overline{\mathcal
M}_{g_1,m_1}\times\overline{\mathcal M}_{g_2,m_2}$,
$${\mathcal C}^{\bf
t}_{(g_1,g_2)}(K_0;\{\alpha_i\}^{m_1}_{i=1};\{\gamma_j\}^{m_2}_{j=1};A_1,
A_2) \eqno(5.15)
$$
of dimension
$2(m_1+m_2)+2c_1(A_1+A_2)+2(3-n)(g_1+g_2-2)-\sum^{m_1}_{i=1}\alpha_i-
\sum^{m_2}_{j=1}\gamma_j$, and that
$$
{\mathcal C}^{\bf
t}_{g_1}(K_0;\{\alpha_i\}^{m_1}_{i=1};\omega,\mu, J, A_1)$$
 in $\overline{\mathcal
M}_{g_1,m_1}$ of dimension $2m_1+
2c_1(A_1)+2(3-n)(g_1-1)-\sum^{m_1}_{i=1}\alpha_i$, and  that
$${\mathcal C}^{\bf
t}_{g_2}(K_0;\{\gamma_j\}^{m_2}_{j=1};\omega,\mu, J, A_2)$$
 in $\overline{\mathcal
M}_{g_2,m_2}$ of dimension $2m_2+
2c_1(A_2)+2(3-n)(g_2-1)-\sum^{m_2}_{j=1}\gamma_j$, such that
\begin{eqnarray*}
&&\quad {\mathcal C}^{\bf
t}_{(g_1,g_2)}(K_0;\{\alpha_i\}^{m_1}_{i=1};\{\gamma_j\}^{m_2}_{j=1};A_1,
A_2)\\
&&={\mathcal C}^{\bf
t}_{g_1}(K_0;\{\alpha_i\}^{m_1}_{i=1};\omega,\mu, J, A_1)\times
{\mathcal C}^{\bf
t}_{g_2}(K_0;\{\gamma_j\}^{m_2}_{j=1};\omega,\mu, J, A_2).
\end{eqnarray*}
As in $\S4$ we may also prove that the homology class
$$\bigl[{\mathcal C}^{\bf
t}_{(g_1,g_2)}(K_0;\{\alpha_i\}^{m_1}_{i=1};\{\gamma_j\}^{m_2}_{j=1};A_1,
A_2)\bigr]\in H_\ast(\overline{\mathcal
M}_{g_1,m_1}\times\overline{\mathcal M}_{g_2,m_2},\mbox{\Bb Q})
\eqno(5.16)
$$
 is independent of all related choices.

 Let $P_i$ (resp. $\rho_i$) be the projections of ${\mathcal
B}^M_{A_1,g_1,m_1}\times{\mathcal B}^M_{A_2,g_2,m_2}$ (resp.
$\overline{\mathcal M}_{g_1,m_1}\times\overline{\mathcal
M}_{g_2,m_2}$) to the $i$-th factor, $i=1,2$. For $i=1,2$ let
${\rm ev}^{(i)}_k$ denote the evaluation at the $k$-th marked
point of ${\mathcal B}^M_{A_i,g_i,m_i}$, $k=1,\cdots, m_i$. If
$$\sum^{m_1}_{s=1}\deg\alpha_s+\!\sum^{m_2}_{r=1}\deg\gamma_r\!=
\!\sum^2_{i=1}\Bigl(2c_1(M)(A_i)+2(3-n)(g_i-1)+2m_i\Bigr)
$$
 we define
\begin{eqnarray*}
\hspace{10mm}&&\quad {\mathcal G}{\mathcal W}^{(\omega,\mu,
J)}_{\{A_i,g_i,m_i\}}(\kappa_1,\kappa_2;\alpha_1,\cdots,\alpha_{m_1};
\gamma_1,\cdots,\gamma_{m_2})\hspace{45mm}(5.17)\\
&&=(-1)^{\dim\kappa_1\dim\kappa_2}\bigl\langle
PD(\kappa_1\times\kappa_2), \bigl[{\mathcal C}^{\bf
t}_{(g_1,g_2)}(K_0;\{\alpha_i\}^{m_1}_{i=1};\{\gamma_j\}^{m_2}_{j=1};A_1,
A_2)\bigr]\bigr\rangle
\end{eqnarray*}
Otherwise we define ${\mathcal G}{\mathcal W}^{(\omega,\mu,
J)}_{\{A_i,g_i,m_i\}}(\kappa_1,\kappa_2;\alpha_1,\cdots,\alpha_{m_1};
\gamma_1,\cdots,\gamma_{m_2})=0$.
 Note that
 \begin{eqnarray*}
 PD(\kappa_1\times\kappa_2)\!\!\!\!&&=(-1)^{(6g_1-6+2m_1-\dim\kappa_1)\dim\kappa_2}
 PD(\kappa_1)\wedge PD(\kappa_2)\\
 &&=(-1)^{\dim\kappa_1\dim\kappa_2}
 PD(\kappa_1)\wedge PD(\kappa_2).
 \end{eqnarray*}
 Using (5.17) we get: \vspace{2mm}

\noindent{\bf Lemma 5.6.}\quad{\it Under the above assumptions it
holds that
\begin{eqnarray*}
&&\quad{\mathcal G}{\mathcal W}^{(\omega,\mu,
J)}_{\{A_i,g_i,m_i\}}(\kappa_1,\kappa_2;\alpha_1,\cdots,\alpha_{m_1};
\gamma_1,\cdots,\gamma_{m_2})\\
&&={\mathcal G}{\mathcal W}^{(\omega,\mu,
J)}_{A_1,g_1,m_1}(\kappa_1;\alpha_1,\cdots,\alpha_{m_1})\cdot{\mathcal
G}{\mathcal W}^{(\omega,\mu,
J)}_{A_2,g_2,m_2}(\kappa_2;\gamma_1,\cdots,\gamma_{m_2}).
\end{eqnarray*}}

 \noindent{\bf Theorem 5.7.}\quad{\it Assume that $\dim H^\ast(M)<\infty$. Let
$\kappa_i\in H_\ast(\overline{\mathcal M}_{g_i, m_i+1},\mbox{\Bb
Q})$, $i=1,2$,  and $\alpha_k\in H^\ast(M,\mbox{\Bb Q})$ for
$k=1,\cdots, m$. Suppose that
 $\alpha_s, \alpha_t\in H_c^\ast(M,\mbox{\Bb Q})$ for some $s\in Q_1$ and $t\in Q_2$.
 Then
\begin{eqnarray*}
\!\!\!{\mathcal G}{\mathcal W}^{(\omega, \mu,
J)}_{A,g,m}(\varphi_{Q\ast}(\kappa_1\times\kappa_2);\alpha_1,\cdots\!,\!\alpha_m)
\!=\epsilon(Q)(-1)^{\deg\kappa_2\sum_{i\in Q_1}\!{\deg\alpha_i}}
\!\!\sum_{A=A_1+A_2}\\
\sum_{k,l} \eta^{kl}\cdot {\mathcal G}{\mathcal W}^{(\omega, \mu,
J)}_{A_1,g_1,m_1+1}(\kappa_1;\{\alpha_i\}_{i\in
Q_2},\beta_k)\cdot{\mathcal G}{\mathcal W}^{(\omega, \mu,
J)}_{A_2,g_2,m_2+1}(\kappa_2; \beta_l, \{\alpha_i\}_{i\in Q_2})
\end{eqnarray*}
Here $\eta^{kl}$ and $\{\beta_i\}$ are  as in Lemma 5.4,
$g=g_1+g_2$, and $\epsilon(Q)$ is the sign of permutation
$Q=Q_1\cup Q_2$ of $\{1,\cdots,m\}$ with $|Q_1|=m_1$ and
$|Q_2|=m_2$.}\vspace{2mm}

By the explanations below Lemma 5.4 this formula is exactly the
ordinary composition law if $M$ is a closed symplectic
manifold.\vspace{2mm}

\noindent{\it Proof of Theorem 5.7.}\quad Without loss of
generality one may assume that $Q_1=\{1,\cdots, m_1\}$,
$Q_2=\{m_1+1,\cdots, m\}$ with $m_2=m-m_1$, and $s=1$, $t=m$.  Let
$A_i\in H_2(M,\mbox{\Bb Z})$, $i=1,2$ and $A=A_1+ A_2$. Consider
the evaluation
$$\Upsilon:\bigcup_{A_1+A_2=A}{\mathcal B}^M_{A_1, g_1, m_1+1}
\times{\mathcal B}^M_{A_2, g_2, m_2+1}\to {\mathcal B}^M_{A, g, m}
$$
 given by $([\Sigma^{(1)},\bar{\bf z}^{(1)}, f^{(1)}],
[\Sigma^{(2)},\bar{\bf z}^{(2)}, f^{(2)}])\mapsto
(f^{(1)}(z_{m_1+1}^{(1)}), f^{(2)}(z_1^{(1)}))$. There is a map
$$\widetilde\varphi_Q: \Upsilon^{-1}(\Delta_M)\to {\mathcal
B}^M_{A, g, m},\;([{\bf f}^{(1)}], [{\bf f}^{(2)}])\mapsto [{\bf
f}^{(1)}\sharp_{f^{(2)}(z_1^{(2)})}{\bf f}^{(2)}],$$
 which is a
lifting of $\varphi_Q$,  such that the following commutative
diagram holds:
\begin{center}\setlength{\unitlength}{1mm}
\begin{picture}(80,30)
\thinlines \put(10,25){$\Upsilon^{-1}(\Delta_M)$}
\put(35,25){\vector (1,0){34}} \put(45,27){$\widetilde\varphi_Q$}
\put(15,22){\vector(0,-1){16}} \put(10,13){$\Pi_r$}
\put(75,25){${\mathcal B}^M_{A, g, m}$}
\put(80,22){\vector(0,-1){16}} \put(71,13){$\Pi^A_{g,m}$}
\put(-5,0){$\overline{\mathcal M}_{g_1,m_1+1}\times
\overline{\mathcal M}_{g_2,m_2+1}$} \put(35,1){\vector(1,0){34}}
\put(45,3){$\varphi_Q$} \put(75,0){$\overline{\mathcal M}_{g,m}$}
\end{picture}\end{center}
Here $\Pi_r$ is the restriction of $\Pi:=\cup_{A_1+
A_2=A}\Pi^{A_1}_{g_1, m_1+1}\times\Pi^{A_2}_{g_2, m_2+1}$ to
$\Upsilon^{-1}(\Delta_M)$.
 Setting
$$S:=\varphi_Q(\overline{\mathcal M}_{g_1,m_1+1}\times
\overline{\mathcal M}_{g_2,m_2+1}),$$
 it is a compact complex
suborbifold of $\overline{\mathcal M}_{g,m}$, and
$\Upsilon^{-1}(\Delta_M)$ may be identified with ${\mathcal B}^{M,
S}_{A, g, m}=(\Pi^A_{g,m})^{-1}(S)$ via the map
$\widetilde\varphi_Q$ naturally.

Take a compact subset $K_0\subset M$ containing the supports of
$\alpha_1$ and $\alpha_m$. If
$$([{\bf f}^{(1)}], [{\bf
f}^{(2)}])\in\Upsilon^{-1}(\Delta_M)\cap\bigl(\overline{\mathcal
M}_{g_1, m_1+1}(M, J, A_1)\times\overline{\mathcal M}_{g_2,
m_2+1}(M, J, A_2)\bigr)$$
  such that $\widetilde\varphi_Q([{\bf f}^{(1)}], [{\bf
f}^{(2)}])$ belongs to  $\overline{\mathcal M}_{g, m}(M,  J, A;
K_3, S)$, then $A=A_1+ A_2$, and it follows from (3.17) that
$[{\bf f}^{(1)}]\in\overline{\mathcal M}_{g_1, m_1+1}(M, J,
A_1;K_4)$ and $[{\bf f}^{(2)}]\in\overline{\mathcal M}_{g_2,
m_2+1}(M, J, A_2; K_4 )$. Conversely, if $([{\bf f}^{(1)}], [{\bf
f}^{(2)}])$ belongs to
 $$
 \bigl(\overline{\mathcal M}_{g_1, m_1+1}(M, J,
A_1;K_3)\times\overline{\mathcal M}_{g_2, m_2+1}(M, J, A_2;
K_3)\bigr)\cap\Upsilon^{-1}(\Delta_M)
$$
 then $\widetilde\varphi_Q([{\bf f}^{(1)}], [{\bf
f}^{(2)}])\in\overline{\mathcal M}_{g, m}(M, J, A; K_3, S)$ for
$A=A_1+ A_2$. Note that there only exist finitely many different
pairs $(A_1, A_2)\in H_2(M,\mbox{\Bb Z})\times H_2(M,\mbox{\Bb
Z})$ such that $A_1+ A_2=A$ and $\overline{\mathcal M}_{g_1,
m_1+1}(M, J, A_1;K_4) \times\overline{\mathcal M}_{g_2, m_2+1}(M,
J, A_2; K_4)\ne\emptyset$. Let  $(A_1^{(i)}, A_2^{(i)})$,
$i=1,\cdots,q$, be these pairs. Set
$$\Omega:=\bigcup^q_{i=1}\Upsilon^{-1}(\Delta_M)\cap(\overline{\mathcal
M}_{g_1, m_1+1}(M, J, A_1^{(i)};K_4)\times\overline{\mathcal
M}_{g_2, m_2+1}(M, J, A_2^{(i)}; K_4)),$$
 then $\overline{\mathcal M}_{g,
m}(M, J, A; K_3,S)\subset\widetilde\varphi_Q(\Omega)\subset
\overline{\mathcal M}_{g, m}(M, J, A; K_4,S)$. As in the proof of
Theorem 5.5 we first choose finitely many points in $\Omega$,
$$[{\bf f}_i]=([{\bf f}_i^{(1)}], [{\bf f}_i^{(2)}]),\,
i=1=k_0+1, \cdots, k_2,\cdots, k_{q-1},\cdots, k_q=k$$
 such that
 $$
 [{\bf f}_i]\in\Upsilon^{-1}(\Delta_M)\cap\bigl(\overline{\mathcal
M}_{g_1, m_1+1}(M, J, A_1^{(t)};K_4)\times\overline{\mathcal
M}_{g_2, m_2+1}(M, J, A_2^{(t)}; K_4)\bigr)
$$
for $i=k_t+1,\cdots,
k_{t+1},\,t=0,\cdots,q-1$. Then we again take finitely many points
$$[{\bf f}_i]\in\bigcup^p_{t=1}\Bigl(\overline{\mathcal M}_{g_1,
m_1+1}(M, J, A_1^{(t)};K_4)\times\overline{\mathcal M}_{g_2,
m_2+1}(M, J, A_2^{(t)}; K_4)\Bigr)\setminus\Omega$$
 for $i=k+1,\cdots,l$. Clearly,
$$f_i^{(1)}(z^{(1)}_{m_1+1})\ne
f_i^{(2)}(z^{(2)}_1),\;i=k+1,\cdots,l.\eqno(5.18)
$$
 If $[{\bf f}_i]=([{\bf f}_i^{(1)}], [{\bf f}_i^{(2)}])\in
\overline{\mathcal M}_{g_1, m_1+1}(M, J,
A_1^{(t)};K_4)\times\overline{\mathcal M}_{g_2, m_2+1}(M, J,
A_2^{(t)}; K_4)$ we can, as before, construct the uniformizers
$\pi_{si}:\widetilde W_{si}\to W_{si}$ of neighborhoods of $[{\bf
f}_i^{(s)}]$ in ${\mathcal B}^M_{A_s^{(t)}, g_s, m_s+1,}$,
$s=1,2$. Then
$$\pi_i:=\pi_{1i}\times\pi_{2i}:\widetilde W_{i}:=
\widetilde W_{1i}\times\widetilde W_{2i}\to W_{i}:=W_{1i}\times
W_{2i}
$$
 is a uniformizer of the neighborhood $W_i$ of $[{\bf f}_i]$ in
${\mathcal B}^M_{A_1^{(t)}, g_1, m_1+1}\times{\mathcal
B}^M_{A_2^{(t)}, g_2, m_2+1}$ with the uniformization group
$\Gamma_i:=\Gamma_{1i}\times\Gamma_{2i}$, where $\Gamma_{si}={\rm
Aut}({\bf f}_i^{(s)})$ for $s=1,2$. By (5.18) we may require
$$W_i\cap\Omega=\emptyset,\,i=k+1,\cdots,l.
$$
Moreover, for each $i=1,\cdots,k$ let us denote
$$\widetilde\Upsilon_i:=\Upsilon\circ\pi_i:\widetilde
W_i:=\widetilde W_{1i}\times\widetilde W_{2i}\to M\times M.
$$
 It is not hard to check that it is transversal to the diagonal
 $\Delta_M$, and that the stratified Banach submanifold
$\widetilde\Upsilon_i^{-1}(\Delta_M)$
 is also invariant under the action of the group $\Gamma_i$.
By increasing $k$ and $l$ if necessary we can use
$\{(\pi_i,\widetilde W_i, W_i,\Gamma_i)\}_{i=1}^l$ to construct a
virtual moduli chain
$$
{\mathcal C}^{\bf t}(K_0)=\sum_{I\in{\mathcal
 N}_l}\frac{1}{|\Gamma_I|}\{\hat\pi_I:{\mathcal M}^{\bf t}_I(K_0)\to{\mathcal
 W}^\ast\}
 $$
  of dimension $2(m+2)+ 2c_1(A)+
2(3-n)(g-2)$ in ${\mathcal W}^\ast\subset {\mathcal
W}=\cup^l_{i=1}W_i\subset\cup_{i=1}^l{\mathcal B}^M_{A_1^{(i)},
g_1, m_1+1}\times{\mathcal B}^M_{A_2^{(i)}, g_2, m_2+1}$
associated with
$$\bigcup^q_{i=1}\bigl(\overline{\mathcal M}_{g_1,
m_1+1}(M, J, A_1^{(i)};K_4)\times\overline{\mathcal M}_{g_2,
m_2+1}(M, J, A_2^{(i)}; K_4)\bigr).
$$

Let $P_i$ be the projections of the product ${\mathcal
B}^M_{A_1,g_1, m_1+1}\times {\mathcal B}^M_{A_2,g_2, m_2+1}$ to
the $i$-th factor, $i=1,2$. Also denote by the evaluations
\begin{eqnarray*}
{\rm ev}_j^{(1)}:{\mathcal B}^M_{A_1,g_1, m_1+1}\to M,\; [{\bf
f}^{(1)}]\mapsto f^{(1)}(z^{(1)}_j),\;j=1,\cdots,m_1+1,\\
{\rm ev}_j^{(2)}:{\mathcal B}^M_{A_2,g_2, m_2+1}\to M,\; [{\bf
f}^{(2)}]\mapsto f^{(2)}(z^{(2)}_j),\;j=1,\cdots,m_2+1.
\end{eqnarray*}

As in (5.12) we can use the virtual moduli chain ${\mathcal
C}^{\bf t}(K_0)$ and the evaluation $\prod^{m_1}_{j=1}{\rm
ev}_j^{(1)}\times({\rm ev}_{m_1+1}^{(1)}\times{\rm
ev}^{(2)}_1)\times\prod^{m_2+1}_{j=2}{\rm ev}_j^{(2)}$ to
construct a rational cycle in $\overline{\mathcal
M}_{g_1,m_1+1}\times\overline{\mathcal M}_{g_2,m_2+1}$,
$${\mathcal C}^{\bf
t}_{(g_1,g_2)}(K_0;\{\bar\alpha_i\}^{m_1}_{i=1},
PD^{II}(\triangle_M),\{\bar\alpha_i\}^m_{i=m_1+1};\omega,\mu, J,
A) \eqno(5.19)
$$
 of dimension $r=2m+2c_1(A)+2(3-n)(g-1)-2-\sum^{m}_{i=1}\deg\alpha_i$, and
prove that its homology class
$$\bigl[{\mathcal C}^{\bf
t}_{(g_1,g_2)}(K_0;\{\alpha_i\}^{m_1}_{i=1},
PD^{II}(\triangle_M),\{\alpha_i\}^m_{i=m_1+1};\omega,\mu, J,
A)\bigr]
$$ in
 $H_r(\overline{\mathcal M}_{g_1,m_1+1}\times\overline{\mathcal M}_{g_2,m_2+1},
 \mbox{\Bb Q})$
 is independent of all related choices. The important is that
  the class of the cycle in (5.19) can also be constructed from
  another direction. In fact we only use
 the virtual moduli chain ${\mathcal C}^{\bf t}(K_0)$ and the
evaluation $\prod^{m_1}_{j=1}{\rm
ev}_j^{(1)}\times\prod^{m_2+1}_{j=2}{\rm ev}_j^{(2)}$ to construct
a rational cycle in $\overline{\mathcal
M}_{g_1,m_1+1}\times\overline{\mathcal M}_{g_2,m_2+1}$,
$${\mathcal C}^{\bf
t}_{(g_1,g_2)}(K_0;\{\bar\alpha_i\}^{m_1}_{i=1},
\{\bar\alpha_i\}^m_{i=m_1+1};\omega,\mu, J, A)
$$
 of dimension $r=2m+2c_1(A)+2(3-n)(g-1)+ 2n-2-\sum^{m}_{i=1}\alpha_i$,
and can also require  that the restriction of $({\rm
ev}_{m_1+1}^{(1)}\times{\rm ev}^{(2)}_1)$  to this cycle (or its
domain) is transversal to $\triangle_M\subset M\times M$.  So we
get a rational cycle in $\overline{\mathcal
M}_{g_1,m_1+1}\times\overline{\mathcal M}_{g_2,m_2+1}$, denoted by
$${\mathcal C}^{\bf
t}_{(g_1,g_2)}(K_0;\{\bar\alpha_i\}^{m_1}_{i=1},
\{\bar\alpha_i\}^m_{i=m_1+1};\omega,\mu, J, A)\cap ({\rm
ev}_{m_1+1}^{(1)}\times{\rm ev}^{(2)}_1)^{-1}(\triangle_M)
$$
  of
dimension $r=2m+2c_1(A)+2(3-n)(g-1)-2-\sum^{m}_{i=1}\alpha_i$.
 It is not hard to check that this cycle gives the same class as one in
 (5.19). Namely, we have
 \begin{eqnarray*}
 \hspace{8mm}&&\quad
 \bigl[{\mathcal C}^{\bf
t}_{(g_1,g_2)}(K_0;\{\alpha_i\}^{m_1}_{i=1},
PD^{II}(\triangle_M),\{\alpha_i\}^m_{i=m_1+1};\omega,\mu, J,
A)\bigr]\hspace{28mm}(5.20)\\
&& =\bigl[{\mathcal C}^{\bf
t}_{(g_1,g_2)}(K_0;\{\alpha_i\}^{m_1}_{i=1},
\{\alpha_i\}^m_{i=m_1+1};\omega,\mu, J, A)\cap ({\rm
ev}_{m_1+1}^{(1)}\times{\rm ev}^{(2)}_1)^{-1}(\triangle_M)\bigr].
\end{eqnarray*}

By the constructions of the virtual moduli chains in $\S3$, it is
not hard using our above arguments to show that for ${\mathcal
N}_k:=\{I\in{\mathcal N}_l\,|\,\max(I)\le k\}$,
$${\mathcal C}^{\bf t}_{\Upsilon^{-1}(\Delta_M)}(K_0):=\sum_{I\in{\mathcal
 N}_k}\frac{1}{|\Gamma_I|}\{\hat\pi_{Ir}:{\mathcal M}^{\bf t}_I(K_0)_r\to{\mathcal
 W}^\ast\cap\Xi^{-1}(\Delta_M)\}
$$
is a virtual moduli chain of dimension
 $2(m+2)+ 2c_1(A)+ 2(3-n)(g-2)-2n$ in ${\mathcal
 W}^\ast\cap\Upsilon^{-1}(\Delta_M)$ associated with
 $\Omega$,
 where ${\mathcal M}^{\bf t}_I(K_0)_r=(\Upsilon\circ\hat\pi_I)^{-1}(\Delta_M)$
 and $\hat\pi_{Ir}$ denotes the restriction of $\hat\pi_I$ to ${\mathcal M}^{\bf
t}_I(K_0)_r$. (This actually implies that $\Upsilon\circ\hat\pi_I:
{\mathcal M}^{\bf t}_I(K_0)\to M\times M$ is transversal to
$\Delta_M$.) As in $\S4.1$ we use the composition
$(\prod^{m_1}_{j=1}{\rm ev}_j^{(1)}\times\prod^{m_2+1}_{j=2}{\rm
ev}_j^{(2)})\circ\hat\pi_{Ir}$ of the evaluation
$\prod^{m_1}_{j=1}{\rm ev}_j^{(1)}\times\prod^{m_2+1}_{j=2}{\rm
ev}_j^{(2)}$ and $\hat\pi_{Ir}$ to construct a rational cycle in
$\overline{\mathcal M}_{g_1,m_1+1}\times\overline{\mathcal
M}_{g_2,m_2+1}$,
$${\mathcal C}^{\bf
t}_{(g_1,g_2)}(K_0;\{\bar\alpha_i\}^{m}_{i=1}, \omega,\mu, J,
A;\Upsilon^{-1}(\triangle_M))
$$
 of dimension $r=2m+2c_1(A)+2(3-n)(g-1)-2-\sum^{m}_{i=1}\alpha_i$, and
prove that its homology class
$$\bigl[{\mathcal C}^{\bf
t}_{(g_1,g_2)}(K_0;\{\bar\alpha_i\}^{m}_{i=1};\omega,\mu, J,
A;\Upsilon^{-1}(\triangle_M))\bigr]
$$
 in $H_r(\overline{\mathcal
M}_{g_1,m_1}\times\overline{\mathcal M}_{g_2,m_2},\mbox{\Bb Q})$
is independent of all related choices. It is not hard to prove
that
\begin{eqnarray*}
&&\quad{\mathcal C}^{\bf
t}_{(g_1,g_2)}(K_0;\{\bar\alpha_i\}^{m}_{i=1},
\omega,\mu, J, A;\Upsilon^{-1}(\triangle_M))\\
&&={\mathcal C}^{\bf
t}_{(g_1,g_2)}(K_0;\{\bar\alpha_i\}^{m_1}_{i=1},
\{\bar\alpha_i\}^m_{i=m_1+1};\omega,\mu, J, A)\cap ({\rm
ev}_{m_1+1}^{(1)}\times{\rm ev}^{(2)}_1)^{-1}(\triangle_M).
\end{eqnarray*}
This and (5.20) give rise to
$$
 \bigl[{\mathcal C}^{\bf
t}_{(g_1,g_2)}(K_0;\{\alpha_i\}^{m_1}_{i=1},
PD^{II}(\triangle_M),\{\alpha_i\}^m_{i=m_1+1};\omega,\mu, J,
A)\bigr]\eqno(5.21)$$
$$\hspace{-22mm} =[{\mathcal C}^{\bf
t}_{(g_1,g_2)}(K_0;\{\alpha_i\}^{m}_{i=1}, \omega,\mu, J,
A;\Upsilon^{-1}(\triangle_M))].
$$

 Moreover using
$\widetilde\varphi_Q$ we can identify
 $\Upsilon^{-1}(\Delta_M)$ with ${\mathcal
 B}^{M,S}_{A,g,m}=(\Pi_{g,m}^A)^{-1}(S)$. (In this case
 $\kappa_1\times\kappa_2\in H_\ast(\overline{\mathcal
 M}_{g_1,m_1+1}\times\overline{\mathcal
 M}_{g_2,m_2+1},\mbox{\Bb Q})$ is identified with $\varphi_{Q\ast}
(\kappa_1\times\kappa_2)\in H_\ast(S,\mbox{\Bb Q})$.)
 So we may consider ${\mathcal C}^{\bf t}_{\Upsilon^{-1}(\Delta_M)}(K_0)$
 as a virtual moduli chain in ${\mathcal  B}^M_{A,g,m}$ associated with
$\overline{\mathcal M}_{g,m}(M,A,J;K_3)$. In this case we have
$$
\varphi_{Q\ast}\bigl[{\mathcal C}^{\bf
t}_{(g_1,g_2)}(K_0;\{\alpha_i\}^{m}_{i=1}, \omega,\mu, J,
A;\Upsilon^{-1}(\triangle_M))\bigr] =[{\mathcal C}_g^{\bf
t}(K_0;\{\alpha_i\}^m_{i=1} ;\omega,\mu,J,A)_S].\eqno(5.22)
$$
Then by Theorem 5.3, (5.22) and (5.21) we get
\begin{eqnarray*}
\hspace{27mm}&&\quad {\mathcal G}{\mathcal W}^{(\omega, \mu,
J)}_{A,g,m}(\varphi_{Q\ast}
(\kappa_1\times\kappa_2);\alpha_1,\cdots,\alpha_m)\hspace{43mm}(5.23)\\
&&={\mathcal G}{\mathcal W}^{(\omega, \mu,
J)}_{A,g,m,S}(\varphi_{Q\ast}
(\kappa_1\times\kappa_2);\alpha_1,\cdots,\alpha_m)\\
&&=\langle PD_S(\varphi_{Q\ast} (\kappa_1\times\kappa_2)),
[{\mathcal C}_g^{\bf
t}(K_0;\{\alpha_i\}^m_{i=1};\omega,\mu,J,A)_S]\rangle\\
&&=\langle PD(\kappa_1\times\kappa_2), [{\mathcal C}^{\bf
t}_{(g_1,g_2)}(K_0;\{\alpha_i\}^{m}_{i=1}, \omega,\mu, J,
A;\Upsilon^{-1}(\triangle_M))]\rangle\\
&&=\langle PD(\kappa_1\times\kappa_2), [{\mathcal C}^{\bf
t}_{(g_1,g_2)}(K_0;\{\alpha_i\}^{m_1}_{i=1},
PD^{II}(\triangle_M),\{\alpha_i\}^m_{i=m_1+1};\omega,\mu, J,
A)]\rangle.
\end{eqnarray*}

Note that the restriction
 of ${\mathcal C}^{\bf t}(K_0)$ to each open
stratified Banach orbifold
$${\mathcal B}^M_{A_1^{(s)},
g_1, m_1+1}\times{\mathcal B}^M_{A_2^{(s)}, g_2, m_2+1}\subset
\bigcup_{t=1}^l{\mathcal B}^M_{A_1^{(t)}, g_1,
m_1+1}\times{\mathcal B}^M_{A_2^{(t)}, g_2, m_2+1}$$
 gives a virtual moduli chain ${\mathcal C}^{\bf t}(K_0)_s$
 of dimension $2(m+2)+ 2c_1(A)+
2(3-n)(g-2)$ associated with
$$\overline{\mathcal M}_{g_1,
m_1+1}(M, J, A_1^{(s)};K_4)\times\overline{\mathcal M}_{g_2,
m_2+1}(M, J, A_2^{(s)}; K_4).$$
 Indeed, let
 \begin{eqnarray*}
&&{\mathcal W}^\ast_s:={\mathcal W}^\ast\cap({\mathcal
B}^M_{A_1^{(s)}, g_1,
m_1+1}\times{\mathcal B}^M_{A_2^{(s)}, g_2, m_2+1}),\\
&& {\mathcal M}^{\bf t}_I(K_0)_s:=(\hat\pi_I)^{-1}({\mathcal
W}^\ast_s)\cap{\mathcal M}^{\bf t}_I(K_0)
\end{eqnarray*}
 and $\hat\pi^s_I$
denote the restriction of $\hat\pi_I$ to ${\mathcal M}^{\bf
t}_I(K_0)_s$. Then we may take
$${\mathcal C}^{\bf t}(K_0)_s=\sum_{I\in{\mathcal
 N}_l}\frac{1}{|\Gamma_I|}\{\hat\pi^s_I:{\mathcal M}^{\bf t}_I(K_0)_s\to{\mathcal
 W}^\ast_s\}.
$$
In constructing rational cycle of (5.19), if we replace ${\mathcal
C}^{\bf t}(K_0)$ by ${\mathcal C}^{\bf t}(K_0)_s$ then we get a
rational cycle in $\overline{\mathcal
M}_{g_1,m_1}\times\overline{\mathcal M}_{g_2,m_2}$,
$${\mathcal C}^{\bf
t}_{(g_1,g_2)}(K_0;\{\bar\alpha_i\}^{m_1}_{i=1},
PD^{II}(\triangle_M);\{\bar\alpha_j\}^{m}_{j=m_1+1};A_1^{(s)},
A_2^{(s)})
$$
 as in (5.15). Note that  ${\mathcal
C}^{\bf t}(K_0)_s$ and ${\mathcal C}^{\bf t}(K_0)_j$ might have
nonempty intersection for $s\ne j$. However their top strata are
disjoint each other, i.e.,
$$T{\mathcal C}^{\bf t}(K_0)_s\cap T{\mathcal C}^{\bf t}(K_0)_j=
\emptyset,\;s\ne j,\,s, j=1,\cdots,l;
$$
 other strata have at least codimension two. It is easily seen
 that these conclusions also hold for ${\mathcal C}^{\bf
t}_{(g_1,g_2)}(K_0;\{\bar\alpha_i\}^{m_1}_{i=1},
PD^{II}(\triangle_M);\{\bar\alpha_j\}^{m}_{j=m_1+1};A_1^{(s)},
A_2^{(s)})$, $s=1,\cdots,q$. It follows that
\begin{eqnarray*}
&&\quad \bigl[{\mathcal C}^{\bf
t}_{(g_1,g_2)}(K_0;\{\alpha_i\}^{m_1}_{i=1},
PD^{II}(\triangle_M),\{\alpha_i\}^m_{i=m_1+1};\omega,\mu, J, A).
\bigr]\\
&&=\sum^q_{s=1}\bigl[{\mathcal C}^{\bf
t}_{(g_1,g_2)}(K_0;\{\alpha_i\}^{m_1}_{i=1},
PD^{II}(\triangle_M);\{\alpha_j\}^{m}_{j=m_1+1};A_1^{(s)},
A_2^{(s)})\bigr].
\end{eqnarray*}
 This and (5.23)  give\vspace{1mm}
$$ \hspace{-40mm}{\mathcal G}{\mathcal W}^{(\omega, \mu,
J)}_{A,g,m}(\varphi_{Q\ast}
(\kappa_1\times\kappa_2);\alpha_1,\cdots,\alpha_m)\eqno(5.24)$$
$$=\sum^q_{s=1}\bigl\langle PD(\kappa_1\times\kappa_2),
\bigl[{\mathcal C}^{\bf
t}_{(g_1,g_2)}(K_0;\{\alpha_i\}^{m_1}_{i=1},
PD^{II}(\triangle_M);\{\alpha_j\}^{m}_{j=m_1+1};A_1^{(s)},
A_2^{(s)})\bigr]\bigr\rangle.
$$
 By Lemma 5.4 we have
\begin{eqnarray*}
&&\quad\bigl[{\mathcal C}^{\bf
t}_{(g_1,g_2)}(K_0;\{\alpha_i\}^{m_1}_{i=1},
PD^{II}(\triangle_M);\{\alpha_j\}^{m}_{j=m_1+1};A_1^{(s)},
A_2^{(s)})\bigr]\\
 &&=\sum_{i,j}c_{ij}[{\mathcal C}^{\bf
t}_{(g_1,g_2)}(K_0;\{\alpha_i\}^{m_1}_{i=1},\beta_i\otimes\beta_j,
\{\alpha_i\}^m_{i=m_1+1};A_1^{(s)}, A_2^{(s)})]\\
&&=\sum_{i,j}c_{ij}[{\mathcal C}^{\bf
t}_{(g_1,g_2)}(K_0;(\otimes^{m_1}_{i=1}\alpha_i)\otimes\beta_i,\beta_j\otimes
(\otimes^m_{i=m_1+1}\alpha_i); A_1^{(s)}, A_2^{(s)})].
\end{eqnarray*}
Moreover, by the definition in (5.17),
$${\mathcal G}{\mathcal
W}^{(\omega,\mu,
J)}_{\{A^{(s)}_i,g_i,m_i\}}(\kappa_1,\kappa_2;\alpha_1,\cdots,\alpha_{m_1},\beta_i;
\beta_j,\alpha_{m_1+1},\cdots,\alpha_m)
$$
is equal to
$$(-1)^{\dim\kappa_1\dim\kappa_2}\langle PD(\kappa_1\times\kappa_2), [{\mathcal C}^{\bf
t}_{(g_1,g_2)}(K_0;(\otimes^{m_1}_{i=1}\alpha_i)\otimes\beta_i,\beta_j\otimes
(\otimes^m_{i=m_1+1}\alpha_i); A_1^{(s)}, A_2^{(s)})]\rangle.$$
From (5.24) and Lemma 5  It follows that
\begin{eqnarray*}
&&\quad {\mathcal G}{\mathcal W}^{(\omega, \mu,
J)}_{A,g,m}(\varphi_{Q\ast}
(\kappa_1\times\kappa_2);\alpha_1,\cdots,\alpha_m)\\
&&=(-1)^\epsilon\sum^q_{s=1}\sum_{i,j}c_{ij}{\mathcal G}{\mathcal
W}^{(\omega,\mu,
J)}_{\{A^{(s)}_i,g_i,m_i\}}(\kappa_1,\kappa_2;\alpha_1,\cdots,\alpha_{m_1},\beta_i;
\beta_j,\alpha_{m_1+1},\cdots,\alpha_m)\\
&&=(-1)^\epsilon\sum^q_{s=1}\sum_{i,j}c_{ij}{\mathcal G}{\mathcal
W}^{(\omega,\mu,
J)}_{A^{(s)}_1,g_1,m_1}(\kappa_1;\alpha_1,\cdots,\alpha_{m_1},\beta_i)\cdot{\mathcal
G}{\mathcal W}^{(\omega,\mu,
J)}_{A^{(s)}_2,g_2,m_2}(\kappa_2;\beta_j,\alpha_{m_1+1},\cdots,\alpha_m).
\end{eqnarray*}
Here $\epsilon=\dim\kappa_1\dim\kappa_2$. Note that the dimension
condition (1.5) implies
$$\sum^{m_1}_{i=1}\deg\alpha_i+ \deg\kappa_1\in 2\mbox{\Bb Z}\quad{\rm and}\quad
\sum^m_{i=m_1+1}\deg\alpha_i+ \kappa_2\in 2\mbox{\Bb Z}.$$
 We get
 $$(-1)^\epsilon=(-1)^{\dim\kappa_1\dim\kappa_2}=
 (-1)^{\deg\kappa_1\deg\kappa_2}=(-1)^{\deg\kappa_2
 \sum^{m_1}_{i=1}\alpha_i}$$
 because $\dim\kappa_j+\deg\kappa_j\in 2\mbox{\Bb Z}$ for $j=1,2$.
This completes the proof of Theorem 5.7. \hfill$\Box$\vspace{2mm}

\subsection{Reduction formulas}

Suppose that $2g+m\ge 3$ and $(g,m)\ne (0,3), (1,1)$.
 For a $m$-pointed stable curve $(\Sigma^\prime,\bar{\bf z}^\prime)$
 of genus $g$ one may obtain a
$(m-1)$-pointed stable curve $(\Sigma,\bar{\bf z})$ of genus $g$
by forgetting the last marked point and contracting the unstable
rational component. This yields a forgetful map ${\mathcal F}_m:
\overline{\mathcal M}_{g,m}\to \overline{\mathcal M}_{g,m-1},
\;[\Sigma,\bar{\bf z}]\mapsto [\Sigma^\prime,\bar{\bf z}^\prime]$.
It is a Lefschetz fibration, whose fibre at $[\Sigma,\bar{\bf z}]$
may be identified with the quotient $\Sigma/{\rm
Aut}(\Sigma,\bar{\bf z})$. But one can still define the
integration along the fibre for it. That is, we have a map
$$({\mathcal F}_m)_\sharp: \Omega^\ast(\overline{\mathcal
M}_{g,m})\to\Omega^{\ast-2}(\overline{\mathcal
M}_{g,m-1})\eqno(5.25)
$$
 that commutes with exterior differentiation $d$ and that
 satisfies the projection formula:
 \begin{description}
 \item[(a)] $({\mathcal F}_m)_\sharp(({\mathcal
 F}_m^\ast\tau)\wedge\omega)=\tau\wedge({\mathcal
 F}_m)_\sharp(\omega)\;\,\forall\tau\in\Omega^\ast(\overline{\mathcal M}_{g,m-1})$ and
  $\omega\in\Omega^\ast(\overline{\mathcal M}_{g,m})$;
 \item[(b)] $\int_{\overline{\mathcal M}_{g,m}}({\mathcal
 F}_m^\ast\tau)\wedge\omega=\int_{\overline{\mathcal
 M}_{g,m-1}}\tau\wedge({\mathcal F}_m)_\sharp(\omega)$ for any
 $\omega\in\Omega^q(\overline{\mathcal  M}_{g,m})$ and
 $\tau\in\Omega^{6g-6+ 2m-q}(\overline{\mathcal
 M}_{g,m-1})$, where the integrations are over the orbifolds.
We still denote by $({\mathcal
F}_m)_\sharp:\;H^\ast(\overline{\mathcal M}_{g,m},\mbox{\Bb Q}
)\to H^{\ast-2}(\overline{\mathcal M}_{g,m-1},\mbox{\Bb Q})$ the
induced map.   ${\mathcal F}_m$ also induces a ``shriek'' map
$$({\mathcal F}_m)_!: H_\ast(\overline{\mathcal M}_{g,m-1};\mbox{\Bb Q})\to
H_{\ast+2}(\overline{\mathcal M}_{g,m};\mbox{\Bb Q})\eqno(5.26)
$$
 given by $\kappa\mapsto PD^{-1}_m\circ{\mathcal
F}_m^\ast\circ PD_{m-1}(\kappa)$, where $PD_m:
H_\ast(\overline{\mathcal M}_{g,m};\mbox{\Bb Q})\to
H^{6g-6+2m-\ast}(\overline{\mathcal M}_{g,m};\mbox{\Bb Q})$ is the
Poincar\'e duality.
\end{description}
Clearly, from (5.26) and (b) we get
$${\mathcal F}_m^\ast\circ PD_{m-1}(\kappa)=PD_m(({\mathcal
F}_m)_!(\kappa))\quad{\rm and}\quad ({\mathcal F}_m)_\sharp\circ
PD_m=PD_{m-1}\circ({\mathcal F}_m)_\ast\eqno(5.27)
$$
 respectively. So for any $\alpha\in
H^\ast(\overline{\mathcal M}_{g,m};\mbox{\Bb Q})$ and $[c]\in
H_{\ast-2}(\overline{\mathcal M}_{g,m-1};\mbox{\Bb Q})$ we have
$$\langle ({\mathcal F}_m)_\sharp(\alpha), [c]\rangle_{m-1}=
\langle \alpha, ({\mathcal F}_m)_{!}([c])\rangle_{m}.\eqno(5.28)
$$
 The reduction formulas are the following two
 theorems.\vspace{2mm}

\noindent{\bf Theorem 5.8.}\quad{\it If $(g,m)\ne (0,3),(1,1)$,
then for any $\kappa\in H_\ast(\overline{\mathcal
M}_{g,m-1};\mbox{\Bb Q})$, $\alpha_1\in H_c^\ast(M;\mbox{\Bb Q})$,
$\alpha_2, \cdots,\alpha_m\in H^\ast(M;\mbox{\Bb Q})$ with
$\deg\alpha_m=2$ it holds that
$${\mathcal G}{\mathcal W}^{(\omega, \mu,
J)}_{A,g,m}(({\mathcal
F}_m)_!(\kappa);\alpha_1,\cdots,\alpha_m)=\alpha_m(A)\cdot
{\mathcal G}{\mathcal W}^{(\omega, \mu,
J)}_{A,g,m-1}(\kappa;\alpha_1,\cdots,\alpha_{m-1}).$$}

\noindent{\bf Theorem 5.9.}\quad{\it If $(g,m)\ne (0,3),(1,1)$,
$\kappa\in H_\ast(\overline{\mathcal M}_{g,m};\mbox{\Bb Q})$ and
$\alpha_1,\cdots,\alpha_{m-1}$ are as in Theorem 5.8 then it also
holds that
$${\mathcal G}{\mathcal W}^{(\omega, \mu,
J)}_{A,g,m}(\kappa;\alpha_1,\cdots,\alpha_{m-1},{\bf 1})={\mathcal
G}{\mathcal W}^{(\omega, \mu,J)}_{A,g,m-1}(({\mathcal
F}_m)_\ast(\kappa);\alpha_1,\cdots,\alpha_{m-1}).$$
 Here ${\bf 1}\in H^0(M,\mbox{\Bb Q})$ is the identity.}\vspace{2mm}

In order to prove these two theorems we need a similar map to
${\mathcal F}_m$. A stable $L^p_k$-map $(f;\Sigma,\bar{\bf z})$ is
called {\it strong} if for each component $\Sigma_s$ of $\Sigma$
satisfying $(f\circ\pi_{\Sigma_s})_\ast([\widetilde\Sigma_s])=0\in
H_2(M,\mbox{\Bb Z})$ we have: (i) $m_s+ 2g_s\ge 3$, and (ii)
$f|_{\Sigma_s}$ is constant. Correspondingly we call the
isomorphism class $[f;\Sigma,\bar{\bf z}]\in{\mathcal
B}^M_{A,g,m}$ strong. Let us denote by
$$({\mathcal
B}^M_{A,g,m})_s\quad ({\rm resp.} ({\mathcal B}^M_{A,g,m})_0)
$$
  the subset consisting of the strong stable
$[f;\Sigma,\bar{\bf z}]\in{\mathcal B}^M_{A,g,m}$ (resp., the
stable $[f;\Sigma,\bar{\bf z}]\in{\mathcal B}^M_{A,g,m}$ which is
still stable after removing the last marked point.) Then
$\overline{\mathcal M}_{g,m}(M, A, J)\subset({\mathcal
B}^M_{A,g,m})_s$, and $({\mathcal B}^M_{A,g,m})_0$ is an open
subset of ${\mathcal B}^M_{A,g,m}$. Carefully checking the proof
of Lemma 23.2 of [FuO] it is also not hard to see that  the domain
$\Sigma$ of each element $[f;\Sigma,\bar{\bf z}]\in {\mathcal
B}^M_{A,g,m}\setminus({\mathcal B}^M_{A,g,m})_0$ has at least two
components. Thus ${\mathcal B}^M_{A,g,m}\setminus({\mathcal
B}^M_{A,g,m})_0$ has at least codimension two in ${\mathcal
B}^M_{A,g,m}$.  The proof of Lemma 23.2 of [FuO] implies the
following result.\vspace{2mm}

\noindent{\bf Lemma 5.10}\quad{\it  Let $2g+m\ge 3$, if $(g,m)\ne
(0,3), (1,1)$ or $A\ne 0$ then the obvious  map
$\widetilde{\mathcal F}_m: ({\mathcal B}^M_{A, g,
m})_0\to{\mathcal B}^M_{A,g,m-1}$ forgetting the last marked point
is surjective and may be extended to a map
$$\widetilde{\mathcal F}_m:({\mathcal
B}^M_{A,g,m})_{s0}:=({\mathcal B}^M_{A,g,m})_s\cup({\mathcal
B}^M_{A,g,m})_0\to {\mathcal B}^M_{A,g,m-1}
$$
  satisfying
 $$\Pi_{g,m-1}\circ\widetilde{\mathcal F}_m={\mathcal
 F}_m\circ (\Pi_{g,m}|_{({\mathcal B}^M_{A,g,m})_{s0}}).\eqno(5.29)
$$
Moreover, the map  may be regarded as a universal family, i.e.,
the fibre
 $$(\widetilde{\mathcal F}_m)^{-1}([{\bf f}])=\Sigma/{\rm Aut}({\bf
 f}),\;\forall\;[{\bf f}]=[f, \Sigma,\bar{\bf z}]\in {\mathcal
B}^M_{A,g,m-1}.$$ }\vspace{2mm}

Indeed,  if $[f;\Sigma,\bar{\bf z}]\in({\mathcal
B}^M_{A,g,m})_s\setminus({\mathcal B}^M_{A,g,m})_0$ and $\Sigma_0$
is the component of $\Sigma$ containing the $m$-th marked point
$z_m$ then $[f;\Sigma,\bar{\bf z}]\in\overline{\mathcal
M}_{g,m}(M, A, J)$ will become unstable after removing $z_m$. As
in [FuO], $f|_{\Sigma_0}$ has the homology class zero in
$H_2(M,\mbox{\Bb Z})$ and the genus $g_0$ of $\Sigma_0$ is zero.
Consequently either $\Sigma_0$ has one singular point and one
marked point $z_k$ with $k\ne m$ or $\Sigma_0$ has two singular
points and unique marked point $z_m$. In the first case one
removes $\Sigma_0$ from $\Sigma$ and replaces $z_k$ with the point
where $\Sigma_0$ was attached. In the second case one removes
$\Sigma_0$ from $\Sigma$ and glues it at the two points where
$\Sigma_0$ was attached. For each case one can get  a stable curve
$[\Sigma^\prime,\bar{\bf z}^\prime]$ in $\overline{\mathcal
M}_{g,m}$, and $f$ naturally induces a map $f^\prime$ on
$\Sigma^\prime$ such that $(f^\prime;\Sigma^\prime,\bar{\bf
z}^\prime)$ is a stable $L^p_k$-map. Then $\widetilde{\mathcal
F}_m([f;\Sigma,\bar{\bf z}])$ is defined as
$[f^\prime;\Sigma^\prime,\bar{\bf z}^\prime]$. One easily sees
that it satisfies Lemma 5.10. Note here it is very key that $f$ is
{\it constant} on $\Sigma_0$. In general  we do not know whether
or not the map $\widetilde{\mathcal F}_m$ may be extended to the
space ${\mathcal B}^M_{A,g,m}$. For the proof of another claim the
reader may refer to the arguments below the proof of Lemma 23.2 of
[FuO].  In particular, the restriction to
 $\overline{\mathcal M}_{g,m}(M, A, J)$ of the map
$\widetilde{\mathcal F}_m$ may be regarded as a universal family.
\vspace{2mm}

To avoid excessive transversality arguments we shall use the
equivalent definition (4.14) in the proofs of the following two
theorems.\vspace{2mm}

 \noindent{\it Proof of Theorem
5.8.}\quad Fix a compact subset $K_0\subset M$ containing ${\rm
supp}(\wedge^m_{i=1}\alpha_i)$. As before we take $[{\bf
f}_i]\in\overline{\mathcal M}_{g,m}(M, A, J;K_3)$ and sections
$\tilde s_{ij}:\widetilde W_i\to\widetilde E_i$, $i=1,\cdots, n_3$
and $j=1,\cdots,q_i$ to  construct a family of cobordant virtual
moduli chains
$${\mathcal C}^{\bf t}(K_0):=\sum_{I\in{\mathcal
 N}}\frac{1}{|\Gamma_I|}\{\hat\pi_I:{\mathcal M}^{\bf t}_I(K_0)\to{\mathcal
 W}^\ast\}\quad\forall{\bf t}\in{\bf B}_{\varepsilon}^{res}(\mbox{\Bb R}^{q})
 $$
of dimension $2m+ 2c_1(M)(A)+ 2(3-n)(g-1)$ in ${\mathcal
 W}^\ast=\cup_{I\in{\cal N}}V_I^\ast\subset{\mathcal B}^M_{A,g,m}$
 associated with $\overline{\mathcal M}_{g,m}(M, A, J;K_3)$. Here $q=q_1+\cdots+ q_{n_3}$.

Denote by ${\bf f}_i^\prime=\widetilde{\mathcal F}_m({\bf f}_i)$
and $[{\bf f}_i^\prime]=\widetilde{\mathcal F}_m([{\bf f}_i])$,
$i=1,\cdots,n_3$. As above we take sections $\tilde
s'_{ij}:\widetilde W'_i\to\widetilde E'_i$, $i=1,\cdots, n_3$ and
$j=1,\cdots,q_i$ to  construct a family of cobordant virtual
moduli chains
$${\mathcal C}^{\bf t}(K_0)':=\sum_{I\in{\mathcal
 N}}\frac{1}{|\Gamma_I'|}\{\hat\pi'_I:{\mathcal M}^{\bf t}_I(K_0)'\to{\mathcal
 W}{\prime\ast}\}\quad\forall{\bf t}\in{\bf B}_{\varepsilon}^{res}(\mbox{\Bb R}^{q})\eqno(5.30)
 $$
of dimension $2(m-1)+ 2c_1(M)(A)+ 2(3-n)(g-1)$ in  ${\mathcal
 W}^{\prime\ast}=\cup_{I\in{\cal N}}V_I^{\prime\ast}\subset{\mathcal B}^M_{A,g,m-1}$
 associated with $\overline{\mathcal M}_{g,m-1}(M, A, J;K_3)$.
This is possible by increasing $n_3$ if necessary.

We now  show that the sections $\tilde s_{ij}$ and $\tilde
s_{ij}'$ can be suitably chosen so that ${\mathcal C}^{\bf
t}(K_0)$ and ${\mathcal C}^{\bf t}(K_0)'$ may be related. Note
that $\tilde s_{ij}'=\beta_i'\cdot\tilde\nu_{ij}'$ and that
$\nu_{ij}'$ has the support disjoint from marked or singular
points on domain $\Sigma'$ of ${\bf f}_i'$. So $\nu_{ij}'$
 may be naturally regarded as an element of
$L^p_{k-1}(\wedge^{0,1}(f^\ast TM))$, $\nu_{ij}=\widetilde{\cal
F}_m^\ast\nu_{ij}'$. For given  $\widetilde W_i^\prime$ we
construct a $\Gamma_i$-invariant $\widetilde W_i$ such that
$$\widetilde W_i^{s0}:=\pi_i^{-1}(({\mathcal
B}^M_{A,g,m})_{s0})\cap\widetilde W_i=\widetilde{\mathcal
F}_m^{-1}(\widetilde W_i^\prime)\eqno(5.31)
$$
 for the projection
$\pi_i:\widetilde W_i\to W_i\subset{\mathcal B}^M_{A,g,m}$. Then
$\widetilde E_i^{s0}:=\widetilde{\mathcal F}_m^\ast\widetilde
E_i^\prime$ is equal to the restriction of the bundle $\widetilde
E_i$ to $\widetilde W_i^{s0}$. We may take a $\Gamma_i$-invariant
cut-off function $\beta_i$ on $\widetilde W_i^{s0}$ such that
$\beta_i=\widetilde{\mathcal F}_m^\ast\beta_i^\prime$ on
$\widetilde W_i^{s0}$. Set $\tilde
s_{ij}=\beta_i\cdot\tilde\nu_{ij}$ then the virtual moduli chains
in (5.28) and (5.29) satisfy: For each ${\bf t}\in{\bf
B}_{\varepsilon}^{res}(\mbox{\Bb R}^{q})$ the set
$$
{\mathcal M}^{\bf t}_I(K_0)_{s0}:={\mathcal M}^{\bf
t}_I(K_0)\cap(\hat\pi_I)^{-1}(({\mathcal B}^M_{A,g,m})_{s0})
$$
has a complementary of at least codimension two in ${\mathcal
M}^{\bf t}_I(K_0)$, and
$$\Pi:=\widetilde{\mathcal F}_m|_{\hat\pi_I({\mathcal M}^{\bf t}_I(K_0)_{s0})}:
\hat\pi_I({\mathcal M}^{\bf
t}_I(K_0)_{s0})\to\hat\pi_I^\prime({\mathcal M}^{\bf
t}_I(K_0)^\prime)\eqno(5.32)
$$
 is also a Lefschetz fibration.
 So it follows from  (4.14) and (5.27) that
 \begin{eqnarray*}
\hspace{26mm}&&\quad {\mathcal G}{\mathcal W}^{(\omega,\mu,
J)}_{A,g,m}(({\mathcal
F}_m)_!(\kappa);\alpha_1,\cdots,\alpha_m)\hspace{50mm}(5.33)\\
&&=\sum_{I\in{\mathcal N}}\frac{1}{|\Gamma_I|}\int_{{\mathcal
M}^{\bf t}_I(K_0)}({\Pi}_{g,m}\circ\hat\pi_I)^\ast({\mathcal
F}_m^\ast\kappa^\ast)\wedge\bigl(\wedge^m_{i=1}({\rm
ev}_i\circ\hat\pi_I)^\ast\alpha_i^\ast\bigr)\\
&&=\sum_{I\in{\mathcal N}}\frac{1}{|\Gamma_I|}\int_{{\mathcal
M}^{\bf t}_I(K_0)_{s0}}({\Pi}_{g,m}\circ\hat\pi_I)^\ast({\mathcal
F}_m^\ast\kappa^\ast)\wedge\bigl(\wedge^m_{i=1}({\rm
ev}_i\circ\hat\pi_I)^\ast\alpha_i^\ast\bigr)\\
&&=\sum_{I\in{\mathcal N}}\frac{1}{|\Gamma_I|}\int_{{\mathcal
M}^{\bf t}_I(K_0)_{s0}}\hat\pi_I^\ast\Bigl(({\mathcal
F}_m\circ{\Pi}_{g,m})^\ast\kappa^\ast\wedge\bigl(\wedge^m_{i=1}({\rm
ev}_i)^\ast\alpha_i^\ast\bigr)\Bigr)\\
&&=\sum_{I\in{\mathcal N}}\int_{\hat\pi_I({\mathcal M}^{\bf
t}_I(K_0)_{s0})}({\mathcal
F}_m\circ{\Pi}_{g,m})^\ast\kappa^\ast\wedge\bigl(\wedge^m_{i=1}({\rm
ev}_i)^\ast\alpha_i^\ast\bigr).
\end{eqnarray*}
Here the last equality comes from the definition of the
integration over the orbifold $\hat\pi_I({\mathcal M}^{\bf
t}_I(K_0)_{s0})={\mathcal M}^{\bf t}_I(K_0)_{s0}/\Gamma_I$.

  Note that ${\rm ev}_i={\rm
ev}_i^\prime\circ\widetilde{\mathcal F}_m$, $1\le i\le m-1$. By
(5.29) and the projection formula we get:
\begin{eqnarray*}
\hspace{17mm} &&\quad \int_{\hat\pi_I({\mathcal M}^{\bf
t}_I(K_0)_{s0})}({\mathcal
F}_m\circ{\Pi}_{g,m})^\ast\kappa^\ast\wedge\bigl(\wedge^m_{i=1}({\rm
ev}_i)^\ast\alpha_i^\ast\bigr)\hspace{37mm}(5.34)\\
&&\!\!\!\!\!\!=\int_{\hat\pi_I({\mathcal M}^{\bf
t}_I(K_0)_{s0})}({\Pi}_{g,m-1}\circ\Pi)^\ast\kappa^\ast\wedge\bigl(\wedge^{m-1}_{i=1}({\rm
ev}_i^\prime\circ\Pi)^\ast\alpha_i^\ast\bigr)\wedge{\rm
ev}_m^\ast\alpha_m^\ast\\
&&\!\!\!\!\!\!=\int_{\hat\pi_I({\mathcal M}^{\bf
t}_I(K_0)_{s0})}\Pi^\ast\Bigl(({\Pi}_{g,m-1})^\ast
\kappa^\ast\wedge\bigl(\wedge^{m-1}_{i=1}({\rm
ev}_i^\prime)^\ast\alpha_i^\ast\bigr)\Bigl)\wedge{\rm
ev}_m^\ast\alpha_m^\ast\\
&& \!\!\!\!\!\!=\int_{\hat\pi_I^\prime({\mathcal M}^{\bf
t}_I(K_0)^\prime)}({\Pi}_{g,m-1})^\ast\kappa^\ast\wedge\bigl(\wedge^{m-1}_{i=1}({\rm
ev}_i^\prime)^\ast\alpha_i^\ast\bigr)\wedge\Pi_\sharp({\rm
ev}_m^\ast\alpha_m^\ast)\\
&&\!\!\!\!\!\!=\alpha_m(A)\int_{\hat\pi_I^\prime({\mathcal M}^{\bf
t}_I(K_0)^\prime)}({\Pi}_{g,m-1})^\ast\kappa^\ast\wedge\bigl(\wedge^{m-1}_{i=1}({\rm
ev}_i^\prime)^\ast\alpha_i^\ast\bigr).
\end{eqnarray*}
Here the third equality comes from (5.32), and the last step uses
the fact that $\Pi_\sharp({\rm
ev}_m^\ast\alpha_m^\ast)=(\widetilde{\mathcal F}_m)_\sharp({\rm
ev}_m^\ast\alpha_m^\ast)=\alpha_m(A)$. Now (5.33) and (5.34) yield
\begin{eqnarray*}
&&\quad {\mathcal G}{\mathcal W}^{(\omega,\mu,
J)}_{A,g,m}(({\mathcal
F}_m)_!(\kappa);\alpha_1,\cdots,\alpha_m)\\
&&=\sum_{I\in{\mathcal N}}\int_{\hat\pi_I({\mathcal M}^{\bf
t}_I(K_0)_{s0})}({\mathcal
F}_m\circ{\Pi}^A_{g,m})^\ast\kappa^\ast\wedge\bigl(\wedge^m_{i=1}({\rm
ev}_i)^\ast\alpha_i^\ast\bigr)\\
&&=\alpha_m(A)\sum_{I\in{\mathcal N}}
\int_{\hat\pi_I^\prime({\mathcal M}^{\bf
t}_I(K_0)^\prime)}({\Pi}^A_{g,m-1})^\ast\kappa^\ast\wedge\bigl(\wedge^{m-1}_{i=1}({\rm
ev}_i^\prime)^\ast\alpha_i^\ast\bigr)\\
&&=\alpha_m(A)\sum_{I\in{\mathcal
N}}\frac{1}{|\Gamma_I^\prime|}\int_{{\mathcal M}^{\bf
t}_I(K_0)^\prime}(\Pi^A_{g,m-1}\circ\hat\pi_I^\prime)^\ast\kappa^\ast\wedge\bigl(\wedge^{m-1}_{i=1}
({\rm ev}_i^\prime\circ\hat\pi_I^\prime)^\ast\alpha_i^\ast\bigr)\\
&&=\alpha_m(A)\cdot {\mathcal G}{\mathcal W}^{(\omega,\mu,
J)}_{A,g,m-1}(\kappa;\alpha_1,\cdots,\alpha_m).
\end{eqnarray*}
Theorem 5.8 is proved.\hfill$\Box$\vspace{2mm}

 \noindent{\it Proof of Theorem
5.9.}\quad Note that (5.29) implies
$\Pi_\sharp\circ(\Pi_{g,m})^\ast=(\Pi_{g,m-1})^\ast\circ({\mathcal
F}_m)_\sharp$. As in the proof of Theorem 5.8 we have
\begin{eqnarray*}
&&\quad{\mathcal G}{\mathcal W}^{(\omega,\mu, J)}_{A,g,m}(\kappa;
\alpha_1,\cdots,\alpha_{m-1},{\bf 1})\\
&&=\sum_{I\in{\mathcal N}}\frac{1}{|\Gamma_I|}
 \int_{{\mathcal M}^{\bf t}_I(K_0)_{s0}}({\Pi}_{g,m}\circ\hat\pi_I)^\ast\kappa^\ast
 \wedge\bigl(\wedge^m_{i=1}({\rm ev}_i\circ\hat\pi_I)^\ast\alpha_i^\ast\bigr)\\
&&=\sum_{I\in{\mathcal N}}\frac{1}{|\Gamma_I|}
 \int_{{\mathcal M}^{\bf t}_I(K_0)_{s0}}
 \hat\pi_I^\ast\Bigl(({\Pi}_{g,m})^\ast\kappa^\ast\wedge\bigl(\wedge^m_{i=1}({\rm
ev}_i)^\ast\alpha_i^\ast\bigr)\Bigr)\\
&&=\sum_{I\in{\mathcal N}}\int_{\hat\pi_I({\mathcal M}^{\bf
t}_I(K_0)_{s0})}({\Pi}_{g,m})^\ast\kappa^\ast\wedge\bigl(\wedge^m_{i=1}({\rm
ev}_i)^\ast\alpha_i^\ast\bigr)\\
&&=\sum_{I\in{\mathcal N}}\int_{\hat\pi_I({\mathcal M}^{\bf
t}_I(K_0)_{s0})}({\Pi}_{g,m})^\ast\kappa^\ast\wedge\bigl(\wedge^{m-1}_{i=1}({\rm
ev}_i^\prime\circ\Pi)^\ast\alpha_i^\ast\bigr)\wedge{\rm
ev}_m^\ast 1\\
&& =\sum_{I\in{\mathcal N}}\int_{\hat\pi_I^\prime({\mathcal
M}^{\bf
t}_I(K_0)^\prime)}\Pi_\sharp(({\Pi}_{g,m})^\ast\kappa^\ast)\wedge\bigl(\wedge^{m-1}_{i=1}({\rm
ev}_i^\prime)^\ast\alpha_i^\ast\bigr)\\
&&=\sum_{I\in{\mathcal N}}\int_{\hat\pi_I^\prime({\mathcal M}^{\bf
t}_I(K_0)^\prime)}({\Pi}_{g,m-1})^\ast(({\mathcal
F}_m)_\sharp\kappa^\ast)\wedge\bigl(\wedge^{m-1}_{i=1}({\rm
ev}_i^\prime)^\ast\alpha_i^\ast\bigr)\\
&&=\sum_{I\in{\mathcal
N}}\frac{1}{|\Gamma_I^\prime|}\int_{{\mathcal M}^{\bf
t}_I(K_0)^\prime}({\Pi}_{g,m-1}\circ\hat\pi_I^\prime)^\ast(({\mathcal
F}_m)_\sharp\kappa^\ast)\wedge\bigl(\wedge^{m-1}_{i=1}({\rm
ev}_i^\prime\circ\hat\pi_I^\prime)^\ast\alpha_i^\ast\bigr)\\
&&={\mathcal G}{\mathcal W}^{(\omega,\mu, J)}_{A,g,m-1}(({\mathcal
F}_m)_\ast(\kappa);\alpha_1,\cdots,\alpha_{m-1}),
\end{eqnarray*}
where the sixth equality is because
$\Pi_\sharp(({\Pi}_{g,m})^\ast\kappa^\ast)=
({\Pi}_{g,m-1})^\ast(({\mathcal F}_m)_\sharp\kappa^\ast)$.
\hfill$\Box$\vspace{2mm}

Using (5.27) and (5.28) it easily follows that the reduction
formulas in Theorems 5.8 and 5.9 can be written as the following
versions
$$
({\cal F}_m)_\ast([{\mathcal C}_g^{\bf t}(K_0;\{\alpha_i\}^m_{i=1}
;\omega,\mu,J,A)])=\alpha_m(A)\cdot[{\mathcal C}_g^{{\bf t}'}
(K_0;\{\alpha_i\}^{m-1}_{i=1} ;\omega,\mu,J,A)],\eqno(5.35)
$$
$$
({\cal F}_m)_!([{\mathcal C}_g^{{\bf
t}'}(K_0;\{\alpha_i\}^{m-1}_{i=1} ;\omega,\mu,J,A)])=[{\mathcal
C}_g^{\bf t}(K_0;\{\alpha_i\}^{m-1}_{i=1},{\bf 1}
;\omega,\mu,J,A)]\eqno(5.36)
$$
respectively.

\section{Quantum Cohomology, WDVV Equation and String Equation}

In this section we shall apply the invariants constructed in the
previous sections to  the constructions of quantum cohomology,
solutions of WDVV equation and Witten's string equation. We only
give main points because they are imitation for the arguments in
[McSa1, RT1, RT2, T1, W1, W2, CoKa] etc.

\subsection{Generalized string equation and dilation equation}

Let $\overline{\mathcal U}_{g,m}$ be the universal curve over
$\overline{\mathcal M}_{g,m}$. The $i$-th marked point $z_i$
yields a section $\tilde z_i$ of the fibration $\overline{\mathcal
U}_{g,m}\to\overline{\mathcal M}_{g,m}$. Denote by ${\mathcal
K}_{{\mathcal U}|{\mathcal M}}$ the cotangent bundle to fibers of
this fibration, and by ${\mathcal L}_i=\tilde z_i({\mathcal
K}_{{\mathcal U}|{\mathcal M}})$. For nonnegative integers
$d_1,\cdots,d_m$ we also denote by $\kappa_{d_1,\cdots,d_m}$ the
Poincar\'e dual of $c_1({\mathcal L}_1)^{d_1}\cup\cdots\cup
c_1({\mathcal L}_m)^{d_m}$.

Now we need to define (or make conventions)
$$\left\{\begin{array}{lll}
{\mathcal G}{\mathcal W}^{(\omega, \mu,
J)}_{A,g,m}(\kappa_{d_1,\cdots,d_m};\underbrace{{\bf 1},\cdots,
{\bf 1}}_{m\,{\rm times}})=0\\
 {\mathcal G}{\mathcal W}^{(\omega,
\mu, J)}_{A,0,3}([\overline{\mathcal M}_{0,3}]; {\bf 1}, {\bf 1},
{\bf 1})=0\\
 {\mathcal G}{\mathcal W}^{(\omega, \mu,
J)}_{A,1,1}([\overline{\mathcal M}_{1,1}]; {\bf 1})=0\\
{\mathcal G}{\mathcal W}^{(\omega, \mu, J)}_{0,1,1}(\kappa_1;
{\bf 1})=0\\
{\mathcal G}{\mathcal W}^{(\omega, \mu, J)}_{0,1,1}([pt]; {\bf
1})=\chi(M)
\end{array}\right.\eqno(6.1)$$
because these cannot be included in the category of our
definition. Here $\chi(M)$ is the Euler characteristic of $M$.
These definitions or conventions are determined by the conclusions
in the case that $M$ is a closed symplectic manifold.

 Since $M$ is noncompact, $H^0_c(M,\mbox{\Bb
Q})=0$. So the cohomology ring $H^\ast_c(M,\mbox{\Bb Q})$ has no
the unit element. But $H^\ast(M,\mbox{\Bb Q})$ has a unit element
${\bf 1}\in H^0(M,\mbox{\Bb Q})$, which is Poincar\'e dual to the
fundamental class $[M]\in H_{2n}^{\rm II}(M,\mbox{\Bb Q})$. Let us
define
$$\widetilde H^\ast_c(M,\mbox{\Bb Q})=H^\ast_c(M,\mbox{\Bb
Q})\oplus{\bf 1}\mbox{\Bb Q}\eqno(6.2)$$
 so that $\widetilde H^\ast_c(M,\mbox{\Bb Q})$ becomes a ring with the unit.
Note that $H^i_c(M,\mbox{\Bb R})$ for every $i$ is at most
countably generated.  There exist at most countable linearly
independent elements $\{\gamma_i\}_{2\le i<N}$ in $H^{\rm even}
_c(M,\mbox{\Bb Q})$ such that $H^{\rm even} _c(M,\mbox{\Bb
Q})={\rm span}(\{\gamma_i\}_{2\le i<N})$. Here $N$ is a natural
number or $+\infty$. Set  $\gamma_1={\bf 1}$. We get that
 ${\rm span}(\{\gamma_i\}_{1\le
i<N})=\widetilde H^{\rm even}_c(M,\mbox{\Bb Q})$. For $1\le a,b<N$
let
$$\zeta_{ab}=\left\{\begin{array}{lll}
0\;&&{\rm if}\;\deg\gamma_a+\deg\gamma_b\ne 2n,\\
\int_M\gamma^\ast_a\wedge\gamma_b^\ast&&{\rm if}\;
\deg\gamma_a+\deg\gamma_b=2n.\end{array}\right.\eqno(6.3)$$
 Then by the second conclusion of Theorem 4.2 we have
 $${\mathcal G}{\mathcal W}^{(\omega, \mu,
J)}_{A,0,3}([\overline{\mathcal M}_{0,3}]; {\bf 1},
\gamma_a,\gamma_b)=\left\{\begin{array}{lll}
0\;&&{\rm if}\; A\ne 0\\
\zeta_{ab}&&{\rm if}\; A=0.\end{array}\right.\eqno(6.4)$$

 Given an integer $g\ge 0$ and $A\in H_2(M)$,  let one class in
  $\{\alpha_i\}_{1\le i\le m}\subset H_c^\ast(M,\mbox{\Bb
Q})\cup H^\ast(M,\mbox{\Bb Q})$ belong to $H_c^\ast(M,\mbox{\Bb
Q})$. We call the invariant
$$\langle\tau_{d_1,\alpha_1}\tau_{d_2,\alpha_2}\cdots\tau_{d_m,\alpha_m}\rangle_{g, A}=
{\mathcal G}{\mathcal W}^{(\omega, \mu,
J)}_{A,g,m}(\kappa_{d_1,\cdots,d_m};\alpha_1,\cdots,\alpha_m)\eqno(6.5)$$
a {\it gravitational correlator}. The  $m$-point genus-$g$
correlators are  defined by
$$
\langle\tau_{d_1,\alpha_1}\tau_{d_2,\alpha_2}\cdots\tau_{d_m,\alpha_m}\rangle_g(q)
=\sum_{A\in H_2(M)}{\mathcal G}{\mathcal W}^{(\omega, \mu,
J)}_{A,g,m}(\kappa_{d_1,\cdots,d_m};\alpha_1,\cdots,\alpha_m)q^A.
\eqno(6.6)$$
 where $q$ is an element of Novikov ring as before.

With the given base $\{\gamma_i\}_{1\le i<N}$ above  and  formal
variables $t^a_r$,
 $1\le a< N$, $r=0, 1, 2,\cdots$, all genus-$g$
correlators can be assembled into a generating function, called
{\it free energy function}([W2]), as follows:
$$F^M_g(t^a_r;q)=\sum_{n_{r,a}}\prod_{r,a}\frac{(t^a_r)^{n_{r,a}}}{n_{r,a}!}\Bigl\langle
 \prod_{r,a}\tau^{n_{r,a}}_{r,\gamma_a}\Bigr\rangle_g(q),\eqno(6.7)$$
where $n_{r,a}$ are arbitrary collections of nonnegative integers,
almost all zero, labelled by $r,a$. Witten's generating function (
[W2]) is the infinite sum
$$F^M(t^a_r;q)=\sum_{g\ge 0}\lambda^{2g-2}F^M_g(t^a_r;q),\eqno(6.8)$$
 where $\lambda$ is the genus expansion parameter. As in [W2] (referring to the proofs
of Lemma 6.1 and Lemma 6.2 in [RT2] for details) one can easily
derive from Theorem 5.9, (6.1) and (6.4) the following equations:
$$
\Bigl\langle\tau_{0,\gamma_1}\prod^m_{i=1}\tau_{d_i,\gamma_{a_i}}\Bigr\rangle_g(q)=
\sum^m_{j=1}\Bigl\langle\prod^m_{i=1}\tau_{d_i-\delta_{i,j},\gamma_{a_i}}\Bigr\rangle_g(q)+
\delta_{m,2}\delta_{d_1,0}\delta_{d_2,0}\zeta_{a_1,a_2},\eqno(6.9)$$
$$\left.\begin{array}{ll}
\bigl\langle\tau_{1,\gamma_1}\prod^m_{i=1}\tau_{d_i,\gamma_{a_i}}\bigr\rangle_g(q)=(2g-2-m)
\sum^m_{j=1}\bigl\langle\prod^m_{i=1}\tau_{d_i,\gamma_{a_i}}\bigr\rangle_g(q)\\
\hspace{38mm}+\frac{1}{24}\chi(M)\delta_{g,1}\delta_{m,0}.\end{array}\right.
\eqno(6.10)$$
 Here the integers $g\ge 0$, $m>0$ with $2g + m\ge 2$,
 $d_1,\cdots,d_m$ are nonnegative integers, and
 $\tau_{r,\alpha}=0$ if $r<0$.   As in [W2, RT2] it immediately follows from
  (6.9) and (6.10) that\vspace{2mm}

\noindent{\bf Theorem 6.1.}\quad{\it The above functions
$F^M(t^a_r;q)$ and $F^M_g(t^a_r;q)$ satisfy respectively the
following the generalized string equation and the dilation
equation:
\begin{eqnarray*}
&&\frac{\partial F^M}{\partial t_0^1}=\frac{1}{2}\zeta_{ab}t_0^a
t_0^b +\sum^\infty_{i=0}\sum_a t^a_{i+1}\frac{\partial
F^M}{\partial t^a_i},\\
&&\frac{\partial F^M_g}{\partial t_1^1}=\Bigl(2g-2+
\sum^\infty_{i=1}\sum_at^a_i\frac{\partial}{\partial
t^a_i}\Bigr)F^M_g+ \frac{\chi(M)}{24}\delta_{g,1}.\end{eqnarray*}
 Moreover, if $c_1(M)=0$, $F^M$ also satisfies the dilation equation
$$\frac{\partial F^M}{\partial t_1^1}=\sum^\infty_{i=1}\sum_a
\Bigl(\frac{2}{3-n}\Bigl(i-1+\frac{1}{2}\deg\gamma_a\Bigr)
+1\Bigr)t^a_i\frac{\partial F^M}{\partial t^a_i}+
\frac{\chi(M)}{24}.$$}

We may also make a bit generalization. Given a collection of
nonzero homogeneous elements $\underline\xi=\{\xi_i\}_{1\le
  i\le l}$  in
  $H^\ast_c(M,\mbox{\Bb C})\cup H^\ast(M,\mbox{\Bb C})$ we replace (6.6) by
$$
\langle\underline\xi|\tau_{d_1,\alpha_1}\tau_{d_2,\alpha_2}\cdots\tau_{d_m,\alpha_m}\rangle_g(q)
=\sum_{A\in H_2(M)}{\mathcal G}{\mathcal W}^{(\omega, \mu,
J)}_{A,g,m+l}(\kappa_{d_1,\cdots,d_m};\alpha_1,\cdots,\alpha_m,\underline\xi)q^A
\eqno(6.11)$$ and make the following

\noindent{\bf Convention}:\quad ${\mathcal G}{\mathcal
W}^{(\omega, \mu, J)}_{A,g,m+2}(\kappa;{\bf 1},{\bf 1},
\alpha_1,\cdots,\alpha_m)=0$ for any $m\ge 0$ whether or not this
case may be included in our definition category. (This is
reasonable if it may be defined.)

Then the same reason as in (6.9) and (6.10) gives rise to
$$
\langle\underline\xi|\tau_{0,\gamma_1}\prod^m_{i=1}\tau_{d_i,\gamma_{a_i}}\rangle_g(q)=
\sum^m_{j=1}\langle\underline\xi|\prod^m_{i=1}\tau_{d_i-\delta_{i,j},\gamma_{a_i}}\rangle_g(q),
\eqno(6.12)$$
$$
\langle\underline\xi|\tau_{1,\gamma_1}\prod^m_{i=1}\tau_{d_i,\gamma_{a_i}}\rangle_g(q)=(2g-2-m)
\sum^m_{j=1}\langle\underline\xi|\prod^m_{i=1}\tau_{d_i,\gamma_{a_i}}\rangle_g(q).
 \eqno(6.13)$$
It follows from them that\vspace{2mm}

\noindent{\bf Theorem 6.2}.\quad{\it The variants of (6.7) and
(6.8),
$$\left.\begin{array}{ll}F^M_g(\underline\xi|t^a_r;q)=\sum_{n_{r,a}}\prod_{r,a}\frac{(t^a_r)^{n_{r,a}}}{n_{r,a}!}\langle
 \underline\xi|\prod_{r,a}\tau^{n_{r,a}}_{r,\gamma_a}\rangle_g(q),\end{array}\right.\eqno(6.14)$$
$$\left.\begin{array}{ll}F^M(\underline\xi|t^a_r;q)=
\sum_{g\ge
0}\lambda^{2g-2}F^M_g(\underline\xi|t^a_r;q),\end{array}\right.\eqno(6.15)$$
 still called Witten's generating function, respectively satisfy
\begin{eqnarray*}
 &&\frac{\partial F^M(\underline\xi|\cdot)}{\partial
t_0^1}=\sum^\infty_{i=0}\sum_a t^a_{i+1}\frac{\partial
F^M(\underline\xi|\cdot)}{\partial t^a_i},\\
&&\frac{\partial F^M_g(\underline\xi|\cdot)}{\partial
t_1^1}=\Bigl(2g-2+
\sum^\infty_{i=1}\sum_at^a_i\frac{\partial}{\partial
t^a_i}\Bigr)F^M_g(\underline\xi|\cdot).\end{eqnarray*}
 They are still called the
generalized string equation and the dilation equation.
 Moreover, if $c_1(M)=0$, $F^M(\underline\xi|\cdot)$ also satisfies the dilation equation
$$\frac{\partial F^M(\underline\xi|\cdot)}{\partial t_1^1}=\sum^\infty_{i=1}
\sum_a \Bigl(\frac{2}{3-n}\Bigl(i-1+\frac{1}{2}\deg\gamma_a\Bigr)
+1\Bigr)t^a_i\frac{\partial F^M(\underline\xi|\cdot)}{\partial
t^a_i}.$$}

 It follows from the properties of Gromov-Witten invariants in $\S5$ that the
gravitational correlators defined in (6.5) have the following
properties(cf. [CoKa]).\vspace{2mm}

\noindent{\bf Degree Axiom}.\quad The gravitational correlator
$\langle\tau_{d_1,\alpha_1}\tau_{d_2,\alpha_2}\cdots\tau_{d_m,\alpha_m}\rangle_{g,
A}$ can be nonzero only if
$$\left.\begin{array}{ll}\sum^m_{i=1}(2d_i+\deg\alpha_i)=2c_1(M)(A)+ 2(3-n)(g-1)+ 2m.\end{array}\right.$$

\noindent{\bf Fundamental Class Axiom}.\quad Define $\tau_{-1,
\alpha}=0$, then
$$\left.\begin{array}{ll}\langle\tau_{d_1,\alpha_1}\tau_{d_2,\alpha_2}\cdots\tau_{d_{m-1},\alpha_{m-1}}\tau_{0,{\bf
1}}\rangle_{g,
A}=\sum^{m-1}_{j=1}\langle\prod^{m-1}_{i=1}\tau_{d_i-\delta_{i,j},\alpha_i}\rangle_{g,
A}.\end{array}\right.$$

\noindent{\bf Divisor Axiom}.\quad Assume that $m+2g\ge 4$ or
$A\ne 0$ and $m\ge 1$. Then for $\alpha_i\in H^\ast_c(M,\mbox{\Bb
Q})$, $i=1,\cdots, m$ and $\deg\alpha_m=2$, it holds that
\begin{eqnarray*}
&&\langle\tau_{d_1,\alpha_1}\tau_{d_2,\alpha_2}\cdots\tau_{0,\alpha_m}\rangle_{g,A}=
\alpha_m(A)\langle\tau_{d_1,\alpha_1}\tau_{d_2,\alpha_2}\cdots\tau_{d_{m-1},\alpha_{m-1}}\rangle_{g,A}\\
&&\qquad+\sum^{m-1}_{j=1}\langle\tau_{d_1,\alpha_1}\cdots\tau_{d_{j-1},\alpha_{j-1}}
\tau_{d_{j-1},\alpha_{j-1}\cup\alpha_m}\tau_{d_{j+1},\alpha_{j+1}}\cdots\tau_{d_{m-1},\alpha_{m-1}}\rangle_{g,A}.
\end{eqnarray*}

\noindent{\bf Dilation Axiom}.\quad If $2g+m\ge 3$ then
$$\langle\tau_{1,{\bf 1}}\tau_{d_1,\alpha_1}\tau_{d_2,\alpha_2}\cdots\tau_{d_m,\alpha_m}\rangle_{g,
A}=(2g-2+m)\langle\tau_{d_1,\alpha_1}\tau_{d_2,\alpha_2}\cdots\tau_{d_m,\alpha_m}\rangle_{g,A}.$$

\noindent{\bf Splitting Axiom}.\quad If $\dim H^\ast(M)<+\infty$
and $\alpha_1,\alpha_m\in H^\ast_c(M,\mbox{\Bb Q})$ then
\begin{eqnarray*}
&&\langle\tau_{d_1,\alpha_1}\tau_{d_2,\alpha_2}\cdots\tau_{d_m,\alpha_m}\rangle_{g,A}
=\sum_{i,j}\sum_{A_1+A_2=A}\eta^{ij}\\
&&\qquad\langle\tau_{d_1,\alpha_1}\tau_{d_2,\alpha_2}\cdots\tau_{d_{m_1},\alpha_{m_1}}
\tau_{0,\beta_i}\rangle_{g_1,A_1}
\langle\tau_{0,\beta_j}\tau_{d_{m_1+1},\alpha_{m_1+1}}\cdots\tau_{d_m,
\alpha_m}\rangle_{g_2,A_2}.
\end{eqnarray*}
where $g_k$, $m_k$, $k=1,2$, $\{\beta_i\}$ and $\eta^{ij}$ are as
in Theorem 5.7.

\subsection{The WDVV equation}

In this subsection we assume that
$$\dim H^\ast(M)<\infty.\eqno(6.16)$$
Let $\{\beta_i\}_{1\le i\le L}$ be a basis of $H^\ast(M, \mbox{\Bb
Q})$ consisting of homogeneous elements as in Theorem 5.7.
  We may assume them to satisfy that $\deg\beta_i$ is even if
  and only if $i\le N$. Let $\underline\alpha=\{\alpha_i\}_{1\le
  i\le k}$ be a collection of nonzero homogeneous elements in
  $H^\ast_c(M,\mbox{\Bb C})\cup H^\ast(M,\mbox{\Bb C})$ and at
  least one of them belongs to $H^\ast_c(M,\mbox{\Bb C})$. Putting $w=\sum t_i\beta_i\in
H^\ast(M,\mbox{\Bb C})$ we
  define $\underline\alpha$-{\it Gromov-Witten potential}  by a formal
power series in $q$,
$$
\Phi_{(q,\underline\alpha)}(w)=\!\!\sum_{A\in H_2(M)}\sum_{m\ge
\max(1, 3-k)}\frac{1}{m!}{\mathcal G}{\mathcal W}^{(\omega,\mu,
J)}_{A,0,k+m}([\overline{\mathcal M}_{0,k+m}]; \underline\alpha,
w, \cdots, w) q^A. \eqno(6.17)$$
 To avoid additional technicalities in using superstructures  we only consider the case
$w=\sum^N_{i=1}t_i\beta_i\\ \in W=H^{\rm even}(M,\mbox{\Bb
  C})$. Taking the third derivative on both sides of
 \begin{eqnarray*}
 \hspace{16mm}
 &&\quad  \Phi_{(q,\overline\alpha)}(w)=\sum_{A\in H_2(M)}\sum_{m\ge
\max(1, 3-k)}\sum_{1\le i_1,\cdots,i_m\le N}\frac{1}{m!}\hspace{35mm}(6.18)\\
&&\hspace{20mm}{\mathcal G}{\mathcal W}^{(\omega,\mu,
J)}_{A,0,k+m}([\overline{\mathcal M}_{0,k+m}];\overline\alpha,
\beta_{i_1}, \cdots, \beta_{i_m})t_{i_1}\cdots t_{i_m} q^A,
\end{eqnarray*}
 one gets
\begin{eqnarray*}
&&\frac{\Phi_{(q,\overline\alpha)}(w)} {\partial t_i\partial
t_j\partial t_a}=\sum_{A\in H_2(M)}\sum^\infty_{m=0}
\sum_{1\le i_1,\cdots,i_m\le N}\frac{1}{m!}\\
&&\hspace{15mm}{\mathcal G}{\mathcal W}^{(\omega,\mu,
J)}_{A,0,k+3+ m}([\overline{\mathcal
M}_{0,k+3+m}];\underline\alpha, \beta_i, \beta_j,\beta_a,
\beta_{i_1}, \cdots, \beta_{i_m})t_{i_1}\cdots t_{i_m}
q^A.\end{eqnarray*}
 It immediately follows from this and Theorem 5.7 that\vspace{2mm}

\noindent{\bf Theorem 6.3.}\quad{\it For any specified number $q$
the function $\Phi_{(q,\underline\alpha)}$ satisfies  {\bf
WDVV}-equation of the following form
$$\sum_{r,s}\frac{\partial^3\Phi_{(q,\underline\alpha)}}{\partial t_i\partial
t_j\partial
t_r}\eta^{rs}\frac{\partial^3\Phi_{(q,\underline\alpha)}}{\partial
t_k\partial t_l\partial
t_s}=\sum_{r,s}\frac{\partial^3\Phi_{(q,\underline\alpha)}}{\partial
t_i\partial t_k\partial
t_r}\eta^{rs}\frac{\partial^3\Phi_{(q,\underline\alpha)}}{\partial
t_j\partial t_l\partial t_s}\eqno(6.19)$$
 for $1\le i,j,k,l\le N$, where
 $\eta^{rs}=\int_M\omega_r\wedge\omega_s$ as in Theorem 5.7.}\vspace{2mm}

\noindent{By the statement} below Theorem 5.7,
$\Phi_{(q,\underline\alpha)}$ satisfies the ordinary WDVV equation
if $M$ is a closed symplectic manifold.\vspace{2mm}

\noindent{\bf Remark 6.4.}\quad The idea to prove (6.19) is that
using Theorem 5.7 and the direct computation show that both sides
of (6.19) are equal to
\begin{eqnarray*}
&&\sum_{A\in H_2(M)}\sum^\infty_{n=0}\sum^n_{m=0}
\sum_{1\le j_1,\cdots,j_m\le N}\frac{1}{m!(n-m)!}\\
&&\hspace{15mm}{\mathcal G}{\mathcal W}^{(\omega,\mu,
J)}_{A,0,2k+4+ m}([\overline{\mathcal
M}_{0,2k+4+m}];\underline\alpha,\underline\alpha, \beta_i,
\beta_j,\beta_k,\beta_l, \beta_{j_1}, \cdots,
\beta_{j_m})t_{i_1}\cdots t_{i_m} q^A.\end{eqnarray*}
 By the skew symmetry of Gromov-Witten invariants it is zero if
 some element of $\underline\alpha$ has odd degree. Therefore
 (6.17) is interesting only if $\underline\alpha$ consists of cohomology classes of even
 degree. Later we shall always assume this case without special
 statements. Clearly, (6.17) gives rise to a family of new solutions of the WDVV equation
  even if $M$ is a closed symplectic manifold.
 The reasons to introduce $\underline\alpha$-Gromov-Witten potential are that
on one hand we need to require $w\in H^\ast(M,\mbox{\Bb C})$
because of the composition laws
 and on the other hand the definition category of our Gromov-Witten invariants requires that at least one
 cohomology class  belongs to $H^\ast_c(M,\mbox{\Bb C})$.  Moreover, if we assume $\beta_1={\bf 1}$
 it follows from Theorem 4.2 that all terms containing $t_1$ at the right side of
 (6.19) are all zero.
  \vspace{2mm}

Following [RT1] we define a family of connections on the tangent
bundle $TW$ over $W$ as follows:
$$\nabla^\epsilon v=\sum_{i,j}\Bigl(\frac{\partial v^i}{\partial
t_j}+\epsilon\sum_{k,l}\eta^{il}\frac{\partial^3\Phi_{(q,\underline\alpha)}}{\partial
t_l\partial t_j\partial t_k}v^k\Bigr)\frac{\partial}{\partial
t_i}\otimes dt_j\eqno(6.20)$$
  for a tangent vector field $v=\sum_i v^i\frac{\partial}{\partial
t_i}$ in $TW$, where $\eta^{il}$ is as in Theorem 6.3. It is
easily checked that $\nabla^\epsilon\circ\nabla^\epsilon=0$ is
equivalent to the above WDVV equation (6.20). Therefore we
have\vspace{2mm}

\noindent{\bf Theorem 6.5.}\quad{\it $\{\nabla^\epsilon\}$ is a
family of the flat connections on the tangent bundle $TW$.}
\vspace{2mm}

More generally, let $\{\xi_i\}_{1\le i\le m}\subset
H^\ast(M,\mbox{\Bb C})\cup H_c^\ast(M,\mbox{\Bb C})$ and at least
one of them belong to $H_c^\ast(M,\mbox{\Bb C})$. For nonnegative
integers $d_1,\cdots, d_m$ and for $w=\sum^N_{i=1}t_i\beta_i\in
W=H^{\rm even}(M,\mbox{\Bb C})$ we may define
$\underline\alpha$-{\it genus $g$ couplings}
$$\left.\begin{array}{ll}
\langle\langle\tau_{d_1,\xi_1}\cdots\tau_{d_m,\xi_m}; w
\rangle\rangle_{g,\underline\alpha}=\!\!\sum_{A\in
H_2(M)}\sum_{m\ge
\max(1, 3-2g-k)}\frac{1}{m!}\\
\hspace{20mm}{\mathcal G}{\mathcal W}^{(\omega,\mu,
J)}_{A,g,k+m}(\kappa_{d_1,\cdots,d_m}; \underline\alpha, w,
\cdots, w) q^A\end{array}\right. \eqno(6.21)$$
 and $\underline\alpha$-{\it genus $g$ gravitational Gromov-Witten
 potential}
$$\left.\begin{array}{ll}
\Phi^{g,\{d_i\}}_{(q,\underline\alpha)}(w)=\!\!\sum_{A\in
H_2(M)}\sum_{m\ge \max(1, 3-2g-k)}\frac{1}{m!}{\mathcal
G}{\mathcal W}^{(\omega,\mu, J)}_{A,g,k+m}(\kappa_{d_1,\cdots,
d_m}; \underline\alpha, w, \cdots, w) q^A,\end{array}\right.
\eqno(6.22)$$
 and study their properties. We omit them.

\subsection{Quantum cohomology}

We still assume that (6.16) holds.  Let $H_2(M)$ be the free part
of $H_2(M,\mbox{\Bb Z})$. So far it has a finite integral basis
$A_1,\cdots, A_d$. Thus every $A\in H_2(M)$ may be written as
$A=r_1A_1+\cdots + r_dA_d$ for unique $(r_1,\cdots,
r_d)\in\mbox{\Bb Z}^d$. Denote by $q_j=e^{2\pi iA_j}$,
$j=1,\cdots, d$, and by $q^A=q_1^{r_1}\cdots q_d^{r_d}$. Following
[HS1, RT1, McSa1] an element of the Novikov ring
$\Lambda_\omega(\mbox{\Bb Q})$ over $\mbox{\Bb Q}$
 is the formal sum
$$\lambda=\sum_{A\in H_2(M)}\lambda_A q^A,$$
where the coefficients $\lambda_A\in\mbox{\Bb Q}$ are subject to
the finiteness condition $\sharp\{A\in H_2(M)\,|\,\lambda_A\ne
0,\;\omega(A)\le c\}<\infty$
 for any $c>0$. The multiplication in this ring is defined by
 $\lambda\ast\mu=\sum_{A,B}\lambda_A\mu_B q^{A+B}$. It has a
 natural grading given by $\deg q^A=2c_1(A)$. Denote by
$QH^\ast(M,\mbox{\Bb Q})=H^\ast(M,\mbox{\Bb
Q})\otimes\Lambda_\omega(\mbox{\Bb Q})$.

Let $\{\beta_i\}_{1\le i\le L}$ and $\underline\alpha$  be as in
$\S6.2$. For $\alpha,\beta\in H^\ast(M,\mbox{\Bb Q})$  we define
an element of $QH^\ast(M,\mbox{\Bb Q})$ by
$$
\alpha\star_{\underline\alpha}\beta=\sum_{A\in H_2(M)}\sum_{i,j}
{\mathcal G}{\mathcal W}^{(\omega,\mu,
J)}_{A,0,3+k}([\overline{\mathcal M}_{0,3+k}];\underline\alpha,
\alpha, \beta,\beta_i)\eta^{ij}\beta_j q^A. \eqno(6.23)$$
 More generally, for a given  $w=\sum^L_{i=1} t_i\beta_i\in H^\ast(M,\mbox{\Bb C})$
we  also define another element of $QH^\ast(M,\mbox{\Bb
C})=H^\ast(M,\mbox{\Bb C})\otimes\Lambda_\omega(\mbox{\Bb C})$ by
\begin{eqnarray*}
\hspace{5mm} &&\alpha\star_{(\underline\alpha,w)}\beta=\sum_{A\in
H_2(M)}\sum_{k,l}\sum_{m\ge 0}\frac{\epsilon(\{t_i\})}{m!}\hspace{66mm}(6.24)\\
&&\qquad\qquad{\mathcal G}{\mathcal W}^{(\omega,\mu,
J)}_{A,0,3+k+m}([\overline{\mathcal
M}_{0,3+k+m}];\underline\alpha, \alpha, \beta, \beta_k,
\beta_{i_1},\cdots, \beta_{i_m})\eta^{kl}\beta_l t_{i_1}\cdots
t_{i_m} q^A,
\end{eqnarray*}
 where $\epsilon(\{t_i\})$ is the sign of the
induced permutation on odd dimensional $\beta_i$. Clearly, (6.23)
is the special case of (6.24) at $w=0$. If $w=\sum^L_{i=1}
t_i\beta_i\in H^\ast(M,\mbox{\Bb Q})$, i.e., $t_i\in\mbox{\Bb Q}$,
$i=1,\cdots,L$, then $\alpha\star_w\beta\in QH^\ast(M,\mbox{\Bb
Q})$. We still  call the operations defined by (6.23) and (6.24)
``small quantum product'' and ``big small product'', respectively.
However, it is unpleasant that both
$\alpha\star_{\underline\alpha}{\bf 1}$ and
$\alpha\star_{(\underline\alpha, w)}{\bf 1}$ are always zero by
Theorem 4.1.  After extending it to $QH^\ast(M,\mbox{\Bb
C})=H^\ast(M,\mbox{\Bb C})\otimes\Lambda_\omega(\mbox{\Bb C})$ by
linearity over $\Lambda_\omega(\mbox{\Bb C})$ we as usual can
derive from Theorem 5.7\vspace{2mm}

\noindent{\bf Theorem 6.5.}\quad{\it Let $w\in H^\ast(M,\mbox{\Bb
C})$ and $\underline\alpha$ as above. Then
$$(\alpha\star_{(\underline\alpha, w)}\beta)\star_{(\underline\alpha, w)}\gamma
=\alpha\star_{(\underline\alpha,
w)}(\beta\star_{(\underline\alpha, w)}\gamma)$$
 for any $\alpha,\beta,\gamma\in H^\ast(M,\mbox{\Bb C})$. Consequently,
$QH^\ast(M,\mbox{\Bb C})$ is a supercommutative ring without
identity under the quantum products in (6.24).}

\end{document}